%% file: main.tex
\documentclass[11pt]{article}

\usepackage{microtype}
\usepackage{graphicx}
\usepackage[margin=1in,lmargin=1.25in,rmargin=1.25in]{geometry}
\usepackage[colorlinks,linkcolor=black,citecolor=black,pagebackref]{hyperref}
\renewcommand*{\backref}[1]{}%
\renewcommand*{\backrefalt}[4]{%
  \ifcase #1 %
    No citations.%
  \or
    Cited on page #2.%
  \else
    Cited on pages #2.%
  \fi
}%

\usepackage{xcolor}
\usepackage{listings}
\usepackage{color}
\usepackage{float}
\usepackage{tabularx}
\usepackage{booktabs}
\usepackage{comment}
\usepackage{subcaption}

\usepackage{ltablex} 
\usepackage{footnotehyper}

\usepackage{multirow}

\usepackage{multicol}

\usepackage{amsmath,amsfonts,bm,amssymb}

\usepackage{mathtools}
\usepackage{algorithm}
\usepackage{algpseudocodex}

\usepackage{etoc}

\newlength{\secindent}
\etocsetstyle{section}
  {\vspace{-\baselineskip}\raggedright\small\noindent}      
  {%
    \small\raggedright\parindent0pt
    \hangindent=\secindent  
    \hangafter=1            
    \leftskip0pt\rightskip0pt
  }
  {%
    \makebox[\secindent][l]{\S\etocnumber\quad}
    \textbf{\etocname}\hspace{1em}\textit{\etocpage}\par
  }
  {\vspace{0.5em}}                          

\etocsetstyle{subsection}
  {}
  {%
    \parindent0pt
    \footnotesize\hangindent=\dimexpr1.25\secindent\relax
    \hangafter=1
    \leftskip\secindent\rightskip0pt
  }
  {%
    \makebox[\dimexpr1.25\secindent\relax][l]{\etocnumber\quad}%
    \etocname\hspace{1em}\textit{\etocpage}\par
  }
  {}

\usepackage{natbib}
\renewcommand{\cite}{\citep}

\usepackage[compact,tiny]{titlesec}
\usepackage{enumitem}
\usepackage{indentfirst}

\lstdefinestyle{MATLAB}{
  language=Matlab,
  frame=single,  
  rulecolor=\color{black},
  tabsize=2,
  breaklines=true,  
  breakatwhitespace=false,
basicstyle=\ttfamily\footnotesize,
  keywordstyle=\color{blue},
  stringstyle=\color{purple},
  commentstyle=\color{green},
  morecomment=[l][\color{magenta}]{\#}
}

\usepackage{array}

\input{commands}

\title{An Empirical Study of Conjugate Gradient Preconditioners for Solving Symmetric Positive Definite Systems of Linear Equations}
\date{}

\begin{document}

\etocdepthtag.toc{main}
\etocsettagdepth{main}{subsection}
\etocsettagdepth{appendix}{section}

\makeatletter
\begin{flushleft}
\noindent \textbf{\@title} \\
\end{flushleft}
\begin{flushright}\itshape 
        Marc A. Tunnell, David F.~Gleich \\
        \normalfont \footnotesize \textsc{purdue university}
\par
\end{flushright}
\vspace{0pt}
\hrule 
\vspace{-9pt}
\thispagestyle{empty}
\makeatother


\newcommand{\numberofSPDmatrices}{79}
\newcommand{\numberofgraphs}{200}
\newcommand{\numberofconfigurations}{108}
\newcommand{\numberofpreconditioners}{10} 
\newcommand{\papertolerance}{\ensuremath{10^{-10}}}

\begin{abstract}
\noindent Despite hundreds of papers on preconditioned linear systems of equations, there remains a significant lack of comprehensive performance benchmarks comparing various preconditioners for solving symmetric positive definite (SPD) systems.
In this paper, we present a comparative study of \numberofSPDmatrices{} matrices
using a broad range of preconditioners.
Specifically, we evaluate \numberofpreconditioners{} widely used preconditoners across \numberofconfigurations{} configurations to assess their relative performance against using no preconditioner.
Our focus is on preconditioners that are commonly used in practice, are available in major software packages, and can be utilized as black-box tools without requiring significant \textit{a priori} knowledge. 
In addition, we compare these against a selection of classical methods.
We primarily compare them without regards to effort needed to compute the preconditioner. 
Our results show that symmetric positive definite systems are mostly likely to benefit from incomplete symmetric factorizations, such as incomplete Cholesky (IC). Multigrid methods occasionally do exceptionally well.
Simple classical techniques, symmetric Gauss Seidel and symmetric SOR, are not productive.
We find that including preconditioner construction costs significantly diminishes the advantages of iterative methods compared to direct solvers; although, tuned IC methods often still outperform direct methods.
Additionally, ordering strategies such as approximate minimum degree significantly enhance IC effectiveness.
We plan to expand the benchmark with larger matrices, additional solvers, and detailed metrics to provide actionable information on SPD preconditioning.
\end{abstract}

\newcommand{\transpose}[1]{\ensuremath{#1^{\top}}}
\newcommand{\ntranspose}[1]{\ensuremath{#1^{-\top}}}


\vspace{0pt}
\hrule 
\begin{multicols}{2}
\small\raggedright
\vspace*{-25pt}
\tableofcontents
\end{multicols}

\vspace{3pt}
\hrule 

\section{Introduction}
\label{sec:introduction}

We consider using preconditioned iterative methods to solve $\mA \vx = \vb$ where $\mA$ is a large, sparse, symmetric positive definite (SPD) matrix. We assume the matrix data is available in a standard sparse data structure such as compressed sparse row or similar. Hence, we are not treating \emph{matrix-free} or \emph{operator} methods. It is well known that the convergence rate of iterative solvers  depends heavily on the condition number of the matrix~\cite[for example]{greenbaum_siam_1997}.
To improve the condition number, preconditioners are used to transform the original system into one with more favorable spectral properties, thereby accelerating the rate of convergence as a function of the number of iterations~\cite{greenbaum_siam_1997,saad_siam_2003,golub_jhup_2013_pcg,benzi_jcp_2002}.

Over the past decades, the world of preconditioners has grown. It now includes relatively simple techniques such as symmetric scaling~\cite{saad_siam_2003}, iterative methods such as symmetric successive over-relaxation (SSOR)~\cite{golub_jhup_2013_pcg,young_tams_1954}, and incomplete factorizations~\cite{benzi_jcp_2002,meijerink_moc_1977,saad_nla_1994,li_siam_2003}. On the other end of the spectrum, multigrid methods~\cite{briggs_siam_2000}, algebraic multigrid (AMG)~\cite{brandt_siam_2011}, and domain decomposition based methods~\cite{smith_cup_2004} represent more complicated techniques that promise even faster convergence.

Despite an abundance of research on novel preconditioners for symmetric, positive definite systems---including reviews of their theoretical performance~\cite{benzi_jcp_2002,stuben_elsevier_2001,pearson_gamm_2020}---there is no comprehensive performance benchmark that provides an empirical comparison of the performance of these preconditioners.

Several studies have attempted empirical evaluations of preconditioners, but each targets a distinct class of problems or otherwise have limitations that differentiate them from this study.
We briefly discuss existing efforts next, organized by the types of systems or preconditioners evaluated, in order to highlight the gap we address. 

\paragraph{Benchmarks of SPD systems.}
Those comparisons that exist for SPD matrices are often limited to only a few matrices as in both~\citet{benzi_jcp_2002} and~\citet{stuben_elsevier_2001}. This is also the case in most papers proposing new preconditioners.
Several studies have attempted a comprehensive empirical evaluation of preconditioners for SPD systems, but these have notable limitations or otherwise fill a gap that is different from that addressed in this study.
For instance,~\citet{george_acm_2012} conducted an empirical analysis of preconditioners for SPD systems.
However, their study was limited to just 30 matrices, conducted on a system with significantly restrictive memory compared to modern hardware, limited to just 1000 solver iterations, and aimed for a relative tolerance of only $10^{-5}$ within that iteration limit.
Furthermore, the types of preconditioners tested were not unique across these libraries; for instance, they evaluated ILU preconditioners from each library, resulting in significant overlap in the preconditioner classes examined.

More recently, \citet{kyng_arxiv_2023} presents a benchmark on 28 symmetric diagonally dominant matrices that contain non-negative off-diagonal values from SuiteSparse and a large number of synthetic graphs.
Thus, the scope of their study partially overlaps with ours but does not directly align with our broader investigation of SPD systems.
Similarly,~\citet{inca_ukaea_2021} conducted a benchmarking study on a narrow set of preconditioners and limited their testing set to PDEs.

\paragraph{Benchmarks of non-symmetric matrices.}
In~\citet{benzi_siam_2000}, the authors provide an empirical evaluation of preconditioners for non-symmetric matrices, while~\citet{ghai_fos_2016} provides a more up-to-date empirical evaluation for large-scale non-symmetric systems.
In~\citet{benzi_bit_1998} and~\citet{benzi_anm_1999}, sparse approximate inverse preconditioners are empirically analyzed and validated against the standard preconditioning techniques of the time.

\citet{chen_iclr_2025} presents a benchmark of their neural preconditioner, comparing against supernodal ILU and AMG methods across a large number of matrices.
This analysis is limited to the strictly non-positive definite case and uses a non-symmetric solver, 
which is entirely disjoint with our scenario in this study.


\paragraph{Other problem classes.}
\citet{Orban2014} studies preconditioning for symmetric quasi-definite systems that arise in optimization and utilizes their unique structure for the problem. 
\citet{schork_springer_2020} evaluates preconditioning in the context of linear programming, testing a large number of medium- to large-scale instances.
Their work targets general linear programming rather SPD systems, and their preconditioning methods focuses on iterative linear algebra approaches tailored to the structure of KKT systems. 

\paragraph{Addressing this gap.}
These limitations highlight the need for a comprehensive performance benchmark of SPD systems across a wide range preconditioners and matrices.
Such a benchmark would help establish clearer guidelines on the applicability and efficiency of different preconditioners across a wide range of real-world SPD systems.

In this study, we aim to address this gap by providing a thorough performance evaluation of a wide range of preconditioner classes on \numberofSPDmatrices{} SPD matrices
from the SuiteSparse Matrix Collection~\cite{davis_acm_2011}.\footnote{There are more SPD matrices in SuiteSparse. There were selected by picking all SPD matrices with at least 10,000 rows that converged within $10n$ CG iterations, we continue to study the remaining cohort of large matrices, including those built from graph Laplacians.} 
Our focus is on preconditioners that are widely used in practice, that are available in large software packages, and may be utilized as a black-box without requiring significant \textit{a priori} knowledge of the underlying properties of the system.
We choose these types of preconditioners to ensure the relevance of this study to the broadest possible community.

Specifically, we test preconditioners from the classes outlined in \autoref{tab:preconditioners} and compare their performance against the unpreconditioned but symmetric diagonally scaled case.\footnote{This symmetric diagonal scaling is often called a Jacobi preconditioner. We treat this as the minimum preconditioning for any problem. See \citet[Section 2.3]{Bradley-thesis} for discussion of how to do this in the matrix-free case.} Note that there are many individual configurations for each of these preconditioners.

\begin{table}[b]
    \centering
    \begin{tabularx}{0.9\linewidth}{>{\raggedleft\arraybackslash}X>{\raggedright\arraybackslash}X}
        \toprule 
        \textbf{Preconditioner Class} & \textbf{Preconditioners} \\ 
        \midrule 
        Classical Methods & 
                Symmetric Gauss-Seidel (SGS), SSOR, 
                Truncated Neumann
            \\ \midrule
        Sparse Approximate Inverse & SPAI \\  \midrule
        Incomplete Factorizations & 
                Incomplete Cholesky (IC), Modified IC,
                SuperLU
            \\  \midrule
        Algebraic Multigrid & Ruge-Stuben, Smoothed Aggregation \\ \midrule
        Graph Laplacian Preconditioners & 
        Approximate Cholesky \\ 
        \bottomrule 
    \end{tabularx}
    \caption{Classes of preconditioners evaluated in this study.}
    \label{tab:preconditioners}
\end{table}


Later, in \autoref{sec:preconditioners}, we provide further details on each preconditioner, including their properties, computational workload, implementation software, optimal parameter configurations, recommended orderings, and additional practical guidance.
We utilize the preconditioned CG (PCG) method as our solver, specifically the Hestenes-Stiefel formulation~\cite{hestenes_jrnbs_1952,golub_jhup_2013_pcg}.
In a future iteration of this manuscript, we will include the results of the minimum residual, \textsc{minres}, method.

We aim to eliminate as many extraneous variables as possible and focus primarily on comparing the effectiveness of the preconditioners in reducing the overall computational workload.
Thus, we study and compare preconditioners in terms of floating point operations needed to converge to a given tolerance.
This approach excludes the time needed to construct the preconditioner, which can be significant, but crucially, it removes variance from differences in implementation quality.
While floating point operations are only a proxy for the more desirable measure of \textit{computation time}, the idea here is to study and assess preconditioner effectiveness independent of hardware and software details and optimizations. 

In further analyses of the data, we include a proxy measure for preconditioner construction time based on the fill-in of an incomplete factorization. This enables a rough comparison against direct methods using the same fill-in construction estimate (\autoref{sec:direct}).

\section{Results}\label{sec:results}

We measure the computational work for the relative residual to fall below $10^{-10}$. For an $n \times n$ matrix, we use the value 
\[ \text{work} = (5n + \text{nonzeros} + \text{preconditioner}) \cdot \text{iterations} + \text{preconditioner} \]
as the work for a preconditioned method. 
This reflects the 2 inner-products and 3 scaling operations in the PCG iteration combined with the work involved in the matrix-vector product as well as any work involved in applying the preconditioner. Note that there is one additional application of the preconditioner to startup the method.
We provide detail on how the cost of each preconditioner is measured in the relevant subsection of \autoref{sec:preconditioners} and note that in many cases, such as classical methods and incomplete factorizations, this has a simple cost in terms of the non-zero structure. 

We discuss key results in four scenarios: the first group reflects the raw results of our measurements for black-box preconditioners on  our test matrices (\autoref{subsec:single_matrix_results}), the second provided performance profiles that summarize aggregate benefits (\autoref{subsec:performance_profiles}),  the third analyzes the impact of reordering schemes on performance (\autoref{sec:ordering}), and the fourth compares against direct methods (\autoref{sec:direct}), both including and excluding preconditioner setup cost.




\subsection{Black-box preconditioner performance}\label{subsec:single_matrix_results}

\begin{figure}[tbp]
    \centering
    \includegraphics[width=\textwidth]{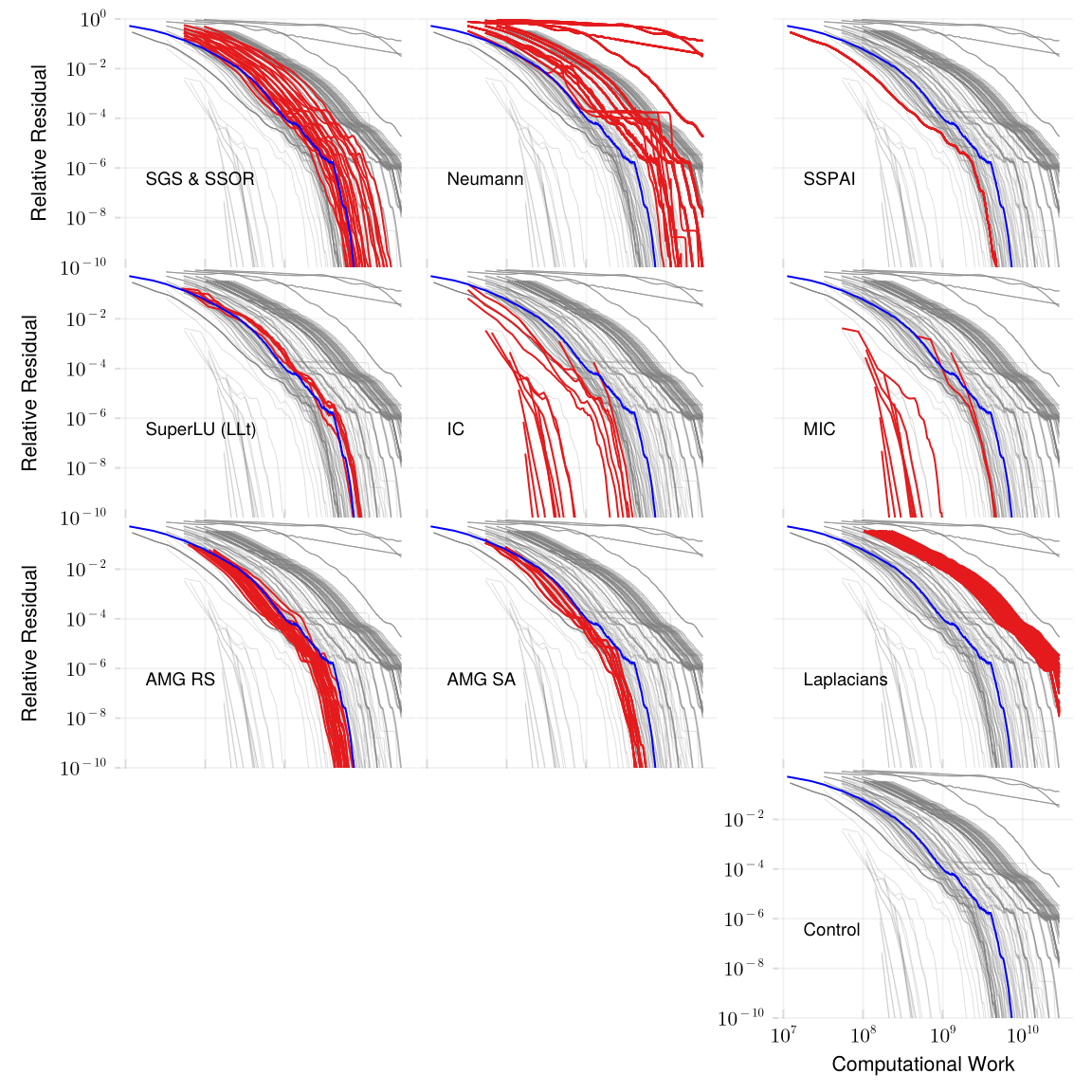}
    \caption{Convergence of the PCG method with various preconditioners applied to the \texttt{crankseg\_1} matrix (53k rows, 10m non-zeros, condition number $2 \cdot 10^4$). The plots have a log-log scale. The control case (diagonal scaling only) is shown in blue for reference, the preconditioner of focus is highlighted in red, and all other preconditioners are shown in gray.}
    \label{fig:crankseg_1}
\end{figure}

\begin{figure}[tbp]
    \centering
    \includegraphics[width=\textwidth]{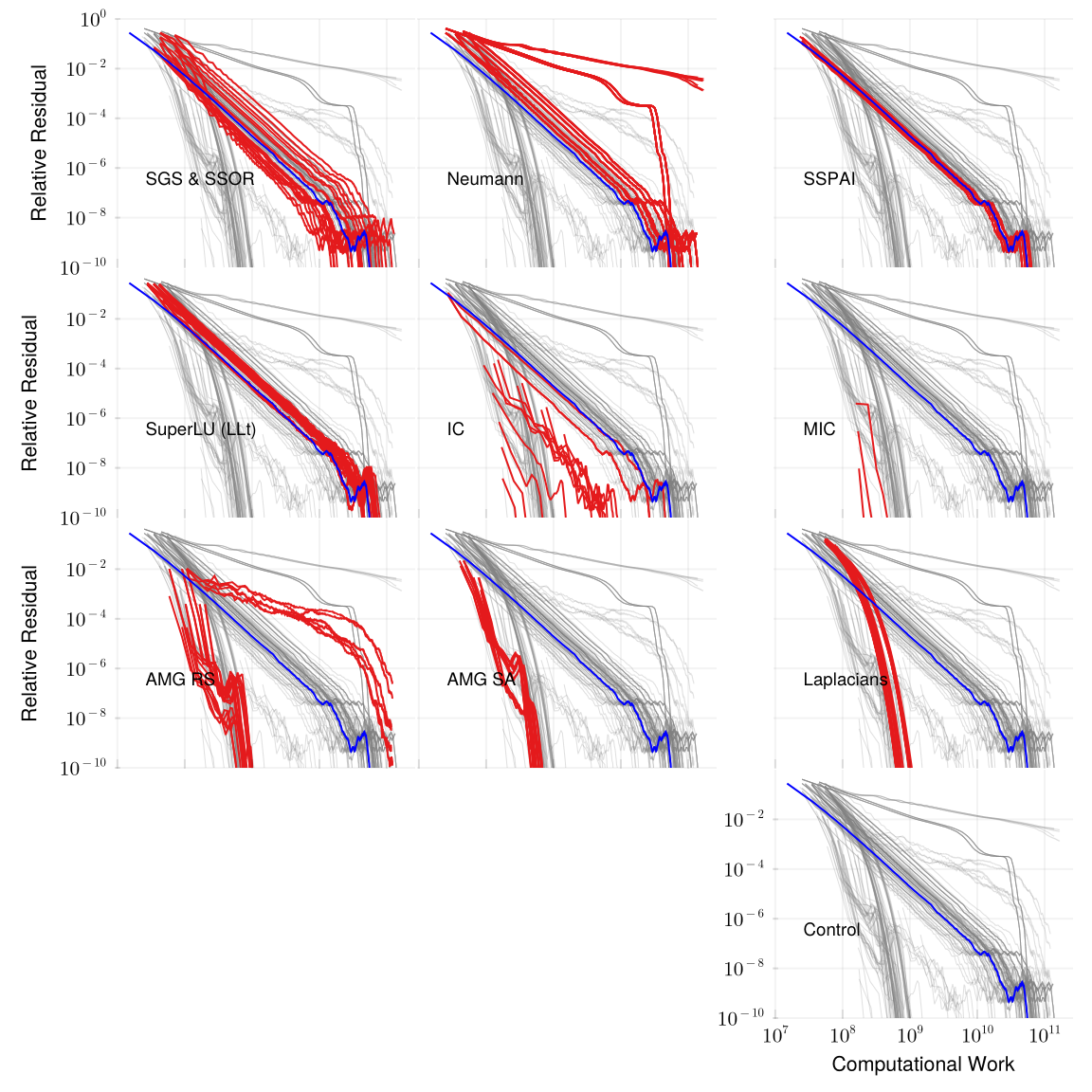}
    \caption{Convergence of the PCG method with various preconditioners applied to the \texttt{ecology2} matrix (1m rows, 5m non-zeros, condition number $6 \cdot 10^7$). The plots have a log-log scale. The control case (diagonal scaling only) is shown in blue for reference, the preconditioner of focus is highlighted in red, and all other preconditioners are shown in gray.}
    \label{fig:ecology2}
\end{figure}

We first present the  measurements from our study. These figures illustrate the convergence of the PCG method with various preconditioners applied to individual matrices.
Each figure displays a grid of subplots, one for each preconditioner class, where the relative residual is plotted against the work.
The control case (diagonal scaling only) is shown in blue for reference, the preconditioner of focus is highlighted in red, and all other preconditioners are shown in gray.

We start by showing a selection of results for \texttt{crankseg\_1}, a matrix from structural analysis that models a crankshaft~\cite{crankseg_1}.
This matrix has dimension $53$k and has $10.6$m non-zero elements.
We compute the condition number of the symmetrically scaled matrix to a tolerance of $10^{-8}$ and find it to be $2.2707 \cdot 10^{4}$.
This matrix is not diagonally dominant and has an average degree of $201.01$, which is relatively high.
The results for this matrix can be seen in \autoref{fig:crankseg_1}.

From the figure, we see that the IC and Modified IC preconditioners are the most effective, reducing the computational cost in the best case by over one order of magnitude.
The AMG and Symmetric SPAI preconditioners show minor improvements in decrease of total work, while the classical methods offer no improvement over the control case.
Because the matrix is not diagonally dominant, the graph Laplacian methods are not effective, and the SuperLU preconditioner is similarly ineffective.

We next show a selection of results for \texttt{ecology2}, a matrix using a 5-point grid stencil on a 1000-by-1000 mesh that is used to model animal movement.
This matrix has dimension $1$m and has $5$m non-zero elements.
We compute the condition number of the symmetrically scaled matrix to a tolerance of $10^{-8}$ and find it to be $6.3183 \cdot 10^7$.
This matrix is diagonally dominant and almost all rows have 4 neighbors.
The results for this matrix can be seen in \autoref{fig:ecology2}.

This is a problem where we would expect multigrid to work well. From the figure, we see that this is indeed the case. The methods based on algebraic multigrid are also effective, reducing the computational cost similarly to the IC preconditioners.
Additionally, the graph-Laplacian-based Cholesky preconditioner is effective on this matrix, reducing the computational cost by a similar margin. Finally, we see that the IC and Modified IC preconditioners are the most effective, reducing the computational cost in the best case by almost two orders of magnitude.
Again, the classical preconditioners are largely ineffective here, although the SGS and SSOR preconditioners show minor improvements in the best case.
The SuperLU and Symmetric SPAI preconditioners are largely ineffective on this matrix.

The performance of the graph Laplacian preconditioner on \texttt{ecology2} is typical of the performance of this preconditioner on approximately half of the diagonally dominant matrices in the test set.
The algebraic multigrid preconditioners have similar performance to the results shown in \autoref{fig:ecology2} on the majority of the diagonally dominant matrices in the test set.

The rest of the raw data is shown in this form in \autoref{sec:appendix1}.

\begin{figure}[t!]
    \centering
    \includegraphics[width=\textwidth]{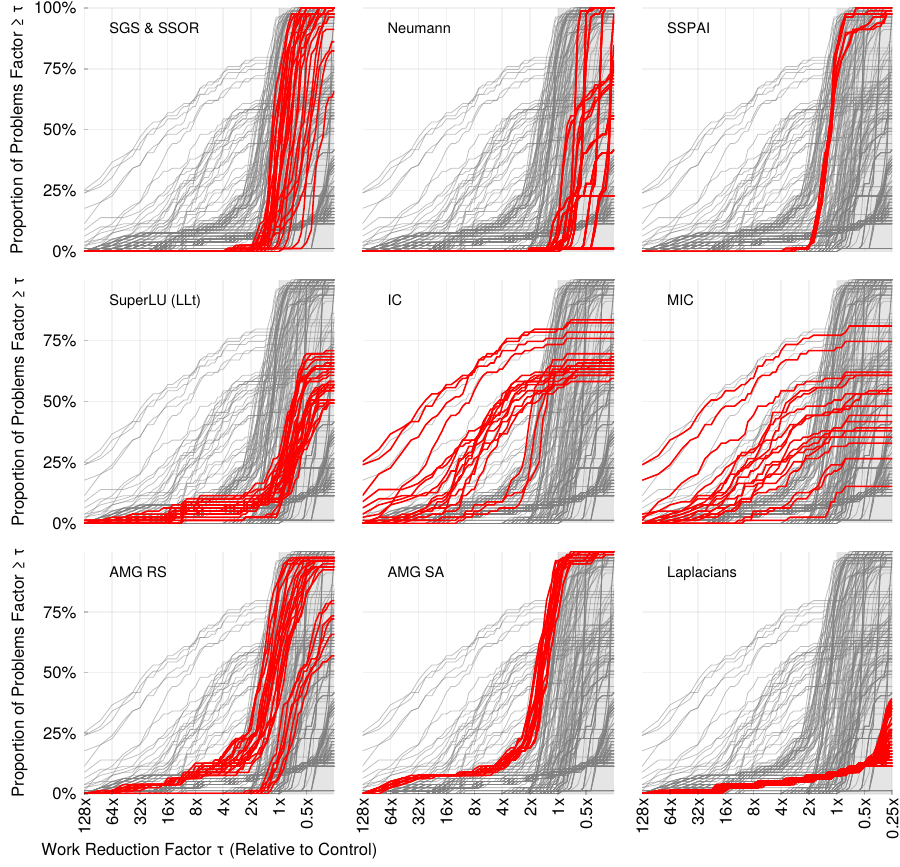}
    \caption{Performance profiles of the PCG method with various preconditioners applied to all matrices in the test set. The red lines show the behavior of multiple configurations of the same preconditioner; the grey lines show the entirely of the results. These show that the IC and MIC have the highest aggregate behavior, but only succeed in slightly over 75\% of problems.}
    \label{fig:all_problems_profile}
\end{figure}

\begin{figure}[t!]
  \centering
  \begin{subfigure}[t]{0.48\textwidth}
    \includegraphics[width=\textwidth]{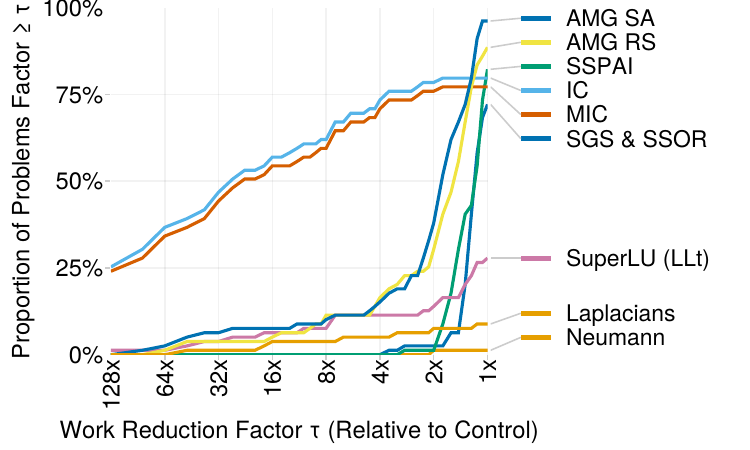}
    \caption{Single‐best configuration}%
    \label{fig:single_best_config}
  \end{subfigure}
  \quad
  \begin{subfigure}[t]{0.48\textwidth}
    \includegraphics[width=\textwidth]{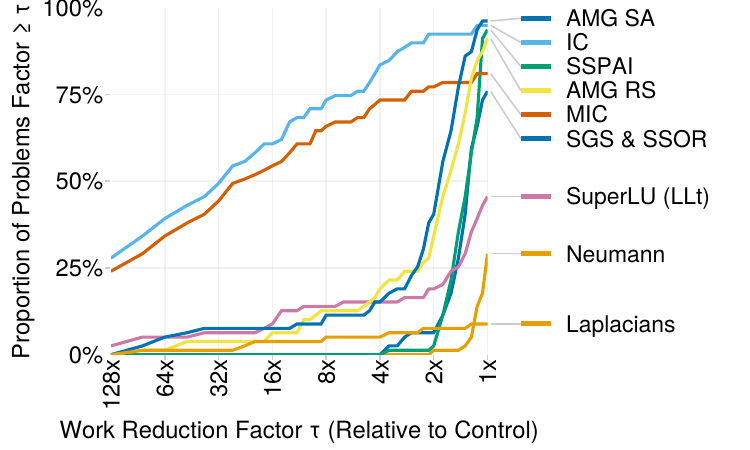}
    \caption{ Tuned‐best configuration.}%
    \label{fig:tuned_best_config}
  \end{subfigure}
  \caption{%
    The single-best configuration selects the best configuration over all matrices whereas the tuned-best configuration selects the best run for each matrix within a given preconditioner group (log‐linear). While there are differences, the single best configuration retains the majority of the performance across a wide range of problems. 
  }
  \label{fig:all_problems_profile_best}
\end{figure}

\subsection{Preconditioner performance profiles}\label{subsec:performance_profiles}

We now study the preconditioners across problems. This is done using performance profiles. A performance profile indicates the fraction of problems for which a preconditioner achieves a certain factor of reduction in computational effort to reduce the relative residual to \papertolerance{} compared with the diagonally scaled \emph{control} case. 

In each figure, there are multiple subplots, displaying one graph per preconditioner group reflecting all the different configurations of that preconditioner. 
The graph is lightly shaded in the region where the preconditioner is more costly than the control, starting at the horizontal axis value of $1\times$.
The preconditioner group in focus is shown in red, while all other preconditioners are shown in gray.
For each of these plots, the horizontal axis represents the factor by which the computational cost is decreased (or increased) relative to the computational cost of convergence to \papertolerance{} with no preconditioner.
The vertical axes show the fraction of problems where the preconditioner achieves that performance.

We show the performance profile over all matrices in the test set in \autoref{fig:all_problems_profile}.
This figure shows that the IC and modified IC preconditioners are generally the most effective at reducing work in solving the preconditioned system once given the preconditioner. Both IC and modified IC reduce the computational cost by one order of magnitude on a significant fraction of the test set. A downside of these methods is that they fail on around 20\% of problems. 
The AMG preconditioners are also at least slightly effective on a significant fraction of the set while succeeding on virtually all problems. 
We found that the classical methods were largely ineffective overall, though they do reduce the computational cost on a small fraction of the test set. 
The AMG methods and SSPAI show more effectiveness on larger problem sizes and we explore this in more depth in \autoref{sec:large}. Note that we would not expect SuperLU or Laplacians to succeed on a large fraction (see discussion in \autoref{sec:caveats}).

Seeing this information for all configurations makes direct comparisons among preconditioner types difficult. 
In \autoref{fig:single_best_config}, we select the single-best configuration applied uniformly across all matrices.
In contrast, \autoref{fig:tuned_best_config} shows the best configuration selected individually for each matrix, simulating per-matrix parameter tuning.
The single-best configuration is determined by computing the area under the performance profile curve (AUC).\footnote{See \autoref{sec:auc} for details on how to compute the AUC.}

This figure shows that, under both the single-best configuration and per-matrix tuning, factorization-based preconditioners remain the most effective. Moreover, there is less difference between the single best choice and a tuned configuration than one might expect. For IC, the biggest difference arises for problems that are between 1 and 4 times faster. For MIC, there is only a small difference. The best configurations for each method are summarized in~\autoref{tab:preconditioners}.


\subsection{Impact of generation cost on preconditioner performance}

\begin{figure}[tbp]
    \centering
    \includegraphics[width=\textwidth]{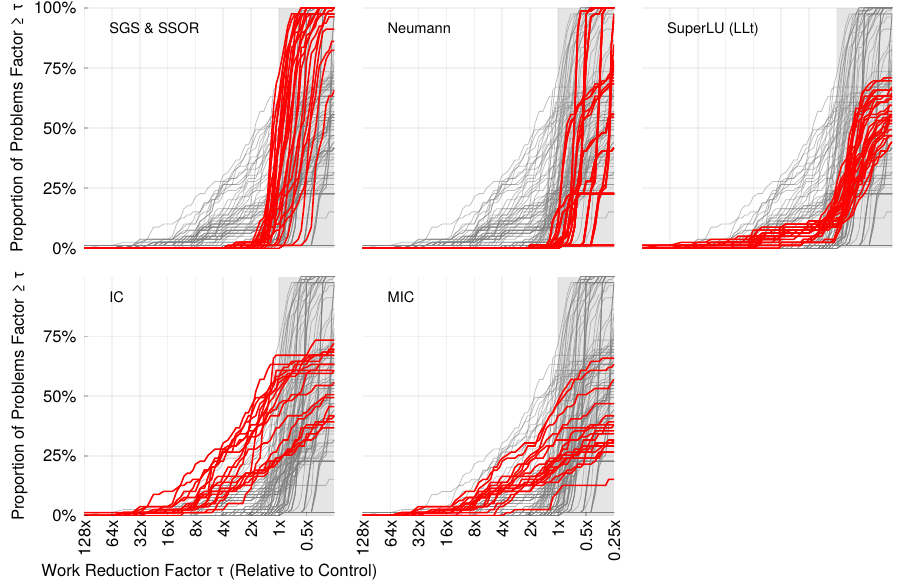}
    \caption{\emph{Including generation work.} Performance profiles of the PCG method with various preconditioners applied to all matrices in the test set where work includes an approximation to the amount of work required to generate the preconditioner. }
    \label{fig:all_problems_profile_and_work}
\end{figure}

\begin{figure}[tbp]
  \centering
  \begin{subfigure}[t]{0.48\textwidth}
    \includegraphics[width=\textwidth]{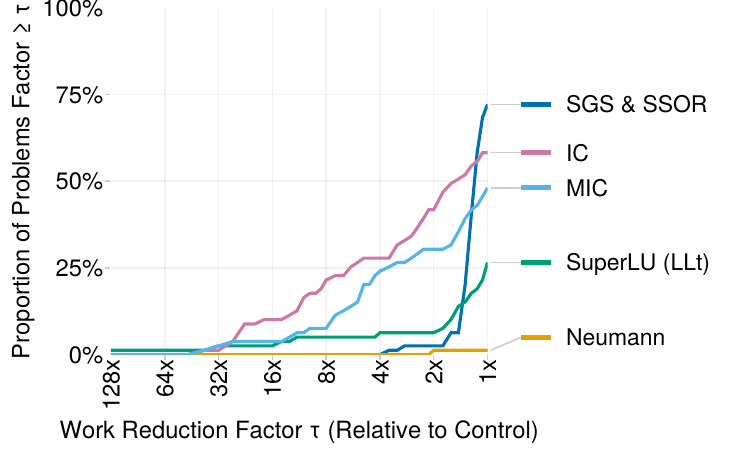}
    \caption{Single‐best configuration.}%
    \label{fig:single_best_config_and_work}
  \end{subfigure}
  \quad
  \begin{subfigure}[t]{0.48\textwidth}
    \includegraphics[width=\textwidth]{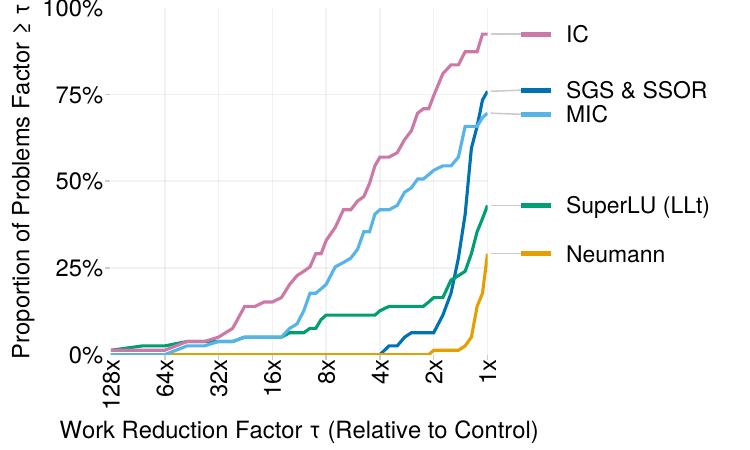}
    \caption{Tuned‐best configuration.}%
    \label{fig:tuned_best_config_and_work}
  \end{subfigure}
  \caption{%
    \emph{Including generation work.}
    The single-best configuration selects the best configuration over all matrices whereas the tuned-best configuration selects the best run for each matrix within a given preconditioner group. Both include an approximation the the amount of work required to generate the preconditioner. In comparison with \autoref{fig:all_problems_profile_best}, we see that tuning the preconditioner now has a large impact on the performance. 
  }
  \label{fig:all_problems_profile_best_and_work}
\end{figure}

For preconditioners where the construction cost can be reasonably estimated (classical and incomplete factorization methods), we now consider the total computational effort required, including both the preconditioner construction and iterative solve.
We discuss how we compute these cost estimates in \autoref{sec:generation-cost}.
We show the same set of figures for this new idea of work in~\autoref{fig:all_problems_profile_and_work} (all individual configurations) \autoref{fig:single_best_config_and_work} (single best configuration), and \autoref{fig:tuned_best_config_and_work} (best configuration for each matrix). As expected, these show that accounting for factorization cost reduces the overall effectiveness of preconditioning. 
Even a modest $4\times$ gains occur on less than $28\%$ of problems for any single-best preconditioner in the best case (IC).
What is most interesting is that the \emph{tuned} configuration now shows a substantial difference from the single-best configuration. 
However, this assumes \emph{a priori} knowledge of the optimal configuration for each matrix: exhaustively trying multiple variants would itself incur significant work and largely erode these apparent gains.

\subsection{Impact of ordering on preconditioner performance}
\label{sec:ordering}

\begin{figure}[tbp]
    \centering
    \includegraphics[width=\textwidth]{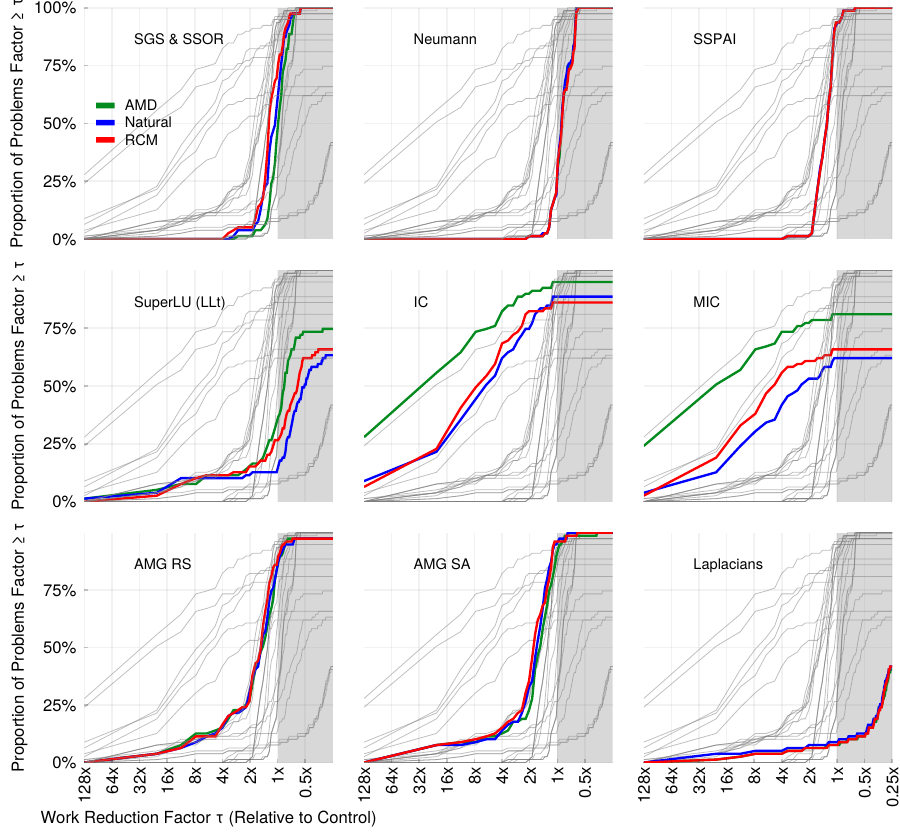}
    \caption{\emph{Focusing on ordering.} Performance profiles of CG with various preconditioners applied to all matrices in the test set, grouped by ordering. We expect no impact to Neumann, SSPAI, AMG RS, AMG SA, and Laplacians. This largely follows although we see some modest effects for the AMG methods. Overall, AMD seems superior and there is further analysis of fill-in in \autoref{sec:fill-in} }
    \label{fig:ordering_profile}
\end{figure}

For sparse factorization methods, such as the incomplete factorizations we study, the matrix ordering has a large impact on overall performance. We now evaluate preconditioner performance by selecting the best run for each matrix (across all parameter settings) within a given ordering and then compare the resulting performance profiles.
This approach simulates tuning the preconditioners parameters for each matrix while isolating the effect of reordering. The results are shown in \autoref{fig:ordering_profile}.


A number of the methods should be insensitive to order including Neumann series, SSPAI, and the AMG methods. This mirrors what we found. There were slight differences with the AMG  methods that we continue to study. The graph Laplacians solver does its own ordering and so is insensitive to the input ordering.

There was a small increase in performance for SSOR and SGS with the RCM ordering compared to natural ordering, while AMD ordering saw a slight degradation in performance for these  preconditioners.
For the IC and MIC preconditioners, the AMD ordering was by far the most effective, although the RCM ordering still provided some performance benefit.
The SuperLU preconditioner was largely invariant to the ordering when considering the curve to the left of the $1\times$ work cutoff, though there were still small differences.

\begin{figure}[tbp]
    \centering
    \includegraphics[width=\textwidth]{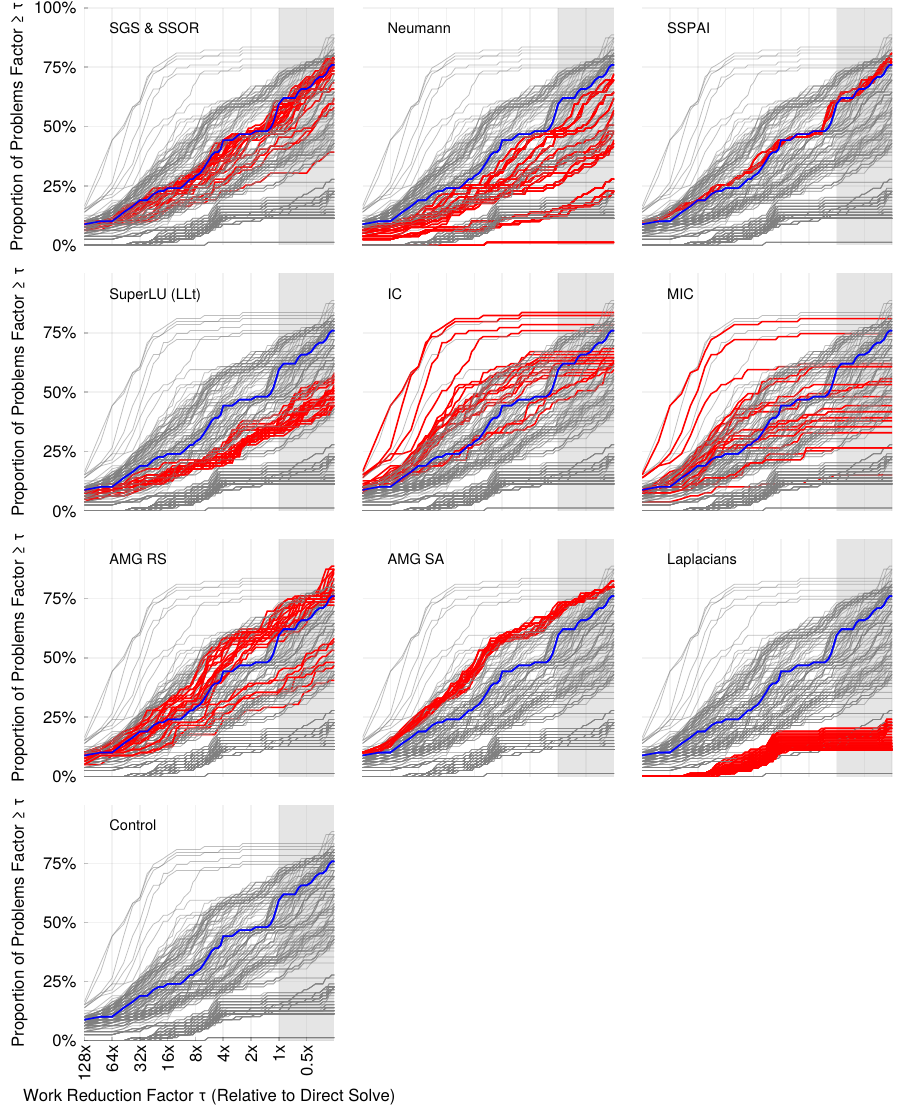}
    \caption{\emph{Comparison to direct solve.} Performance profiles of CG  with various preconditioners relative to the work required to solve the system with a sparse direct method. The blue line shows the unpreconditioned baseline whereas the   }
    \label{fig:performance_profile_vs_direct}
\end{figure}

\subsection{Comparison to direct methods}\label{sec:direct}
Direct methods for sparse systems have been honed over decades to be efficient, reliable, and scalable. 
This subsection presents the same performance profiles as in previous sections, but with computational effort measured relative to solving the system directly.
We use SuiteSparse to compute the factorization used as the baseline for directly solving \citep{davis_acm_2011}.
The cost of direct methods is estimated by the number of nonzero elements in the factorization plus an approximation of the computational work required, computed as the sum of squares of the nonzero counts per row in the factorization. See additional discussion in \autoref{sec:generation-cost}. We also compared this to a measure of factorization flops previously reported on the SuiteSparse website (e.g.~\url{https://www.cise.ufl.edu/research/sparse/matrices/AMD/G3_circuit.html}), and we found that it largely captured the costs although sometimes there were significant differences in either direction.
As before, the performance profiles indicate the fraction of problems for which a preconditioner achieves a certain factor of reduction in computational effort to reduce the relative residual to $10^{-10}$.

\begin{figure}[t]
  \centering
  \begin{subfigure}[t]{0.48\textwidth}
    \includegraphics[width=\textwidth]{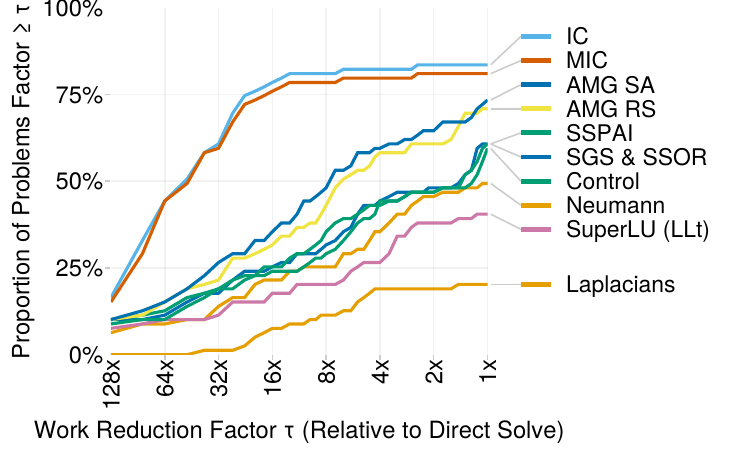}
    \caption{Single‐best configuration; }%
    \label{fig:single_best_config_direct}
  \end{subfigure}
  \quad
  \begin{subfigure}[t]{0.48\textwidth}
    \includegraphics[width=\textwidth]{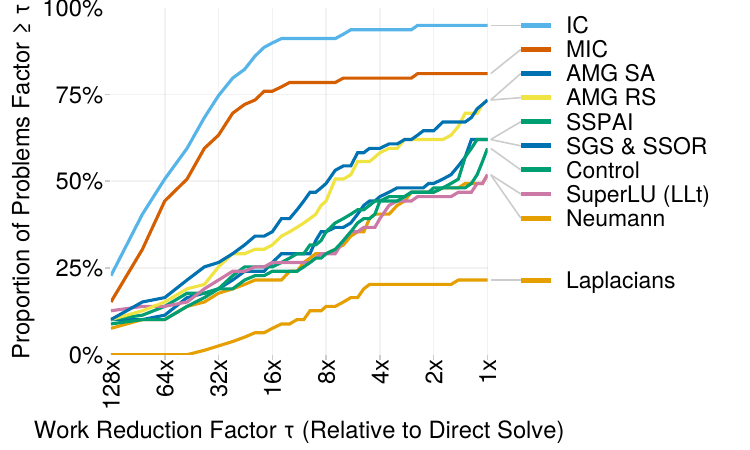}
    \caption{Tuned‐best configuration.}%
    \label{fig:tuned_best_config_direct}
  \end{subfigure}
  \caption{%
    \emph{Comparison to direct solve.}
    The Single-best configuration selects the best configuration over all matrices whereas the tuned-best configuration selects the best run for each matrix within a given preconditioner group. Plotted relative to solving directly. (log‐linear scale)
  }
  \label{fig:all_problems_best_direct}
\end{figure}

We present the aggregated performance profile across all matrices in \autoref{fig:performance_profile_vs_direct} and show the best single configuration and tuned configurations in Figures~\ref{fig:single_best_config_direct},~\ref{fig:tuned_best_config_direct}. 
Consistent with previous observations, IC and MIC preconditioners remain the most effective overall.
Notably, AMG methods also reduce work compared with direct methods for the majority of matrices.
Interestingly, using unpreconditioned CG alone already provides at least a twofold improvement over direct solution on $46.84\%$ of problems, and at least parity ($1\times$) on $59.49\%$ of problems.



\subsection{Comparison to direct methods with generation cost}

\begin{figure}[tpb]
    \centering
    \includegraphics[width=\textwidth]{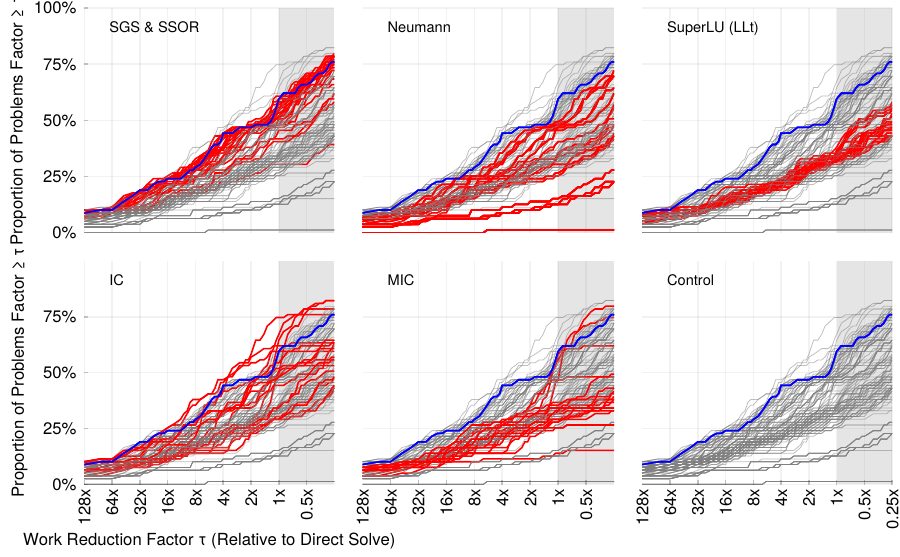}
    \caption{\emph{Comparison to direct solve including generation cost.} Performance profiles of the PCG method with various preconditioners applied to all matrices in the test set, grouped by ordering. This is relative to the work required to solve the system with a direct method and includes an estimate to the amount of work required to generate the preconditioner.
    The graphs in this figure are plotted on a log-linear scale.}
    \label{fig:all_problems_direct_and_work}

\vspace{\baselineskip}
\nextfloat
  \centering
  \begin{subfigure}[t]{0.48\textwidth}
    \includegraphics[width=\textwidth]{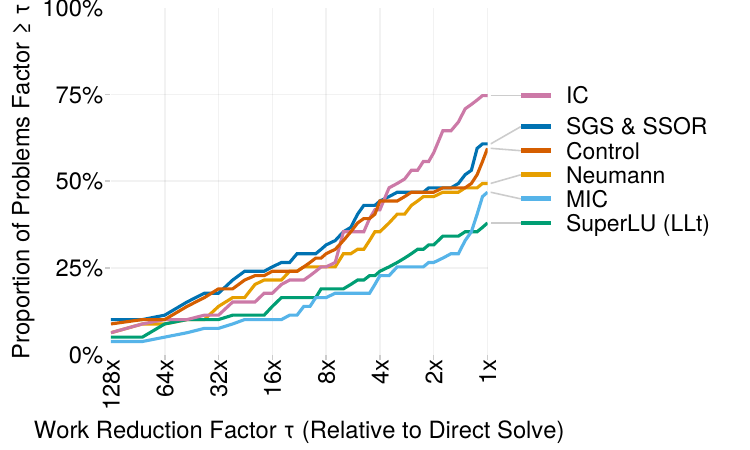}
    \caption{Single‐best configuration; }%
    \label{fig:single_best_config_direct_and_work}
  \end{subfigure}
  \quad
  \begin{subfigure}[t]{0.48\textwidth}
    \includegraphics[width=\textwidth]{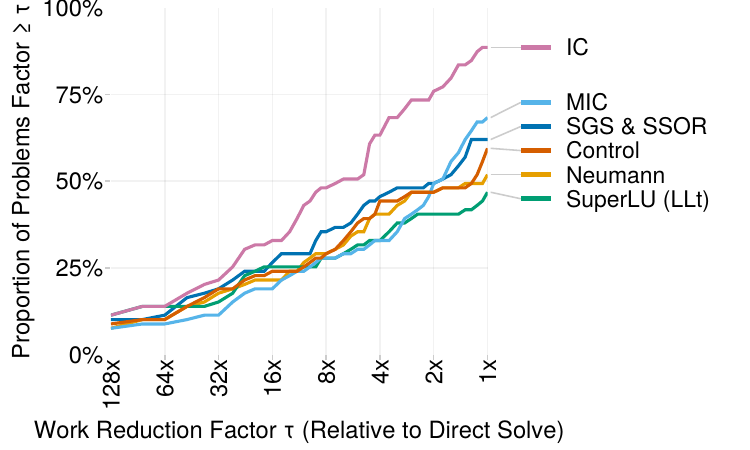}
    \caption{Tuned‐best configuration.}%
    \label{fig:tuned_best_config_direct_and_work}
  \end{subfigure}
  \caption{%
  \emph{Comparison to direct solve including generation cost.}
    The Single-best configuration selects the best configuration over all matrices whereas the tuned-best configuration selects the best run for each matrix within a given preconditioner group. Plotted relative to solving directly. Both include an approximation the the amount of work required to generate the preconditioner. (log‐linear scale)
  }
  \label{fig:all_problems_profile_best_direct_and_work}
\end{figure}

Our final comparison is, perhaps, the most relevant to end-users. This shows our best attempt at comparing the work needed to solve $\mA \vx = \vb$ between a direct method and an iterative method. We include the generation work with the work from the iterations and compare to the reference of an estimate of the operations needed for a sparse direct method. (See \autoref{sec:generation-cost} for information on the generation cost and the sparse direct method cost estimates.) We show the same plots as before: \autoref{fig:all_problems_direct_and_work} shows performance profiles normalized by a direct factorization (including an approximation to each preconditioner's build cost). We summarize these with the single best configuration comparison in \autoref{fig:single_best_config_direct_and_work} and the tuned best configuration comparison in \autoref{fig:tuned_best_config_direct_and_work}. 

These show that preconditioning is successful at reducing total work in many cases, however the difference between preconditioning and no preconditioning shrinks dramatically. In fact, the difference between incomplete Cholesky and the baseline control case of no preconditioner is only around problems that are accelerated by about a factor 2. Of note, the optimal configuration for modified incomplete Cholesky (MIC) changes substantially when we \emph{tune} for the best result for any matrix. 


\section{Key Findings \emph{\&} Discussion }
\label{sec:key-findings}

We use a lengthy passage from \citet{Parlett-1978-progress} to contextualize our discussion.
\begin{quote}\color{darkgray}
\emph{Iterative Methods for $\mA \vx = \vb$}. This was the most fashionable research
topic from 1950-1965. Has the point of diminishing returns been reached? Here is a
quotation from a sophisticated user. Italics are mine. 
\begin{quote}\footnotesize 
    Iterative techniques for processing large sparse linear systems were popular in the late 1950’s and early 1960’s (\emph{and their decaying remains still pollute some computational circles}). When iterative methods finally departed from the finite element scene in the mid 1960’s--having been replaced by direct sparse-matrix methods--the result was a quantum leap in the reliability of linear analysis packages, which contributed significantly to the rapid acceptance of FE analysis at the engineering group level. (This effect, it should be noted, had nothing to do with the relative computational efficiency; in fact iterative methods can run faster on many problems if the user happens to know the optimal acceleration parameters.) Presently, linear FE analyzers are routinely exercised as black box devices;...
\end{quote}
Our own view of the situation is different. By their training, the experts in
iterative methods expect to collaborate with users. Indeed the combination of user,
numerical analyst and iterative method can be incredibly effective. Of course, by the
same token, inept use can make any iterative method not only slow but prone to
failure. Gaussian elimination, in contrast, is a classical black box algorithm demanding
no cooperation from the user. As the tower of scientific computation grew so did the
value of a reliable black box program. In the 1950’s it did not seem possible that
Gaussian elimination could ever be adapted efficiently for large systems of orders
exceeding 200. So the alacrity with which serious users switched to direct methods in
the 1960’s must have been a painful surprise to the adepts of relaxation methods. As
computer systems developed one saw time and again that as soon as direct methods
became feasible for a particular problem they were preferred to their more computer
efficient iterative rivals. Surely the moral of the story is not that iterative methods are dead but that too little attention has been paid to the user’s current needs?
\end{quote}


In this context our results attempt to evaluate the current state of \emph{black box} iterative methods and how well they serve the needs of a user -- almost 50 years after this passage was written. 

While it may be tempting to conclude that IC and MIC are the preferred preconditioners from our study, this conclusion is premature. Both of these preconditioners perform exceptionally well \emph{when they are able to be generated}. Unfortunately, they fail to generate surprisingly often. This is reported on further in~\autoref{tab:best_configurations}. In the best case, we see around 70-80\% success rate with the MATLAB implementation of these preconditioners. This is even worse when we look at the results on the largest problems. 
For the single-best configuration, IC and MIC both failed to generate a preconditioner of 12 of the 16 largest (by number of non-zero values) matrices in the test set.
The large number of failure cases ($\approx 15.19\%)$ of matrices means that software needs more graceful fallbacks to be a true black-box preconditioner. 
Additionally, once preconditioner setup costs are accounted for, significant benefits become less common under a single-best configuration, with IC providing at least factor 4 total-work improvement on roughly $28\%$ of matrices.

Perhaps surprisingly, comparing iterative approaches against factorization-based direct methods, unpreconditioned CG alone already provides at least a twofold computational improvement on nearly half of the test matrices. We caution that these computational estimates for the work involved in generation of the factorizations are only an approximation. Nevertheless, this is encouraging.
Furthermore, IC remains competitive even when including setup costs, surpassing direct methods in total computational effort for approximately $74\%$ of the matrices with a fixed, single-best configuration, increasing to approximately $88\%$ under ideal tuning.
A key difference, however, is that much of the work in direct methods is amenable to parallel computing whereas the iterative methods have a fundamental serial bottleneck. 

In terms of reliability, we echo the note from Parlett's \emph{sophisticated user}. The reliability of modern direct methods is remarkable and the authors of these packages have set a very high bar for iterative methods to clear. 

Returning to our results, the algebraic multigrid preconditioners were reliable, if not the fastest. Unfortunately, we do not have consistent methods to analyze the construction costs for these problems. Anecdotally, they were long -- indicating expensive construction. 
They do seem to be more consistent on larger problems (see \autoref{sec:large}), although they only achieve high work-reduction rates ($64\times$ or greater) on relatively few problems. 



Classical preconditioners, SGS and SSOR, demonstrate only modest improvements even when tuned, with fewer than $10\%$ of problems experiencing significant speedups.  The SSPAI preconditioner showed small, but consistent, benefits on larger problems (see \autoref{sec:large}).


One area where software needs improvement is in terms of matrix ordering. This should not be left up to the user except in unusual circumstances. The reason why is clear. We found that the choice of ordering has a substantial impact on the quality of the preconditioner.
Surprisingly, the fill-reducing ordering AMD consistently outperformed both RCM and natural orderings, even when accounting for the cost of constructing the preconditioner.
While it is known that AMD can be effective when large fill is permitted \citep{benzi_jcp_2002}, prior empirical work has generally favored RCM for preconditioning \citep{saad_siam_2003,duff_bit_1989}.
Our findings challenge this conventional wisdom by demonstrating that, under moderate to high fill levels, AMD often yields more effective preconditioners than RCM -- suggesting a reevaluation of the standard guidance may be warranted.
We explore this behavior further in \autoref{sec:incomplete-cholesky}.


\subsection{Ongoing work}

We are actively expanding our test suite by more than a factor of three by incorporating a large number of graph-based problems that naturally yield symmetric positive definite graph Laplacian matrices.
Many of these graphs are significantly larger than our current examples, allowing us to study preconditioner performance across a wide range of problem sizes -- from thousands to billions of non-zeros. 
In our current results, we noted that the behavior of IC and MIC was more muted on larger problems, especially because the preconditioners often failed to generate on these problems. 
This expanded dataset will help us determine whether certain preconditioners excel only on smaller problems or maintain their advantages as problem size increases.
We additionally seek to add a greater number of diagonally dominant matrices, which will allow us to more thoroughly evaluate preconditioner performance in settings where strong diagonal dominance plays a role (see \autoref{sec:sdd} for our current study of this special cases).

In addition to expanding our test set, we are integrating additional reordering strategies into our framework.
Beyond our current use of AMD and RCM orderings, we plan to incorporate nested dissection and METIS orderings.
These additional ordering schemes will allow us to further explore how structural reordering impacts preconditioner efficiency on large-scale graphs and SPD systems, and to provide clearer guidance on best practices for different problem classes.
Finally, we are also actively testing the preconditioned \textsc{minres} method to compare its performance with that of PCG, thereby extending our analysis across different power subspace solvers.

\subsection{Limitations}
\label{sec:limitations}

Our study has several limitations that we are actively addressing in ongoing work.
First, the current test suite does not yet include a sufficiently diverse set of large-scale problems; most matrices are what we would consider small- to medium-scale for modern scientific computing.
This limitation restricts our ability to fully assess how preconditioner performance interacts with problem size, particularly in the context of very large graphs and systems.

Another limitation arises from our current focus on floating point operations as a proxy for computational cost.
While this measure provides a machine-independent benchmark, it does not capture all aspects of actual runtime performance, such as memory access patterns or the impact of hardware-specific optimizations.
While we developed a means of comparing construction work among factorization based methods and direct methods, we should note that these are likely underestimates. They do not take into account time in computing orderings or other similar types of work. Generalizing this activity across the full range of preconditioners is needed. 

\begin{table}[tb]
    \centering
    \begin{tabularx}{0.75\linewidth}{>{\centering\arraybackslash}X>{\centering\arraybackslash}X>{\centering\arraybackslash}X}
        \toprule
        \textbf{Matrix} &
            \mbox{\textbf{CG Iterations}} \textbf{(Seeded RHS)}
         & 
            \mbox{\textbf{CG Iterations}} \textbf{(Provided RHS)}
         \\
         \midrule 
        2cubes\_sphere & 26 & 28 \\
        af\_0\_k101 & 18416 & 46588 \\
        af\_1\_k101 & 18401 & 45186 \\
        af\_2\_k101 & 17574 & 40810 \\
        af\_3\_k101 & 16342 & 29932 \\
        af\_4\_k101 & 17930 & 36687 \\
        af\_5\_k101 & 17783 & 34254 \\
        af\_shell3 & 2464 & 139 \\
        af\_shell7 & 2463 & 757 \\
        BenElechi1 & 32283 & 33932 \\
        bone010 & 8124 & 10001 \\
        boneS01 & 2137 & 2401 \\
        boneS10 & 23163 & 37350 \\
        nasasrb & 9633 & 20888 \\
        offshore & 300 & 1226 \\
        parabolic\_fem & 2376 & 2568 \\
        Pres\_Poisson & 729 & 812 \\
        smt & 3371 & 3615 \\
        thermal1 & 1283 & 1511 \\
        thermal2 & 4506 & 5392 \\
        \bottomrule
    \end{tabularx}
    \caption{Comparison of the number of unpreconditioned CG iterations to converge with our right-hand side vector compared to the SuiteSparse Matrix Collection right-hand side vector. In this table, ``Seeded'' refers to the right-hand side vector generated by our method, and ``Provided'' refers to the right-hand side vector provided by the SuiteSparse Matrix Collection. In the case where multiple right-hand side vectors were provided, we used the first one.}
    \label{tab:iterations_comparison}
\end{table}

Finally, in this study we focus on the relative performance of preconditioners across a diverse set of matrices.
In performing so many runs, we are limited in the number of right-hand side vectors we can generate and test on.
We recognize that the choice of the right-hand side vector can cause variations in the reported performance of the preconditioners, or even the convergence of the solver.
This is a fundamental limitation to our study. We propose a reproducible method for generating right hand sides in \autoref{subsec:right_hand_side}. We have found that using our method of generating the right-hand side vector provides similar results to those provided by the SuiteSparse Matrix Collection. 
In \autoref{tab:iterations_comparison}, we show a table which shows the number of unpreconditioned CG iterations to converge with our right-hand side vector compared to the SuiteSparse Matrix Collection right-hand side vector (when these are available).

This shows that there is some variance in the number of iterations required to converge, but that the variance is always within a reasonable range.
Note that there are relatively few matrices in the test set which come with provided right-hand side vectors, but this provides evidence that our method of generating the right-hand side vector is reasonable.

\subsection{Remainder of the paper}
In the remainder of the paper, we provide additional details for those who wish to understand the details of our work. The next sections review the technical background that we expect many readers to be familiar with, but reiterate it to explain the precise details of our experiments. We explain the details of each preconditioner in \autoref{sec:preconditioners}. Finally, we present a few additional snapshots of our results in \autoref{sec:extra-results}.

\subsection{Data and reproducibility} 
We host the results of our study online in two places. The results are available from 
\begin{center}
    \url{https://tunnellm.github.io/CAPPA-repository}
\end{center}
and the code to reproduce it is availalbe from 
\begin{center}
    \url{https://github.com/tunnellm/CAPPABenchmarkTools}
\end{center}

\section{Background}\label{sec:background}

We first discuss the notation and terminology used in this paper.
We utilize boldface lowercase letters to denote vectors (e.g., $\vx$), boldface uppercase letters to denote matrices (e.g., $\mA$), and lowercase letters to denote scalars (e.g., $\gamma$).
The $i^\text{th}$ element of a vector $\vx$ is denoted by $x_i$, and the $(i,j)^\text{th}$ element of a matrix $\mA$ is denoted by $a_{i,j}$.
We  $\vx^{(k)}$ to denote the $k^\text{th}$ iterate of some vector $\vx$, whereas $\transpose{\mA}$,  $\mA^{-1}$, and $\mA^{k}$ refers to the transpose, the inverse, and the $k^\text{th}$ power of some matrix $\mA$, respectively.

We refer to some arbitrary preconditioner as $\mM$, where its application to a vector $\vr$ is denoted as $\vz = \mM^{-1}\vr$.
We make frequent reference to linear operators in this paper and often denote some arbitrary linear operator as $\gL$.
The application of a linear operator $\gL$ is equivalent to applying the inverse of $\mM$ and is denoted as $\vz = \gL(\vr) = \mM^{-1}\vr$.


\subsection{Krylov subspace methods}\label{sec:solvers}

In this preliminary study, we focus on Krylov subspace methods for solving large, sparse linear systems of equations.
Specifically, we utilize the Preconditioned Conjugate Gradient (PCG) method for SPD systems.
We recognize that the performance of preconditioners may be different when applied with different solvers and continue to study other methods including \textsc{minres}~\cite{paige_siam_1975}.

\subsection{Preconditioned conjugate gradient}\label{subsec:pcg}

\paragraph{Implementation details}

We implemented the formulation of PCG due to~\citet{hestenes_jrnbs_1952}, which can be found in Algorithm~11.5.1 of~\citet{golub_jhup_2013_pcg}, in Julia~\cite{bezanson_siam_2017}.
We use one matrix-vector product per iteration to compute the update and one preconditioner application per iteration.
Our implementation stores information on most variables after each step of computation, allowing us to track the performance of the solver with high granularity. We also use one additional matrix-vector product to recompute the exact residual at each iteration but this is not included in the computational cost of the method. 
We abstract the application of the preconditioner $\mM$ as the linear operator $\gL(\vr) = \mM^{-1}\vr$, allowing us to wrap all of the external preconditioners in a common interface that conforms to the implementation.
Our implementation checks for convergence based on both the relative residual and the normwise relative backward error at each step. For normwise relative backward error, we computed the operator induced matrix $2$-norm (which is just the largest eigenvalue of $\mA$ for an SPD matrix) via the {GenericArpack} package~\cite{genericarpack}, a pure Julia translation of Symmetric ARPACK~\cite{lehoucq_siam_1998}. 

\paragraph{Computational cost}

The dominant computational costs per iteration are the matrix-vector product $\mA\vp^{(k)}$ and the preconditioner application $\vz^{(k)} = \gL(\vr^{(k)})$.
The cost of the matrix-vector product is proportional to the number of non-zero values in the matrix, and the cost of the preconditioner application is dependent on the type of preconditioner -- ranging from the number of non-zero values in the preconditioner to the complexity of the cycle/application in the case of algebraic multigrid.
We provide more detail on the computational cost of each preconditioner in \autoref{sec:preconditioners}.

\section{Systems of Linear Equations}\label{sec:systems}

In this section, we describe the selection of symmetric positive definite matrices used in our study.
Our goal is to evaluate the performance of various preconditioners on large-scale SPD systems that are representative of real-world problems.
We are actively adding to to the selection of matrices.

\subsection{Symmetric positive definite matrices}\label{subsec:spd_matrices}


We utilize matrices from the SuiteSparse Matrix Collection~\cite{davis_acm_2011}, a widely recognized repository of sparse matrix benchmarks that are representative of real-world problems.
To ensure our study focuses on large-scale and practically relevant problems, we filter matrices based on the following criteria:
\begin{itemize}[noitemsep,topsep=0pt]
    \item \textbf{Matrix Dimension}: Only matrices of dimension $n \geq 10{,}000$ are considered. This ensures that the matrices are sufficiently large to reflect computational challenges encountered in real-world applications.
    \item \textbf{Symmetric Positive Definite}: We consider only those matrices that are real-valued and symmetric positive definite in this study.
    \item \textbf{Convergence without Preconditioning}: The unpreconditioned Conjugate Gradient method should converge within $10 \cdot n$ iterations, where $n$ is the matrix dimension. This allows us to compare the relative performance of preconditioners in aggregate against the unpreconditioned case.
\end{itemize}

By limiting our study to a test set of matrices that fit these criteria, we ensure that we are evaluating preconditioners on those matrices that present a computational challenge that is typical in large-scale scientific and engineering problems.
Additionally, requiring that the unpreconditioned CG method converges within a reasonable number iterations enables us to compare the relative performance of the preconditioners both in aggregate and on a per-matrix basis.

The selected matrices come from a diverse set of applications, including structural analysis, fluid dynamics, and optimization problems.
These matrices comprise a wide range of sizes and structures, reflecting the diversity of problems encountered in practice.
This diversity allows us to evaluate the robustness and effectiveness of preconditioners across a broad spectrum of real-world problems.
The matrices used in this study are categorized in \autoref{tab:matrix_categories}.

\begin{table}[tbp]
    \centering
    \begin{tabularx}{\linewidth}{>{\centering\arraybackslash}X>{\centering\arraybackslash}X>{\centering\arraybackslash}c}
        \toprule
        \textbf{Matrix Category} & \textbf{Citations} & \textbf{Count} \\
        \midrule
        Circuit Simulation & \citet{amd_group,mcrae_group} & 3  \\
        \midrule 
        Combinatorial & \citet{trefethen_group} & 2  \\
        \midrule 
        Computational Fluid Dynamics & \citet{acusim_group,janna_siam_2011,um_group,wisgott_group,maxplanck_group} & 6 \\
        \midrule 
        Computer Vision & \citet{lourakis_group,mazaheri_group} & 2 \\
        \midrule 
        Electromagnetics & \citet{um_group,cemw_group} & 2 \\
        \midrule 
        Materials & \citet{boeing_group} & 2 \\
        \midrule 
        Model Reduction & \citet{rudnyi_proceeding_2005} & 6 \\
        \midrule 
        Optimization & \citet{ghs_psdef_group,berger_sorrr_1995} & 5 \\
        \midrule 
        Miscellaneous PDE & \citet{benelechi_group,janna_sjsc_2015,ferronato_ijnme_2011,utep_group,williams_group,cunningham_group,ghs_psdef_group} & 12 \\
        \midrule 
        Structural & \citet{janna_sjsc_2010,janna_aes_2009,ferronato_ijnamg_2007,ferronato_cmame_2008,janna_ijnme_2009,ferronato_ecm_2010,schenk_afe_group,ghs_psdef_group,duff_acm_1989,pothen_group,tkk_group,ghs_psdef_group,inpro_group,dnvs_group} & 34 \\
        \midrule 
        Thermal & \citet{bindel_group,dabrowski_geochemistry_2008,botonakis_group} & 5 \\
        \bottomrule
    \end{tabularx}
    \caption{Matrix categories}
    \label{tab:matrix_categories}
\end{table}

\phantomsection
\subsection{Matrix scaling}
\label{sec:scaling}

To isolate potentially poor scaling of certain matrices, we scale with symmetric diagonal scaling (symmetric Jacobi) to improve the conditioning of the system~\cite{golub_jhup_2013_pcg,saad_siam_2003}.
Explicitly, we scale $\mA$ on the left and right by $\mD^{-\frac{1}{2}}$, the right-hand side vector by $\mD^{-\frac{1}{2}}$, and the solution vector by $\mD^{\frac{1}{2}}$.
We symmetrize the matrix by averaging every off diagonal entry and explicitly set the diagonal entries of the matrix equal to the identity. Formally, if $\mA'$ is the result of scaling and scaling 
\[ \mA' \gets \mD^{-\frac{1}{2}} \mP\mA\transpose{\mP} \mD^{-\frac{1}{2}} \] 
then we use 
\[ 
    \mA'' \gets \frac{1}{2}(\mA' + \transpose{(\mA')})
\]    
and then explicitly set the diagonal of $\mA''$ to the identity,
\[ a_{i,i}'' \gets 1 \quad \text{ for all $i$}. \]

\subsection{Matrix ordering}\label{subsec:matrix_ordering}

The ordering of a sparse matrix is well known to be crucial for the performance of a large number of preconditioners~\cite{benzi_jcp_2002,saad_siam_2003}.
In the symmetric case, the rows and columns of the matrix must be permuted in the same way to preserve the symmetry of the matrix, so we only consider symmetric reorderings in this study.
In the context of direct solvers such as Gaussian elimination or those preconditioners that compute an incomplete factorization, the ordering of the matrix is known to have a significant impact on the fill-in and the quality of the factorization~\cite{benzi_jcp_2002,davis_siam_2006}.


\paragraph{Natural}

The natural ordering of a matrix is the order in which it is obtained from the application that generated it.
This serves as a baseline for comparison with other ordering strategies.

\paragraph{Reverse Cuthill-McKee}

The Reverse Cuthill-McKee (RCM) algorithm was originally proposed for reducing the bandwidth of a matrix~\cite{cuthill_acm_1969}.
By minimizing the bandwidth of the matric, RCM aims to cluster the non-zero elements closer to the diagonal and has been shown in practice to reduce the fill-in of the matrix during factorization~\cite[Chapter 7]{davis_siam_2006}.
We utilize the \texttt{SymRCM} package \citep{symrcm} in Julia to compute the RCM ordering of the matrices.


\paragraph{Symmetric approximate minimum degree}

The Symmetric Approximate Minimum Degree (SYMAMD) ordering is an approximate version of the minimum degree algorithm that is designed to preserve the symmetry of the matrix~\cite{amestoy_siam_1996} and minimize the fill-in of a Cholesky factorization of a sparse matrix~\cite[Chapter 7]{davis_siam_2006}.
We utilize the \texttt{AMD} package in Julia~\cite{dominique_zenodo_2023}.



\paragraph{Slight differences with the baseline case}
It is known that Krylov methods (including CG) are unaffected by matrix ordering in exact arithmetic, so the control case should remain unchanged.
However, we found extremely slight differences when evaluated. 
For this reason, we run the control for each ordering to establish the baseline of comparison against each preconditioner.

\subsection{Right-hand side vector}\label{subsec:right_hand_side}

We generate a random solution vector $\vx^{(*)} \in \sR^{n}$ using the \texttt{StableRNGs} package in Julia, which reproduces a random number generator (RNG) suitable for scientific computing.
Specifically, we use the seed $123456789$ to initialize the random number generator, allowing us to reproduce the same sequence of random number across different runs and systems.


Out of an abundance of caution for the statistical properties of the generated random numbers, we ``warm up'' the RNG by generating and discarding $\log_2(n)$ random vectors of length $n$ (rounded to the nearest integer), where $n$ is the dimension of the matrix $\mA$.
This practice reduces the likelihood of any biases or correlations in the sequence that may be present in the initial random numbers generated by the RNG~\cite{gentle_springer_2003}.


Using the last generated random vector as the solution vector $\vx^{(*)}$, we construct the right-hand side vector $\vb$ in the following manner:
We locate the largest in magnitude $\log_2(n) + 1$ elements of $\vx^{(*)}$ (rounded to the nearest integer).
For each of these elements, we set its value to either $1$ or $-1$, chosen at random with equal probability.
We set the remaining elements of $\vx^{(*)}$ to zero, resulting in a sparse vector.


Finally, we compute $\vb = \mA \vx^{(*)}$, which ensures that the right-hand side vector is consistent with the solution vector and the matrix $\mA$.
Performing this procedure allows us to evaluate the performance of preconditioners in terms of its error $\normof{\vx^{(k)} - \vx^{(*)}}$ during iterations of the solver.

As noted in~\autoref{tab:iterations_comparison}, we saw no characteristic differences in terms of iterations to convergence between our choice of right hand side and those right hand side vectors that were provided in SuiteSparse.


\section{Preconditioners}\label{sec:preconditioners}

In this section, we describe the preconditioners that are evaluated in this preliminary study.
In the following subsections, we provide an overview of each preconditioner, discussing the software packages used and the expected computational costs associated with their use.



\subsection{Caveats with SuperLU and graph Laplacian preconditioners}
\label{sec:caveats}
Both the SuperLU and graph Laplacian preconditioners have important caveats. First, SuperLU is not designed to produce the \emph{symmetric} preconditioners we need for CG. We are able to adapt its strategy to evaluate the method as described below. However, this means that our evaluation represents an evaluation of our adaptation instead of the direct software. Second, graph Laplacian preconditioners are designed for symmetric diagonally dominant systems. We are evaluating them on the superclass of symmetric positive definite systems. Both of these studies are done to provide a reference and should not be interpreted as a critique on the design ideas of the underlying preconditioners or software. Both meet the high bar of being able to complete our evaluation.

\subsection{Truncated Neumann series}\label{sec:truncated_neumann_series}

The Truncated Neumann Series (TNS) preconditioner is based on approximating the inverse of a matrix using a finite number of terms in the Neumann series expansion~\cite{dubois_springer_1979}.
Given a unit diagonally scaled $\mA$, the Truncated Neumann series can be expressed as
\begin{align*}
    \mM^{-1} &= \alpha \sum_{k=0}^m \left(\mI -\alpha \mA \right)^k,
\end{align*}
where $\alpha$ is chosen such that the spectral radius of $\alpha\mA$ is less than 1.


\phantomsection
\paragraph{Theoretical considerations}

Consider the linear system $\mA \vx = \vb$ where $\mA$ has unit diagonal and $\mR$ is the remainder matrix.
Without loss of generalization for symmetric positive definite matrices, assume $\mA$ has eigenvalues strictly inside the region $(0,1)$, then the TNS for this system is given by
\begin{align*}
    \mM^{-1} &= \sum_{k=0}^{m} \left( -\mR  \right)^k,
\intertext{where $\mR = \mI - \mA$ is the remainder matrix. Suppose we let $m = 1$, we have}
    \mM^{-1} &= \mI - \mR.
\intertext{Because $\mR$ is the remainder matrix, we have}
    \mM^{-1} &= \mA.
\end{align*}

The preconditioned Krylov subspace is given by
\begin{align*}
    \mathcal{K}^{(k+1)} &= \text{span}\{\vz^{(0)}, \mA\vz^{(0)}, \mA^2\vz^{(0)}, \dots, \mA^k\vz^{(0)}\},
\intertext{where}
    \vz^{(0)} &= \mM^{-1} \vr^{(0)}.
\intertext{Thus the TNS preconditioner with $m = 1$ gives}
    \vz^{(0)} &= \mA \vr^{(0)},
\intertext{and the Krylov subspace is}
    \mathcal{K}^{(k+1)} &= \text{span}\{\mA\vr^{(0)}, \mA^2\vr^{(0)}, \mA^3\vr^{(0)}, \dots, \mA^{k+1}\vr^{(0)}\}.
\end{align*}

In this scenario, the preconditioned residual $\vz^{(0)}$ is $\mA$ times the original residual $\vr^{(0)}$.
The Krylov subspace is effectively built from higher powers of $\mA$ applied to $\vr^{(0)}$, starting from $\mA \vr^{(0)}$ instead of $\vr^{(0)}$.
As a function of the work performed, the TNS preconditioner builds the Krylov subspace at approximately the same cost as the unpreconditioned system.
Therefore for matrices that are scaled to have unit diagonal, there should be no significant improvement in the convergence of the system with the TNS preconditioner as a function of the work performed.

\paragraph{Computational cost}

The computational cost of the TNS preconditioner is equal to the number of non-zero values in the matrix $\mA$ multiplied by the number of terms in the Neumann series expansion.

\paragraph{Preconditioner settings used}

The only setting to tune for the TNS preconditioner are the number of iterations to apply and the scaling factor.
We run this preconditioner for all combinations of iterations in $\{1,2,3,4\}$ and
\[\frac{1}{\alpha} \in \{\onormof{\mA}{F},\onormof{\mA}{\infty}, \onormof{\mA}{1}, \frac{\onormof{\mA}{2}}{2}, 1\}.\]

\subsection{Symmetric Gauss-Seidel and successive over-relaxation}\label{sec:sgs-ssor}

Gauss-Seidel (GS) and Successive Over-Relaxation (SOR) are classical iterative methods that update the solution vector component-wise~\cite{golub_jhup_2013_pcg,saad_siam_2003}.
These methods improve over the Jacobi method by using the most recent values available to update each component of the solution vector.
The SOR method accelerates the convergence of GS by introducing a relaxation parameter $\omega$ that over-relaxes the update in the direction of the residual, which can improve the convergence rate of the solver if $\omega$ is chosen appropriately.
Consequently, the GS method can be viewed as a special case of SOR with $\omega = 1$.

Given a linear system $\mA \vx = \vb$, where $\mA \in \mathbb{R}^{n \times n}$ is a nonsingular matrix, $\vx \in \mathbb{R}^n$ is the solution vector, and $\vb \in \mathbb{R}^n$ is the right-hand side vector, the SOR method partitions $\mA$ into its diagonal, strictly lower triangular, and strictly upper triangular components:
\begin{equation*}
\mA = \mD + \mL + \mU,
\end{equation*}
where $\mD$ are the diagonal elements of $\mA$, $\mL$ contains the entries strictly below the diagonal, and $\mU$ contains the entries strictly above the diagonal~\cite{golub_jhup_2013_pcg}.
These methods are not suitable for use in solvers that require symmetric matrices, such as PCG and \textsc{minres}, as the preconditioner is not symmetric.


\paragraph*{Symmetric GS and SOR}
Symmetric Successive Over-Relaxation (SSOR) and Symmetric Gauss-Seidel (SGS) are variants of the SOR and GS methods that update the solution vector in a symmetric fashion, alternating between forward and backward sweeps~\cite{golub_jhup_2013_pcg,saad_siam_2003,young_tams_1954,young_ams_1950}.
Performing the update symmetrically smooths the error in the solution vector and can improve the convergence rate of the solver.
Just as before, the SGS method can be viewed as a special case of SSOR with $\omega = 1$.

Given the same linear system $\mA \vx = \vb$ and splitting defined for the SOR method, the SSOR method performs a forward sweep followed by a backward sweep to update the solution vector.
Formally, the forward sweep is given by:
\begin{align*}
    \vx^{(k+\frac{1}{2})} &= (\mD + \omega \mL)^{-1} \left[ \omega \vb - (\omega \mU + (1 - \omega) \mD) \vx^{(k)} \right],
\intertext{and the backward sweep is given by}
    \vx^{(k+1)} &= (\mD + \omega \mU)^{-1} \left[ \omega \vb - (\omega \mL + (1 - \omega) \mD) \vx^{(k+\frac{1}{2})} \right].
\end{align*}

\paragraph{Computational cost}

The computational cost per sweep of GS, SGS, SOR and SSOR is equal to the number of non-zero values in the matrix $\mA$.
There are an additional $2n$ operations per sweep to compute the proportional update to the solution vector $\vx$ for SOR and SSOR.
The cost of the symmetric variants are exactly twice the cost of the non-symmetric variants, as they perform both a forward and backward sweep.
These costs are given for a single iteration of the method, and the total cost of the preconditioner is proportional to the number of iterations applied.

\paragraph{Preconditioner settings used}

We run these preconditioners for all combinations of iterations in $\{1,2\}$ and relaxation parameters in $\{1.0, 1.2, 1.5, 1.8\}$.
When running with a relaxation parameter of $1.0$, we use a standard SGS implementation that does not perform the additional computations to compute the relaxation.

Additionally, it was shown in~\citet{young_tams_1954} that, for certain matrices, there is a known optimal value of $\omega$.
Those matrices that possess ``Property-A'' due to~\citet{young_tams_1954} have an optimal $\omega$ parameter described by
\begin{align*}
    \omega = 1 + \wpar{\frac{\normof{\mJ}}{1 + \sqrt{\wpar{1 - \normof{\mJ}^2}}}}^2,
\end{align*}
where $\normof{\mJ}$ is the spectral norm of the Jacobi iteration matrix of $\mA$.
This is well-defined for matrices with $\normof{\mJ} \leq 1$, but may or may not correspond with the optimal value of $\omega$ for matrices not possessing Property-A.
We compute this value for all matrices in which it is well defined and run for $1$ and $2$ iterations.

\subsection{Symmetric sparse approximate inverse}\label{sec:ssai}

Sparse Approximate Inverse (SAI) preconditioners construct an explicit sparse approximation to $\mA^{-1}$ for efficient application during iterative methods~\citep{benzi_siam_1998,grote_siam_1997,regev_stanford_2020}.
The symmetrized version, such as SSAI described in~\citep{regev_stanford_2020,regev_stanford_dissertation_2022}, extends this idea by adding an explicit symmetrization step after computation.
The implementation we use is available at \url{https://stanford.edu/group/SOL/software/minres/minres20.zip}

\phantomsection
\paragraph{Computational cost}

The computational cost of the SSAI preconditioner is equal to the number of non-zero values in the matrix $\mM^{-1}$.

\paragraph{Preconditioner settings used}
The authors suggest using the average number of elements in a column as the level of fill allowed in the preconditioner.
We run the preconditioner with fill levels set to the average number of elements per column multiplied by the values in $\{0.5, 1.0, 2.0, 3.0\}$.


\subsection{Simplicial Incomplete Cholesky and modified incomplete Cholesky}
\label{sec:incomplete-cholesky}
\label{sec:modified-incomplete-cholesky}

Incomplete Cholesky (IC) approximates the Cholesky factorization for SPD matrices by computing a sparse triangular matrix that retains only selected nonzero entries along a sparsity pattern, reducing computational cost while ideally capturing essential features of $\mA$~\citep{meijerink_moc_1977,saad_siam_2003}.
MATLAB offers a widely used implementation with an IC(0) variant, allowing no fill-in, or a threshold-based variant, allowing fill-in only above a specified magnitude.
The Modified IC (MIC) variant extends this approach by adjusting the diagonal entries to compensate for dropped fill-in, which has been show to enhance numerical stability for some poorly conditioned problems~\citep{robert_laa_1982,gustafsson_springer_1978}.
These variants of incomplete factorizations are simplicial, and the systems are solved column-by-column (or row-by-row) and stored in compressed sparse format.
We utilize the MATLAB implementations of these two preconditioners as they are widely used in practice.
Both preconditioners are exposed through various options passed to the \texttt{ichol} function.

\phantomsection
\paragraph{Preconditioner settings used}
We run these preconditioners with the ILU($0$) variant as well as thresholding set to all values in $\{10^{-4}, 10^{-5}, 10^{-6}, 10^{-7}, 10^{-8}\}$.

\phantomsection
\paragraph{Computational cost}

The computational cost of applying the incomplete factorization is equal to the number of non-zero values in each factor that is inverted.
In the IC and MIC cases, this is two times the number of non-zero values in the lower triangular factor.
In the ILU case, this is equal to the number of non-zero values in the lower and upper triangular factors.

\subsection{Supernodal incomplete LU factorization}\label{sec:superlu}

An ILU approximates the LU factorization similarly to IC but is not limited to symmetric positive definite matrices.
For symmetric positive definite matrices, it is generally equivalent to IC.
We use the implementation in the SuperLU library~\citep{li_siam_2003,li_ACM_2010,demmel_siam_1999,li_lbln_1999,li_toms_2005}, which uses a supernodal approach to group columns for improved cache efficiency and reduce operation counts, while offering control over fill-in levels and drop tolerances.
In our study, we use SuperLU's ILU factorization to generate a preconditioner; however the result of this library is not necessarily symmetric even with symmetric input, so we extract the lower triangular portion and construct a symmetric preconditioner. In some preliminary experiments, this was superior to using the upper triangular portion. To be explicit, we scale the lower triangular factor with the square root of the diagonal entries of $\mU$, forming $\mL' = \mL \mD^{\frac{1}{2}}$ and setting $\mM = \mL'\transpose{\mL'}$.

\phantomsection
\paragraph{Preconditioner settings used}
We run this preconditioner with all combinations of a drop tolerance in $\{10^{-4}, 10^{-5}, 10^{-6}\}$ and fill-in levels in $\{1, 2, 3\}$.

\paragraph{Computational cost}
This is the same as the incomplete factorizations discussed in \autoref{sec:incomplete-cholesky}.

\subsection{Algebraic multigrid}\label{subsec:amg}

Algebraic Multigrid (AMG) methods are an advanced class of iterative methods that are designed to solve large and sparse linear systems efficiently.
These methods work by constructing a hierarchy of coarser and coarser grids that approximate the original problem, allowing for the solution to be computed at different levels of resolution~\cite{brandt_elsevier_1986}.
Unlike geometric multigrid methods, which we do not consider in this study, AMG constructs the hierarchy based on the matrix itself, making it more general and applicable to a wider range of problems~\cite{brandt_elsevier_1986,saad_siam_2003,stuben_elsevier_2001,brandt_siam_2011}.

AMG is known to be particularly effective for elliptic PDEs, where the matrix is symmetric and positive definite, and the problem is discretized on a structured or unstructured grid.
For these types of problems, they are known to achieve optimal or nearly-optimal convergence rates~\cite{brandt_elsevier_1986,brandt_siam_2011}.


In this study, we focus on the Ruge-Stuben~\cite{ruge_siam_1987} AMG method and the Smoothed Aggregation~\cite{vanek_nm_2001} AMG method, which are two popular AMG methods that are widely used in practice~\cite{bell_joss_2023}.
We utilize the \texttt{PyAMG} library~\cite{bell_joss_2023} to construct the AMG preconditioners.

\phantomsection
\paragraph{Computational cost}

The PyAmg library provides a function to compute the complexity of a cycle of the AMG method used.
We use this as an approximation to the computational cost of the AMG preconditioner.

\phantomsection
\paragraph{Preconditioner settings used}

The Ruge-Stuben method is provided with \texttt{V}, \texttt{W}, \texttt{F}, and \texttt{AMLI} cycles.
The Smoothed Aggregation method is provided with \texttt{V}, \texttt{W}, and \texttt{F} cycles.
We test both methods and all cycles in our study, running for one and two cycles.
We otherwise use the default settings for the AMG preconditioner.

\subsection{Graph Laplacian preconditioners}\label{subsec:graph_laplacian}

This subsection explores preconditioners designed specifically for graph Laplacian matrices.
Although we understand that these preconditioners are not designed for general SPD matrices, we still run them on general SPD matrices and report their results.

For general SPD matrices that are not SDD, we do not expect these preconditioners to perform well.
For matrices that are SDD when unscaled but not SDD when scaled, we run the preconditioner on the unscaled matrix.
For matrices that are SDD when scaled, we run the preconditioner on the scaled matrix.
For the remaining case in which the matrices are not SDD when scaled or unscaled, we create the preconditioner on a matrix that was modified to be SDD by increasing the diagonal entries of the matrix to be at least equal to the row sum of the absolute values of the off-diagonal entries.

For an SDD matrix of dimension $n$, we can construct a graph Laplacian matrix of dimension $2n$ through an augmentation described in~\cite{kelner_arxiv_2013}, which is an improvement over Gremban's reduction that would produce a matrix of dimension $2n+1$~\cite{gremban_cmu_1996}.
For an SDD matrix $\mA$, this reduction is given by the following equation:
\begin{align*}
    \mL &= \begin{bmatrix}
        \mD + \mN + \frac{1}{2} \mS & -\wpar{\mP + \frac{1}{2} \mS}\\
        -\wpar{\mP + \frac{1}{2} \mS} & \mD + \mN + \frac{1}{2} \mS
        \end{bmatrix},
\end{align*}
where $\mP$ are the positive off-diagonal entries of $\mA$, $\mN$ are the negative off-diagonal entries of $\mA$, and where
\begin{align*}
    \mD &= \text{diag}\wpar{\mP \ve - \mN \ve},\\
\intertext{and}
\mS &= s_{i,i} = a_{i,i} - d_{i,i}.
\end{align*}
The diagonal matrix $\mS$ represents the excess diagonal entries of $\mA$ over the row sum of the absolute values of the off-diagonal entries.
We then distribute half of this slack to the off-diagonal entries of the matrix and half to the diagonal entries.

Given a RHS vector $\vb$ to the original system, the augmented RHS vector is given by:
\begin{equation*}
    \vb' = \begin{bmatrix}
        \vb\\
        -\vb
    \end{bmatrix}.
\end{equation*}
We can then solve the augmented system $\mL \vx' = \vb'$ to obtain the solution to the original system.
The solution to the original system is then recovered by
\begin{equation*}
    2\vx = \wpar{\vx'}_{1:n} - \wpar{\vx'}_{n+1:2n}.
\end{equation*}
We utilize this reduction when applying the graph Laplacian preconditioners to the general SPD matrices.
When solving for the solution to graphs, we use the original graph Laplacian matrix.
Next, we briefly introduce the graph Laplacian preconditioners.








We utilize the \texttt{Laplacians} package in Julia.\footnote{We use the version with SHA hash \texttt{d89118f61fe4346d2401e7595c197ef3ef73d32c}}
We use the \texttt{approxchol\_lap} function to construct the preconditioner, which is a randomized variant that constructs an LD$\transpose{\text{L}}$ decomposition of a subtree of the graph~\cite{kyng_arxiv_2023}.

As of the writing of this study, the \texttt{Laplacians} package does not expose the preconditioner for use in outside solvers.
We worked around this by extracting the relevant portions of the package and creating a function that applies the preconditioner to a matrix.

\phantomsection
\paragraph{Computational cost}

We compute the work expended by the preconditioner as the number of multiplications used in solving the systems in the  LD$\transpose{\text{L}}$ factorization.
Note that for matrices that have multiple connected components in the augmented graph, multiple LD$\transpose{\text{L}}$ factorizations are computed, one for each connected component, and the work expended is the sum of the work expended for each connected component.

\phantomsection
\paragraph{Preconditioner settings used}

There are two parameters that can be tuned: \texttt{split} and \texttt{merge}, which controls the number of multi-edges used.
We run these parameters for all combinations of split and merge values in $\{0,1,2,3,4,5\}$.

\section{Additional Results}
\label{sec:extra-results}

Here, we share additional results that offer insight into preconditioner behavior.
In the final section, we provide a table that highlights notable configurations for each method and include observations that are specific to certain configurations or test cases.

\subsection{Large problems}\label{sec:large}

\begin{figure}[tbp]
  \centering
  \begin{subfigure}[t]{0.48\textwidth}
    \includegraphics[width=\textwidth]{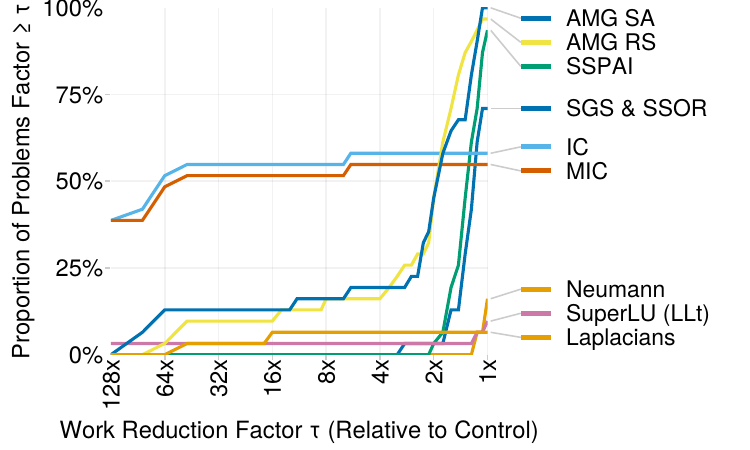}
    \caption{Single‐best configuration; }%
    \label{fig:single_best_config_large_size}
  \end{subfigure}
  \quad
  \begin{subfigure}[t]{0.48\textwidth}
    \includegraphics[width=\textwidth]{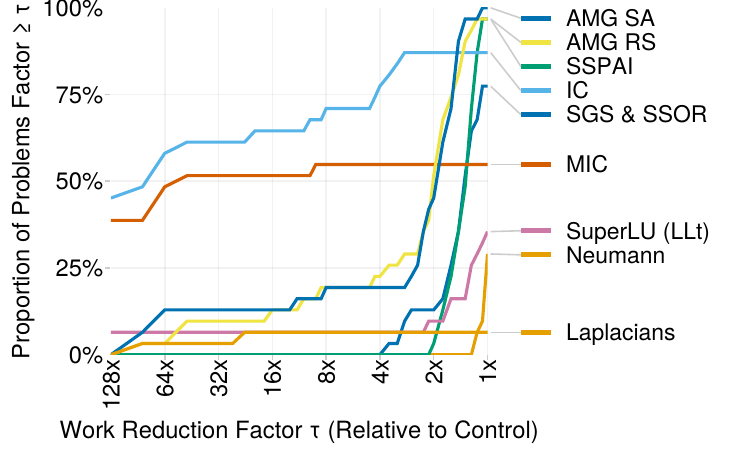}
    \caption{Tuned‐best configuration.}%
    \label{fig:tuned_best_config_large_size}
  \end{subfigure}
  \caption{%
    The single-best configuration selects the best configuration over matrices of dimension $250$k or higher, whereas the tuned-best configuration selects the best run for each matrix within a given preconditioner group. Plotted relative to the unpreconditioned baseline and excluding work involved in generating the preconditioners. 
  }
  \label{fig:large_problems_best_configuration}
\end{figure}

For problems up to around  $100,000$ unknowns, it's reasonable to expect the memory involved in a direct factorization to be available. This changes considerably as the problems become larger. In this study, we wanted to look at what happens on systems with over $250,000$ unknowns. There are $31$ such systems in our set.  \autoref{fig:large_problems_best_configuration} shows the single-best and tuned-best performance profiles over matrices of dimension $250$k and larger, excluding generation work. 
As shown in the figure, we find that classical methods perform similarly to their performance on the entire dataset, the single-best incomplete factorization has a larger proportion of failure cases, and the two AMG methods and SSPAI start to see higher proportions with moderate work reduction gains.

\begin{figure}[tbp]
  \centering
  \begin{subfigure}[t]{0.48\textwidth}
    \includegraphics[width=\textwidth]{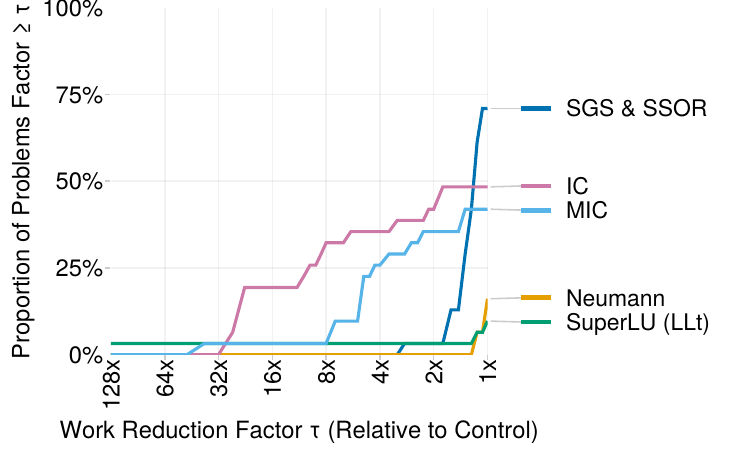}
    \caption{Single‐best configuration; }%
    \label{fig:single_best_config_and_work_large_size}
  \end{subfigure}
  \quad
  \begin{subfigure}[t]{0.48\textwidth}
    \includegraphics[width=\textwidth]{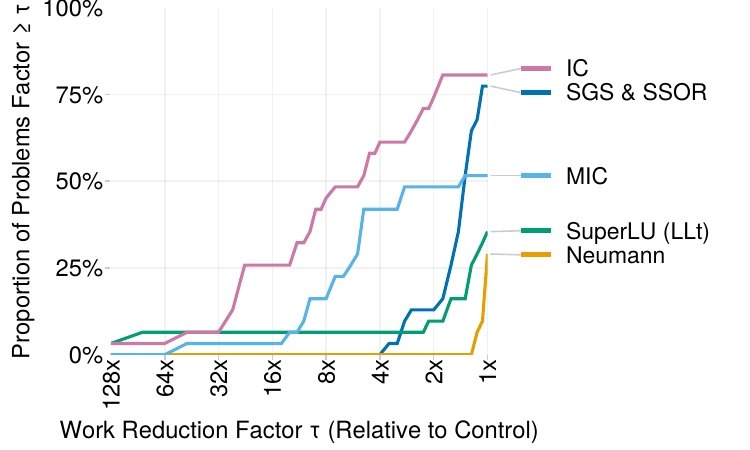}
    \caption{Tuned‐best configuration.}%
    \label{fig:tuned_best_config_and_work_large_size}
  \end{subfigure}
  \caption{%
    \emph{Including generation work.} The single-best configuration selects the best configuration over matrices of dimension $250$k or higher, whereas the tuned-best configuration selects the best run for each matrix within a given preconditioner group. Plotted relative to the unpreconditioned baseline and including an estimate to the work involved in generating the preconditioners. 
  }
  \label{fig:large_problems_and_work_best_configuration}
\end{figure}

In comparing tuned and best configurations, again, we see a marked improvement for IC. Also, Neumann series preconditioning finally appears to offer some limited performance improvements.

\autoref{fig:large_problems_and_work_best_configuration} shows single-best and tuned-best performance plots over matrices of dimension $250$k or larger but includes an estimate to the generation work.
Here, for the single-best configuration (\autoref{fig:single_best_config_and_work_large_size}), we find that IC and MC fail only show any improvement on half the problems.
The SGS \& SSOR single-best configuration was the most robust, leading to minor improvements on a greater proportion that the other preconditioners.
When considering per-matrix tuning in \autoref{fig:tuned_best_config_and_work_large_size}, IC achieves any speedup on over $75\%$ of the subset and leads to moderate to high gains on half of the subset.
The SGS \& SSOR preconditioner improves slightly.

\subsection{Symmetric diagonally dominant}\label{sec:sdd}

\begin{figure}[tbp]
  \centering
  \begin{subfigure}[t]{0.48\textwidth}
    \includegraphics[width=\textwidth]{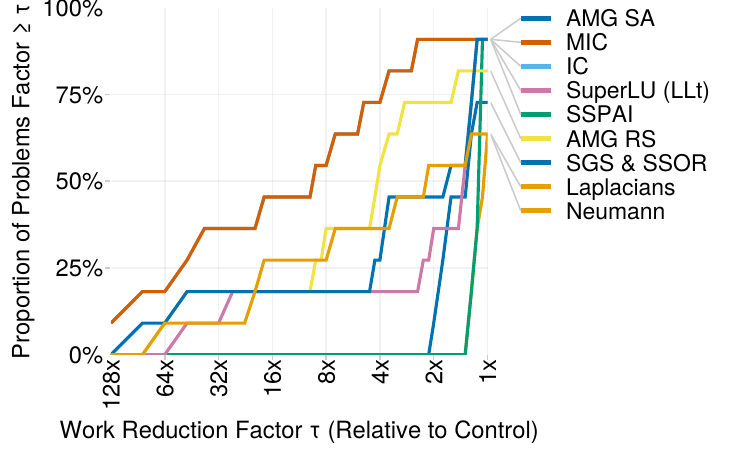}
    \caption{Single‐best configuration; }%
    \label{fig:single_best_config_sdd}
  \end{subfigure}
  \quad
  \begin{subfigure}[t]{0.48\textwidth}
    \includegraphics[width=\textwidth]{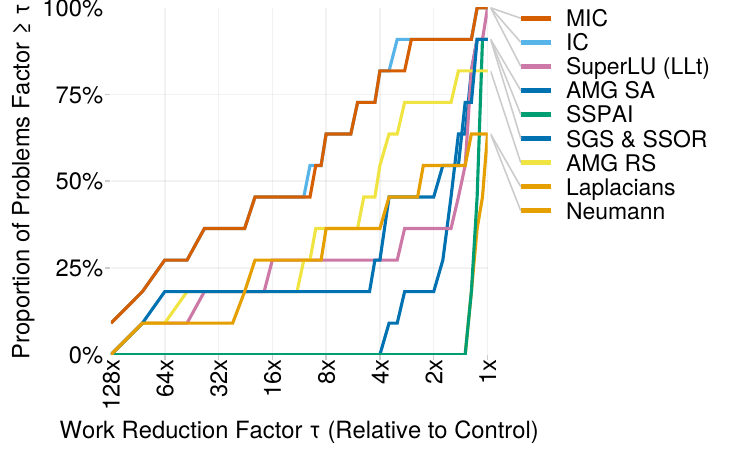}
    \caption{Tuned‐best configuration.}%
    \label{fig:tuned_best_config_sdd}
  \end{subfigure}
  \caption{%
    The Single-best configuration selects the best configuration over matrices that are diagonally dominant, whereas the tuned-best configuration selects the best run for each matrix within a given preconditioner group. Plotted relative to the unpreconditioned baseline and excluding work involved in generating the preconditioners. 
  }
  \label{fig:sdd_problems_best_configuration}
\end{figure}

Laplacian-based preconditioning offers the strongest theoretical guarantees for symmetric diagonally dominant (SDD) systems.
We specifically study this subset to identify any distinct behaviors.
Here, we consider matrices that are diagonally dominant either before or after symmetric scaling.
For matrices that lose diagonal dominance after scaling, we apply the graph Laplacian preconditioner to the original, unscaled matrix.
Otherwise, we use the scaled matrix.

We next isolate the $11$ SDD matrices and plot their profiles in \autoref{fig:sdd_problems_best_configuration}. Again, this shows results without including generation work.
Relative to the full test set, the shape of the curves change appreciably.
Laplacian-based preconditioning becomes relevant and AMD gains ground, but classical methods largely remain the same.
\emph{Note that performance curves for this subset appear choppy due to the limited sample size (11 matrices). Observed trends might not generalize to larger collections of SDD matrices; therefore, we continue to expand this subset to be able to draw more robust conclusions.}

Of note, every method attains parity with the unpreconditioned baseline on a majority of problems in this subset. In fact, the similarity among methods is interesting. The Laplacians preconditions show similar results to the incomplete factorizations.  The AMG methods also do well on this subset. 

%

Considering per-matrix tuning, which assumes \emph{a priori} knowledge of the best matrix configuration or a large setup time to find the best configuration, shows that it's now possible to accelerate 100\% of the solvers with per-matrix tuning for the two incomplete factorizations. 



\subsection{Behavior of fill-in and orderings on incomplete factorizations}
\label{sec:fill-in}

\begin{figure}[tbp]
    \centering
    \includegraphics[width=\textwidth]{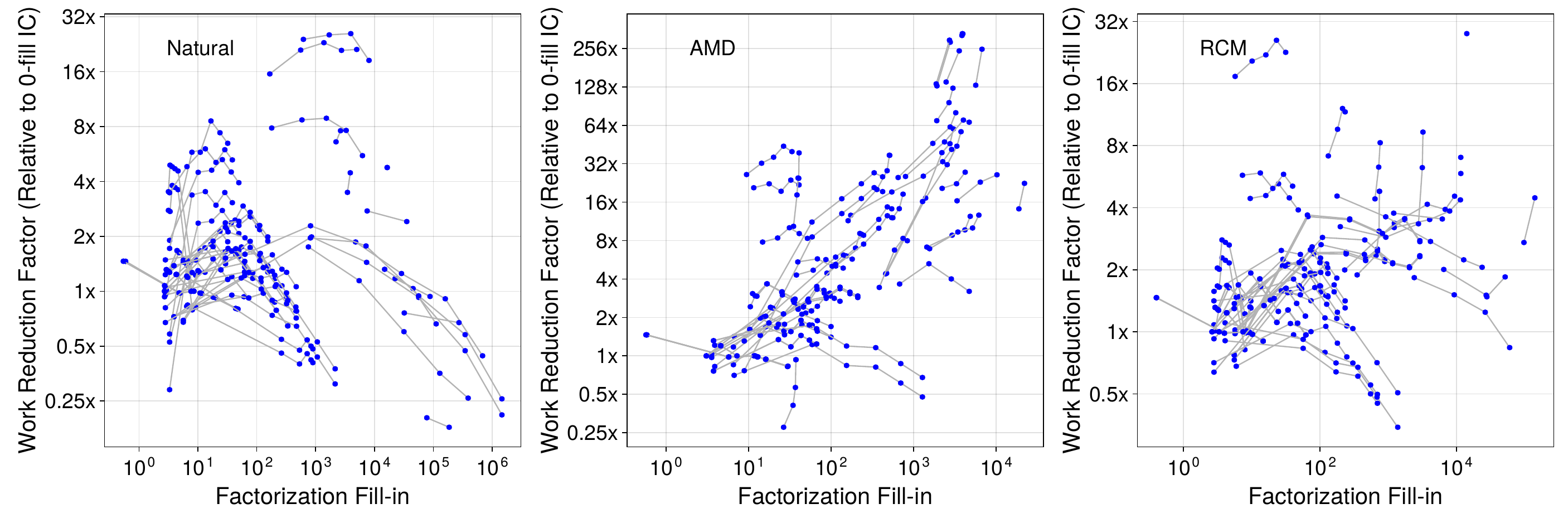}
    \caption{Factorization Fill-in versus Work Reduction Factor relative to the work reduction obtained by the IC(0) factorization. We include an estimate of the work required to generate the preconditioner in these plots.
    Each matrix is grouped and plotted with a thin gray line connecting each dot. All dots are plotted in blue.
    Note that some IC(0) factorizations failed, causing there to be no data for that matrix.}
    \label{fig:cholesky_density_versus_ic0}
\end{figure}

\begin{figure}[tbp]
    \centering
    \includegraphics[width=\textwidth]{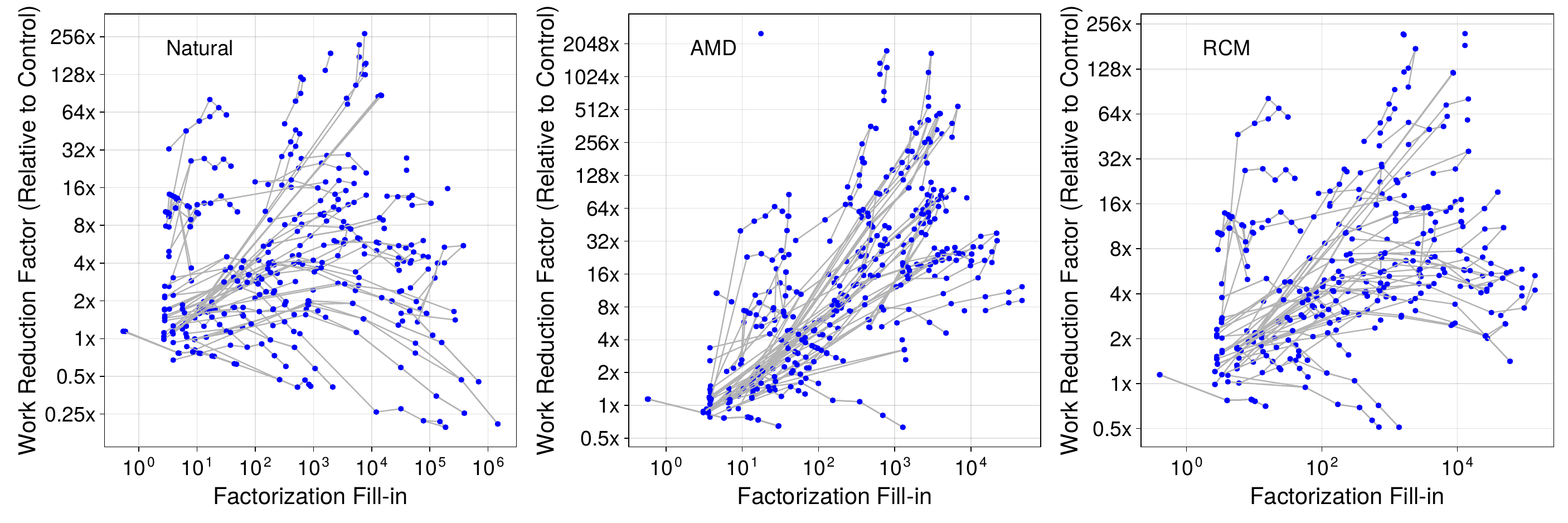}
    \caption{Factorization Fill-in versus Work Reduction Factor relative to the work reduction factor obtains by the unpreconditioned baseline. We include an estimate of the work required to generate the preconditioner in these plots. Each matrix is grouped and plotted with a thin gray line connecting each dot. All dots are plotted in blue.}
    \label{fig:cholesky_density_versus_control}
\end{figure}

Here, we discuss the impact of fill-in and ordering. In \autoref{sec:ordering} we found that AMD had the highest performance out of all the orderings we studied for incomplete factorizations. One might expect that total work -- excluding generation work -- goes down with increased fill-in of the factorization. 
We illustrate the impact of fill-in allowed during IC factorizations on work reduction compared to different baselines: a 0-fill IC preconditioner in~\autoref{fig:cholesky_density_versus_ic0} and the unpreconditioned baseline in~\autoref{fig:cholesky_density_versus_control}.
In each plot, a grouping of dots connected by a gray line represents the work reduction of different configurations of IC applied to the same matrix.\footnote{The outlier points with less than $10^0$ fill-in are associated with \texttt{Trefethen20000b} and \texttt{Trefethen20000}, two matrices that represent a curious combinatorial problem.}


Figures~\ref{fig:cholesky_density_versus_ic0} and~\ref{fig:cholesky_density_versus_control} clearly indicate that ordering strategies significantly influence the extent to which additional fill-in translates to computational savings.
The fill-reducing AMD ordering consistently achieves the greatest computational efficiency, often by orders-of-magnitude relative to Natural and RCM orderings, with substantial gains realized at moderate to high fill-in levels ($10^2$-$10^4$).
In contrast, Natural and RCM orderings exhibit diminishing returns at comparatively lower fill levels, typically saturating before $10^3$.
The qualitative ranking of orderings, AMD $\gg$ RCM $\gg$ Natural, remains unchanged regardless of baseline normalized against.
With IC(0) as the baseline, relative gains appear smaller.
Normalizing to the unpreconditioned baseline inflates those ratios but leaves the ordering curves essentially unchanged and extends coverage to more matrices.

\begin{figure}[tbp]
  \centering
  \begin{subfigure}[t]{0.48\textwidth}
    \includegraphics[width=\textwidth]{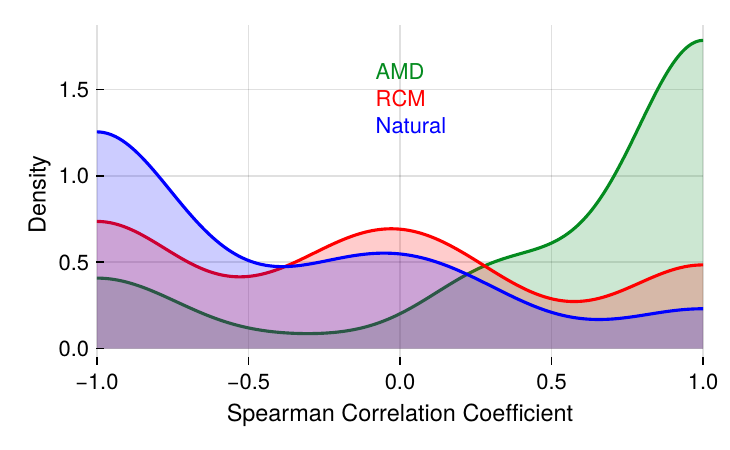}
    \caption{Relative to IC(0); }%
    \label{fig:slope-dist-0fill}
  \end{subfigure}
  \quad
  \begin{subfigure}[t]{0.48\textwidth}
    \includegraphics[width=\textwidth]{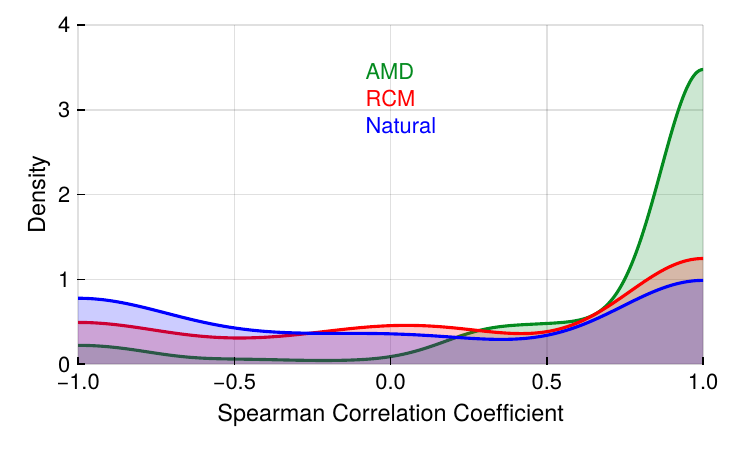}
    \caption{Relative to Unpreconditioned Baseline.}%
    \label{fig:slope-dist-control}
  \end{subfigure}
  \caption{%
    Kernel density estimate plot of the Spearman correlation coefficients of Factorization Fill-in versus Work Reduction factor for each matrix. \autoref{fig:slope-dist-0fill} is relative to the cost of convergence with the IC(0) configuration and \autoref{fig:slope-dist-control} is relative to the cost of convergence with the unpreconditioned baseline. We utilize a bounded KDE on $[-1, 1]$ with boundary conditions corrected using the reflection method~\cite{jones_sac_1993}.
  }
  \label{fig:slope-dist}
\end{figure}

To statistically study this relationship, we show a kernel density estimate of Spearman correlation correlations between fill in and work-reduction in~\autoref{fig:slope-dist}. We do this both relative to the IC(0) baseline as well as the unpreconditioned control. 
Both confirm that the strongest relationship between fill in and work reduction occurs only for AMD. 

These results have implications for practical preconditioner tuning strategies.
Given sufficient memory resources, AMD ordering with moderate-to-high fill-in is clearly advantageous from a computational work reduction viewpoint. However, this also involves substantial generation time at high fill. 
This finding challenges conventional practices in incomplete factorization-based preconditioning \citep{benzi_jcp_2002,saad_siam_2003,duff_bit_1989}, highlighting the need to reconsider ordering strategies, particularly when high fill levels are permissible.\footnote{\emph{This analysis of work reduction does not account for cache utilization, wall-clock time (for generating the factor or applying its inverse), or memory resources required to generate large fill-in factorizations.}}

\subsection{Results for best configuration}

We now present the results that highlight what was selected as the best configuration for each preconditioner when excluding generation work and comparing to the unpreconditioned baseline. 


\input{table-best-configs}


\section*{Acknowledgments}

Tunnell and Gleich acknowledge funding and support from DOE award DE-SC0023162.

\begin{multicols}{2}
\footnotesize
\bibliographystyle{dgleich-bib3}
\bibliography{bib/book_citations,bib/iterative_methods,bib/matrices,bib/preconditioners,bib/related,bib/software,bib/full_matrix,bib/matrices_related,bib/matrix_groups}
\end{multicols}

\clearpage

\appendix 
\section*{Appendix}

\etocsetstyle{subsection}
  {}
  {%
    \parindent0pt
    \footnotesize\hangindent=\dimexpr1.25\secindent\relax
    \hangafter=1
    \leftskip\secindent\rightskip0pt
  }
  {%
    \makebox[\dimexpr1.75\secindent\relax][l]{\etocnumber\quad}%
    \etocname\hspace{1em}\textit{\etocpage}\par
  }
  {}
  
\etocdepthtag.toc{appendix}

\etocsettagdepth{main}{none}
\etocsettagdepth{appendix}{subsubsection}

\begin{multicols}{2}
\tableofcontents    
\end{multicols}

\pagebreak

\section{Supplemental details}

\subsection{AUC calculation for performance profiles}
\label{sec:auc}

We use the area under the curve to summarize each performance-profile curve into a single scalar value, preserving both the proportion of problems solved and the efficiency in terms of work reduction. 
To compute the AUC for a performance profile, we take as input the fraction of problems solved $y_i$ along with the work reduction factor $x_i$ for each datapoint. 
Integration is performed in logarithmic space and normalized so that a method solving all problems within a $128\times$ work reduction of the best achieves an AUC of $1.00$.
The curve with the highest AUC is considered optimal.

We integrate using the composite trapezoidal rule over the domain $[2^{-2},2^7]$.
The $\log_2$-width of this domain is $9$, providing the constant divisor in the Julia code below. When we use this code, we have $x$ and $y$ sampled on a regular grid within this domain instead of looking at values only at the jump points of the performance profile. 

\begin{verbatim}
@inline function auc(x::AbstractVector, y::AbstractVector)
    @assert length(x) == length(y)
    @views 0.5 * sum((y[1:end-1] .+ y[2:end]) .* diff(log2(x))) / 9
end
\end{verbatim}

Due to the piecewise-linear nature of performance profiles, the trapezoidal rule yields exact integration for the plotted values.
Calculating this metric in $\log_2$ space treats each multiplicative factor in work reduction as equal distance, aligning with the visual presentation used throughout this study.

\subsection{Generation cost for incomplete factorizations and sparse direct methods}
\label{sec:generation-cost}
For the incomplete factorizations and 
sparse direct methods, we utilize an approximation to the amount of work required to generate the symmetric factorization. Let $\mL$ be the factor in the approximate factorization $\mL \mL^T \approx \mA$. Then we use 
\begin{align*}
    \texttt{work} &= \sum_{i=1}^{n}{c_i^2},
\end{align*}
where $c_i$ is the number of elements in the $i$-th column of $\mL$. 

SuiteSparse's Cholesky factorization software uses precisely this metric during its symbolic phase to compute factorization costs and select optimal matrix orderings.
For incomplete factorizations, particularly threshold-based approaches, this measure slightly underestimates actual computational work.
In practice, additional entries beyond those retained after thresholding are temporarily computed and subsequently dropped.
Therefore, the metric presented provides a lower bound rather than an exact measure.

\input{all-problems.tex}

\end{document}

%% file: commands.tex
\newcommand{\normof}[1]{\|#1\|}
\newcommand{\onormof}[2]{\|#1\|_{#2}}

\newcommand{\wpar}[1]{\ensuremath{\left(#1\right)}}

\providecommand{\mA}{\ensuremath{\mat{A}}}

\providecommand{\mD}{\ensuremath{\mat{D}}}

\providecommand{\mI}{\ensuremath{\mat{I}}}
\providecommand{\mJ}{\ensuremath{\mat{J}}}

\providecommand{\mL}{\ensuremath{\mat{L}}}
\providecommand{\mM}{\ensuremath{\mat{M}}}
\providecommand{\mN}{\ensuremath{\mat{N}}}

\providecommand{\mP}{\ensuremath{\mat{P}}}

\providecommand{\mR}{\ensuremath{\mat{R}}}
\providecommand{\mS}{\ensuremath{\mat{S}}}

\providecommand{\mU}{\ensuremath{\mat{U}}}

\providecommand{\vb}{\ensuremath{\vec{b}}}

\providecommand{\ve}{\ensuremath{\vec{e}}}

\providecommand{\vp}{\ensuremath{\vec{p}}}

\providecommand{\vr}{\ensuremath{\vec{r}}}

\providecommand{\vx}{\ensuremath{\vec{x}}}

\providecommand{\vz}{\ensuremath{\vec{z}}}

\def\eqref#1{equation~\ref{#1}}

\def\1{\bm{1}}

\def\vb{{\bm{b}}}

\def\ve{{\bm{e}}}

\def\vp{{\bm{p}}}

\def\vr{{\bm{r}}}

\def\vx{{\bm{x}}}

\def\vz{{\bm{z}}}

\def\mA{{\bm{A}}}

\def\mD{{\bm{D}}}

\def\mI{{\bm{I}}}
\def\mJ{{\bm{J}}}

\def\mL{{\bm{L}}}
\def\mM{{\bm{M}}}
\def\mN{{\bm{N}}}

\def\mP{{\bm{P}}}

\def\mR{{\bm{R}}}
\def\mS{{\bm{S}}}

\def\mU{{\bm{U}}}

\DeclareMathAlphabet{\mathsfit}{\encodingdefault}{\sfdefault}{m}{sl}
\SetMathAlphabet{\mathsfit}{bold}{\encodingdefault}{\sfdefault}{bx}{n}

\def\gL{{\mathcal{L}}}

\def\sR{{\mathbb{R}}}



%% file: table-best-configs.tex
\newlength{\mycolwidth}
\setlength{\mycolwidth}{29pt}
{\footnotesize
    \begin{longtable}{
        >{\raggedright\arraybackslash}p{1.8cm}   
        >{\raggedright\arraybackslash}p{1cm} 
        >{\raggedright\arraybackslash}p{\mycolwidth}        
        >{\raggedleft\arraybackslash}p{\mycolwidth}         
        >{\raggedleft\arraybackslash}p{\mycolwidth}         
        >{\raggedleft\arraybackslash}p{\mycolwidth}         
        >{\raggedleft\arraybackslash}p{\mycolwidth}         
        >{\raggedleft\arraybackslash}p{\mycolwidth}         
        >{\raggedleft\arraybackslash}p{\mycolwidth}         
        >{\raggedleft\arraybackslash}p{\mycolwidth}         
    }\caption{Best and worst configurations for each preconditioner. 
    The \textbf{Succ.} column is the fraction of problems in which the preconditioner was successfully generated; 
    \textbf{Parity} is the percentage of problems where the method required no more work than the unpreconditioned baseline (i.e., $1\times$ the work or better), indicating no degradation;
    \textbf{AUC} gives the area under the performance profile curve to  measure overall effectiveness across the test set (higher values indicate better performance), normalized to 1; 
    \textbf{Geo.} measures the geometric mean of the work reduction factor over the entire testset, where failure to converge or generate is defined as $4\times$ control (higher is better, max is 128);
    columns $\bm{\geq2\times}$, $\bm{\geq4\times}$, and $\bm{\geq8\times}$ refer to the percentage of problems solved in, e.g., 4 times less work or better.}
    \label{tab:methods_overview}\label{tab:best_configurations}\\%
        \toprule
        \textbf{Class} & & \textbf{Config} &
        \textbf{AUC} & \textbf{Geo.} & \textbf{Succ.} & \textbf{Parity} &
        $\bm{\geq 2\times}$ & $\bm{\geq 4\times}$ & $\bm{\geq 8\times}$ \\
        \endfirsthead
        \caption[]{\textit{Best and worst configurations for each preconditioner (continued)}}\\
        \toprule
        \textbf{Class} & \textbf{Case} & \textbf{Config} & \textbf{AUC} & \textbf{Geo.} &
        \textbf{Succ.} & \textbf{Parity} & $\bm{\ge 2\times}$ &
        $\bm{\ge 4\times}$ & $\bm{\ge 8\times}$ \\ \midrule
        \endhead
        \addlinespace 
        \multicolumn{10}{r}{\textit{continued on next page}}\\
        \addlinespace
        \midrule
        \endfoot
        \bottomrule
        \endlastfoot
        \midrule
        \multirow{3}{*}{\parbox[t]{2cm}{Truncated\\Neumann%
        \footnote{All configurations were insensitive to matrix ordering. No significant outlier configurations observed; however, many configurations provided no improvement on any matrix. The best configuration by AUC that did not use the spectral norm was 1 iter of $1/\onormof{\mA}{F}$. See \autoref{sec:truncated_neumann_series} for a description of the preconditioner.}%
        }}%
        & Best        & 2~iters $2/\onormof{\mA}{2}$ & 0.16 & 0.69 & 100\% & 1\% & 1\% & 0\% & 0\% \\*
        & Second & 1~iter $2/\onormof{\mA}{2}$ & 0.13 & 0.57 & 100\% & 11\% & 0\% & 0\% & 0\% \\*
        & Worst       & 4~iters $1/\onormof{\mA}{F}$ & 0.00 & 0.25 & 100\% & 0\% & 0\% & 0\% & 0\% \\*
        \midrule
                \multirow{3}{*}{\parbox[t]{2cm}{SGS%
                \footnote{For SGS, with the 1 iteration configuration, RCM was as good or better than natural and AMD orderings on $72\%$ of the test set. AMD ordering lead to a degradation in performance compared to natural ordering on $71\%$ of the test set with this config. See \autoref{sec:sgs-ssor}.}%
                }}%
        & Best        & 1~iter RCM & 0.24 & 1.14 & 100\% & 72\% & 3\% & 0\% & 0\% \\*
        & Second & 1~iter Nat. & 0.23 & 1.07 & 100\% & 62\% & 1\% & 0\% & 0\% \\*
        & Worst       & 2~iters AMD & 0.20 & 0.87 & 100\% & 27\% & 1\% & 0\% & 0\% \\*
        \midrule
                        \multirow{3}{*}{\parbox[t]{2cm}{SSOR%
                        \footnote{For SSOR, similar to SGS, on the 1 iteration and $\omega = 1.2$ configuration, RCM is the best ordering for SSOR on $75\%$ of the test set. AMD performs poorly again, leading to a degradation in performance compared to natural ordering on $75\%$ of the test set. See \autoref{sec:sgs-ssor}.}%
                        }}%
        & Best        & 1~iter RCM $\omega1.2$ & 0.23 & 1.05 & 100\% & 63\% & 1\% & 0\% & 0\% \\*
        & Second & 1~iter Nat. $\omega1.2$ & 0.22 & 0.97 & 100\% & 44\% & 1\% & 0\% & 0\% \\*
        & Worst       & 2~iters AMD $\omega1.8$ & 0.05 & 0.33 & 100\% & 0\% & 0\% & 0\% & 0\% \\*
        \midrule
                \multirow{3}{*}{\parbox[t]{2cm}{Sym.\\Sparse\\ Approx.\\Inverse\\(SSPAI)%
                \footnote{This preconditioner was insensitive to matrix ordering.  See \autoref{sec:ssai} for preconditioner details.}%
                }}%
        & Best        & $\frac{1}{2}\times$~fill & 0.26 & 1.24 & 100\% & 82\% & 1\% & 0\% & 0\% \\*
        & Second & $1\times$~fill & 0.25 & 1.22 & 100\% & 84\% & 3\% & 0\% & 0\% \\*
        & Worst       & $3\times$~fill & 0.23 & 1.05 & 100\% & 71\% & 3\% & 0\% & 0\% \\*
        \addlinespace[.85cm]
        \midrule
                \multirow{3}{*}{\parbox[t]{2cm}{SuperLU%
                \footnote{For SuperLU, AMD was the best ordering for this preconditioner, where the top 6 configurations by AUC were of this ordering. The top configuration was the best single preconditioner on \texttt{bundle1}, \texttt{bundle\_adj}, and \texttt{gyro\_m}. See \autoref{sec:superlu} for details.}}}%
        & Best        & $10^{-4}$ $1\times$~fill AMD & 0.19 & 0.85 & 68\% & 28\% & 14\% & 11\% & 8\% \\*
        & Second & $10^{-5}$ $3\times$~fill AMD & 0.17 & 0.73 & 66\% & 24\% & 10\% & 8\% & 5\% \\*
        & Worst       & $10^{-6}$ $3\times$~fill Nat. & 0.10 & 0.47 & 56\% & 10\% & 8\% & 8\% & 8\% \\*
        \midrule
                \multirow{3}{*}{\parbox[t]{2cm}{Incomplete\\Cholesky%
                \footnote{AMD was the ordering in the top $4$ configurations (by AUC). The top configuration was the single best overall preconditioner (or tied with the best) in terms of work reduction on $63.29\%$ of the test set. 
                See \autoref{sec:incomplete-cholesky}.}%
                }}%
        & Best        & $10^{-8}$ AMD & 0.64 & 17.71 & 84\% & 80\% & 78\% & 73\% & 62\% \\*
        & Second & $10^{-7}$ AMD & 0.62 & 14.77 & 82\% & 80\% & 78\% & 73\% & 63\% \\*
        & Worst       & IC(0) AMD & 0.23 & 1.06 & 70\% & 67\% & 18\% & 6\% & 4\% \\*
        \midrule
                \multirow{3}{*}{\parbox[t]{2cm}{Modified\\Incomplete\\Cholesky%
                \footnote{AMD was the best ordering for the top $3$ configurations for this preconditioner. MIC performed identically to IC on a number of matrices and the top configuration was tied with the best configuration of IC on $60.76\%$ of the test set. See \autoref{sec:modified-incomplete-cholesky}.}
                }}%
        & Best        & $10^{-8}$ AMD & 0.61 & 15.01 & 81\% & 77\% & 77\% & 71\% & 59\% \\*
        & Second & $10^{-7}$ AMD & 0.55 & 8.78 & 75\% & 72\% & 71\% & 67\% & 58\% \\*
        & Worst       & IC(0) AMD & 0.04 & 0.33 & 16\% & 13\% & 3\% & 1\% & 1\% \\*
        \midrule
              \multirow{3}{*}{\parbox[t]{2cm}{AMG\\Ruge-Stuben%
              \footnote{Although AMG is considered to be invariant to ordering, this implementation had differing behavior on small numbers of matrices depending on the ordering. On \texttt{Flan\_1565}, the \texttt{F} and \texttt{W} cycle configurations stalled out with natural ordering but not with other orderings. Other observed slow convergence behavior largely occurred with \texttt{AMLI} cycles. See \autoref{subsec:amg} for preconditioner and setting details.}
              }}%
        & Best        & 1~iter \texttt{V} cyc. (RCM) & 0.34 & 2.05 & 100\% & 89\% & 30\% & 16\% & 11\% \\*
        & Second & 2~iters \texttt{V} cyc. (RCM) & 0.31 & 1.74 & 100\% & 84\% & 24\% & 11\% & 8\% \\*
        & Worst       & 2~iters \texttt{AMLI}~cyc. (AMD) & 0.09 & 0.43 & 100\% & 6\% & 0\% & 0\% & 0\% \\*
        \midrule
        \multirow{3}{*}{\parbox[t]{2cm}{AMG\\Smoothed Agg.%
        \footnote{This preconditioner was invariant to ordering. One iteration of each cycle were each either the best overall preconditioner or tied with the best on $2.5\%$ of the test set. Performing a single iteration of the preconditioner tended to be the best -- 2 iterations of \texttt{V} cycles (the best 2 iteration cycle) had .34 AUC and 2.09 work reduction geometric mean and lead to a $2\times$, $4\times$ or $8\times$ or greater speedup at lower frequencies than any smoothing aggregation cycle with a single iteration. See \autoref{subsec:amg}.}
        }}%
        & Best        & 1~iter \texttt{V} cyc. & 0.36 & 2.39 & 100\% & 96\% & 38\% & 15\% & 10\% \\*
        & Second & 1~iter \texttt{W} cyc. & 0.36 & 2.32 & 100\% & 94\% & 38\% & 14\% & 9\% \\*
        & Worst       & 2~iters \texttt{W} cyc. & 0.31 & 1.72 & 100\% & 73\% & 19\% & 11\% & 9\% \\*
    \end{longtable}%
}

\keepXColumns  

%% file: all-problems.tex
\section{All Problems Convergence Results}\label{sec:appendix1}

Each of the following pages shows a convergence result from one of the problems. 

\clearpage 
\subsection{2cubes\_sphere}

\begin{figure}[!ht]
    \centering
    \includegraphics[width=\textwidth]{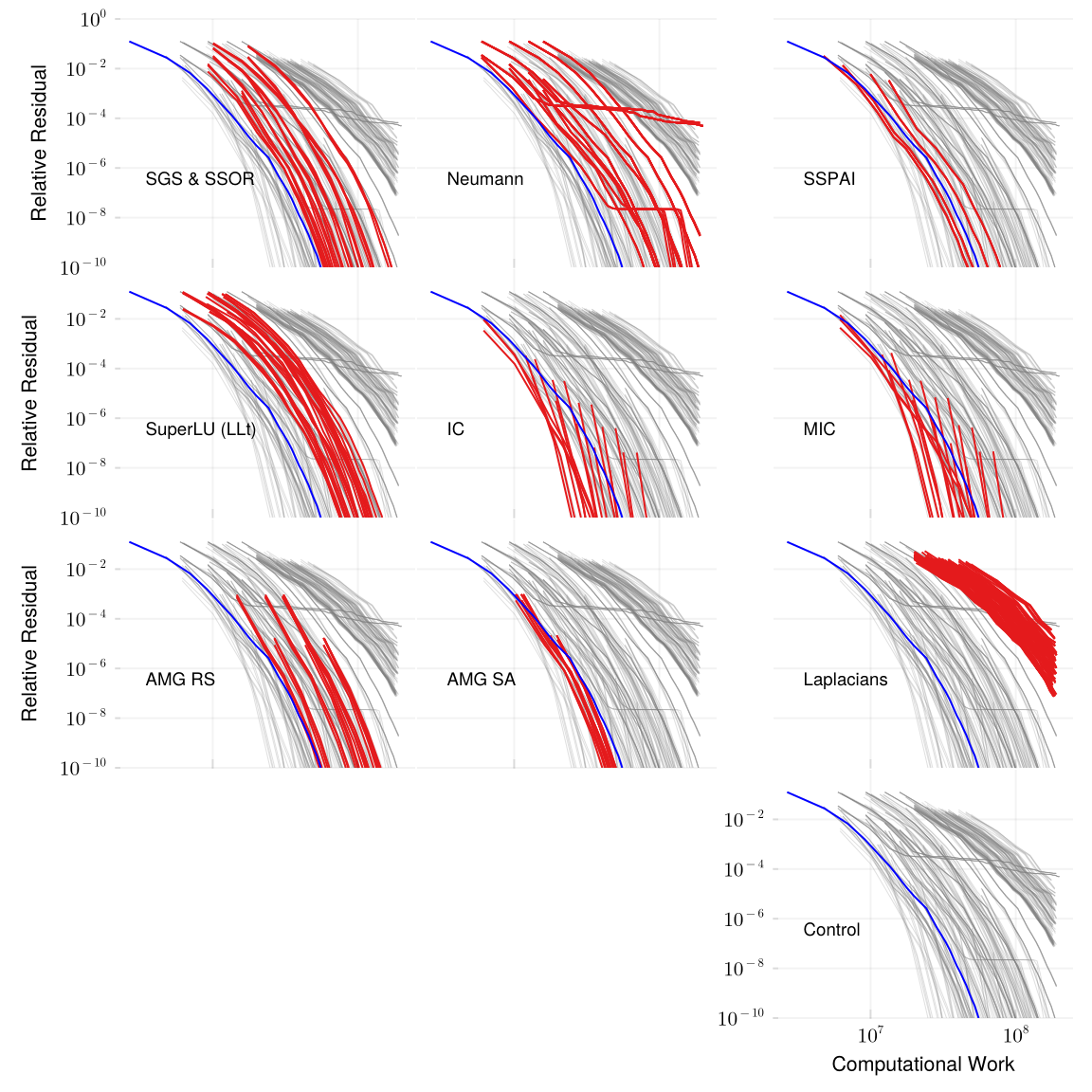}
    \caption{Convergence of the PCG method with various preconditioners applied to the \texttt{2cubes\_sphere} matrix (101k rows, 1.6m non-zeros). The plots have a log-log scale.}
    \label{fig:2cubes_sphere}
\end{figure}
\clearpage 

\subsection{BenElechi1}

\begin{figure}[!ht]
    \centering
    \includegraphics[width=\textwidth]{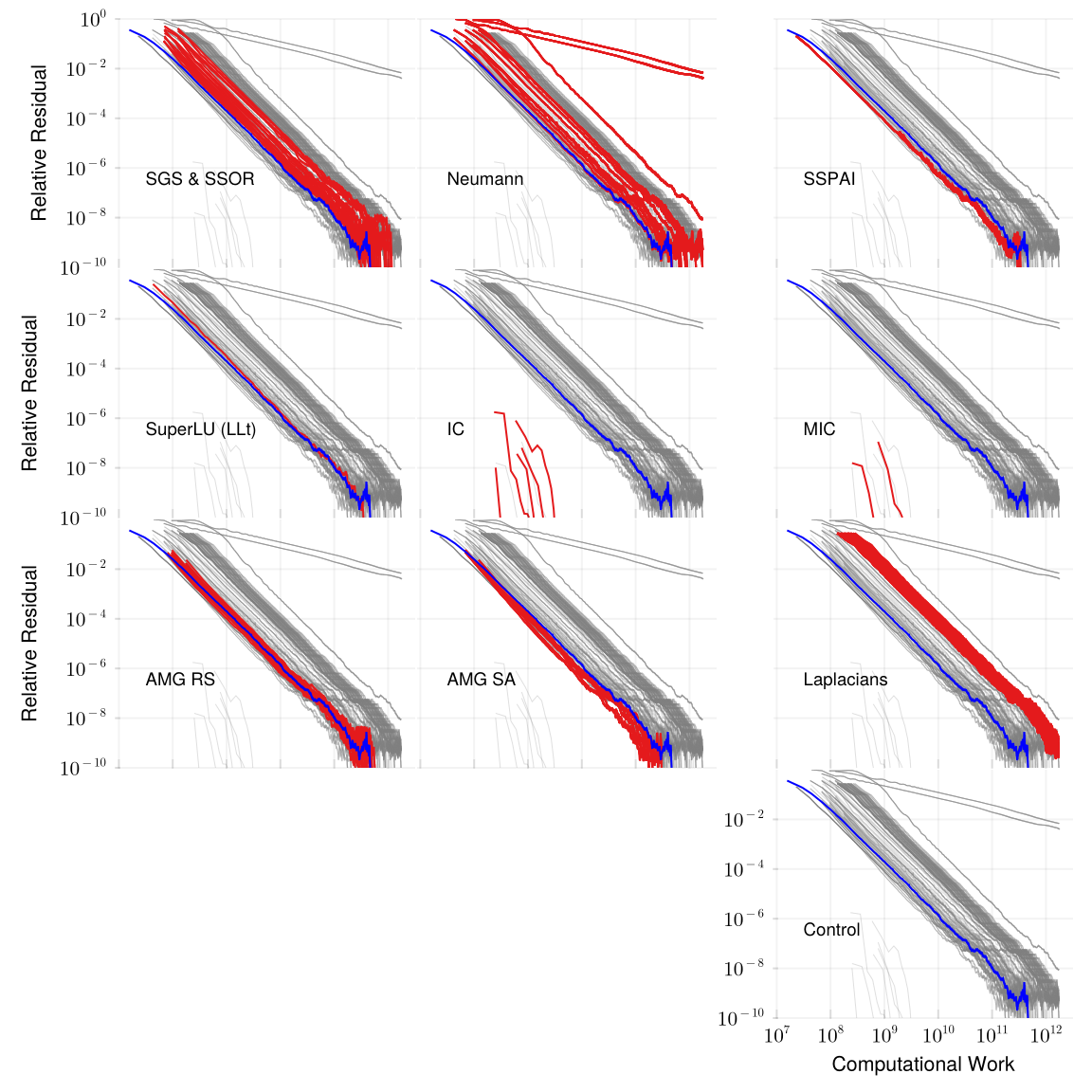}
    \caption{Convergence of the PCG method with various preconditioners applied to the \texttt{BenElechi1} matrix (246k rows, 13.2m non-zeros). The plots have a log-log scale.}
    \label{fig:benElechi1}
\end{figure}
\clearpage 

\subsection{Bump\_2911}

\begin{figure}[!ht]
    \centering
    \includegraphics[width=\textwidth]{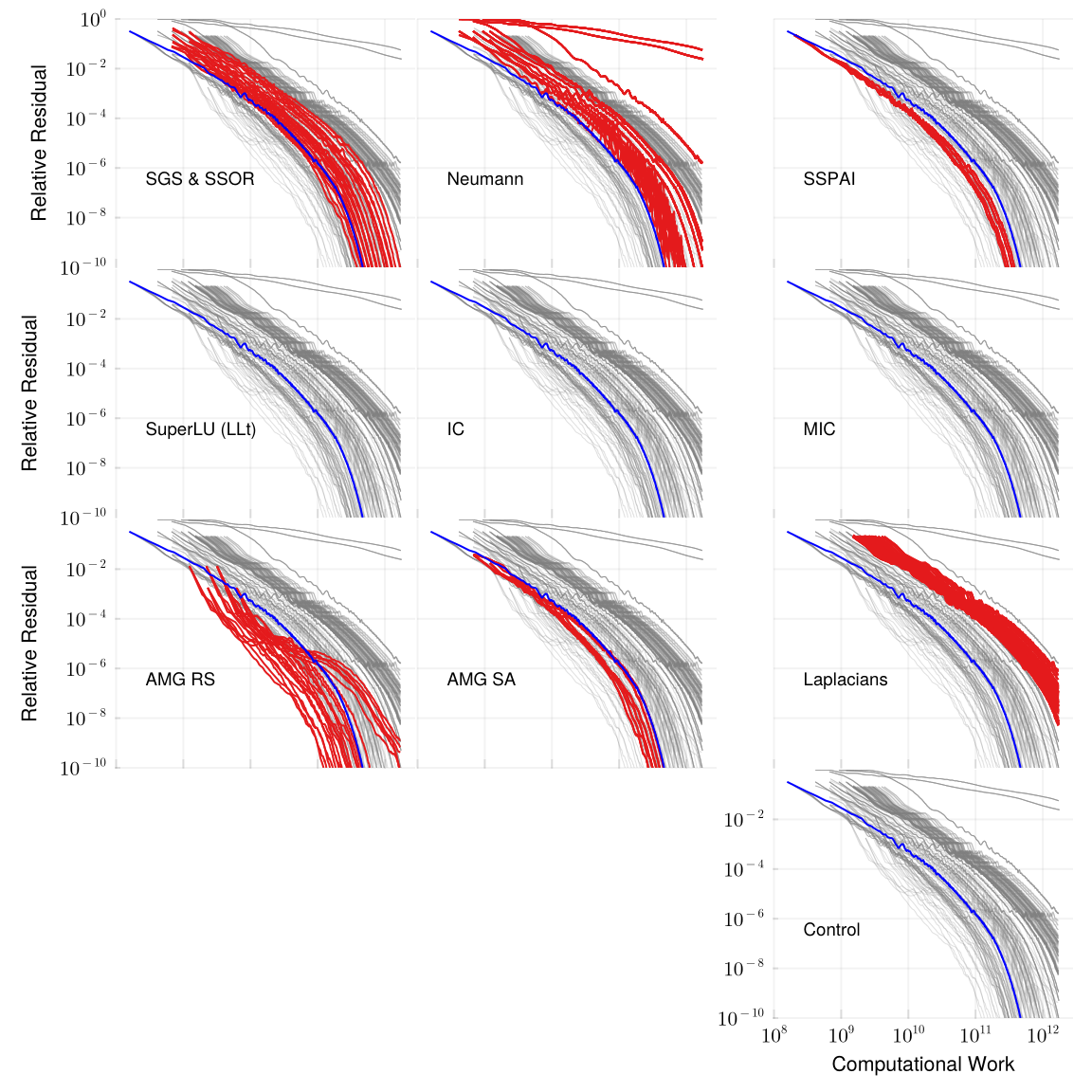}
    \caption{Convergence of the PCG method with various preconditioners applied to the \texttt{Bump\_2911} matrix (2.9m rows, 128m non-zeros). The plots have a log-log scale.}
    \label{fig:Bump_2911}
\end{figure}
\clearpage 

\subsection{Dubcova1}

\begin{figure}[!ht]
    \centering
    \includegraphics[width=\textwidth]{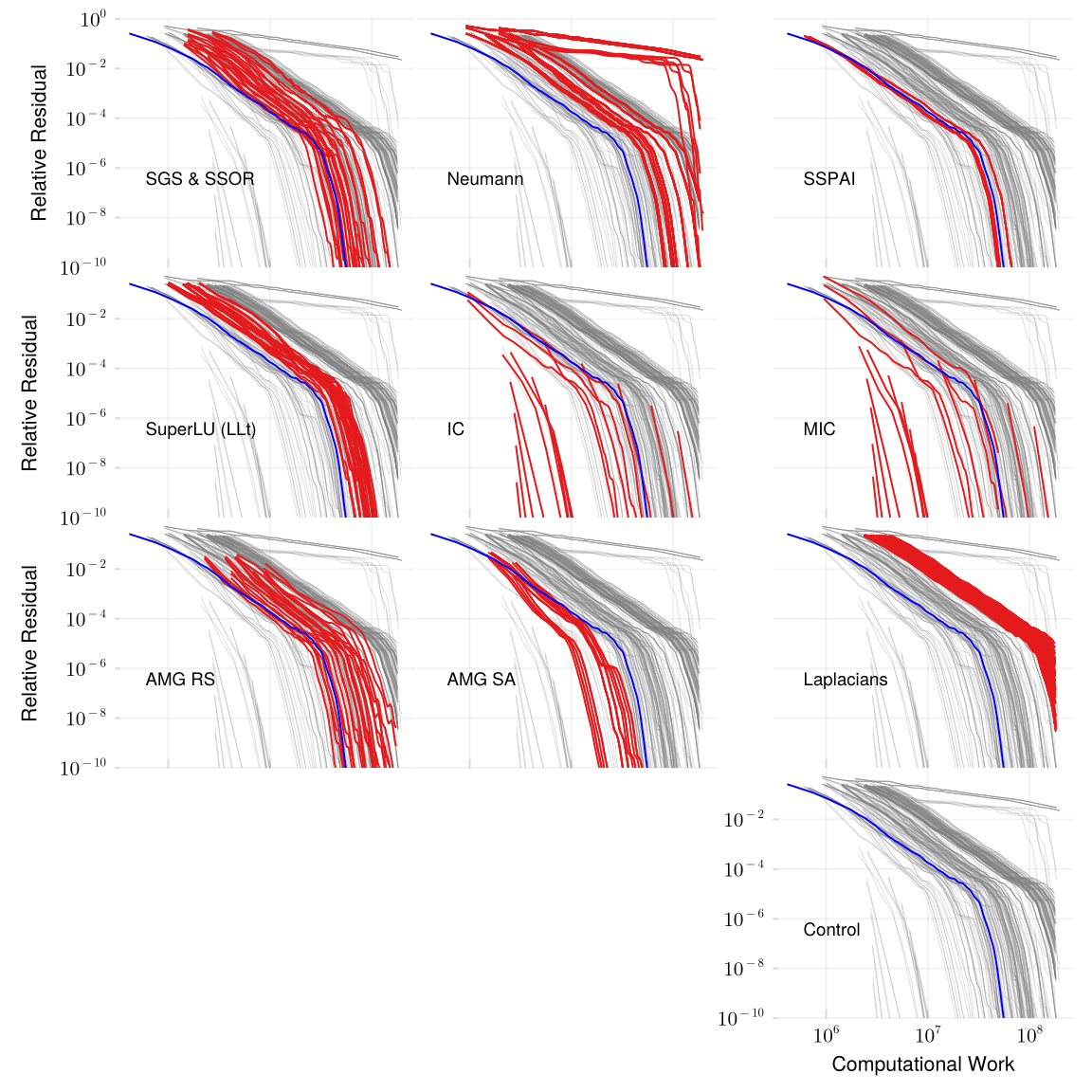}
    \caption{Convergence of the PCG method with various preconditioners applied to the \texttt{Dubcova1} matrix (16k rows, 253k non-zeros). The plots have a log-log scale.}
    \label{fig:Dubcova1}
\end{figure}
\clearpage 

\subsection{Dubcova2}

\begin{figure}[!ht]
    \centering
    \includegraphics[width=\textwidth]{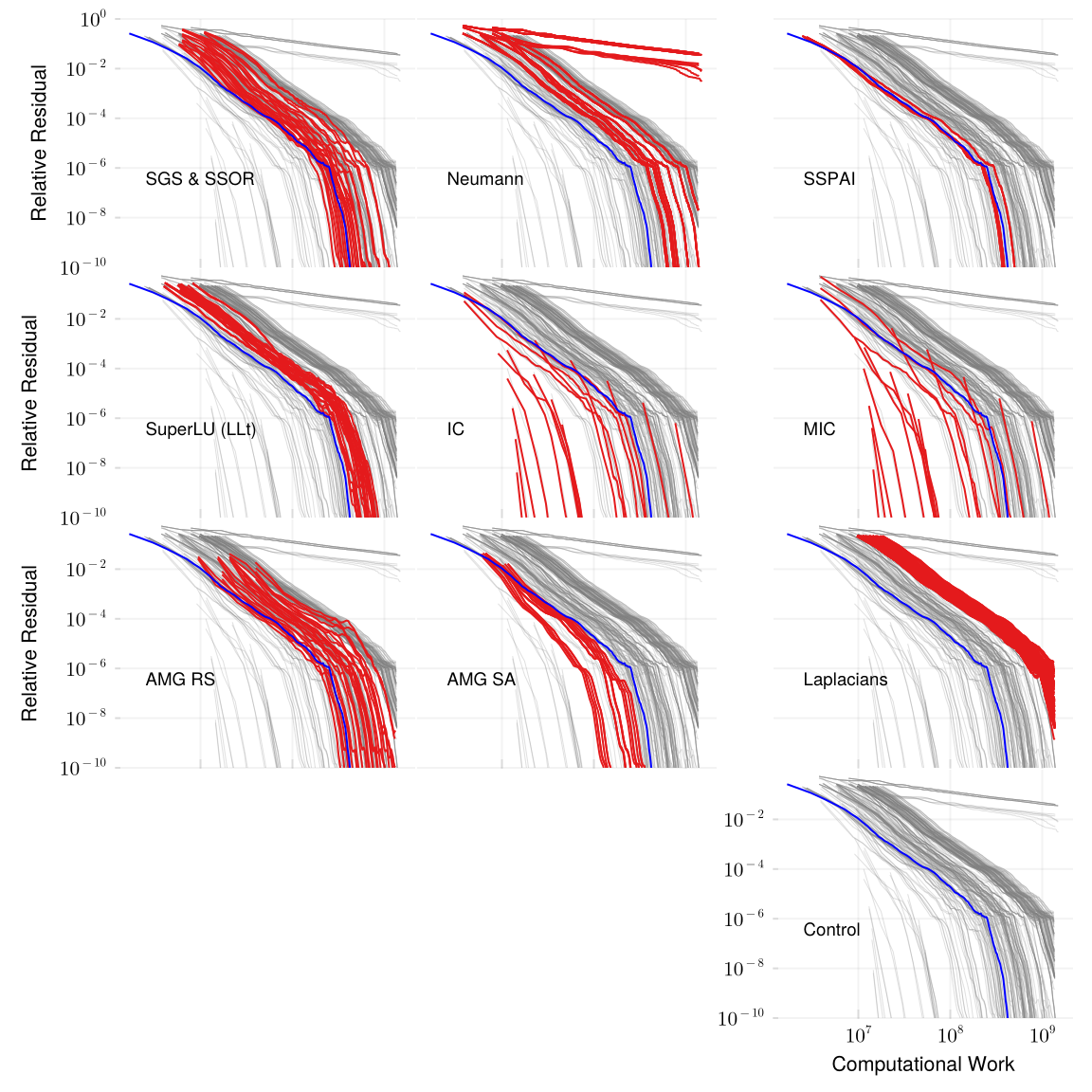}
    \caption{Convergence of the PCG method with various preconditioners applied to the \texttt{Dubcova2} matrix (65k rows, 1m non-zeros). The plots have a log-log scale.}
    \label{fig:Dubcova2}
\end{figure}
\clearpage 

\subsection{Dubcova3}

\begin{figure}[!ht]
    \centering
    \includegraphics[width=\textwidth]{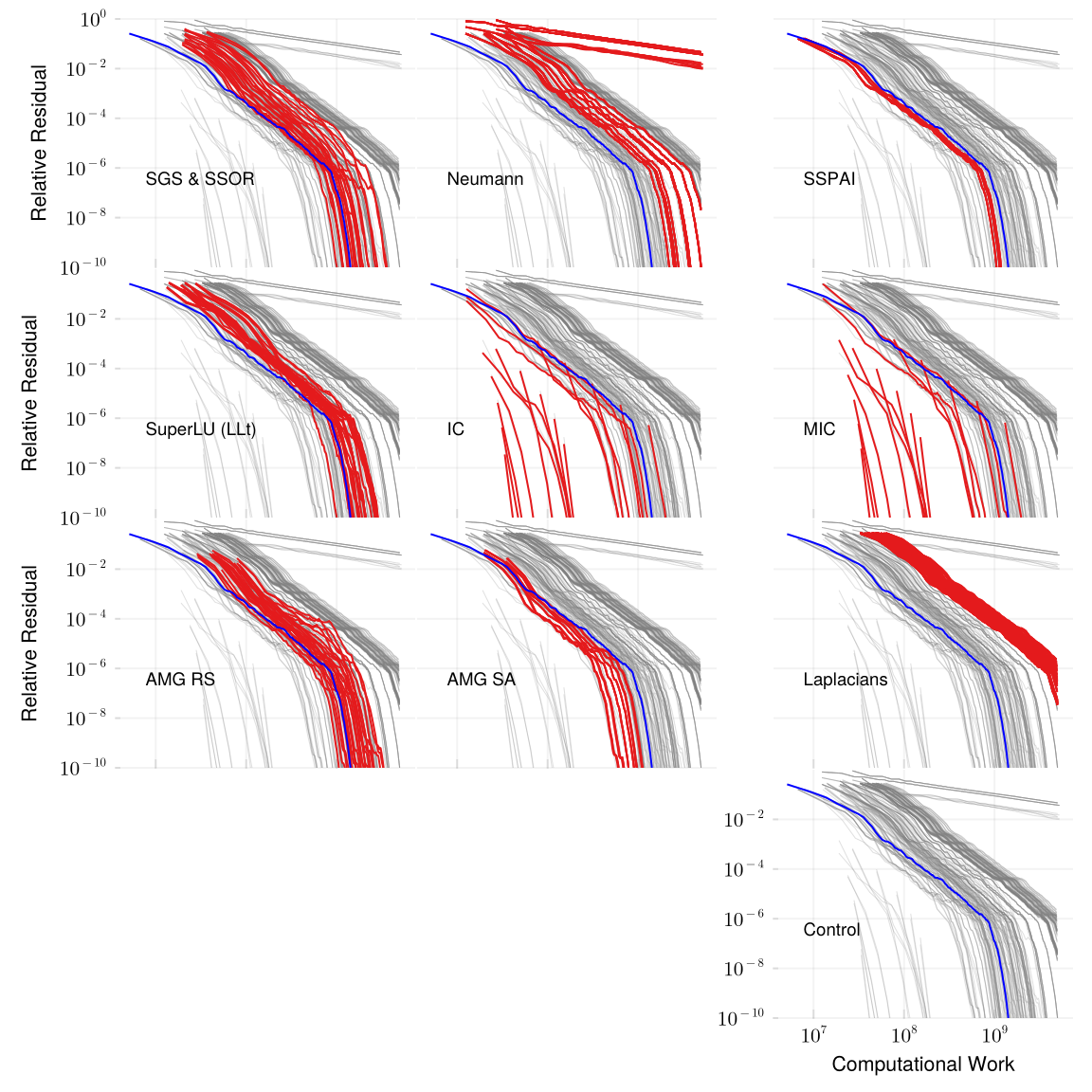}
    \caption{Convergence of the PCG method with various preconditioners applied to the \texttt{Dubcova3} matrix (147k rows, 3.6m non-zeros). The plots have a log-log scale.}
    \label{fig:Dubcova3}
\end{figure}
\clearpage 

\subsection{Emilia\_923}

\begin{figure}[!ht]
    \centering
    \includegraphics[width=\textwidth]{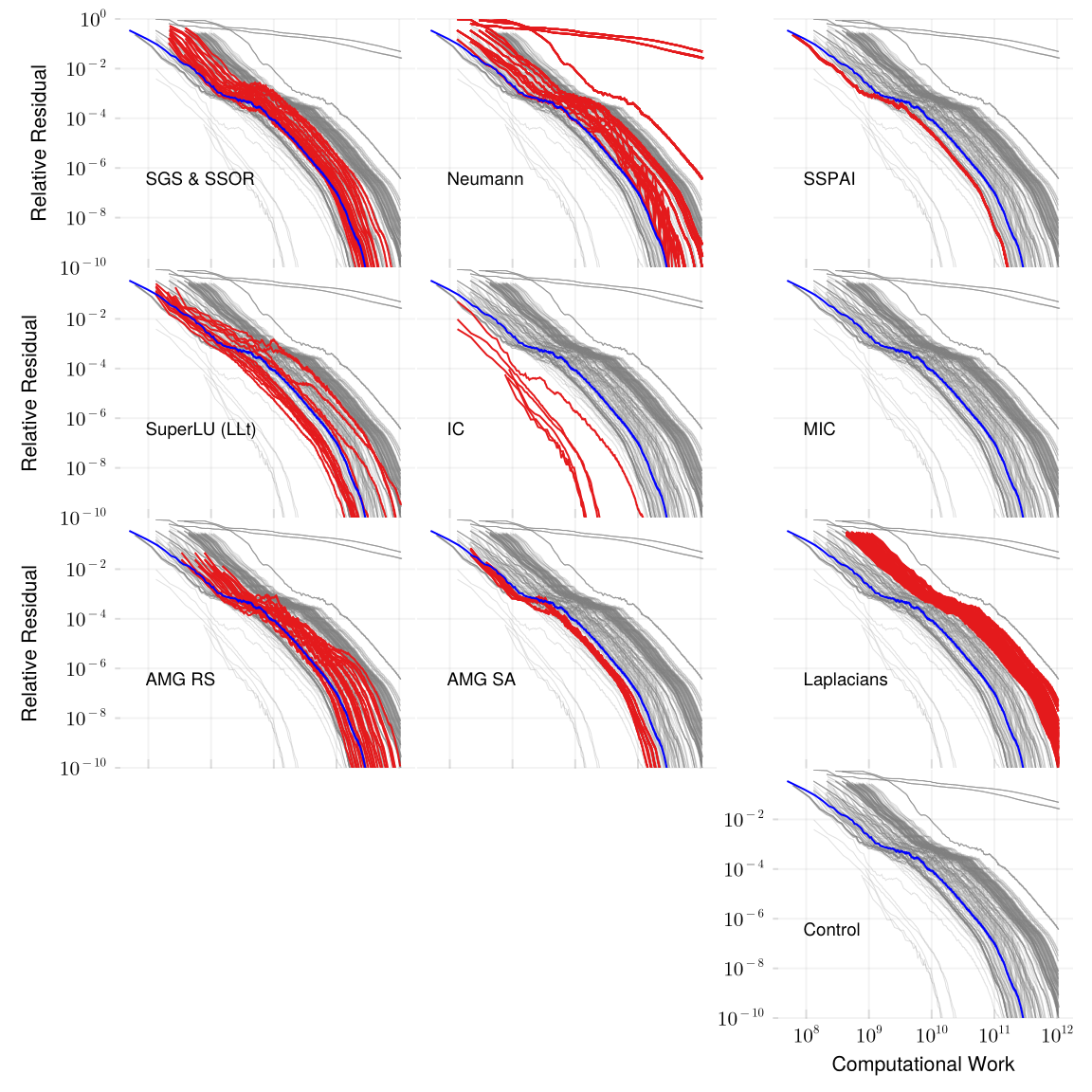}
    \caption{Convergence of the PCG method with various preconditioners applied to the \texttt{Emilia\_923} matrix (923k rows, 40.4m non-zeros). The plots have a log-log scale.}
    \label{fig:Emilia_923}
\end{figure}
\clearpage 

\subsection{Fault\_639}

\begin{figure}[!ht]
    \centering
    \includegraphics[width=\textwidth]{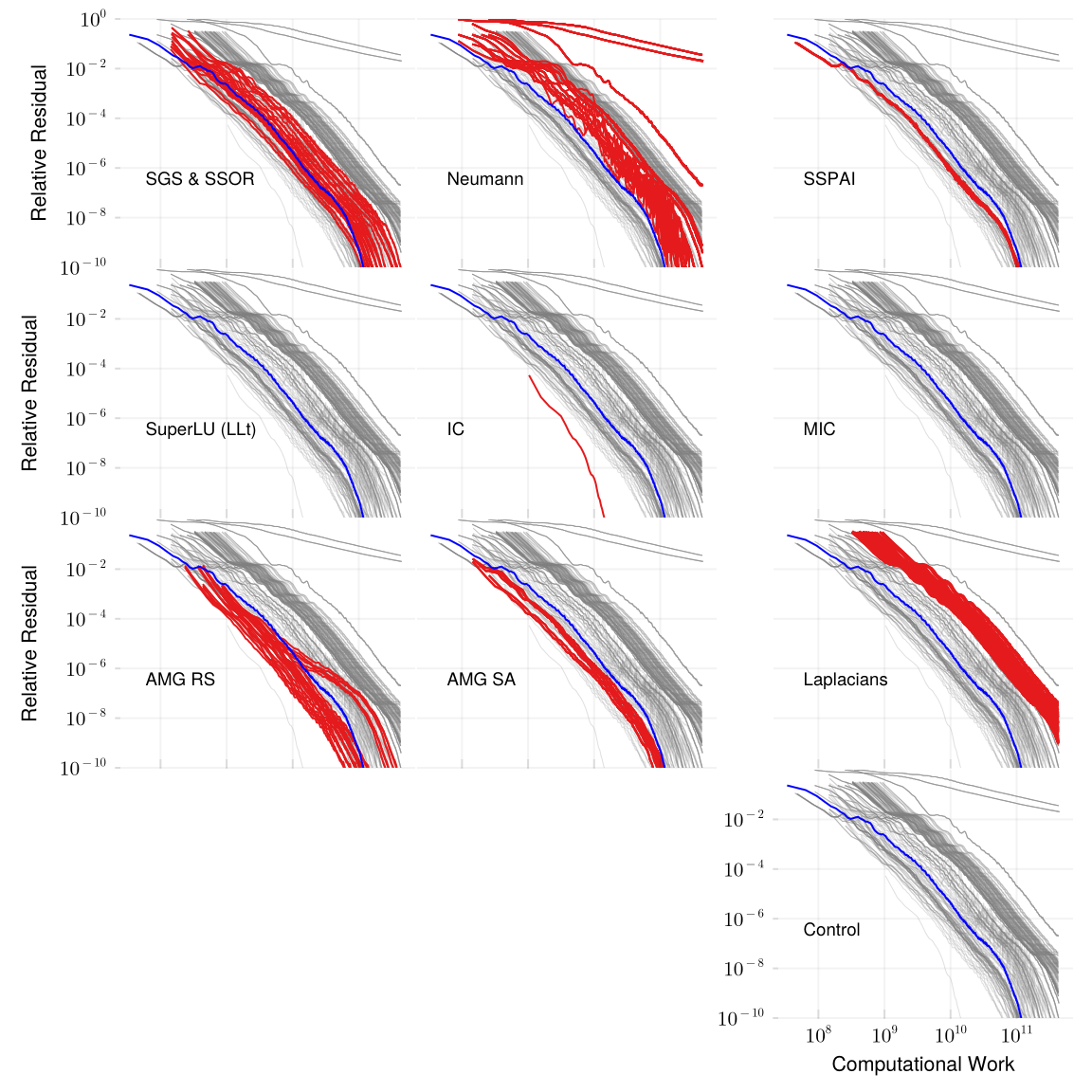}
    \caption{Convergence of the PCG method with various preconditioners applied to the \texttt{Fault\_639} matrix (639k rows, 27.2m non-zeros). The plots have a log-log scale.}
    \label{fig:Fault_639}
\end{figure}
\clearpage 

\subsection{Flan\_1565}

\begin{figure}[!ht]
    \centering
    \includegraphics[width=\textwidth]{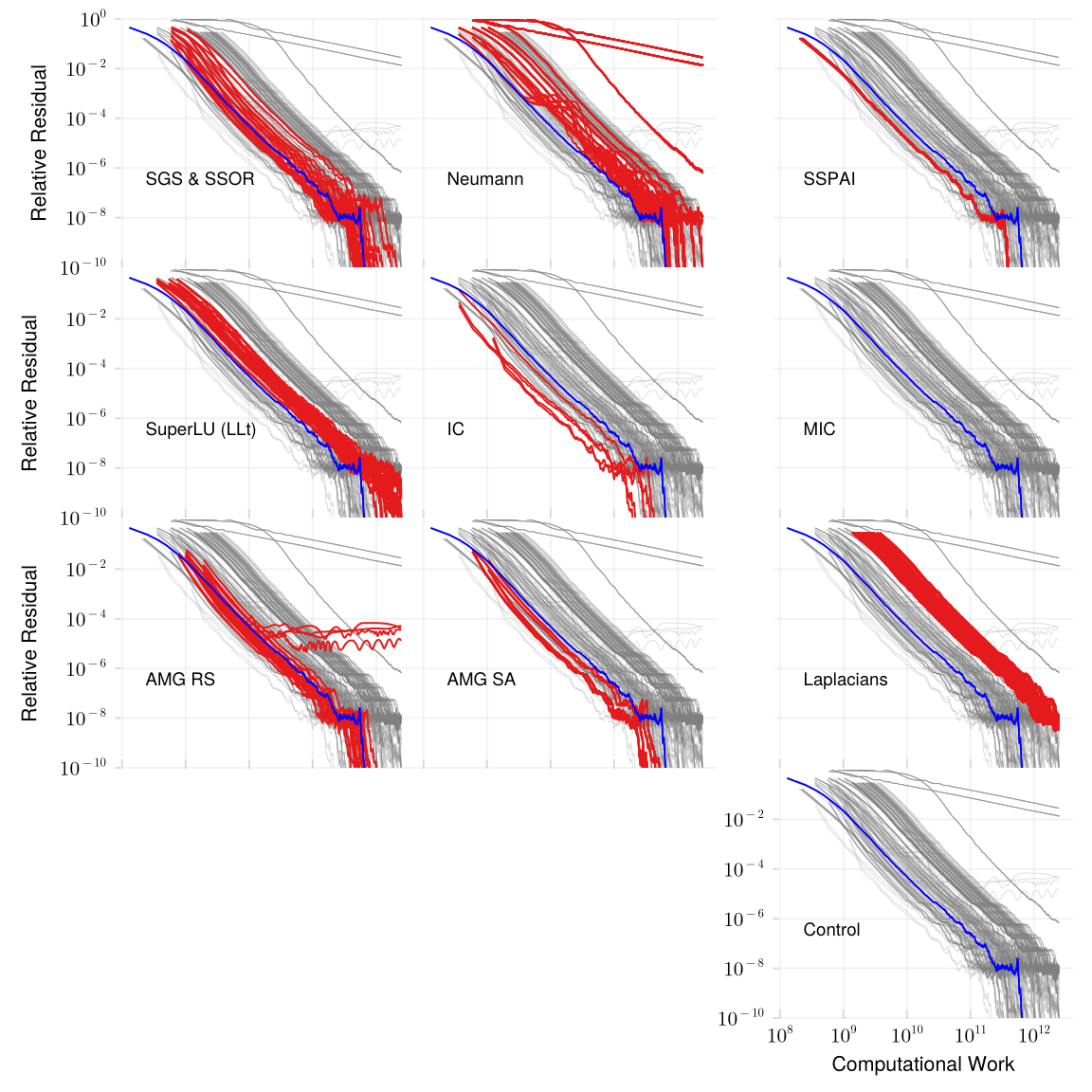}
    \caption{Convergence of the PCG method with various preconditioners applied to the \texttt{Flan\_1565} matrix (1.6m rows, 114.2m non-zeros). The plots have a log-log scale.}
    \label{fig:Flan_1565}
\end{figure}
\clearpage 

\subsection{G2\_circuit}

\begin{figure}[!ht]
    \centering
    \includegraphics[width=\textwidth]{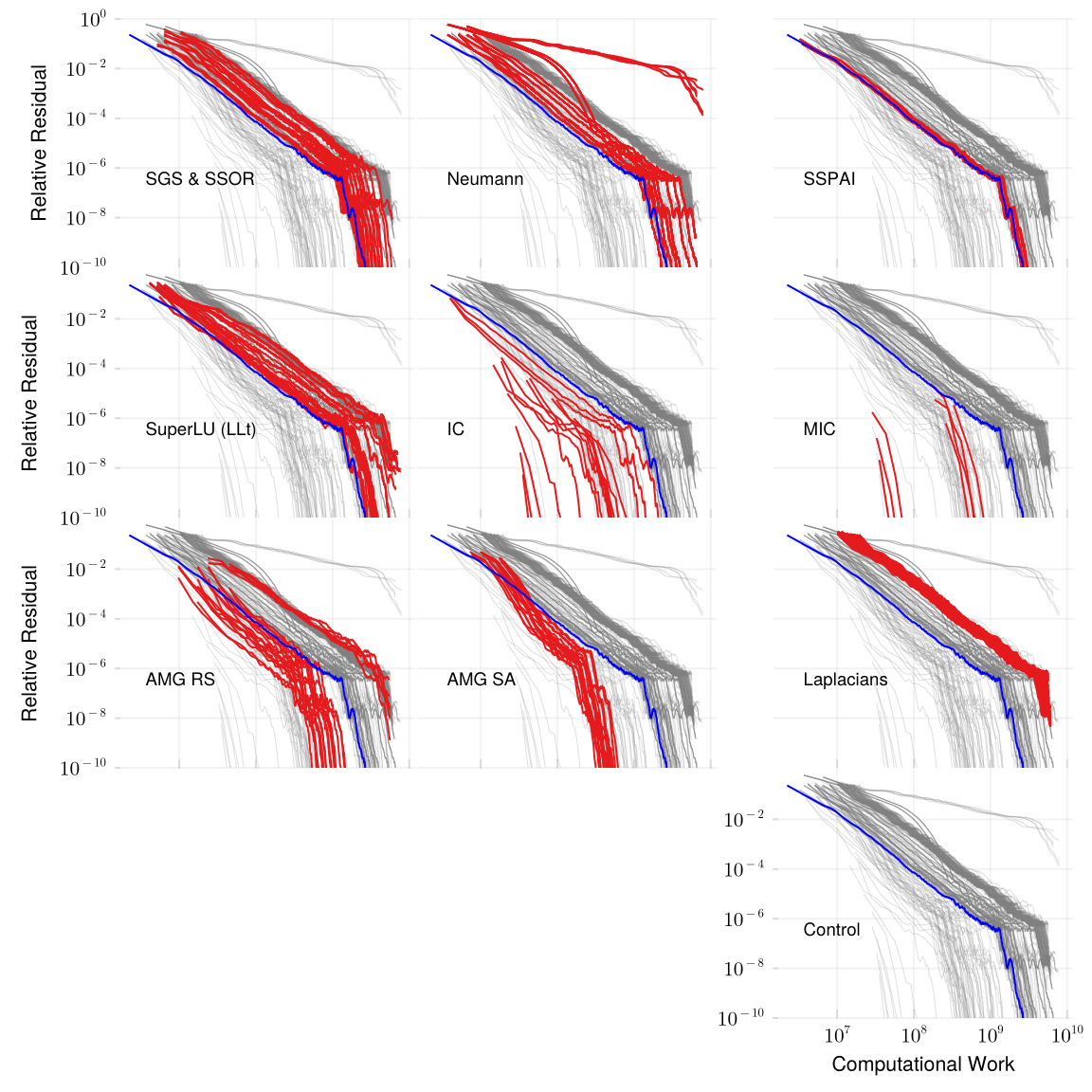}
    \caption{Convergence of the PCG method with various preconditioners applied to the \texttt{G2\_circuit} matrix (150k rows, 727k non-zeros). The plots have a log-log scale.}
    \label{fig:G2_circuit}
\end{figure}
\clearpage 

\subsection{G3\_circuit}

\begin{figure}[!ht]
    \centering
    \includegraphics[width=\textwidth]{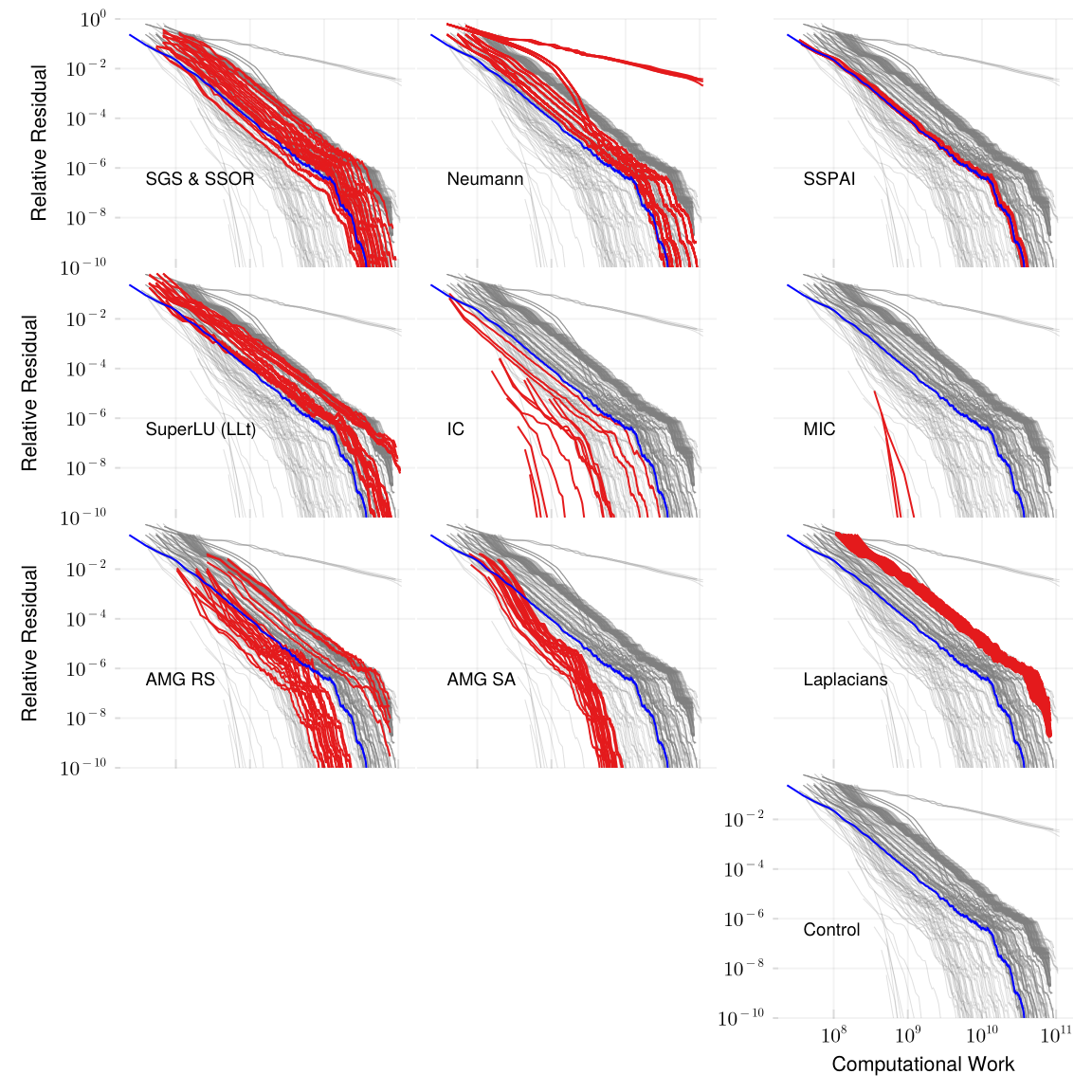}
    \caption{Convergence of the PCG method with various preconditioners applied to the \texttt{G3\_circuit} matrix (1.6m rows, 7.7m non-zeros). The plots have a log-log scale.}
    \label{fig:G3_circuit}
\end{figure}
\clearpage 

\subsection{Geo\_1438}

\begin{figure}[!ht]
    \centering
    \includegraphics[width=\textwidth]{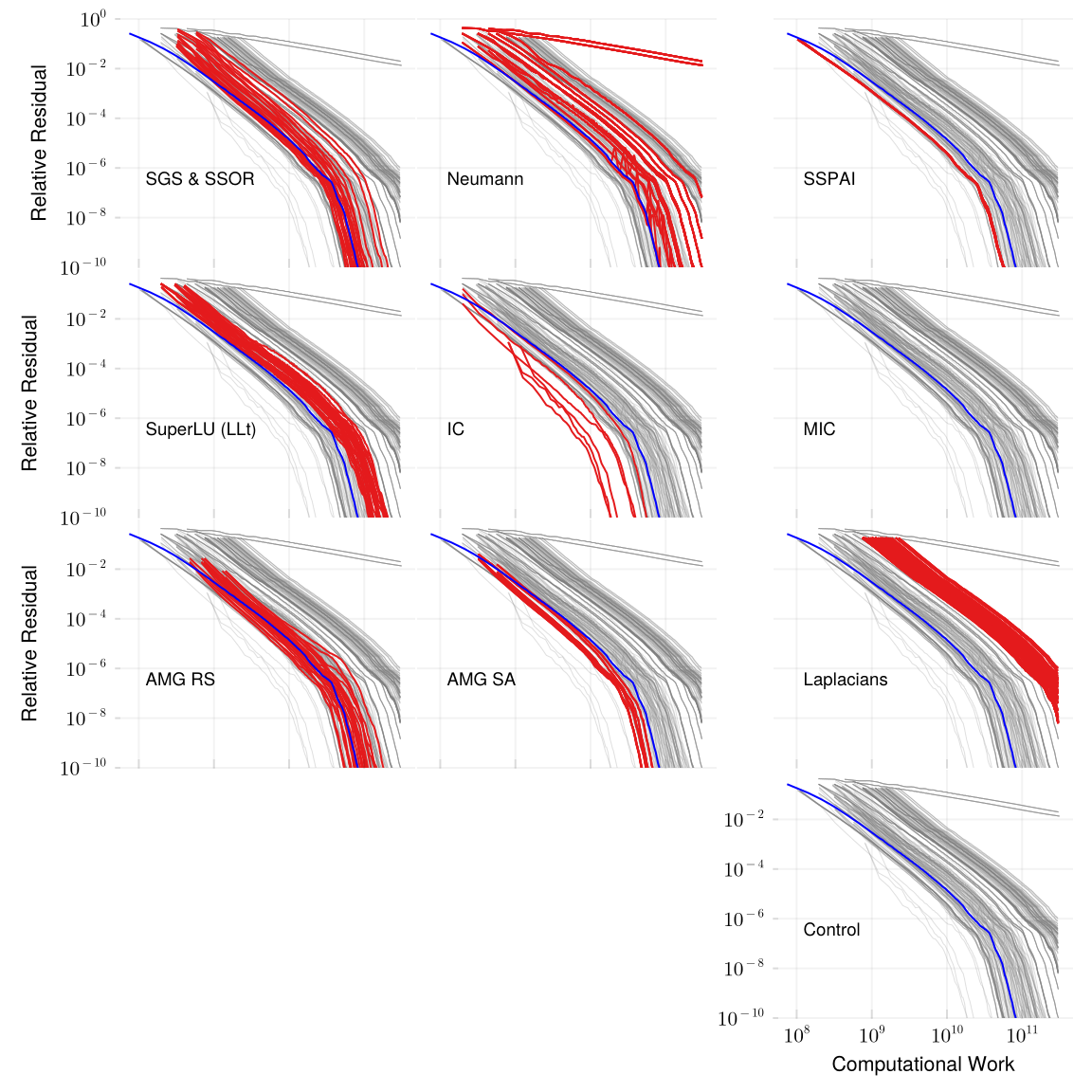}
    \caption{Convergence of the PCG method with various preconditioners applied to the \texttt{Geo\_1438} matrix (1.4m rows, 60.2m non-zeros). The plots have a log-log scale.}
    \label{fig:Geo_1438}
\end{figure}
\clearpage 

\subsection{Hook\_1498}
\begin{figure}[!ht]
    \centering
    \includegraphics[width=\textwidth]{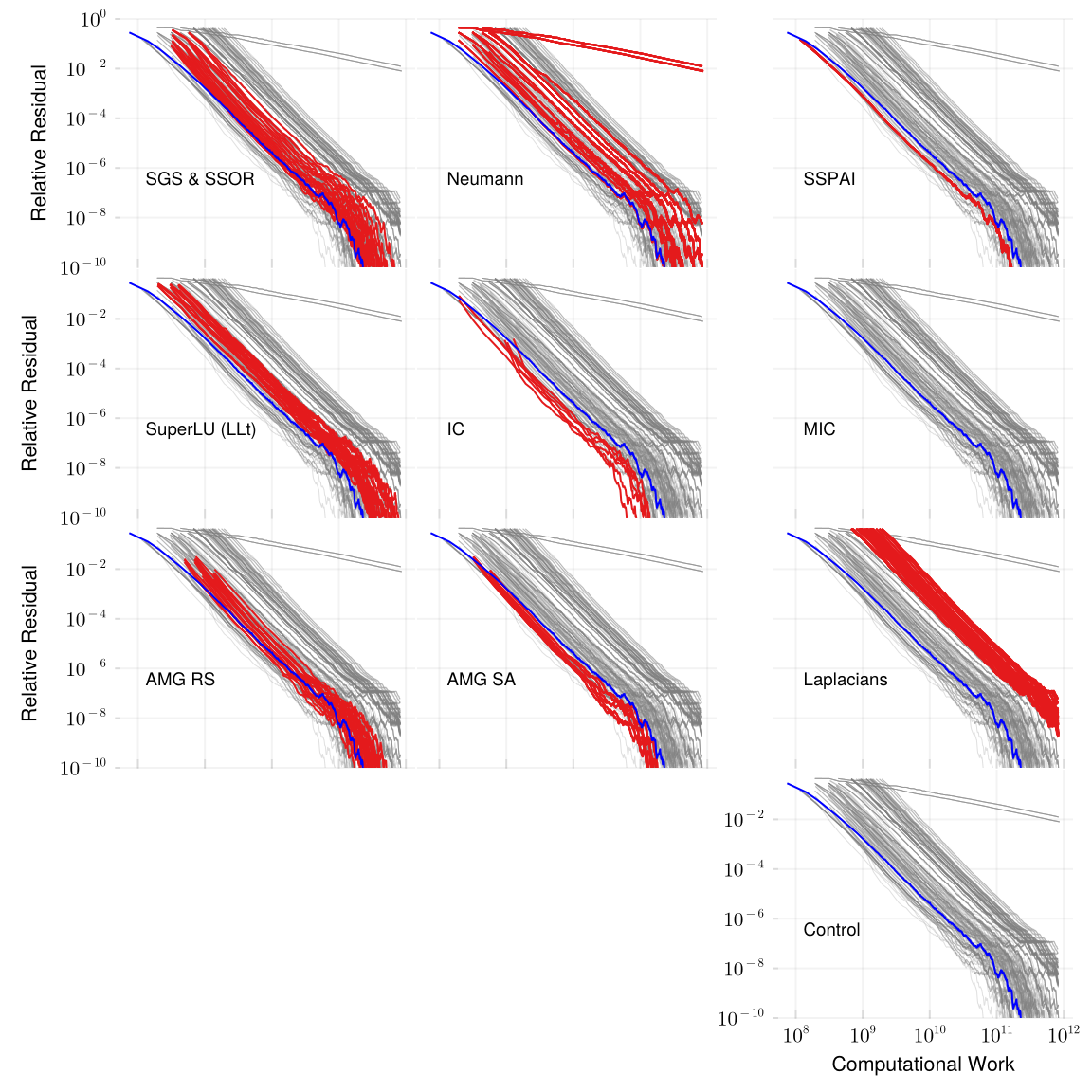}
    \caption{Convergence of the PCG method with various preconditioners applied to the \texttt{Hook\_1498} matrix (1.5m rows, 59.4m non-zeros). The plots have a log-log scale.}
    \label{fig:Hook_1498}
\end{figure}
\clearpage 

\subsection{PFlow\_742}

\begin{figure}[!ht]
    \centering
    \includegraphics[width=\textwidth]{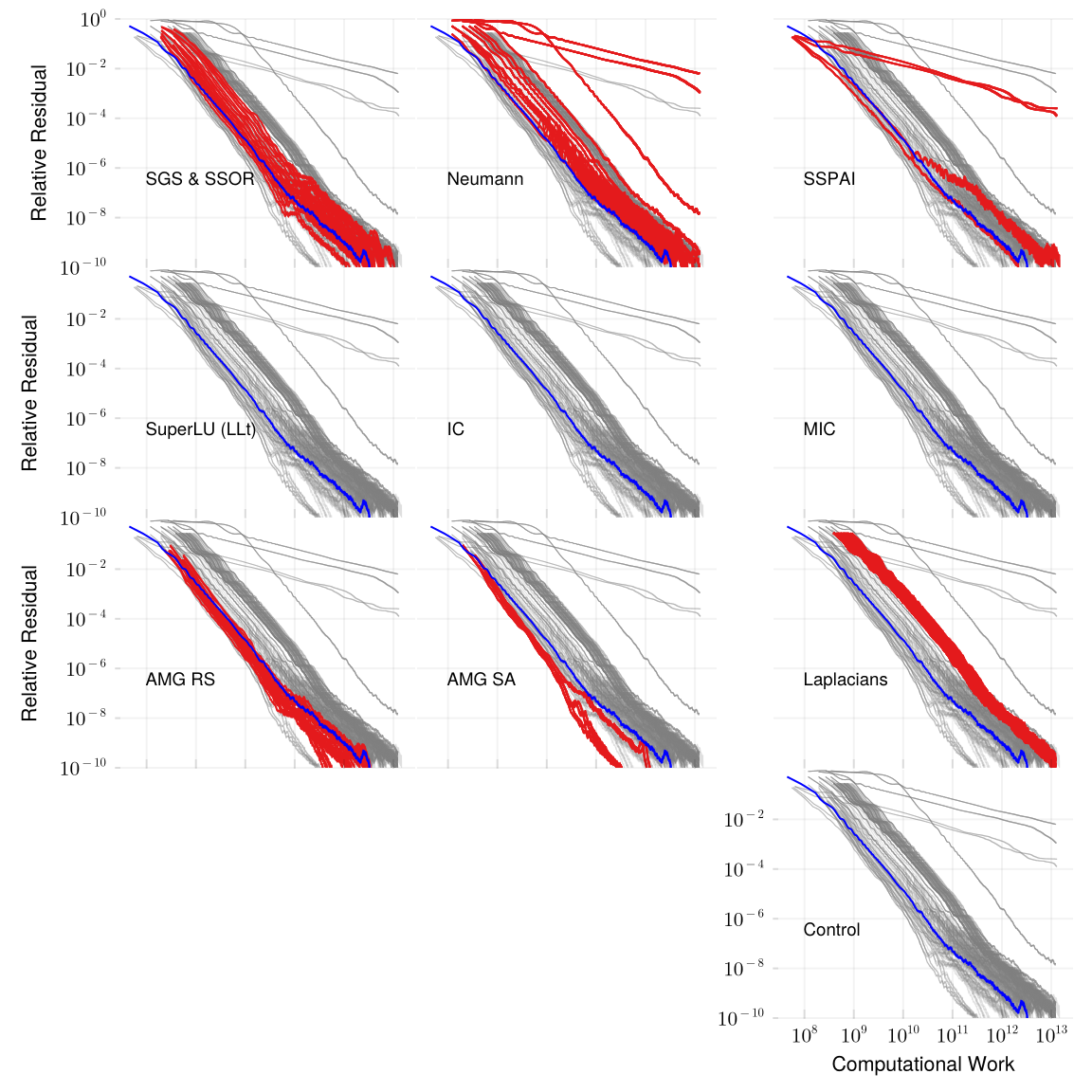}
    \caption{Convergence of the PCG method with various preconditioners applied to the \texttt{PFlow\_742} matrix (743k rows, 37.1m non-zeros). The plots have a log-log scale.}
    \label{fig:PFlow_742}
\end{figure}
\clearpage 

\subsection{Pres\_Poisson}

\begin{figure}[!ht]
    \centering
    \includegraphics[width=\textwidth]{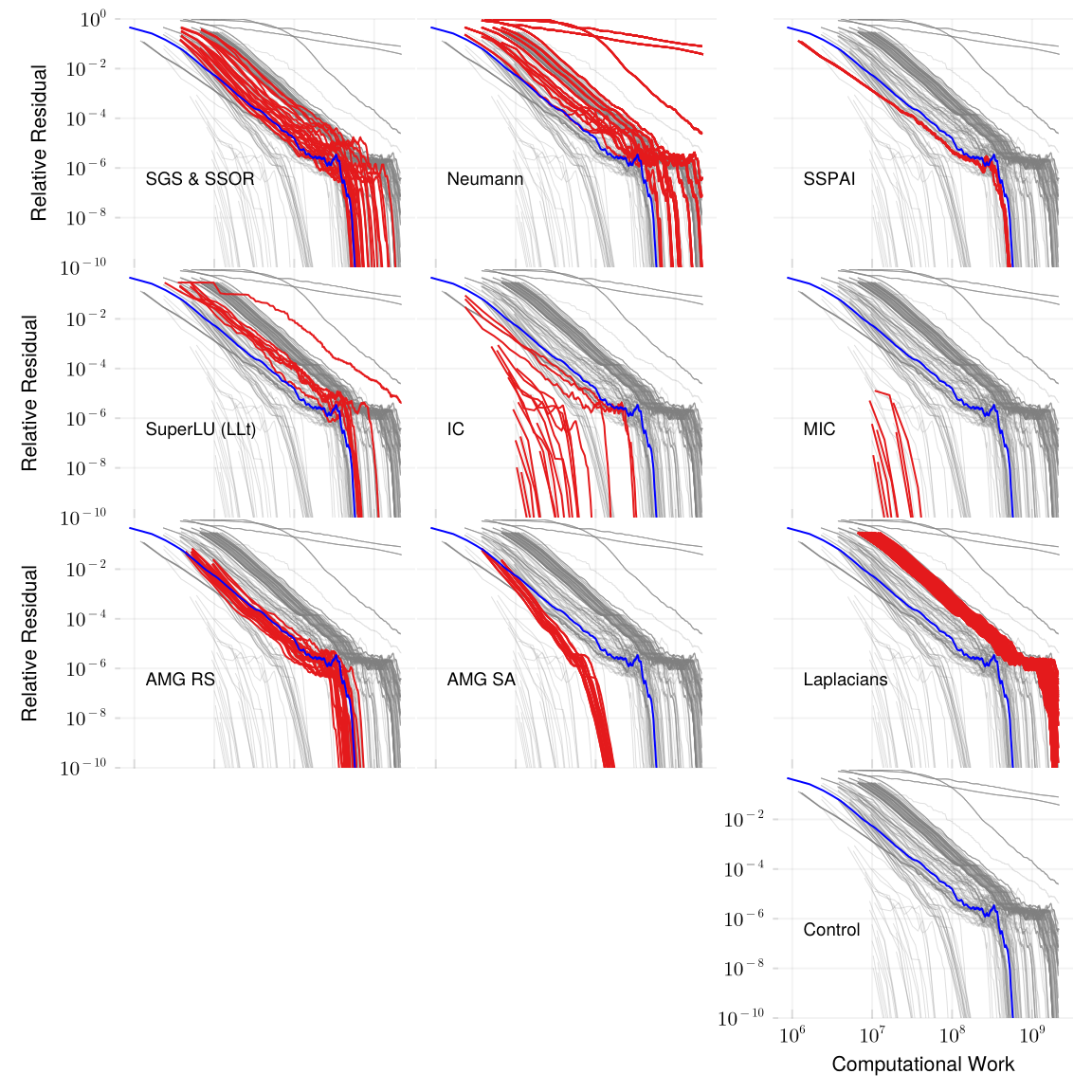}
    \caption{Convergence of the PCG method with various preconditioners applied to the \texttt{Pres\_Poisson} matrix (15k rows, 716k non-zeros). The plots have a log-log scale.}
    \label{fig:Pres_Poisson}
\end{figure}
\clearpage 

\subsection{Queen\_4147}

\begin{figure}[!ht]
    \centering
    \includegraphics[width=\textwidth]{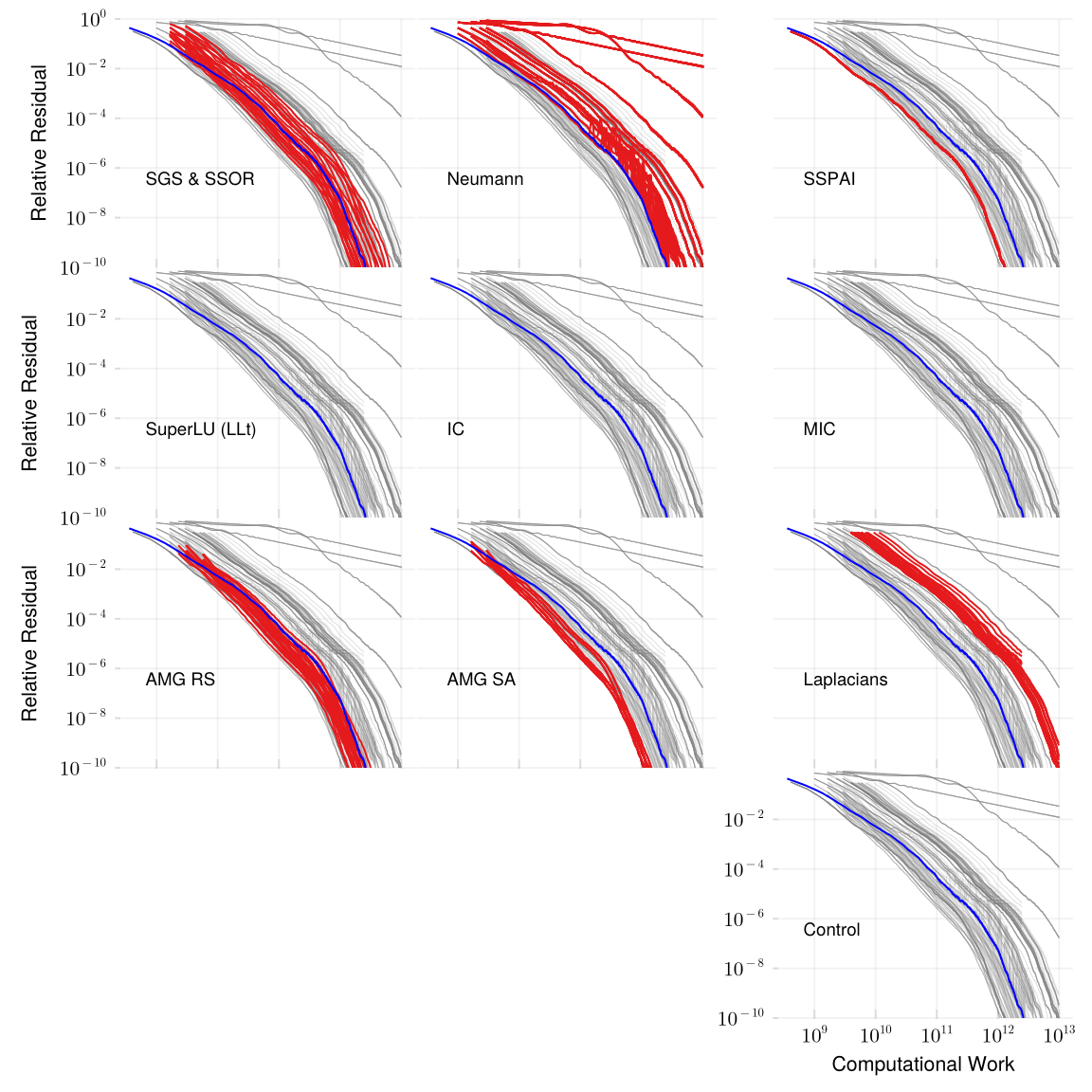}
    \caption{Convergence of the PCG method with various preconditioners applied to the \texttt{Queen\_4147} matrix (4.1m rows, 316.5m non-zeros). The plots have a log-log scale.}
    \label{fig:Queen_4147}
\end{figure}
\clearpage 

\subsection{Serena}

\begin{figure}[!ht]
    \centering
    \includegraphics[width=\textwidth]{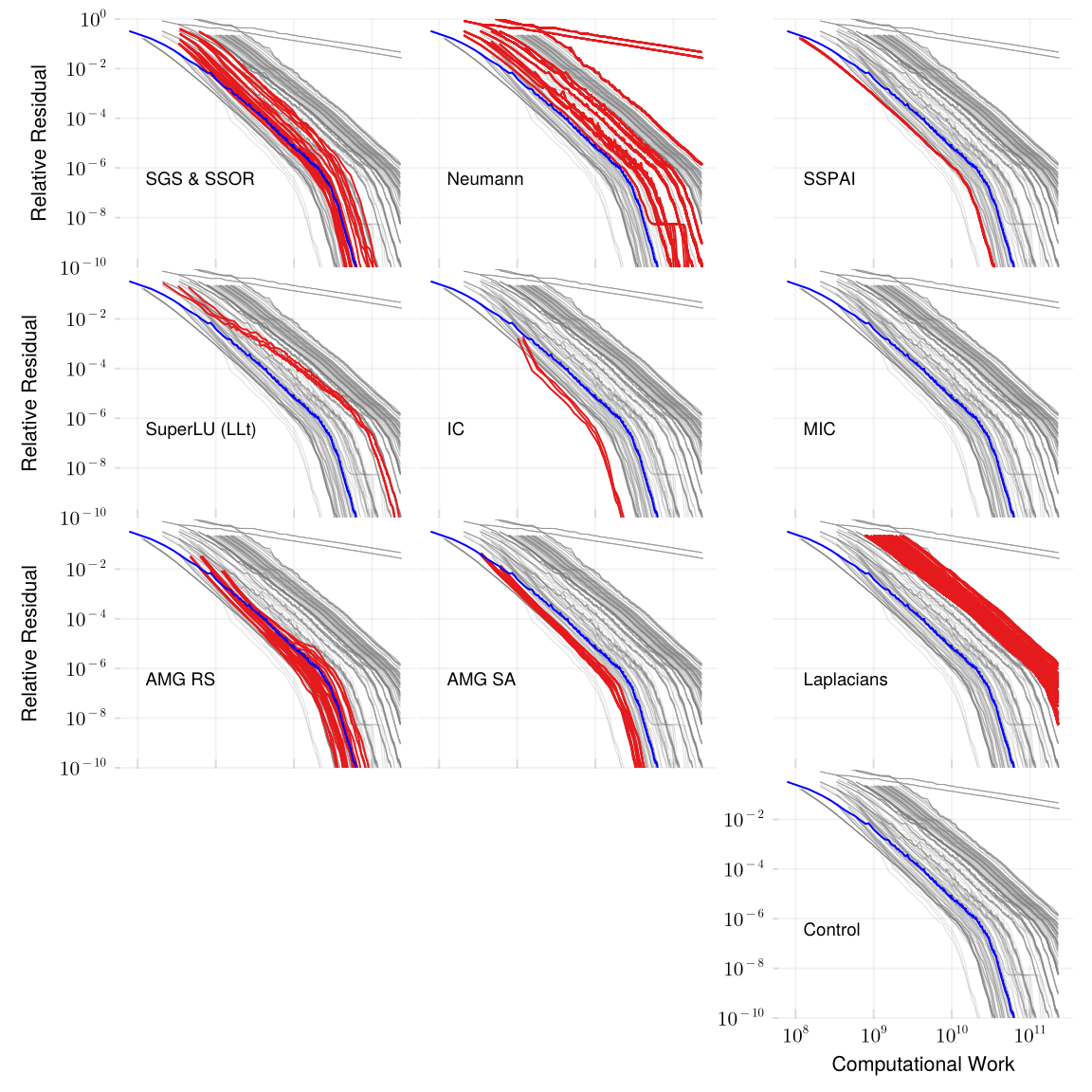}
    \caption{Convergence of the PCG method with various preconditioners applied to the \texttt{Serena} matrix (1.4m rows, 64.1m non-zeros). The plots have a log-log scale.}
    \label{fig:Serena}
\end{figure}
\clearpage 

\subsection{StocF-1465}

\begin{figure}[!ht]
    \centering
    \includegraphics[width=\textwidth]{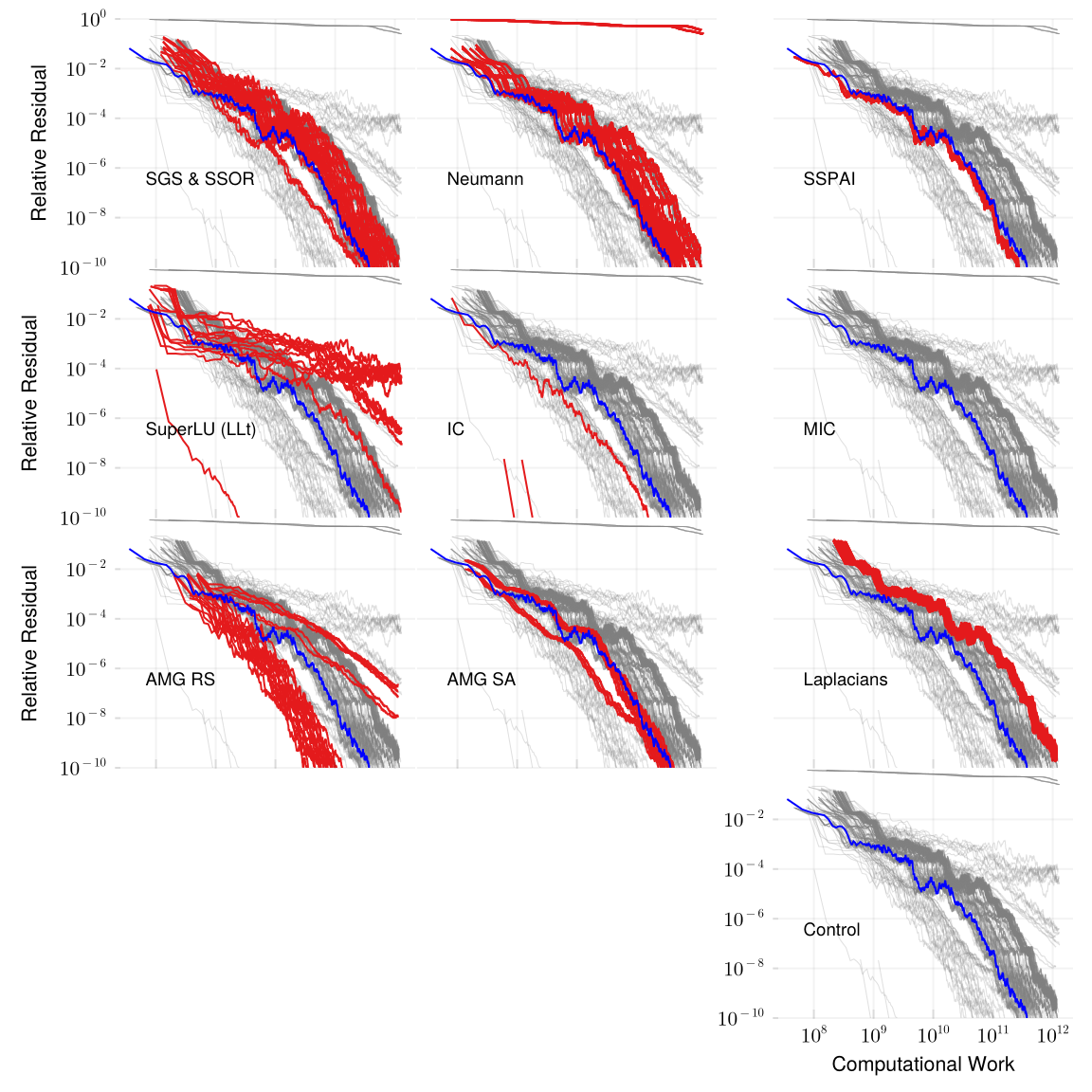}
    \caption{Convergence of the PCG method with various preconditioners applied to the \texttt{StocF-1465} matrix (1.5m rows, 21m non-zeros). The plots have a log-log scale.}
    \label{fig:StocF-1465}
\end{figure}
\clearpage 

\subsection{Trefethen\_20000}
\begin{figure}[!ht]
    \centering
    \includegraphics[width=\textwidth]{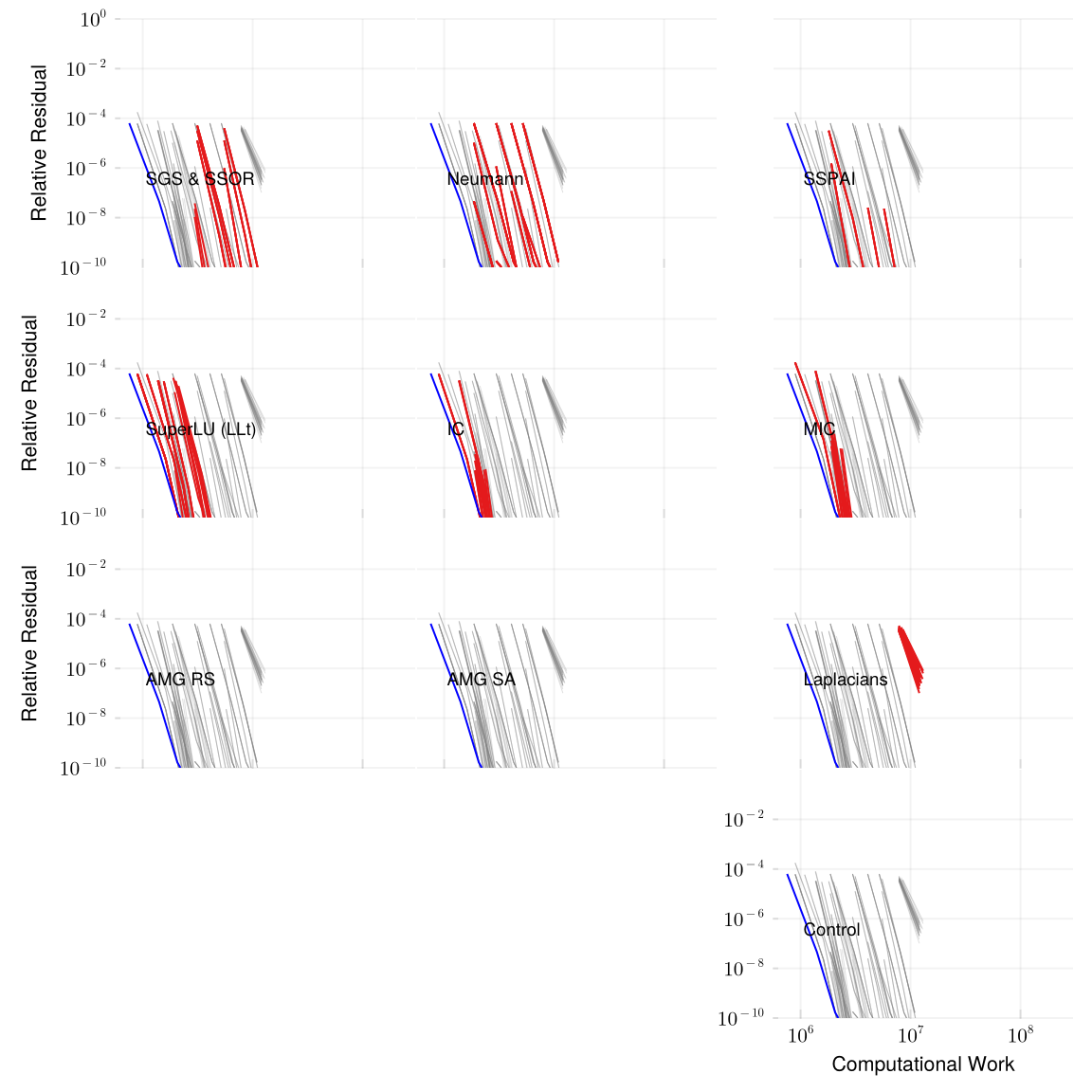}
    \caption{Convergence of the PCG method with various preconditioners applied to the \texttt{Trefethen\_20000} matrix (20k rows, 554k non-zeros). The plots have a log-log scale.}
    \label{fig:Trefethen_20000}
\end{figure}
\clearpage 

\subsection{Trefethen\_20000b}
\begin{figure}[!ht]
    \centering
    \includegraphics[width=\textwidth]{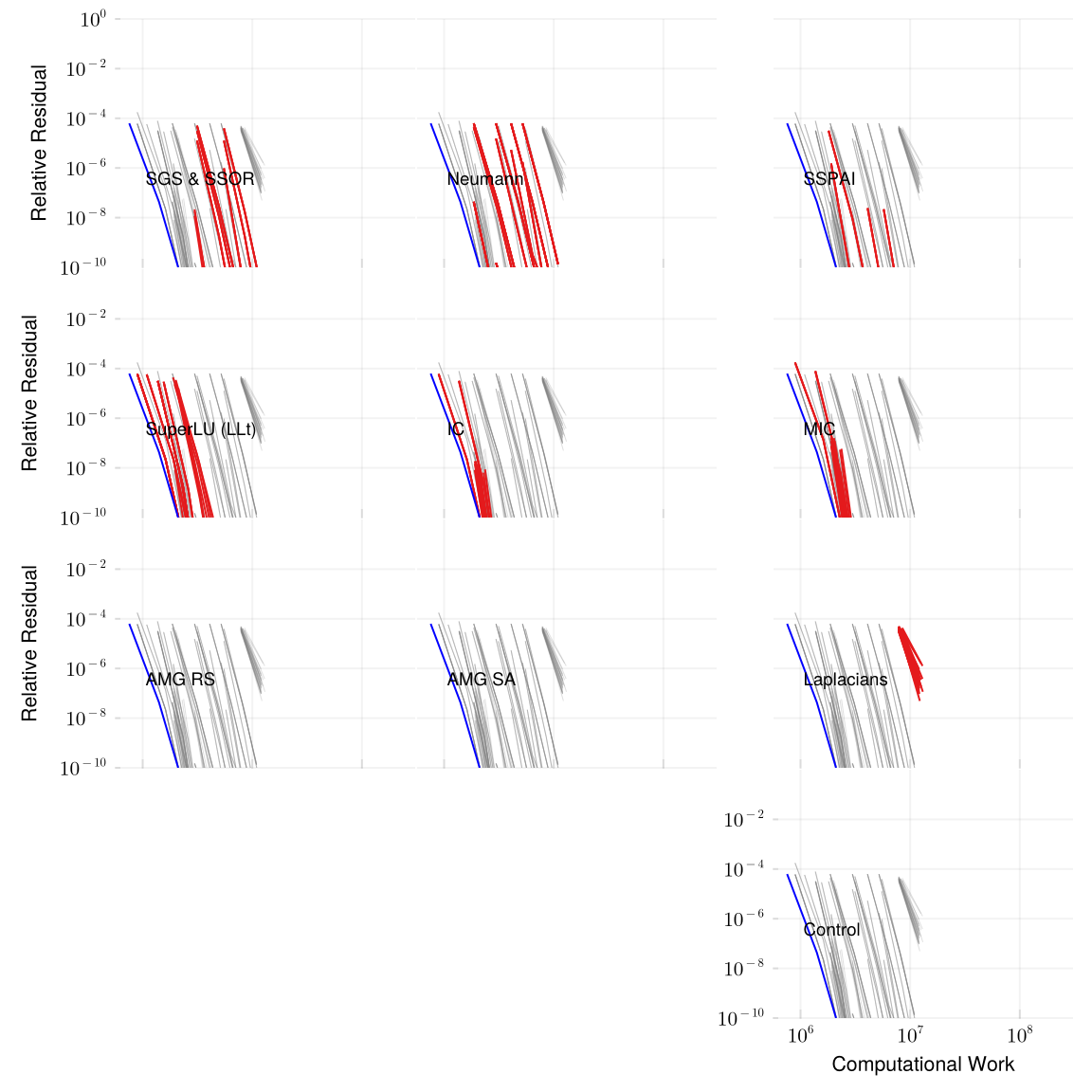}
    \caption{Convergence of the PCG method with various preconditioners applied to the \texttt{Trefethen\_20000b} matrix (20k rows, 554k non-zeros). The plots have a log-log scale.}
    \label{fig:Trefethen_20000b}
\end{figure}
\clearpage 

\subsection{af\_0\_k101}

\begin{figure}[!ht]
    \centering
    \includegraphics[width=\textwidth]{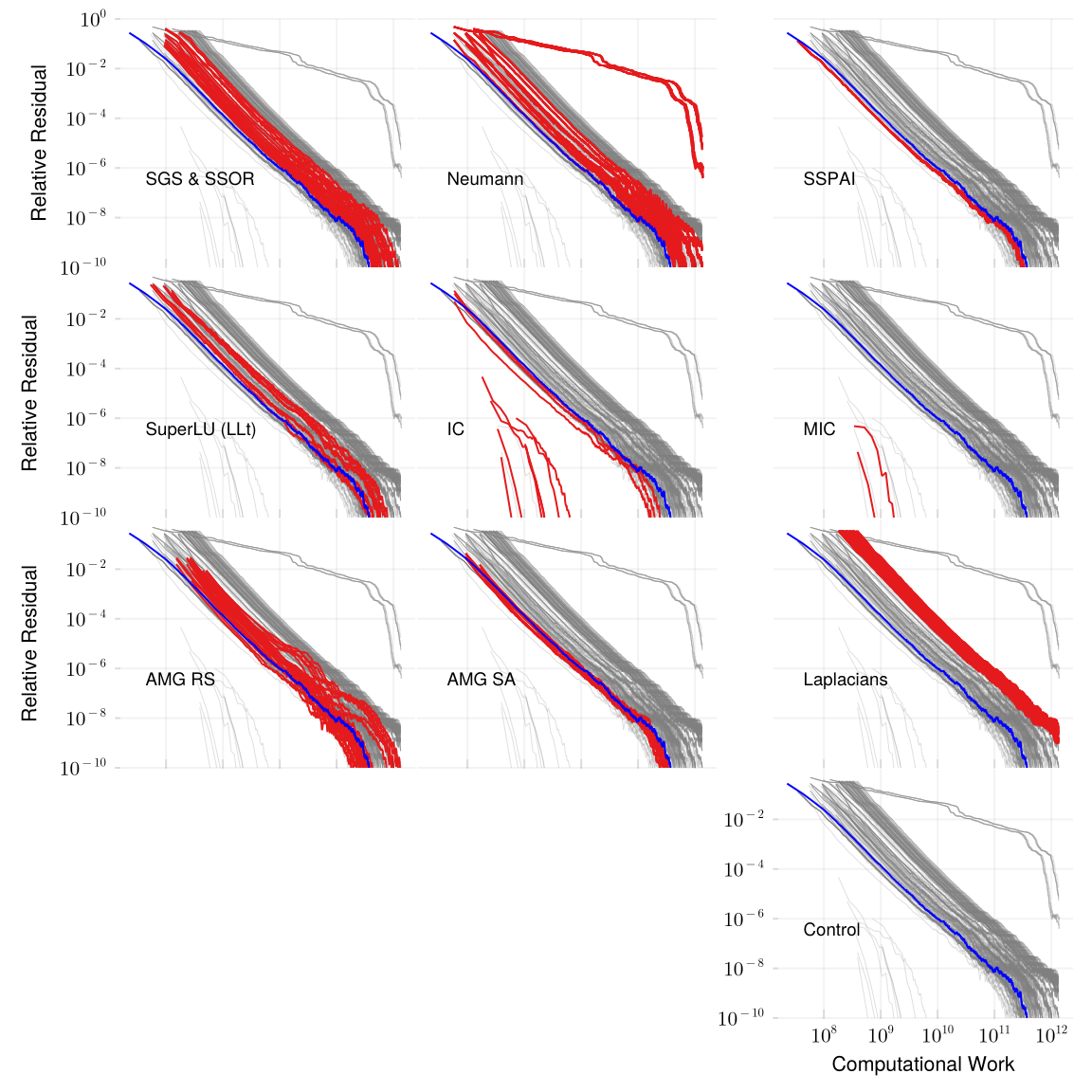}
    \caption{Convergence of the PCG method with various preconditioners applied to the \texttt{af\_0\_k101} matrix (504k rows, 17.6m non-zeros). The plots have a log-log scale.}
    \label{fig:af_0_k101}
\end{figure}
\clearpage 

\subsection{af\_1\_k101}

\begin{figure}[!ht]
    \centering
    \includegraphics[width=\textwidth]{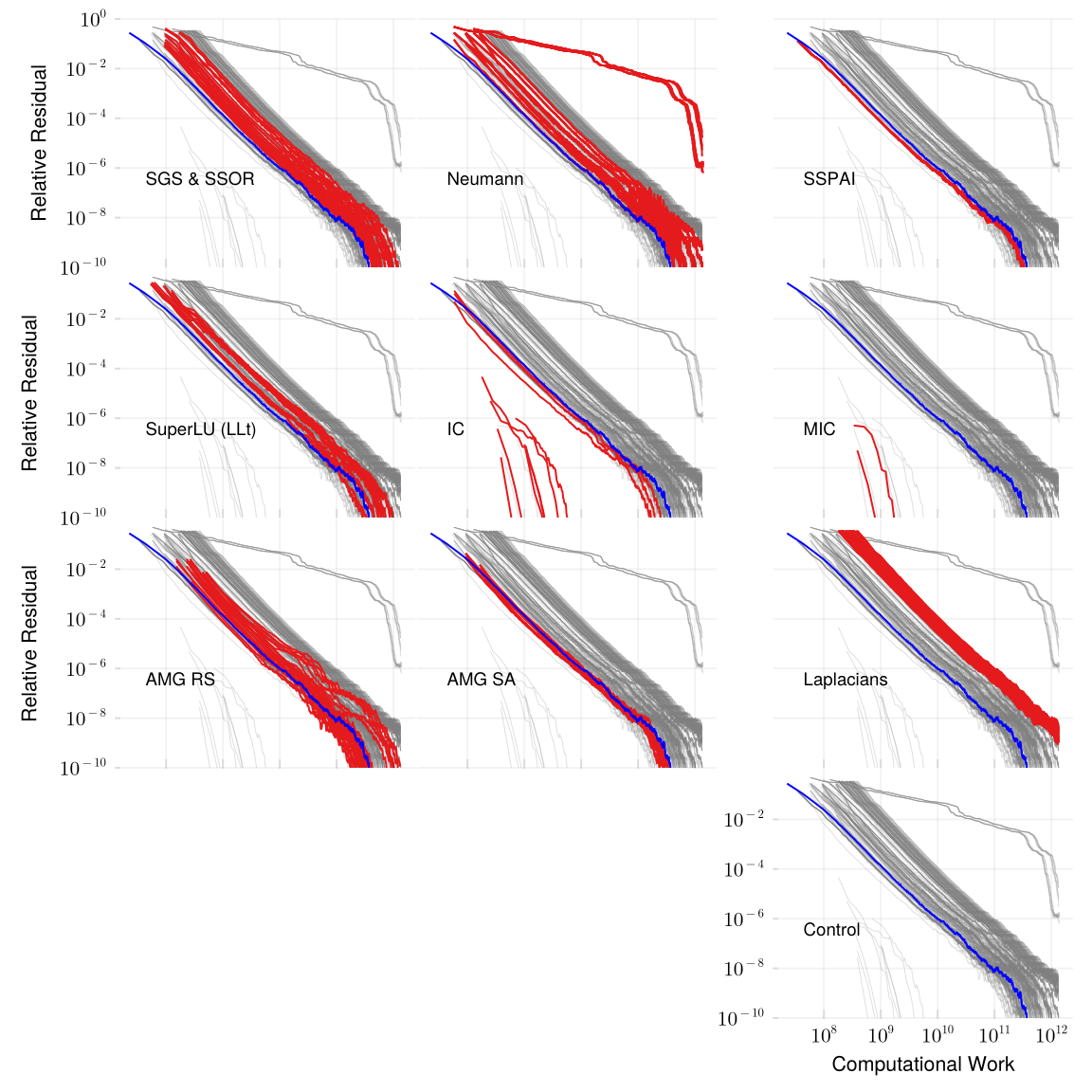}
    \caption{Convergence of the PCG method with various preconditioners applied to the \texttt{af\_1\_k101} matrix (504k rows, 17.6m non-zeros). The plots have a log-log scale.}
    \label{fig:af_1_k101}
\end{figure}
\clearpage 

\subsection{af\_2\_k101}

\begin{figure}[!ht]
    \centering
    \includegraphics[width=\textwidth]{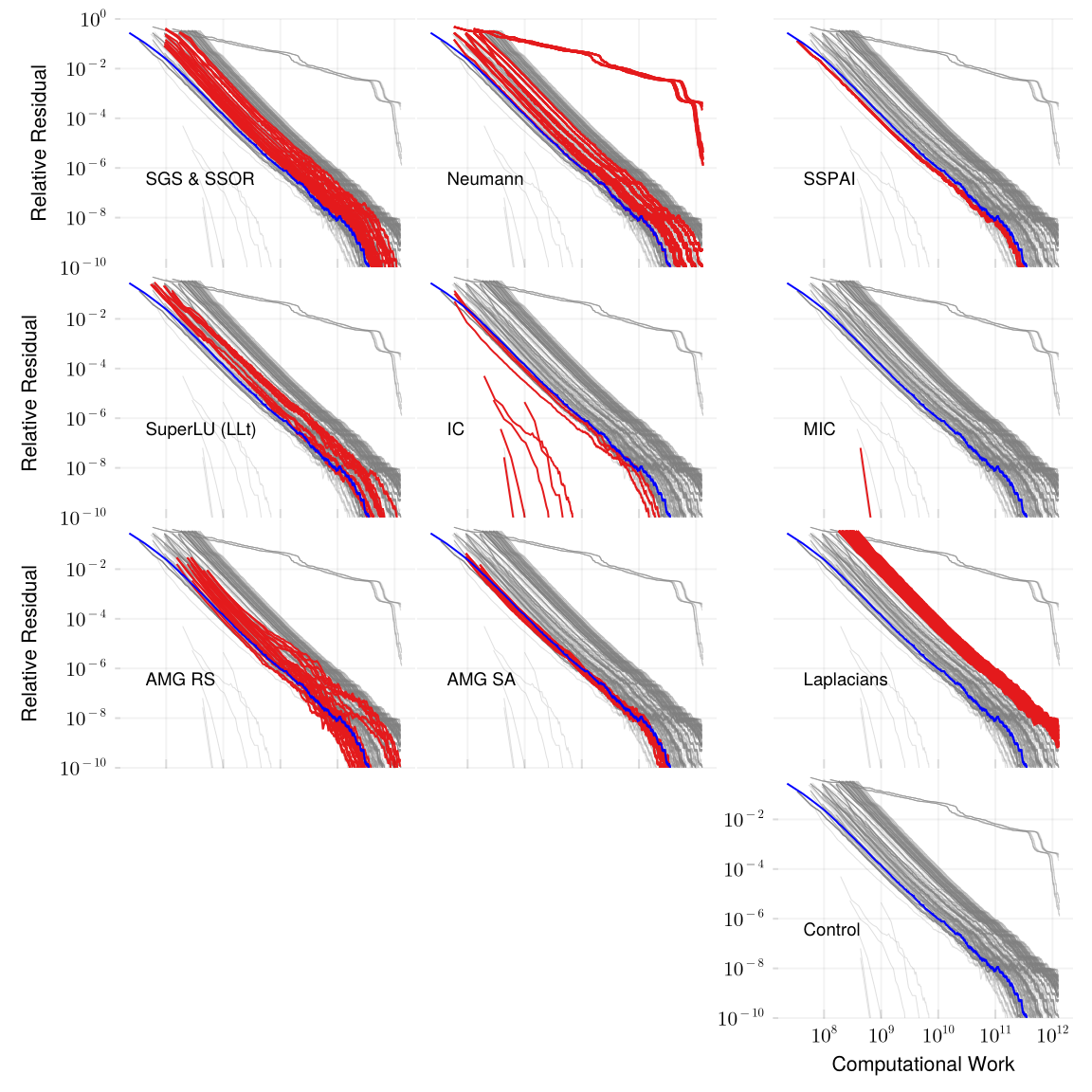}
    \caption{Convergence of the PCG method with various preconditioners applied to the \texttt{af\_2\_k101} matrix (504k rows, 17.6m non-zeros). The plots have a log-log scale.}
    \label{fig:af_2_k101}
\end{figure}
\clearpage 

\subsection{af\_3\_k101}

\begin{figure}[!ht]
    \centering
    \includegraphics[width=\textwidth]{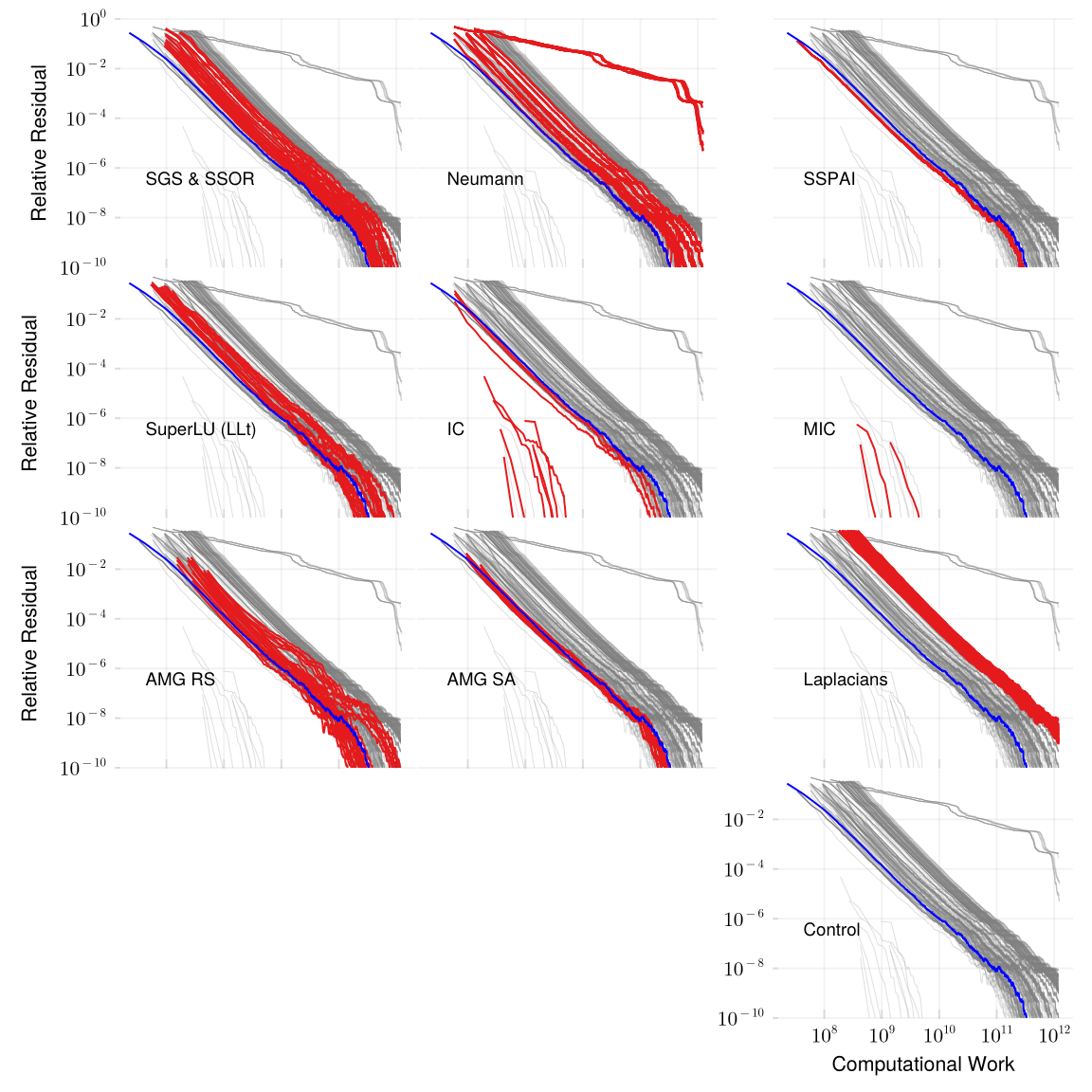}
    \caption{Convergence of the PCG method with various preconditioners applied to the \texttt{af\_3\_k101} matrix (504k rows, 17.6m non-zeros). The plots have a log-log scale.}
    \label{fig:af_3_k101}
\end{figure}
\clearpage 

\subsection{af\_4\_k101}

\begin{figure}[!ht]
    \centering
    \includegraphics[width=\textwidth]{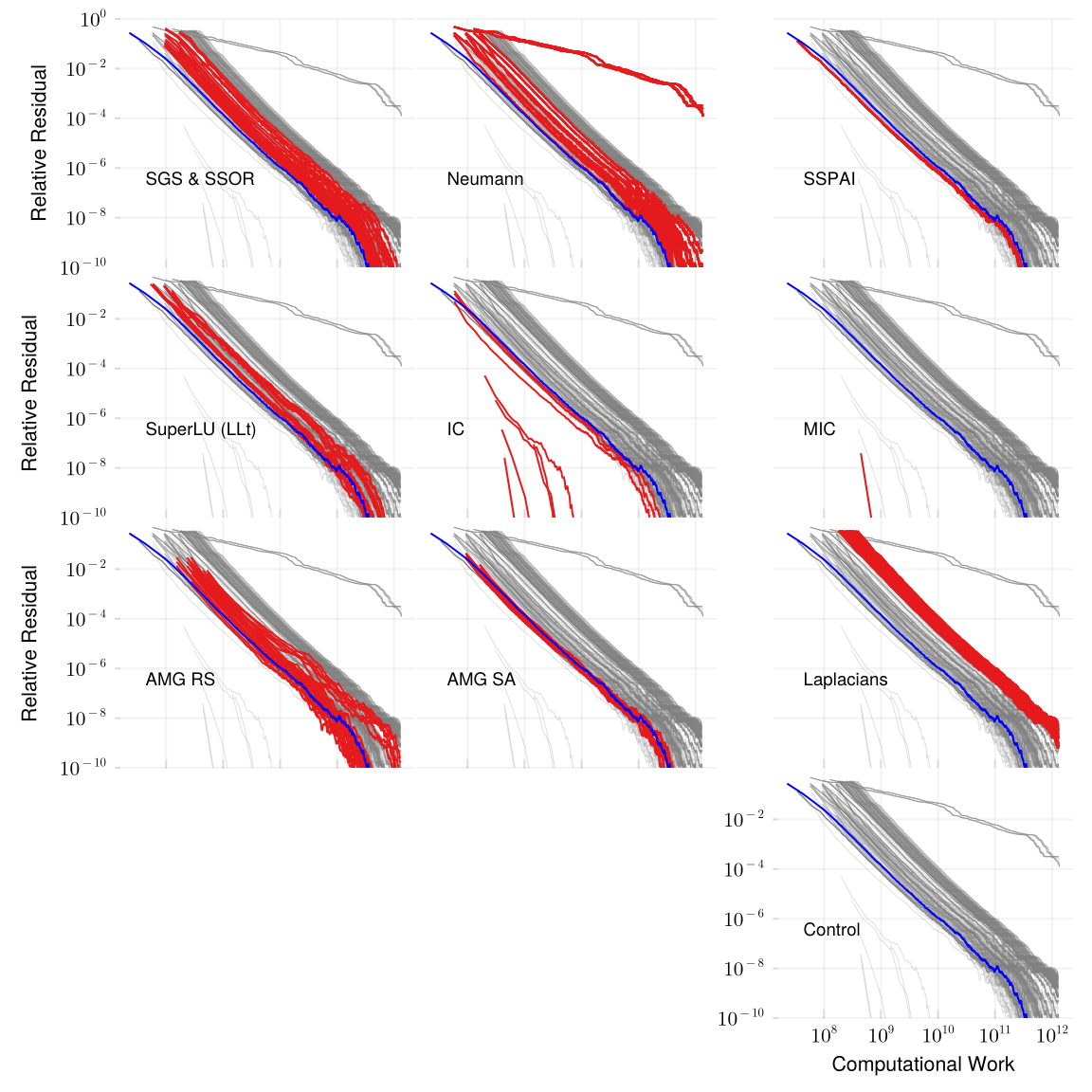}
    \caption{Convergence of the PCG method with various preconditioners applied to the \texttt{af\_4\_k101} matrix (504k rows, 17.6m non-zeros). The plots have a log-log scale.}
    \label{fig:af_4_k101}
\end{figure}
\clearpage 

\subsection{af\_5\_k101}

\begin{figure}[!ht]
    \centering
    \includegraphics[width=\textwidth]{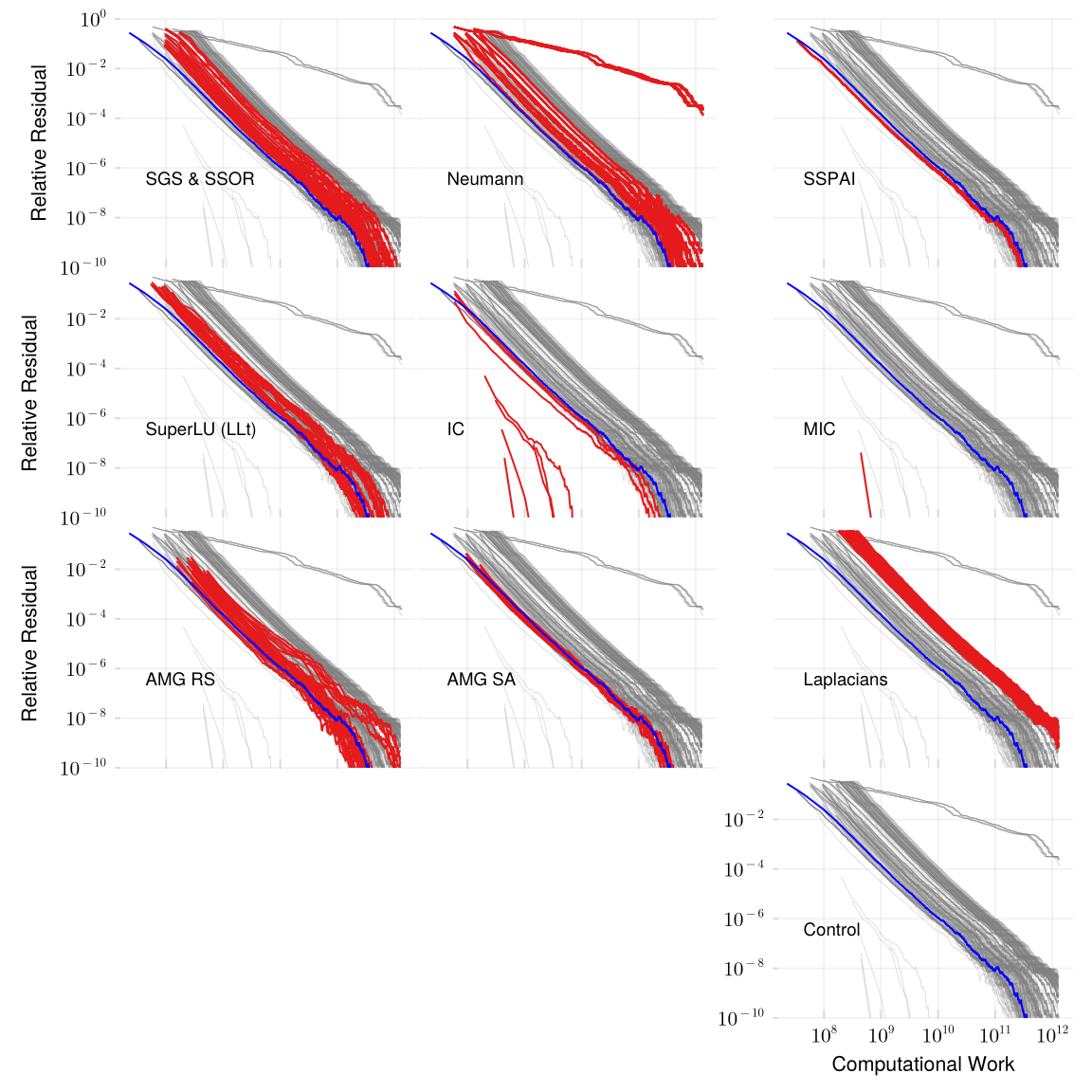}
    \caption{Convergence of the PCG method with various preconditioners applied to the \texttt{af\_5\_k101} matrix (504k rows, 17.6m non-zeros). The plots have a log-log scale.}
    \label{fig:af_5_k101}
\end{figure}
\clearpage 

\subsection{af\_shell3}

\begin{figure}[!ht]
    \centering
    \includegraphics[width=\textwidth]{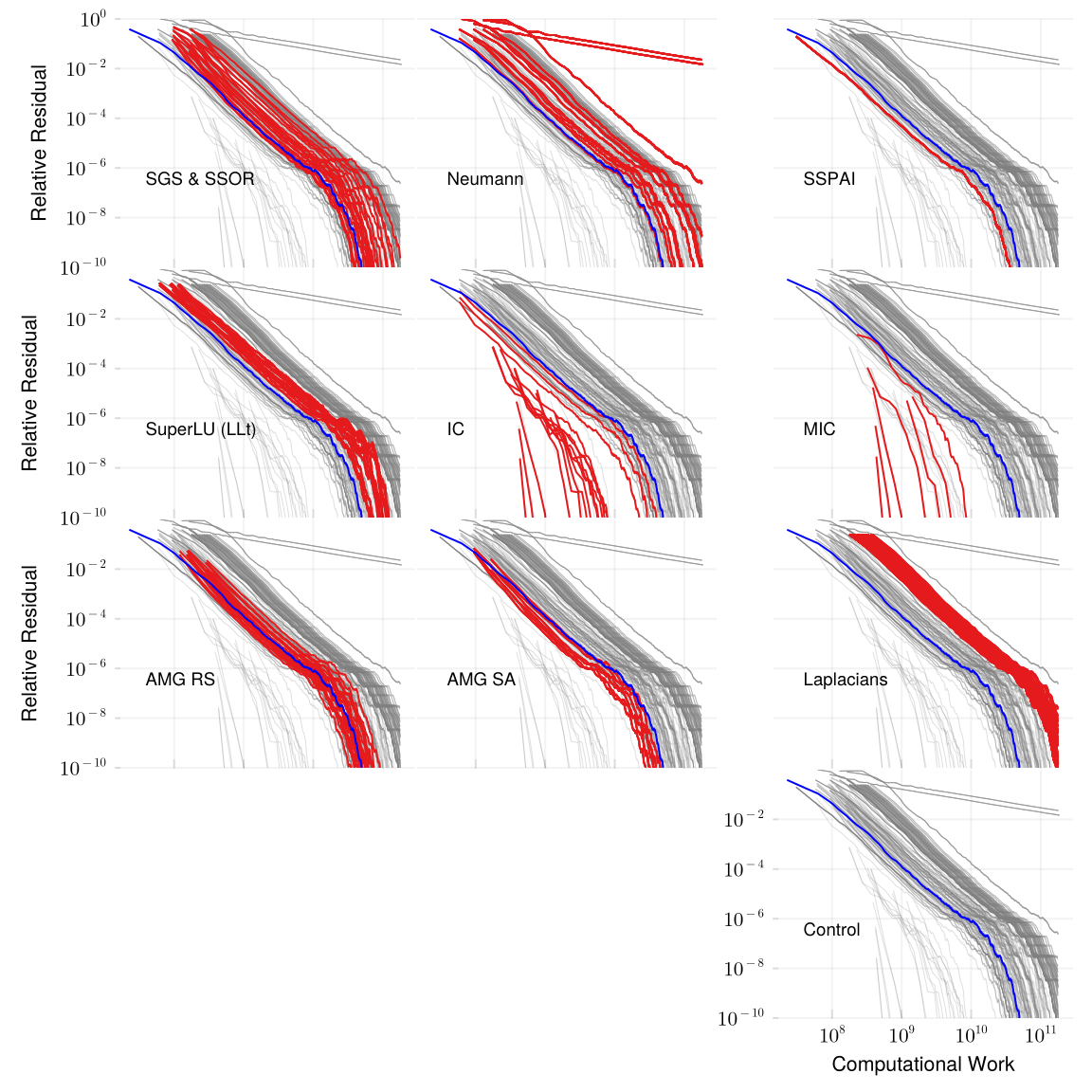}
    \caption{Convergence of the PCG method with various preconditioners applied to the \texttt{af\_shell3} matrix (505k rows, 17.6m non-zeros). The plots have a log-log scale.}
    \label{fig:af_shell3}
\end{figure}
\clearpage 

\subsection{af\_shell7}

\begin{figure}[!ht]
    \centering
    \includegraphics[width=\textwidth]{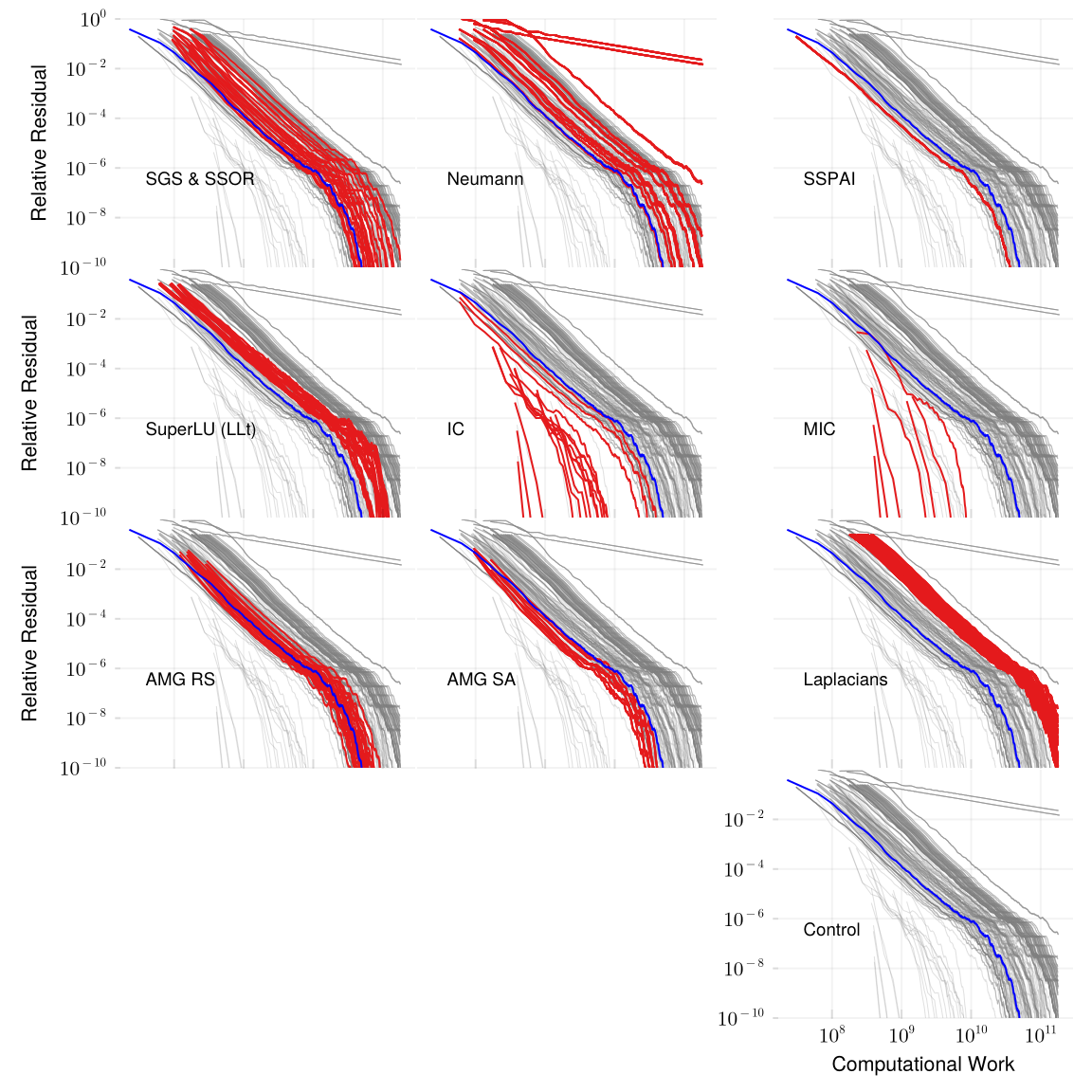}
    \caption{Convergence of the PCG method with various preconditioners applied to the \texttt{af\_shell7} matrix (505k rows, 17.6m non-zeros). The plots have a log-log scale.}
    \label{fig:af_shell7}
\end{figure}
\clearpage 

\subsection{audikw\_1}

\begin{figure}[!ht]
    \centering
    \includegraphics[width=\textwidth]{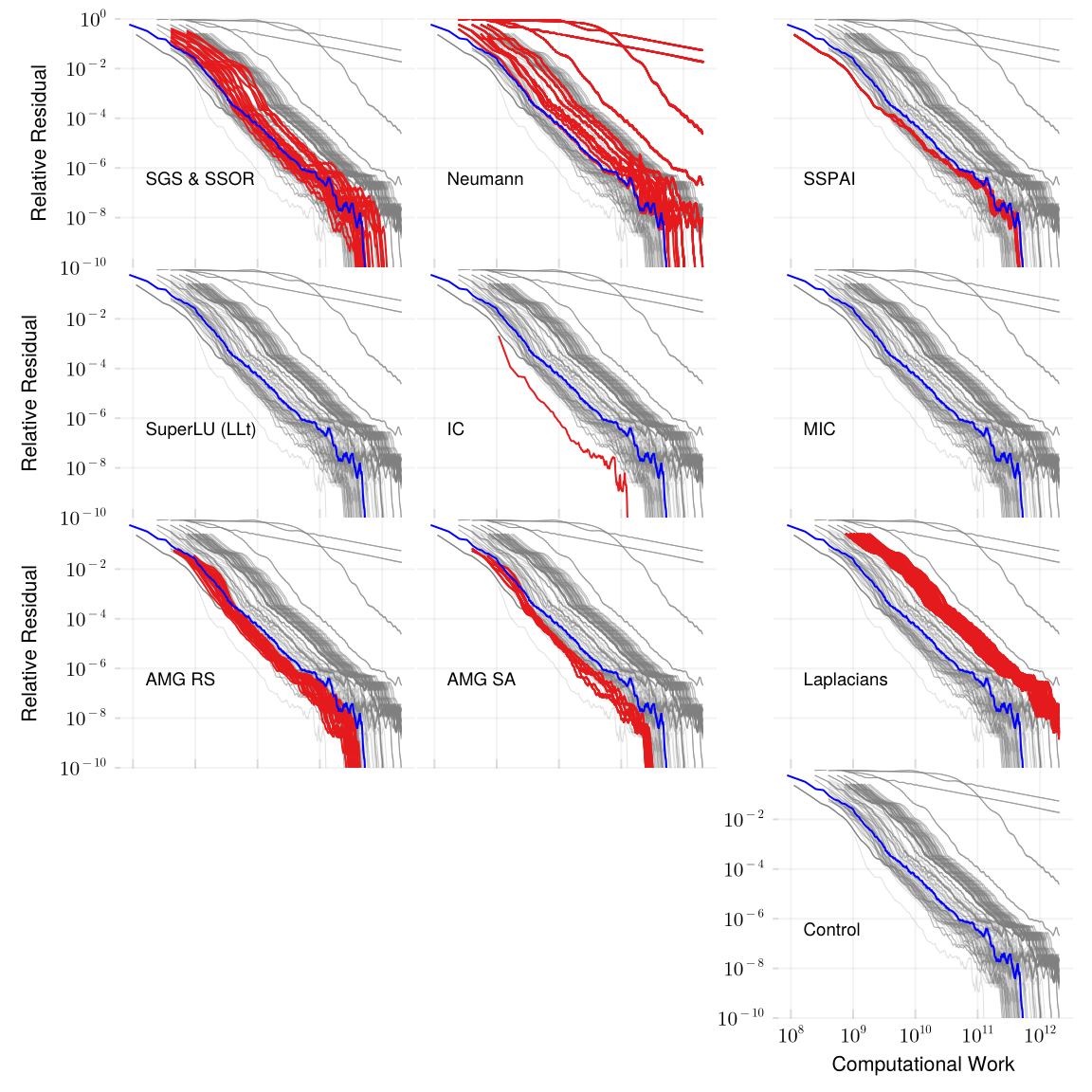}
    \caption{Convergence of the PCG method with various preconditioners applied to the \texttt{audikw\_1} matrix (944k rows, 77.7m non-zeros). The plots have a log-log scale.}
    \label{fig:audikw_1}
\end{figure}
\clearpage 

\subsection{bcsstk17}

\begin{figure}[!ht]
    \centering
    \includegraphics[width=\textwidth]{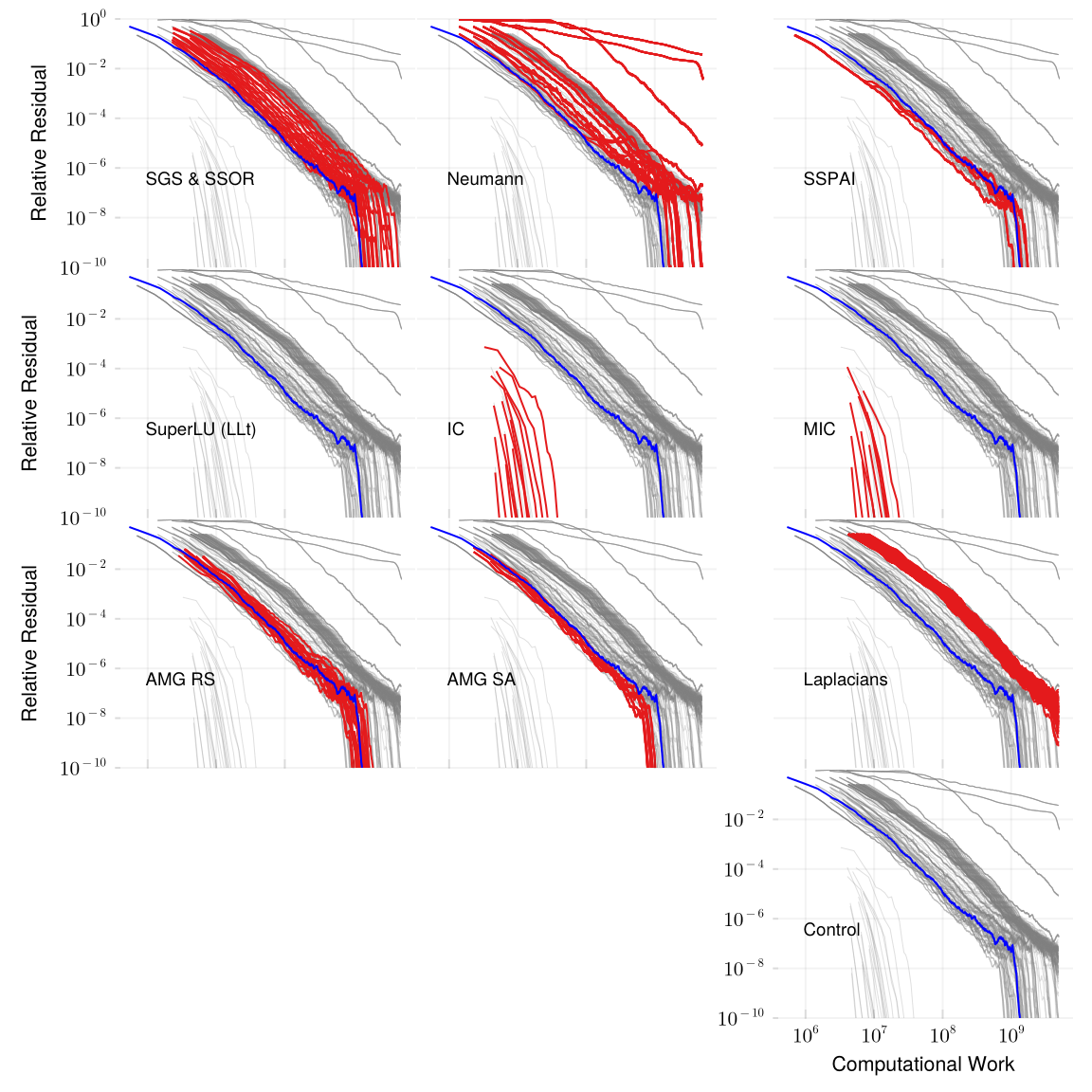}
    \caption{Convergence of the PCG method with various preconditioners applied to the \texttt{bcsstk17} matrix (11k rows, 429k non-zeros). The plots have a log-log scale.}
    \label{fig:bcsstk17}
\end{figure}
\clearpage 

\subsection{bcsstk18}

\begin{figure}[!ht]
    \centering
    \includegraphics[width=\textwidth]{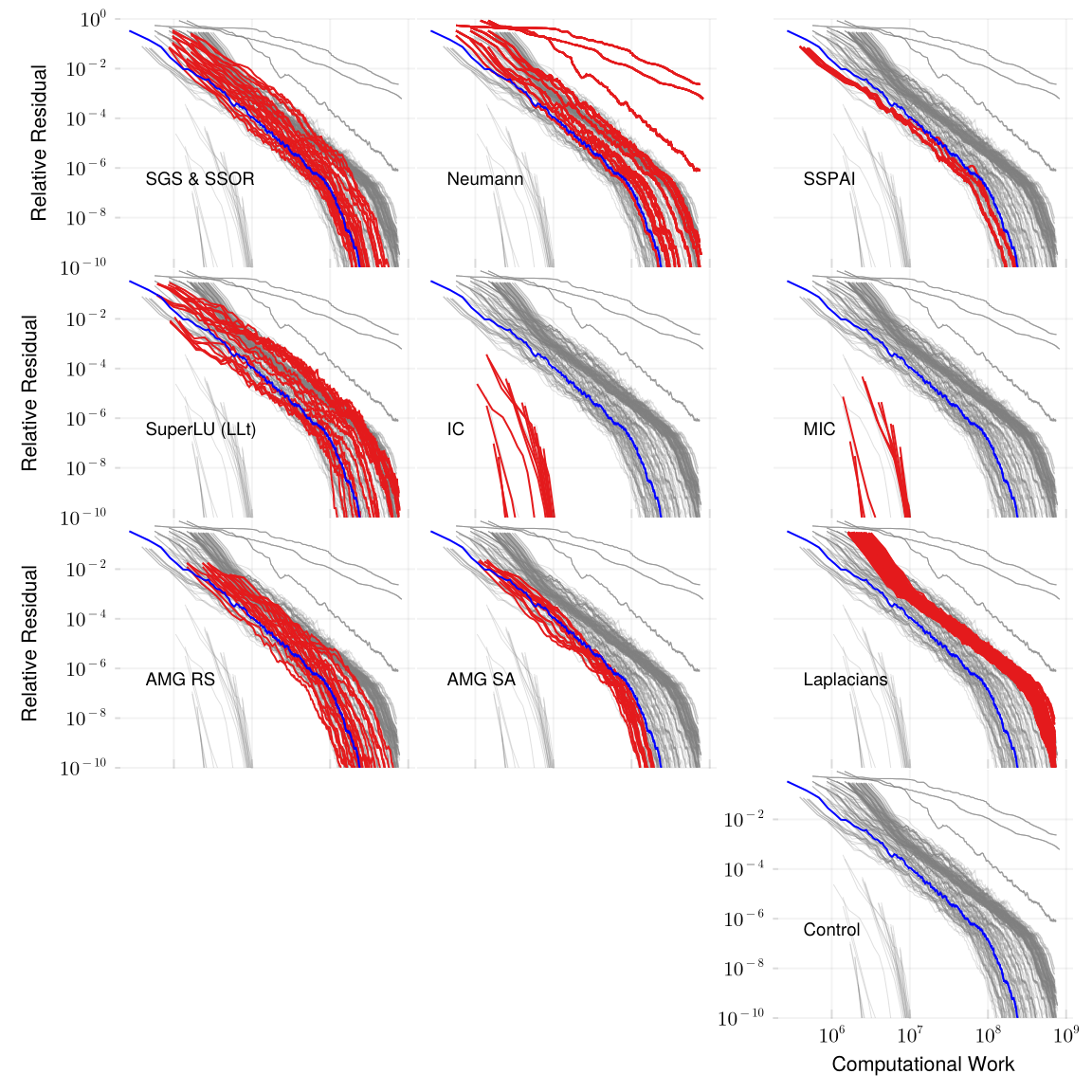}
    \caption{Convergence of the PCG method with various preconditioners applied to the \texttt{bcsstk18} matrix (12k rows, 149k non-zeros). The plots have a log-log scale.}
    \label{fig:bcsstk18}
\end{figure}
\clearpage 

\subsection{bodyy4}

\begin{figure}[!ht]
    \centering
    \includegraphics[width=\textwidth]{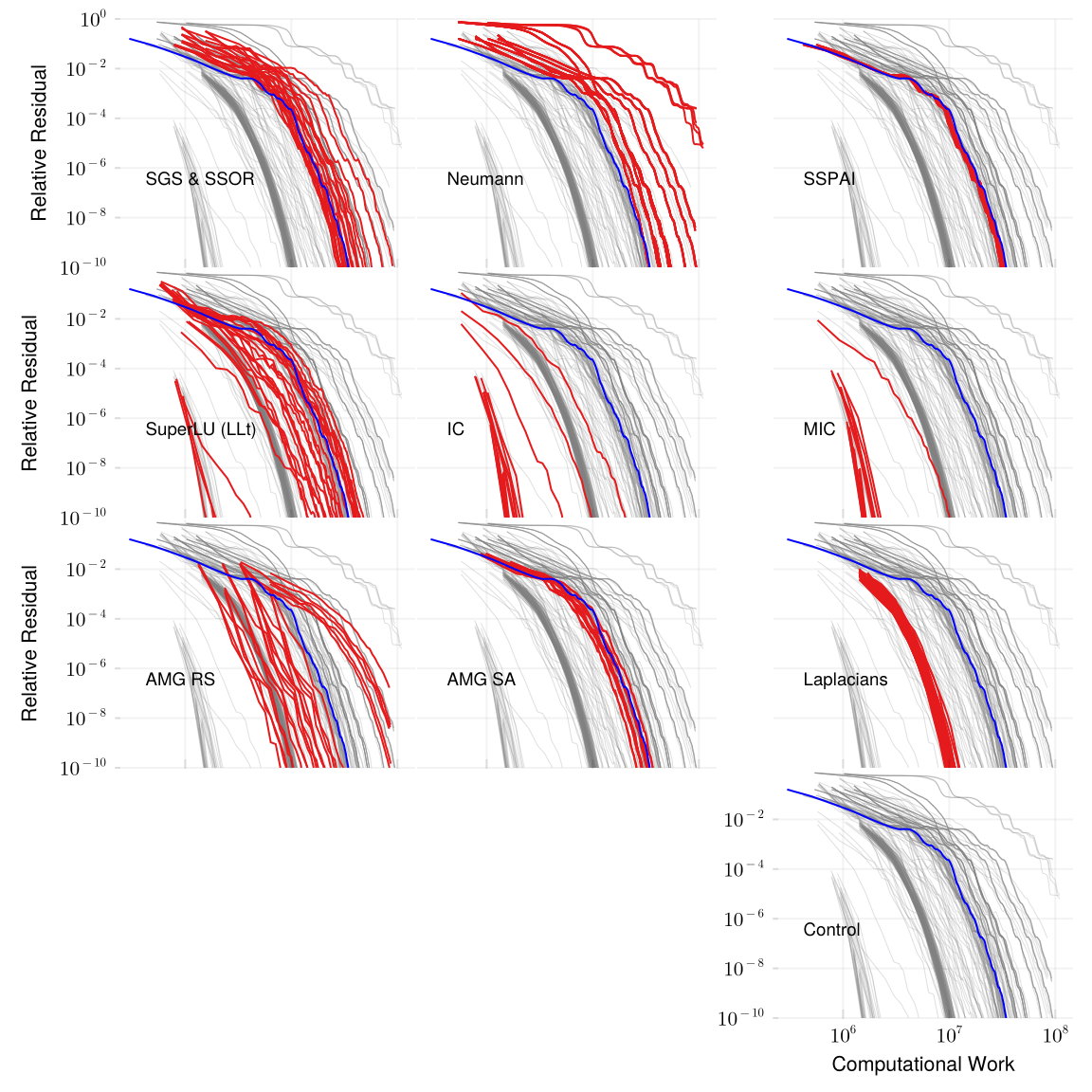}
    \caption{Convergence of the PCG method with various preconditioners applied to the \texttt{bodyy4} matrix (18k rows, 122k non-zeros). The plots have a log-log scale.}
    \label{fig:bodyy4}
\end{figure}
\clearpage 

\subsection{bodyy5}

\begin{figure}[!ht]
    \centering
    \includegraphics[width=\textwidth]{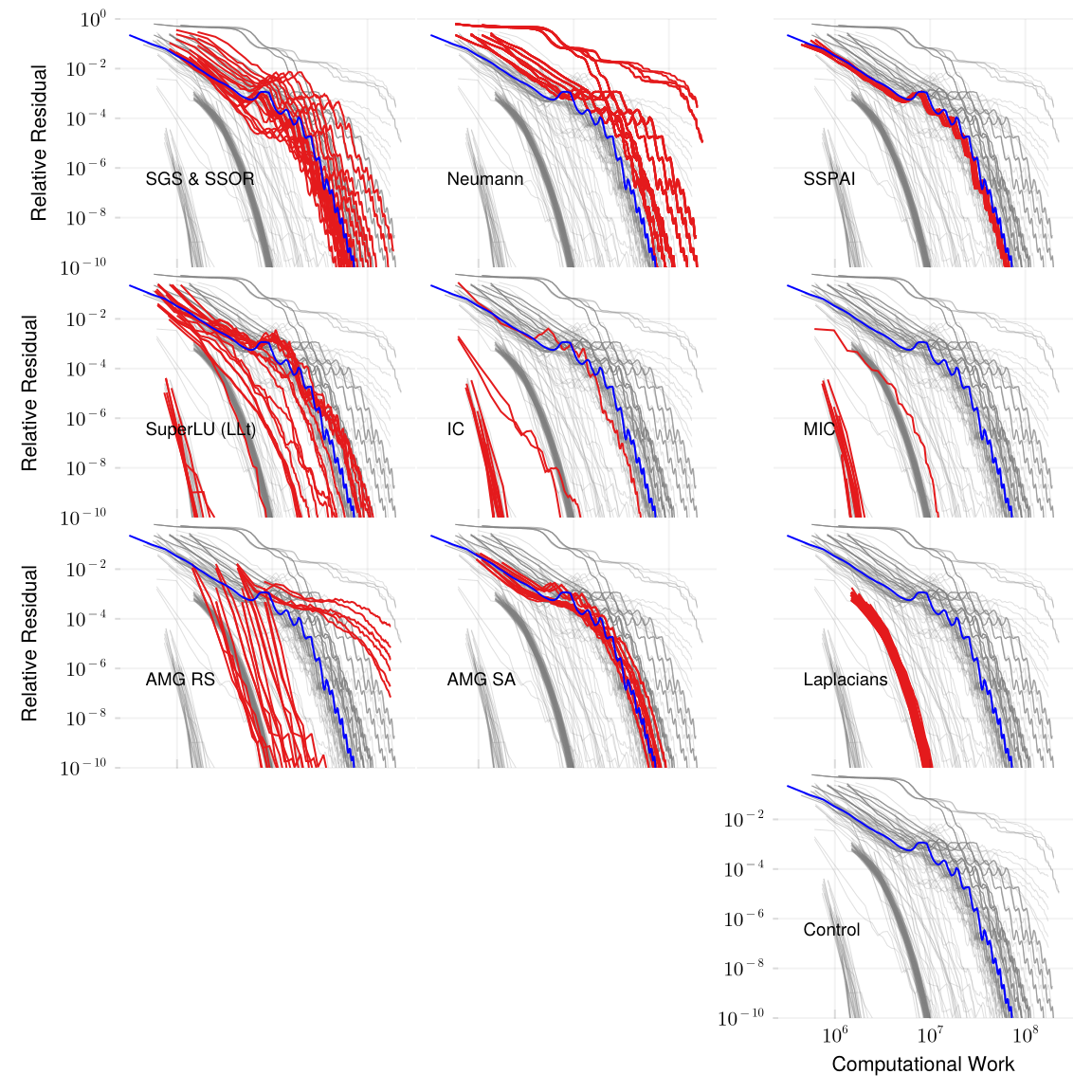}
    \caption{Convergence of the PCG method with various preconditioners applied to the \texttt{bodyy5} matrix (19k rows, 129k non-zeros). The plots have a log-log scale.}
    \label{fig:bodyy5}
\end{figure}
\clearpage 

\subsection{bodyy6}

\begin{figure}[!ht]
    \centering
    \includegraphics[width=\textwidth]{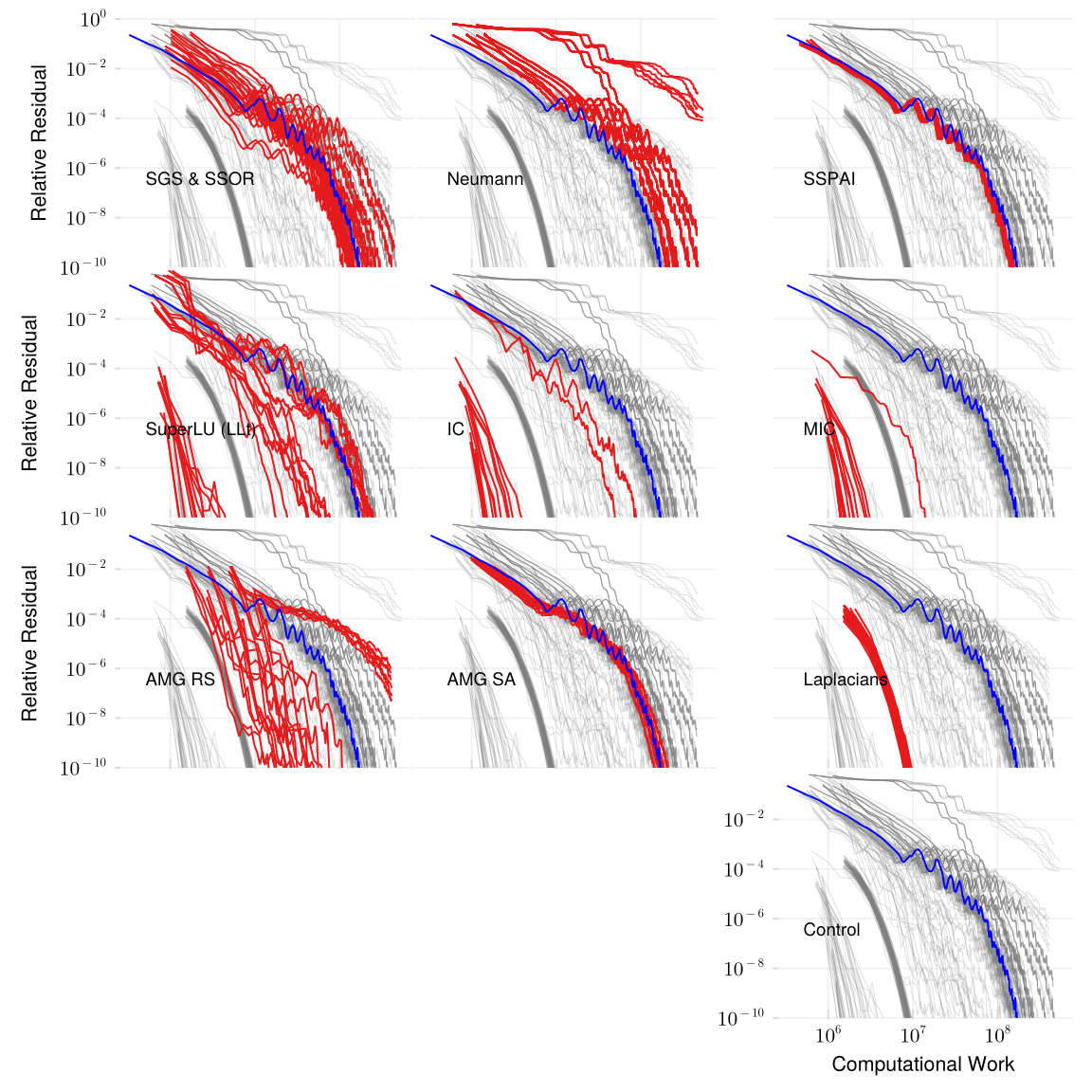}
    \caption{Convergence of the PCG method with various preconditioners applied to the \texttt{bodyy6} matrix (19k rows, 134k non-zeros). The plots have a log-log scale.}
    \label{fig:bodyy6}
\end{figure}
\clearpage 

\subsection{bone010}

\begin{figure}[!ht]
    \centering
    \includegraphics[width=\textwidth]{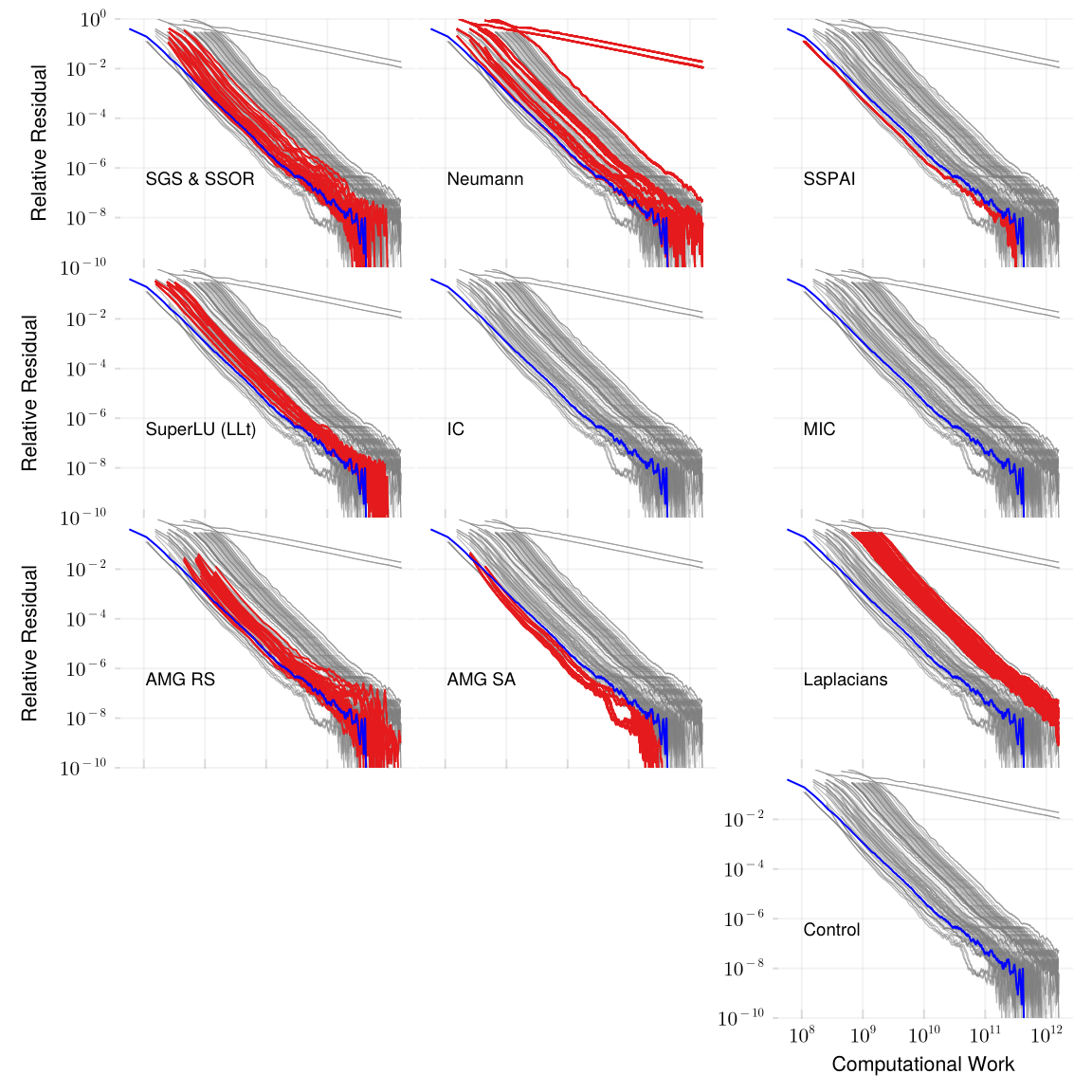}
    \caption{Convergence of the PCG method with various preconditioners applied to the \texttt{bone010} matrix (987k rows, 47.9m non-zeros). The plots have a log-log scale.}
    \label{fig:bone010}
\end{figure}
\clearpage 

\subsection{boneS01}

\begin{figure}[!ht]
    \centering
    \includegraphics[width=\textwidth]{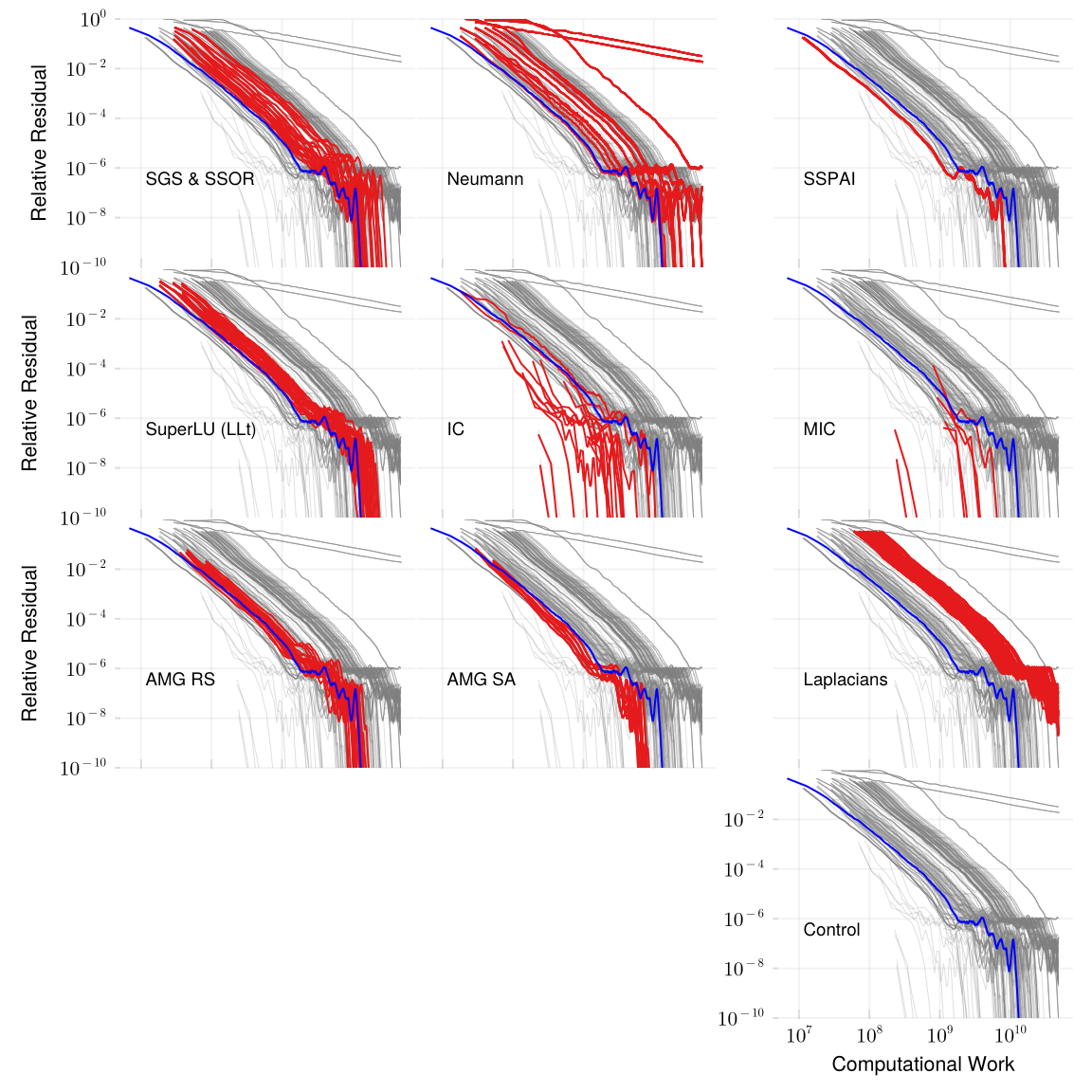}
    \caption{Convergence of the PCG method with various preconditioners applied to the \texttt{boneS01} matrix (127k rows, 5.5m non-zeros). The plots have a log-log scale.}
    \label{fig:boneS01}
\end{figure}
\clearpage 

\subsection{boneS10}

\begin{figure}[!ht]
    \centering
    \includegraphics[width=\textwidth]{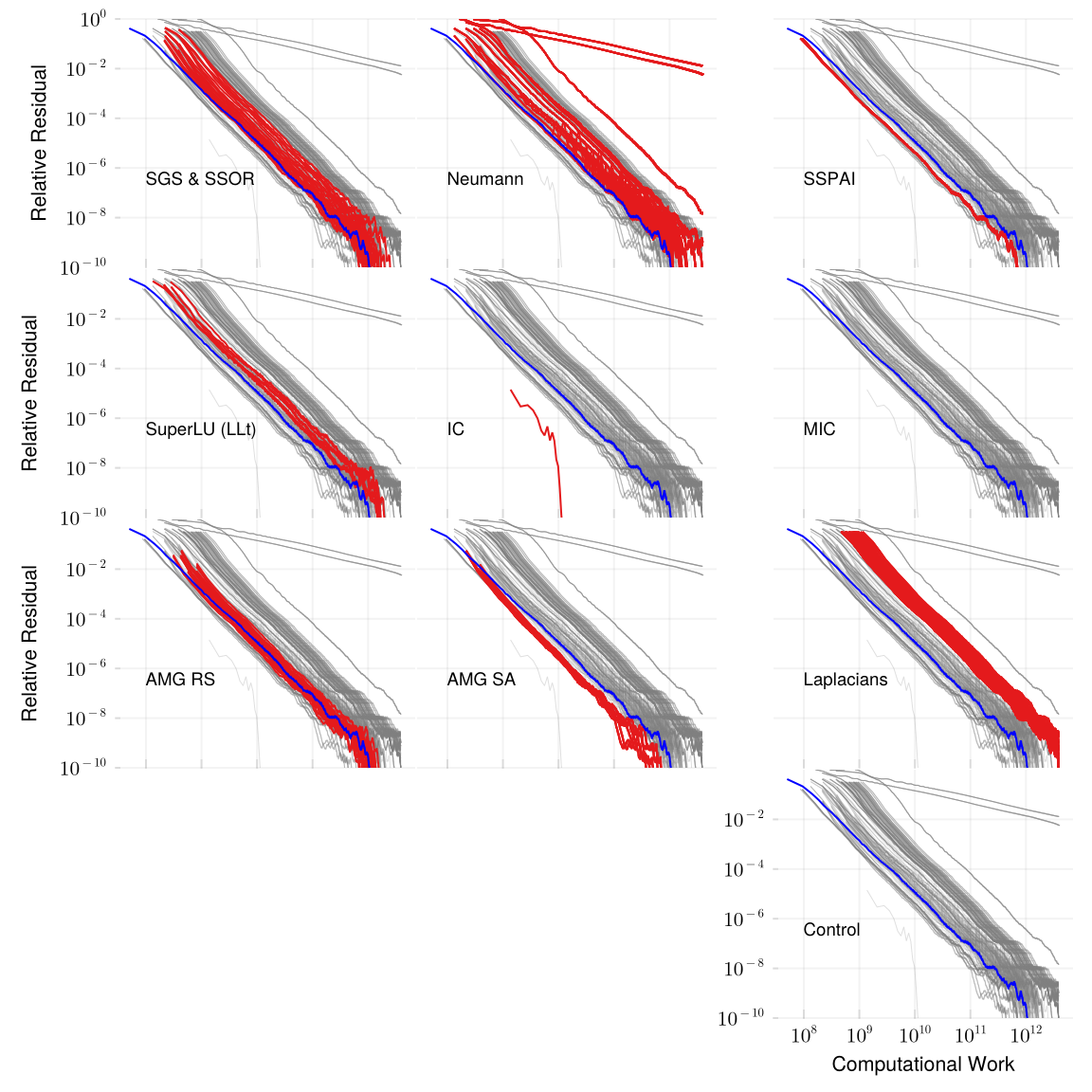}
    \caption{Convergence of the PCG method with various preconditioners applied to the \texttt{boneS10} matrix (915k rows, 40.9m non-zeros). The plots have a log-log scale.}
    \label{fig:boneS10}
\end{figure}
\clearpage 

\subsection{bundle1}

\begin{figure}[!ht]
    \centering
    \includegraphics[width=\textwidth]{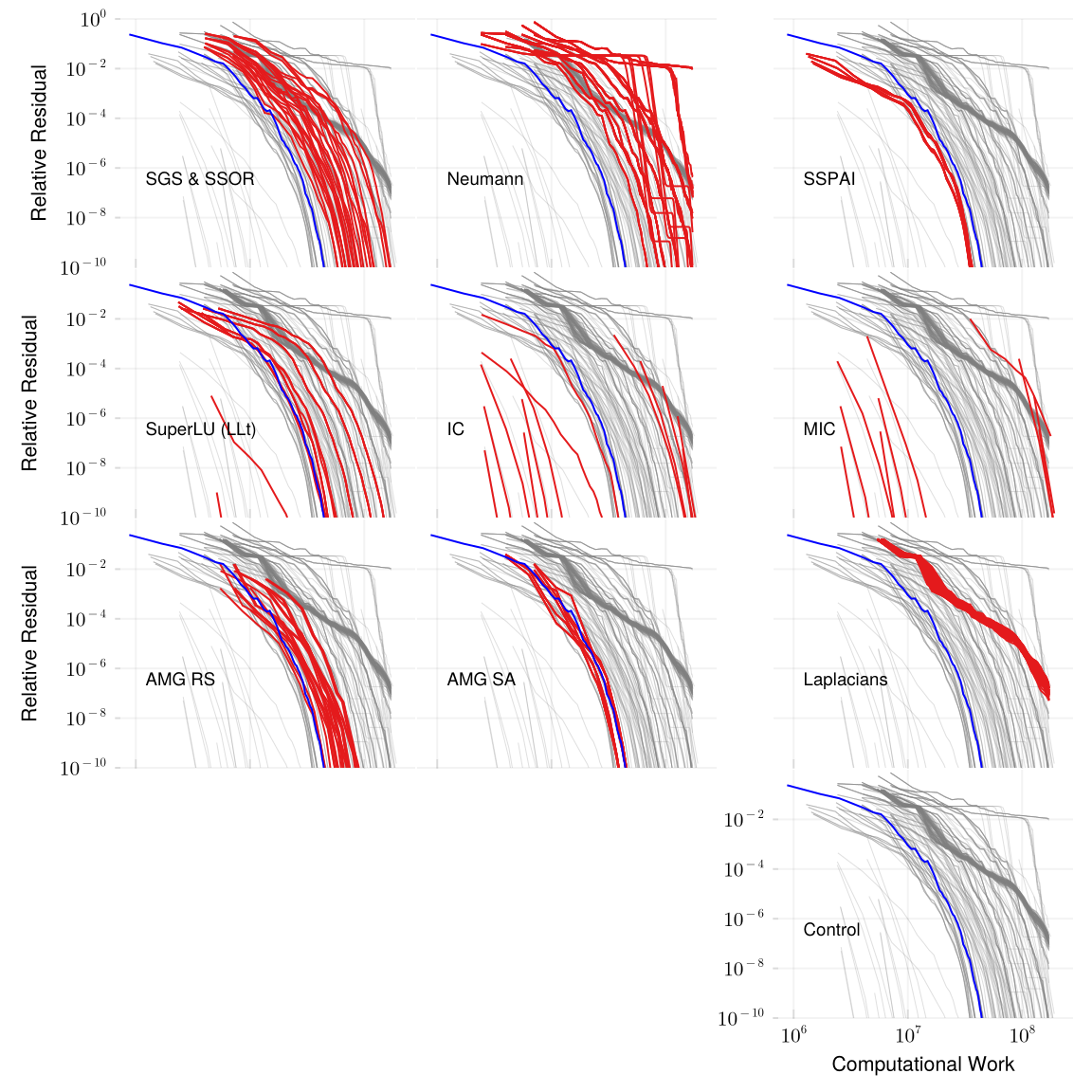}
    \caption{Convergence of the PCG method with various preconditioners applied to the \texttt{bundle1} matrix (11k rows, 771k non-zeros). The plots have a log-log scale.}
    \label{fig:bundle1}
\end{figure}
\clearpage 

\subsection{bundle\_adj}

\begin{figure}[!ht]
    \centering
    \includegraphics[width=\textwidth]{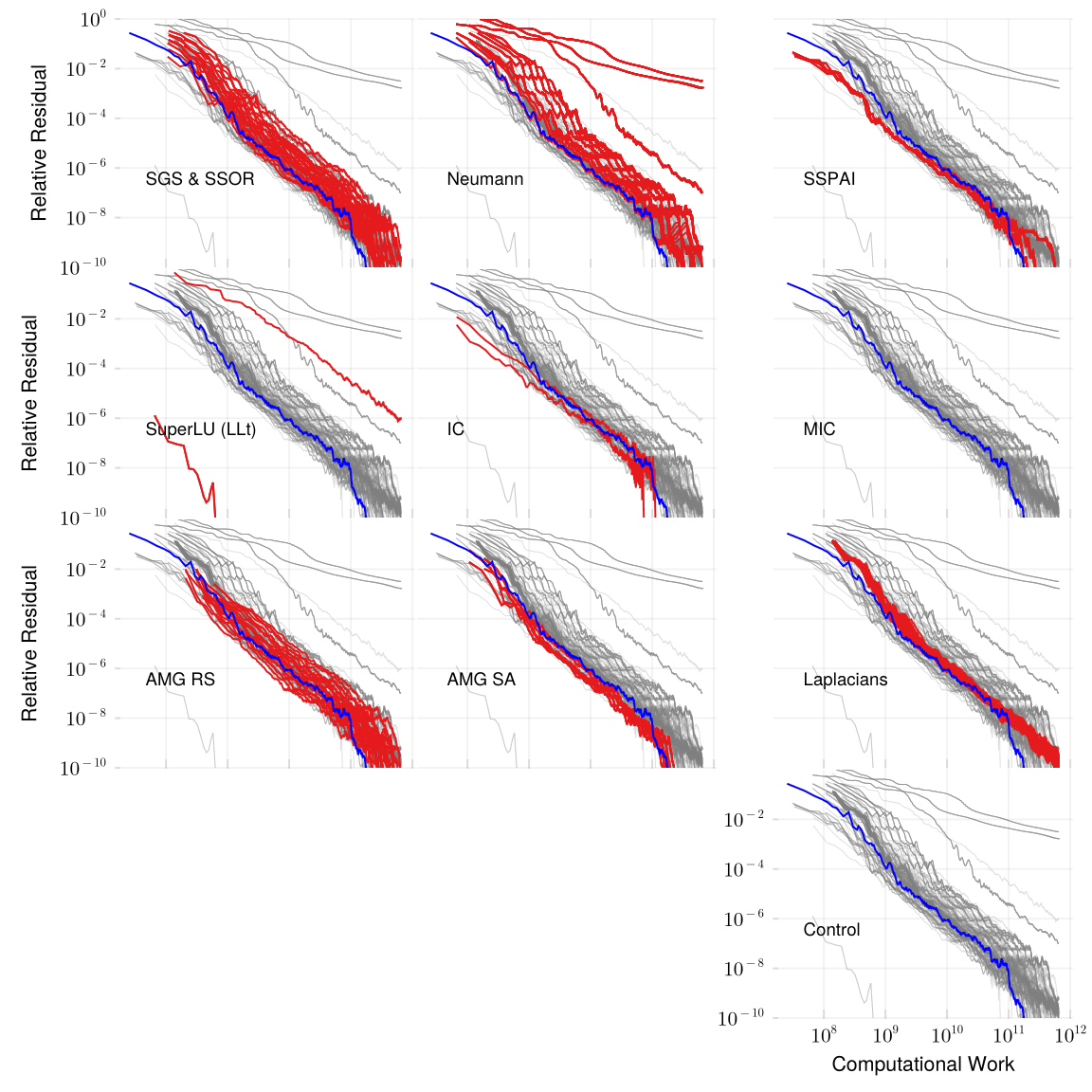}
    \caption{Convergence of the PCG method with various preconditioners applied to the \texttt{bundle\_adj} matrix (513k rows, 20.2m non-zeros). The plots have a log-log scale.}
    \label{fig:bundle_adj}
\end{figure}
\clearpage 

\subsection{cbuckle}

\begin{figure}[!ht]
    \centering
    \includegraphics[width=\textwidth]{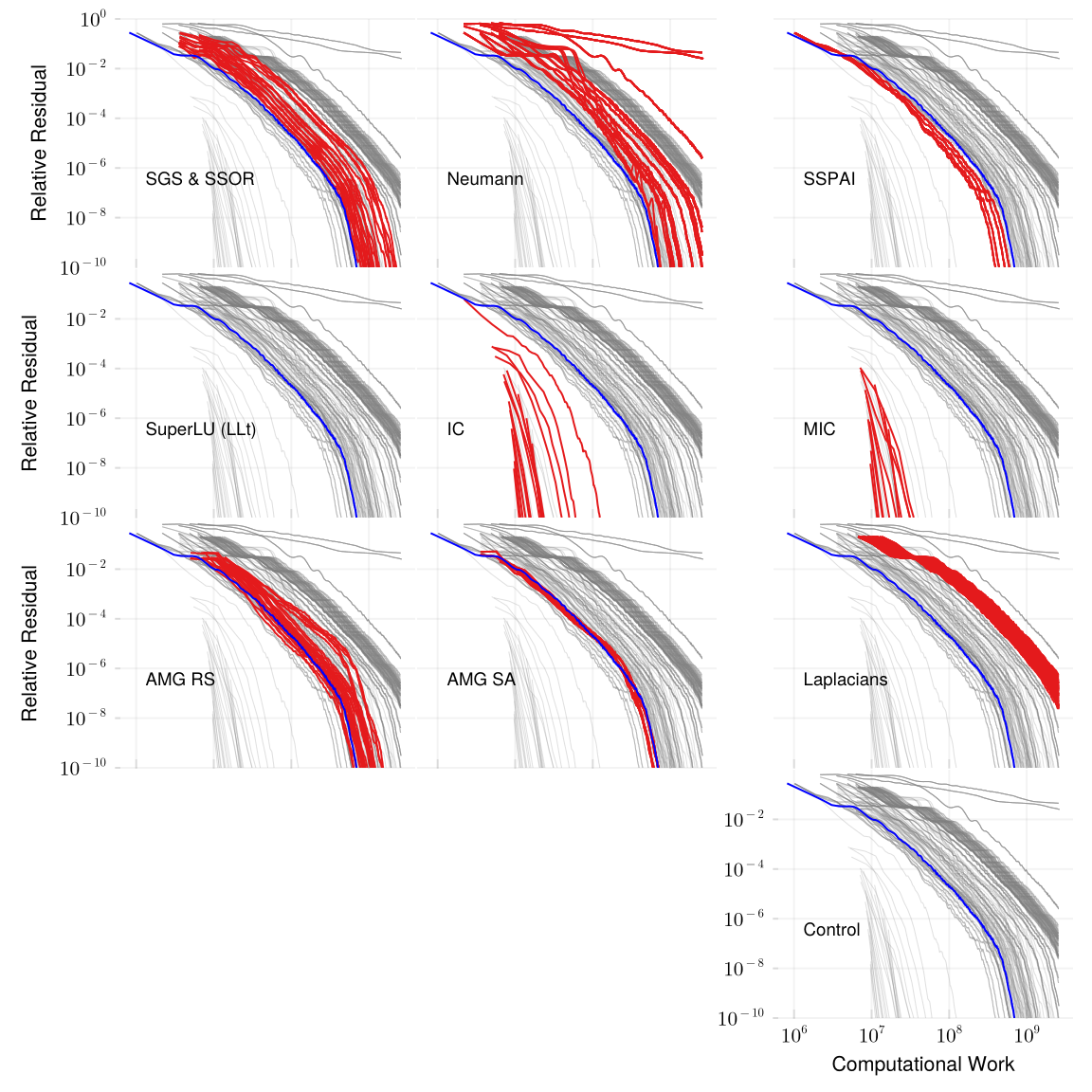}
    \caption{Convergence of the PCG method with various preconditioners applied to the \texttt{cbuckle} matrix (14k rows, 677k non-zeros). The plots have a log-log scale.}
    \label{fig:cbuckle}
\end{figure}
\clearpage 

\subsection{consph}

\begin{figure}[!ht]
    \centering
    \includegraphics[width=\textwidth]{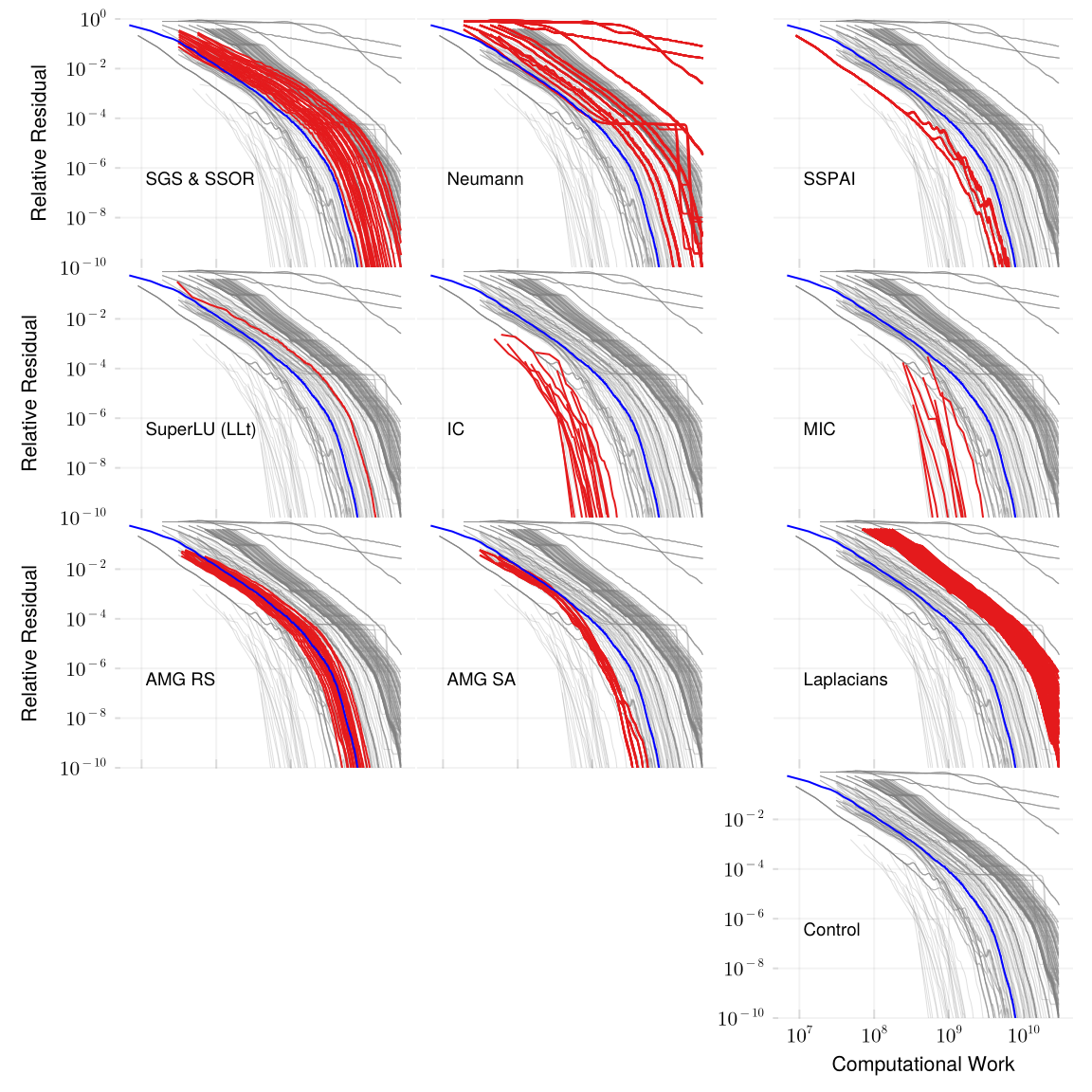}
    \caption{Convergence of the PCG method with various preconditioners applied to the \texttt{consph} matrix (83k rows, 6m non-zeros). The plots have a log-log scale.}
    \label{fig:consph}
\end{figure}
\clearpage 

\subsection{crankseg\_1}

\begin{figure}[!ht]
    \centering
    \includegraphics[width=\textwidth]{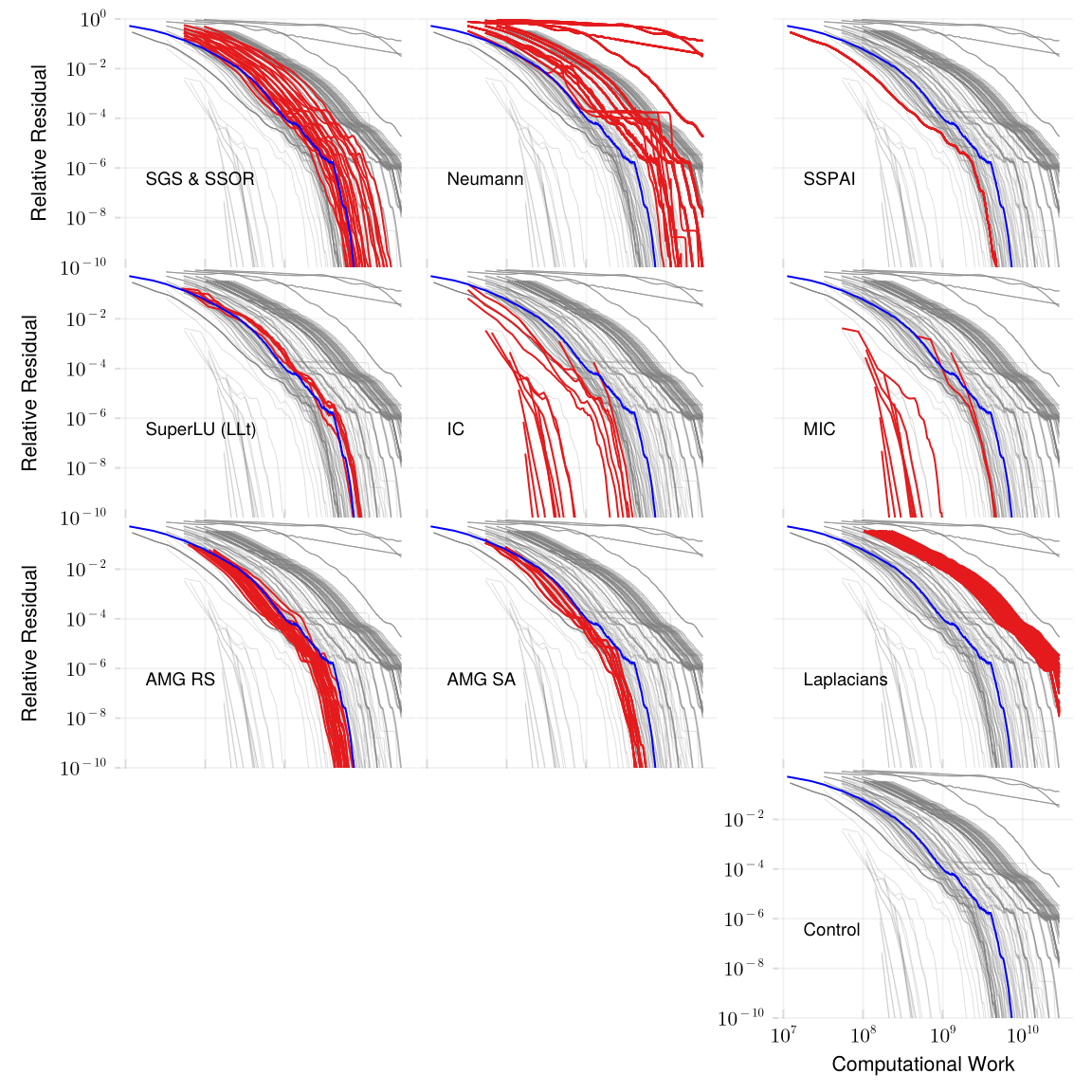}
    \caption{Convergence of the PCG method with various preconditioners applied to the \texttt{crankseg\_1} matrix (53k rows, 10.6m non-zeros). The plots have a log-log scale.}
    \label{fig:crankseg_1_}
\end{figure}
\clearpage 

\subsection{crankseg\_2}

\begin{figure}[!ht]
    \centering
    \includegraphics[width=\textwidth]{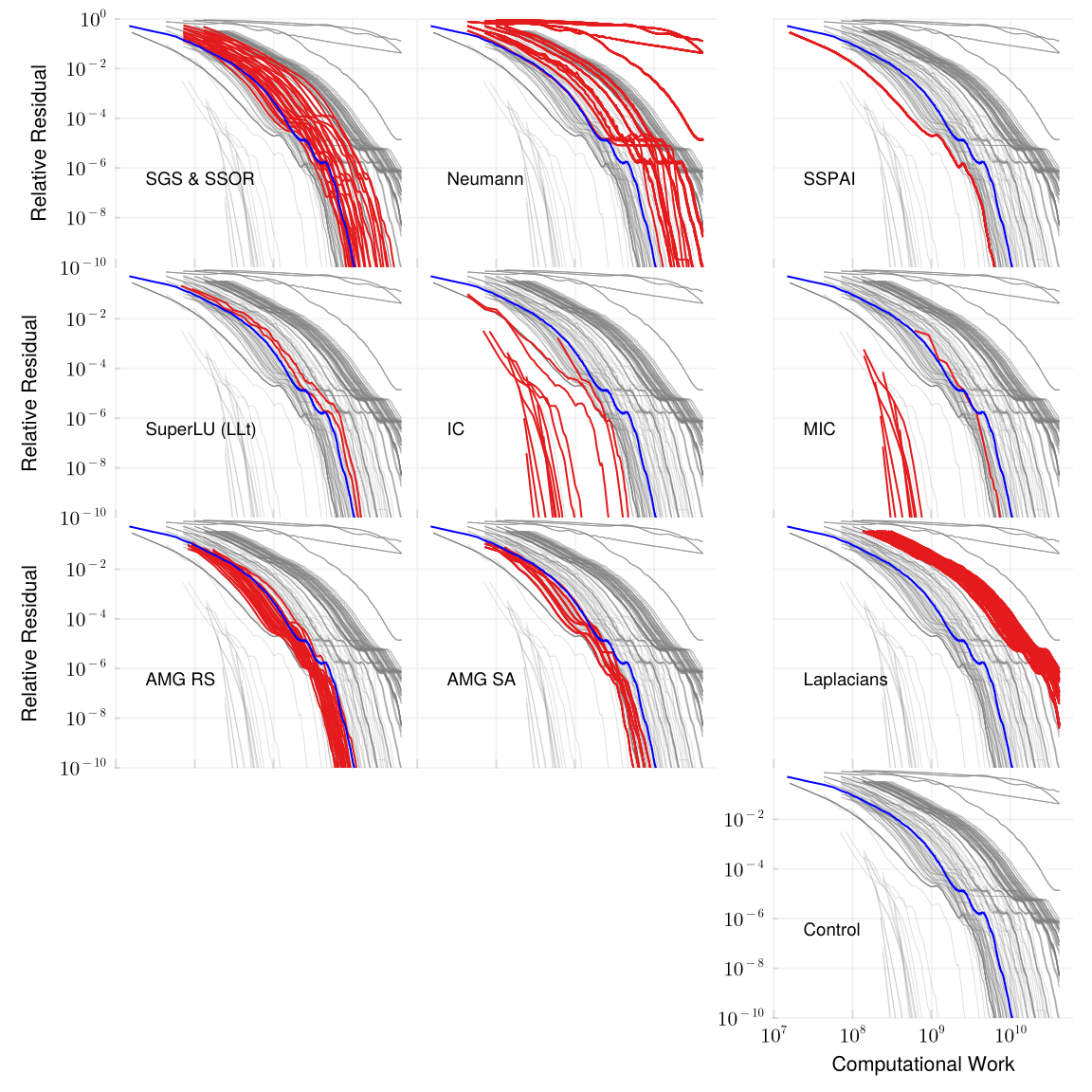}
    \caption{Convergence of the PCG method with various preconditioners applied to the \texttt{crankseg\_2} matrix (64k rows, 14.1m non-zeros). The plots have a log-log scale.}
    \label{fig:crankseg_2_}
\end{figure}
\clearpage 

\subsection{crystm02}

\begin{figure}[!ht]
    \centering
    \includegraphics[width=\textwidth]{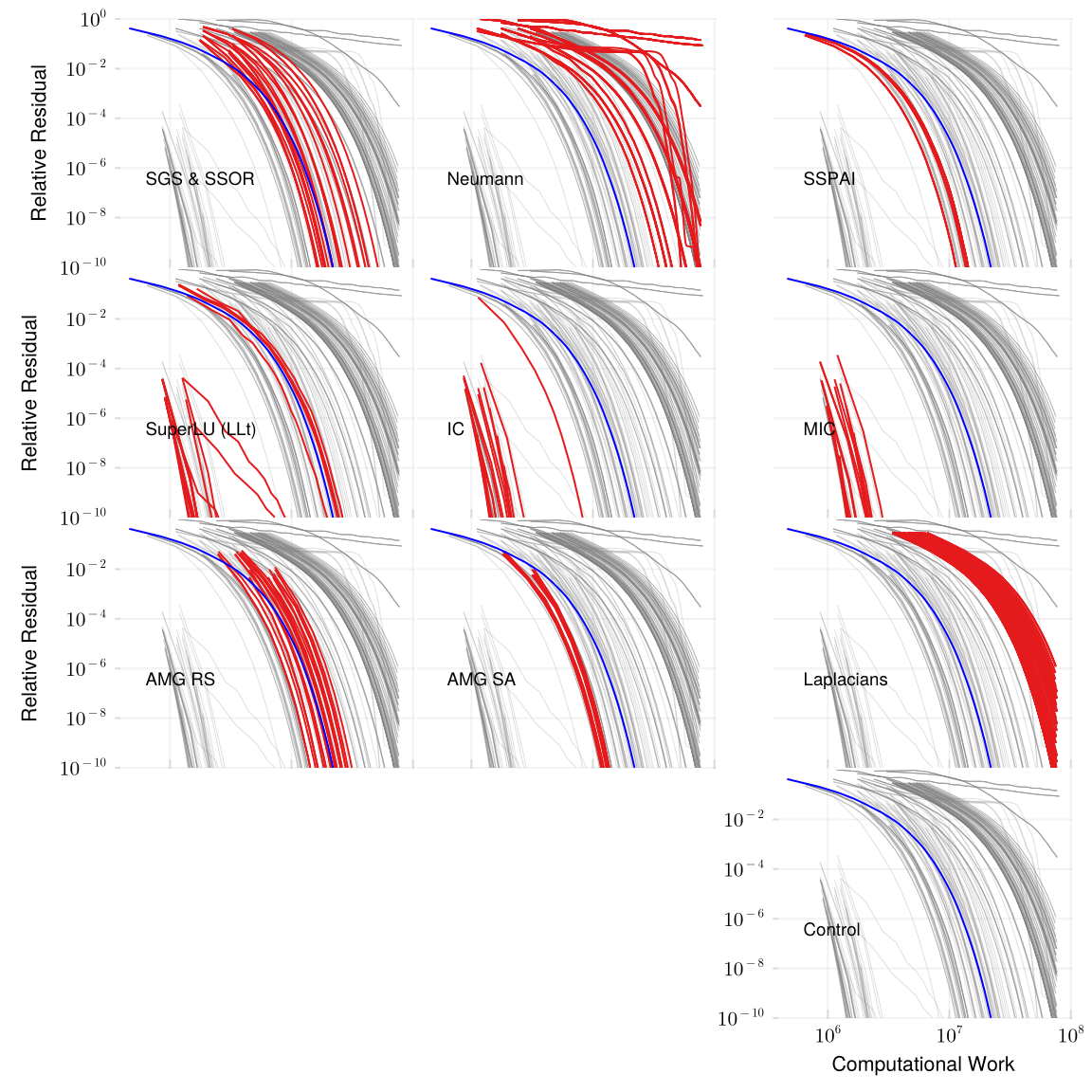}
    \caption{Convergence of the PCG method with various preconditioners applied to the \texttt{crystm02} matrix (14k rows, 333k non-zeros). The plots have a log-log scale.}
    \label{fig:crystm02}
\end{figure}
\clearpage 

\subsection{crystm03}

\begin{figure}[!ht]
    \centering
    \includegraphics[width=\textwidth]{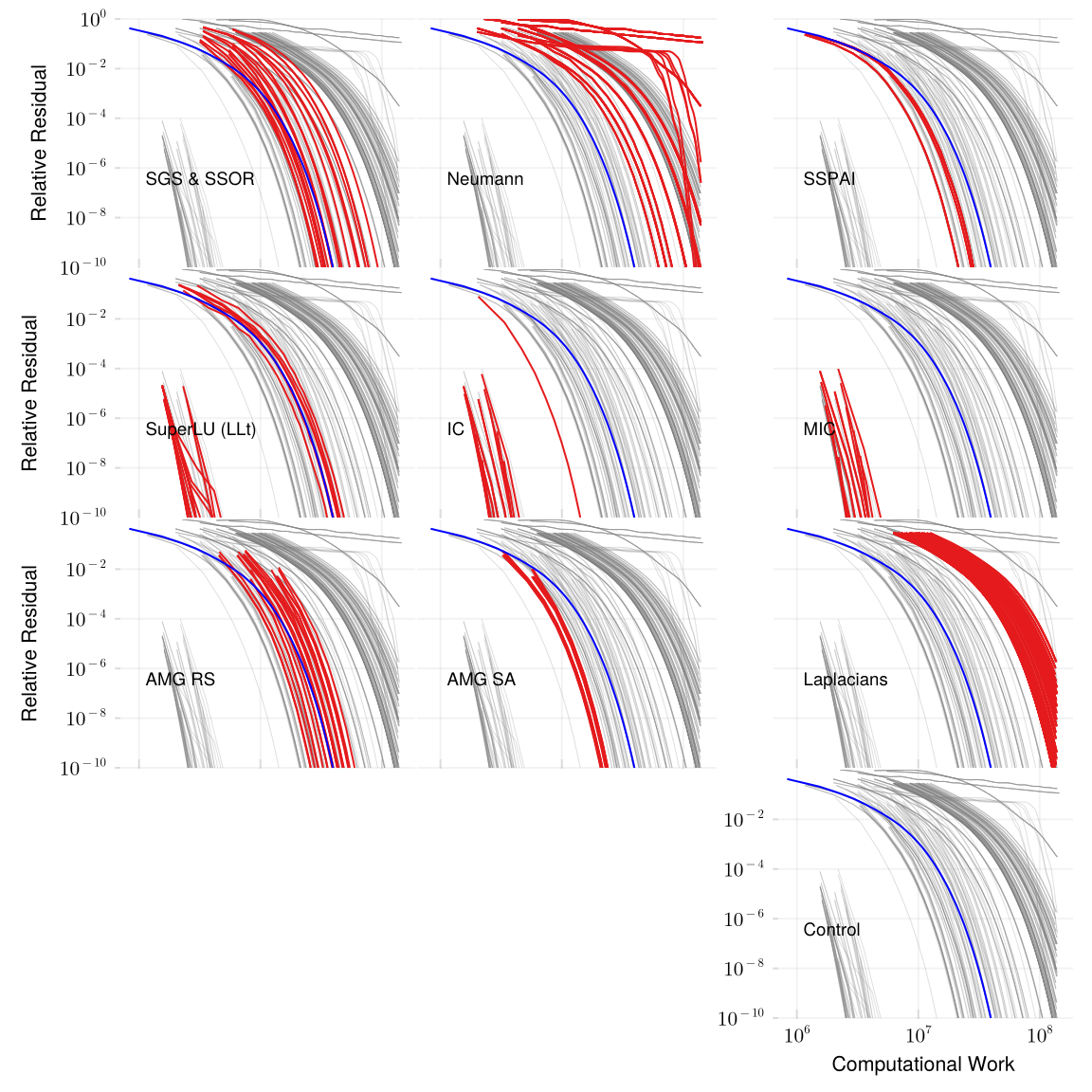}
    \caption{Convergence of the PCG method with various preconditioners applied to the \texttt{crystm03} matrix (25k rows, 584k non-zeros). The plots have a log-log scale.}
    \label{fig:crystm03}
\end{figure}
\clearpage 

\subsection{cvxbqp1}

\begin{figure}[!ht]
    \centering
    \includegraphics[width=\textwidth]{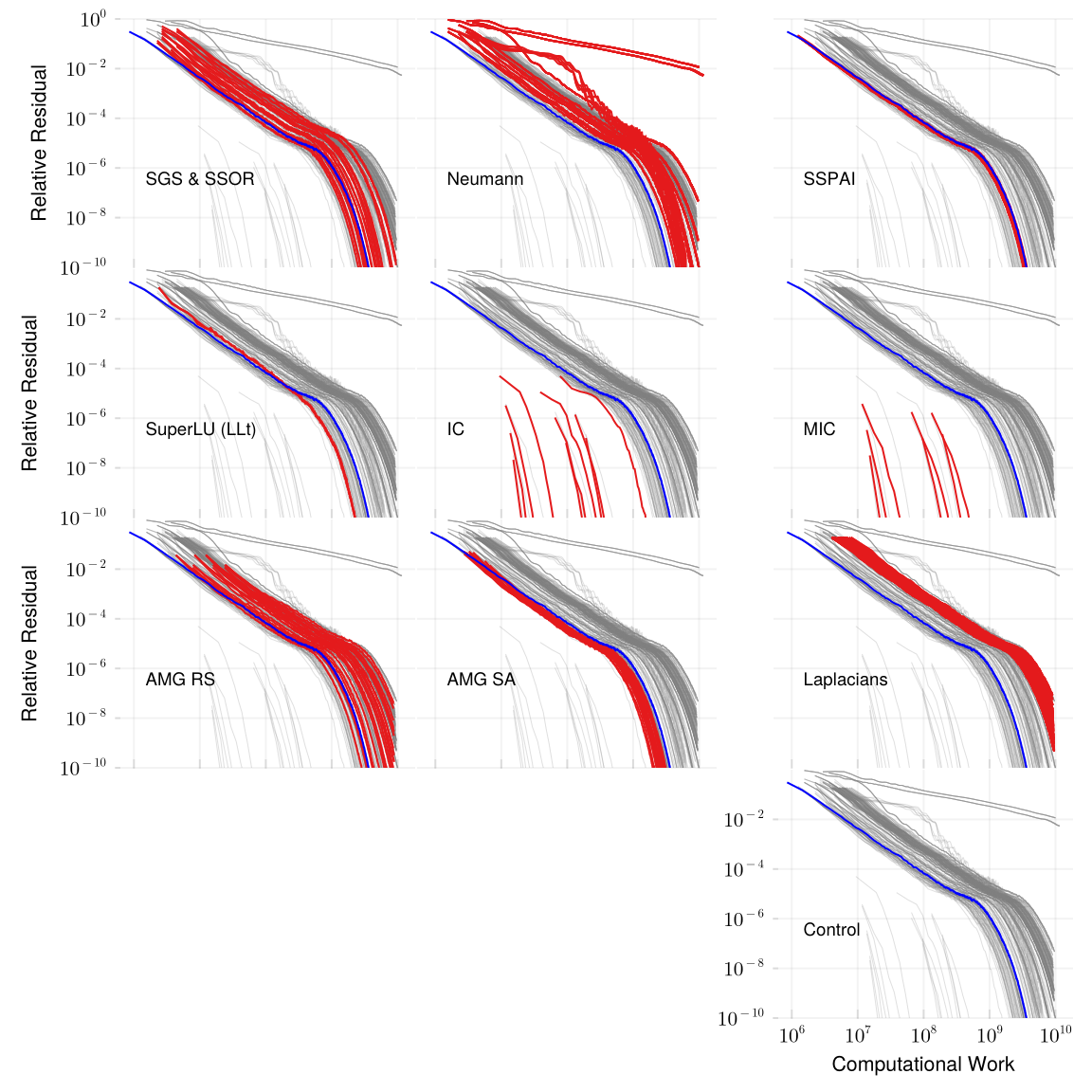}
    \caption{Convergence of the PCG method with various preconditioners applied to the \texttt{cvxbqp1} matrix (50k rows, 350k non-zeros). The plots have a log-log scale.}
    \label{fig:cvxbqp1}
\end{figure}
\clearpage 

\subsection{ecology2}

\begin{figure}[!ht]
    \centering
    \includegraphics[width=\textwidth]{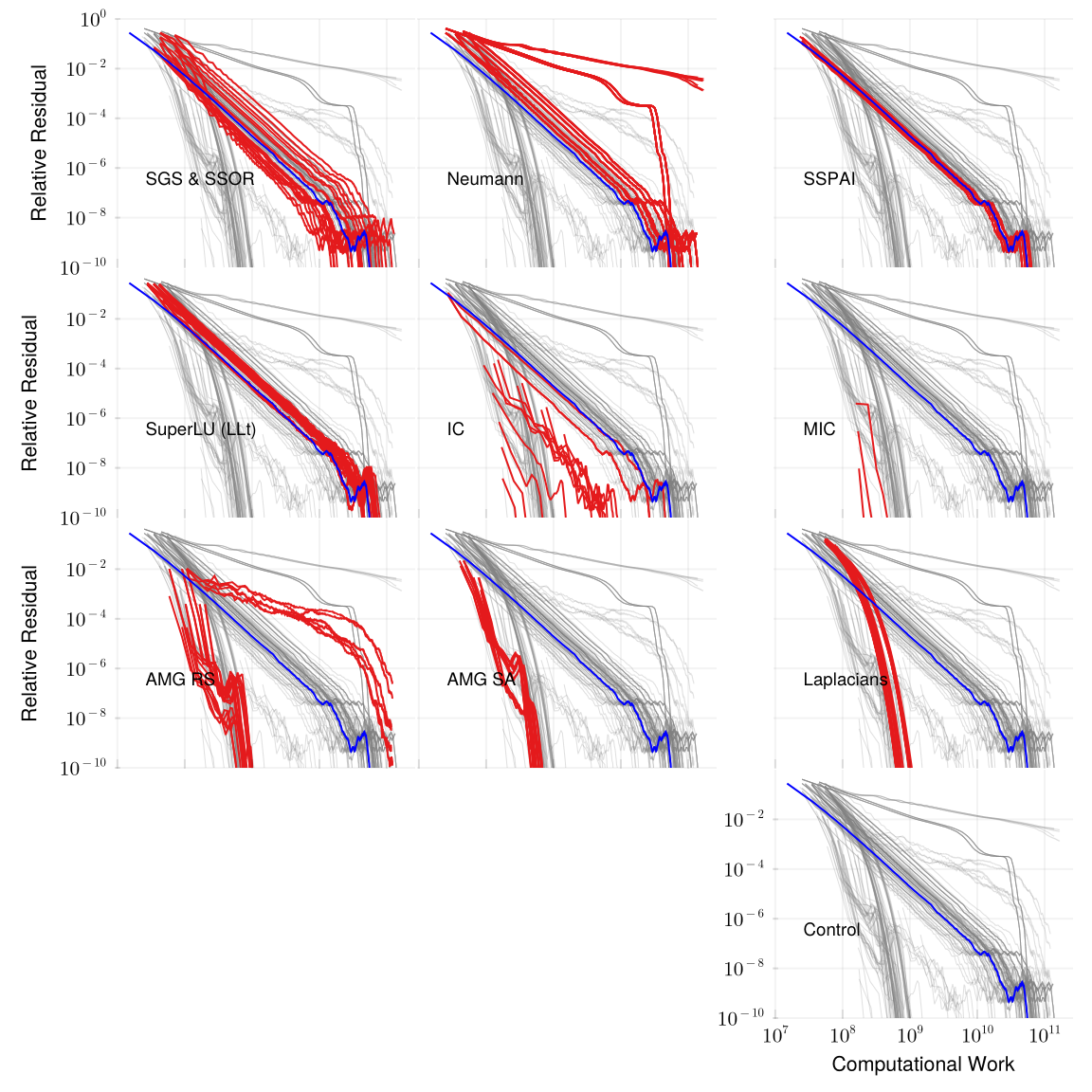}
    \caption{Convergence of the PCG method with various preconditioners applied to the \texttt{ecology2} matrix (1m rows, 5m non-zeros). The plots have a log-log scale.}
    \label{fig:ecology2_}
\end{figure}
\clearpage 

\subsection{finan512}

\begin{figure}[!ht]
    \centering
    \includegraphics[width=\textwidth]{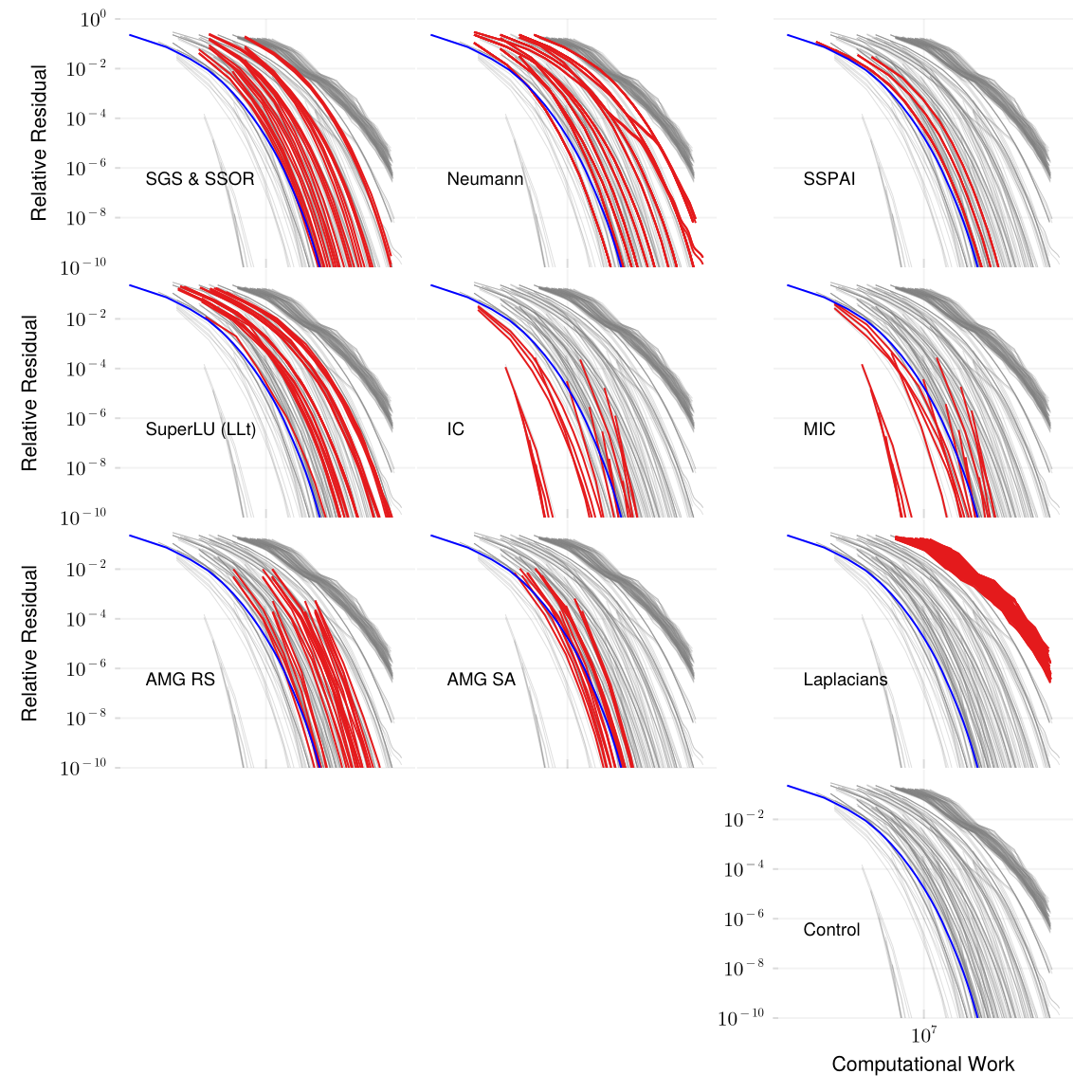}
    \caption{Convergence of the PCG method with various preconditioners applied to the \texttt{finan512} matrix (75k rows, 597k non-zeros). The plots have a log-log scale.}
    \label{fig:finan512}
\end{figure}
\clearpage 

\subsection{gyro\_m}

\begin{figure}[!ht]
    \centering
    \includegraphics[width=\textwidth]{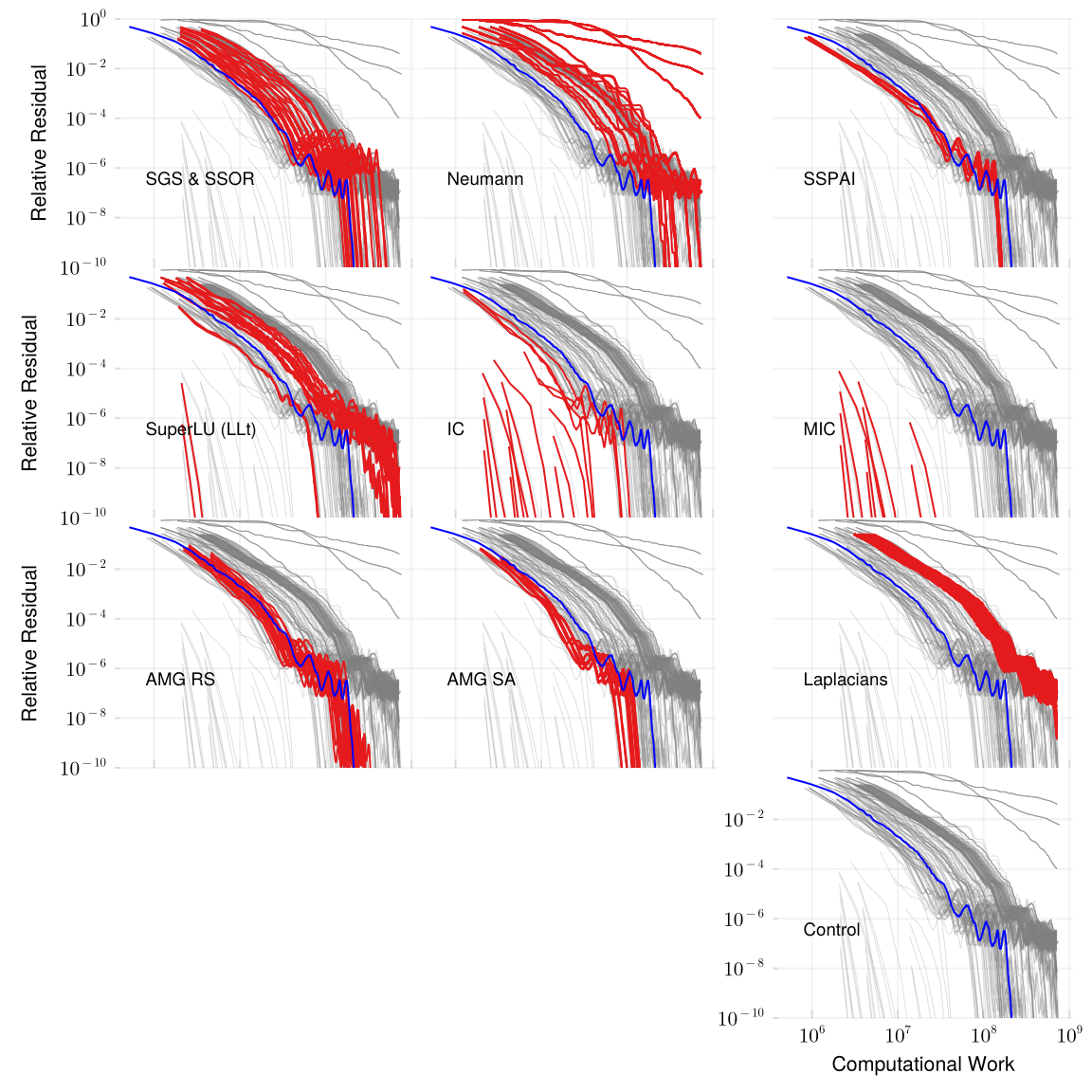}
    \caption{Convergence of the PCG method with various preconditioners applied to the \texttt{gyro\_m} matrix (17k rows, 340k non-zeros). The plots have a log-log scale.}
    \label{fig:gyro_m}
\end{figure}
\clearpage 

\subsection{hood}

\begin{figure}[!ht]
    \centering
    \includegraphics[width=\textwidth]{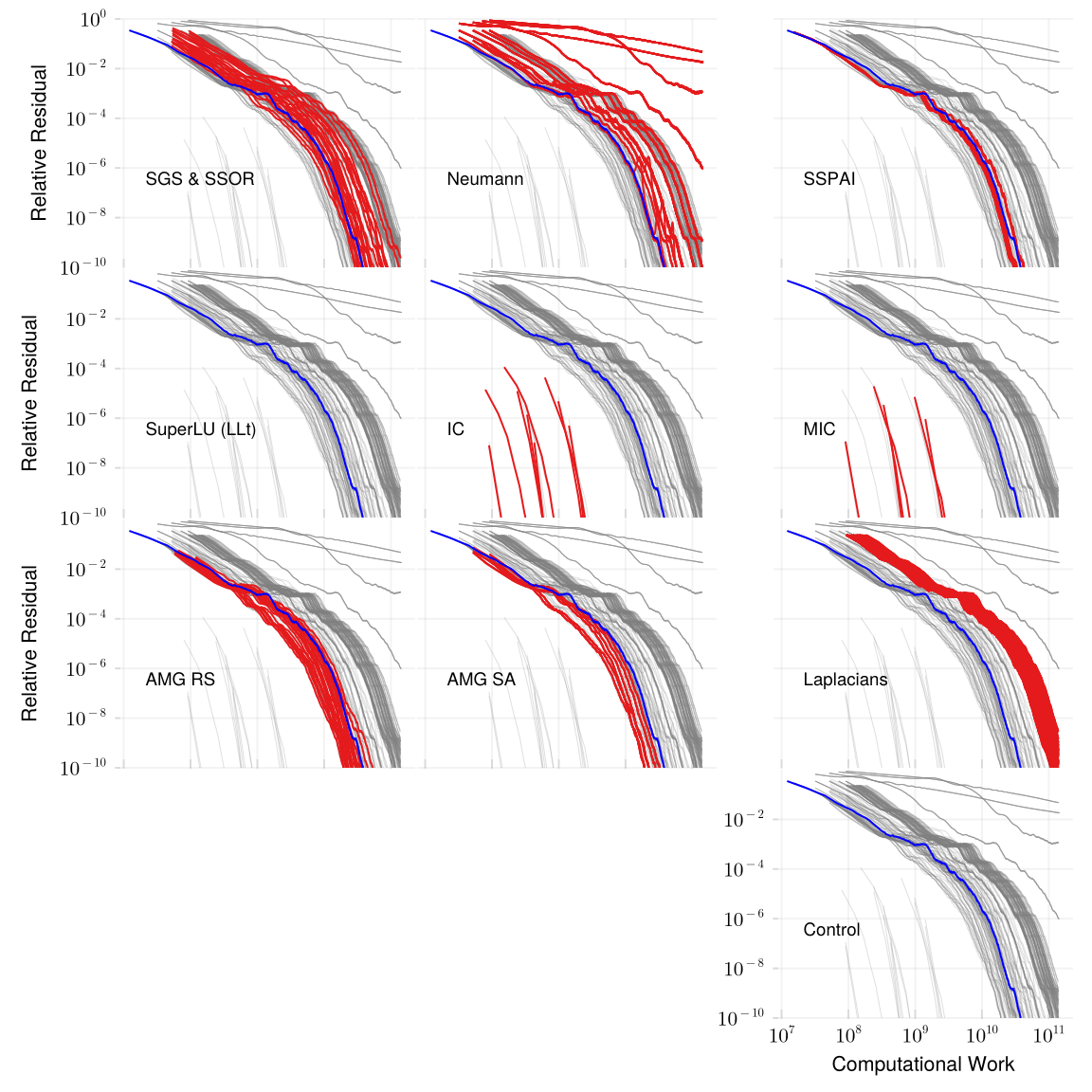}
    \caption{Convergence of the PCG method with various preconditioners applied to the \texttt{hood} matrix (221k rows, 9.9m non-zeros). The plots have a log-log scale.}
    \label{fig:hood}
\end{figure}
\clearpage 

\subsection{inline\_1}

\begin{figure}[!ht]
    \centering
    \includegraphics[width=\textwidth]{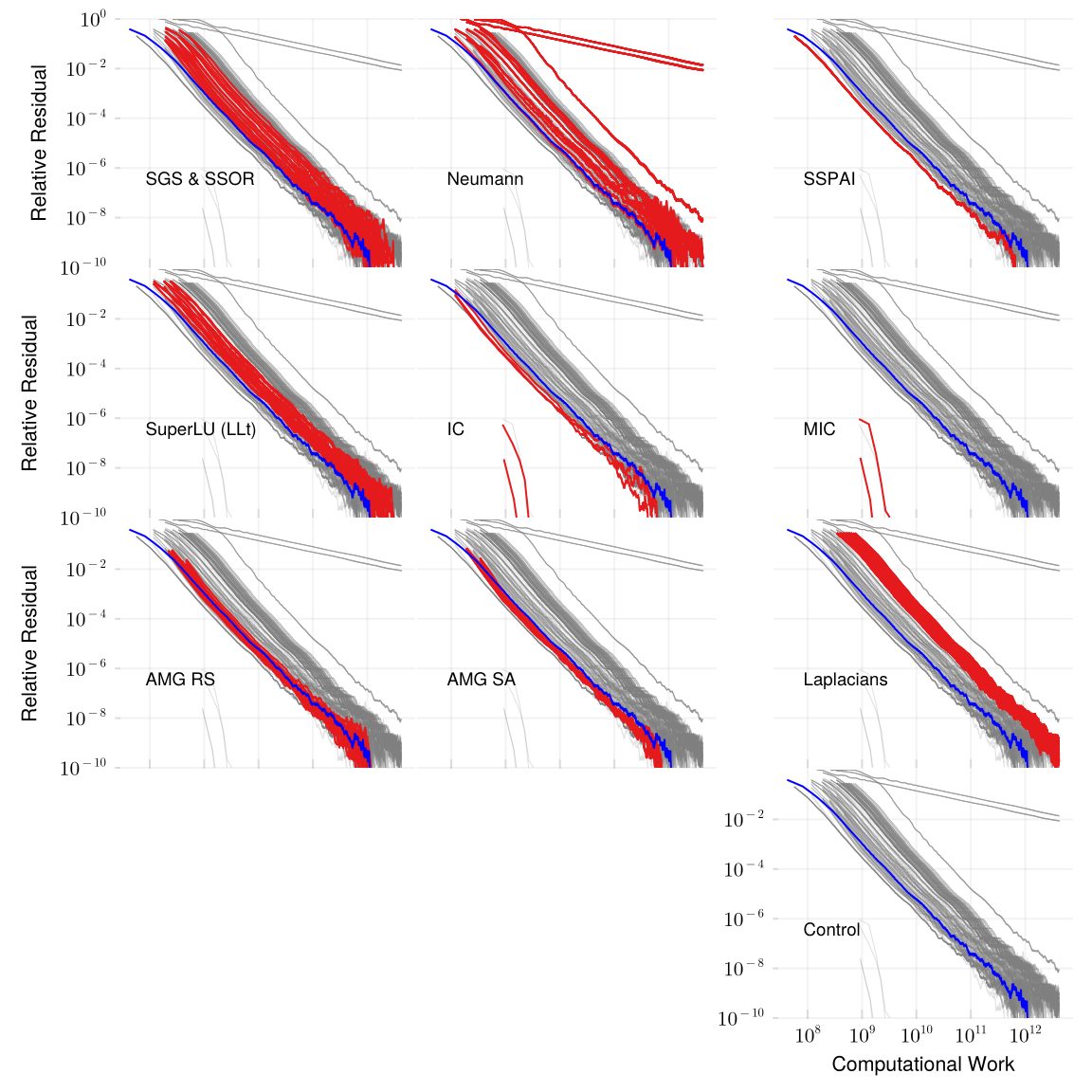}
    \caption{Convergence of the PCG method with various preconditioners applied to the \texttt{inline\_1} matrix (504k rows, 36.8m non-zeros). The plots have a log-log scale.}
    \label{fig:inline_1}
\end{figure}
\clearpage 

\subsection{jnlbrng1}

\begin{figure}[!ht]
    \centering
    \includegraphics[width=\textwidth]{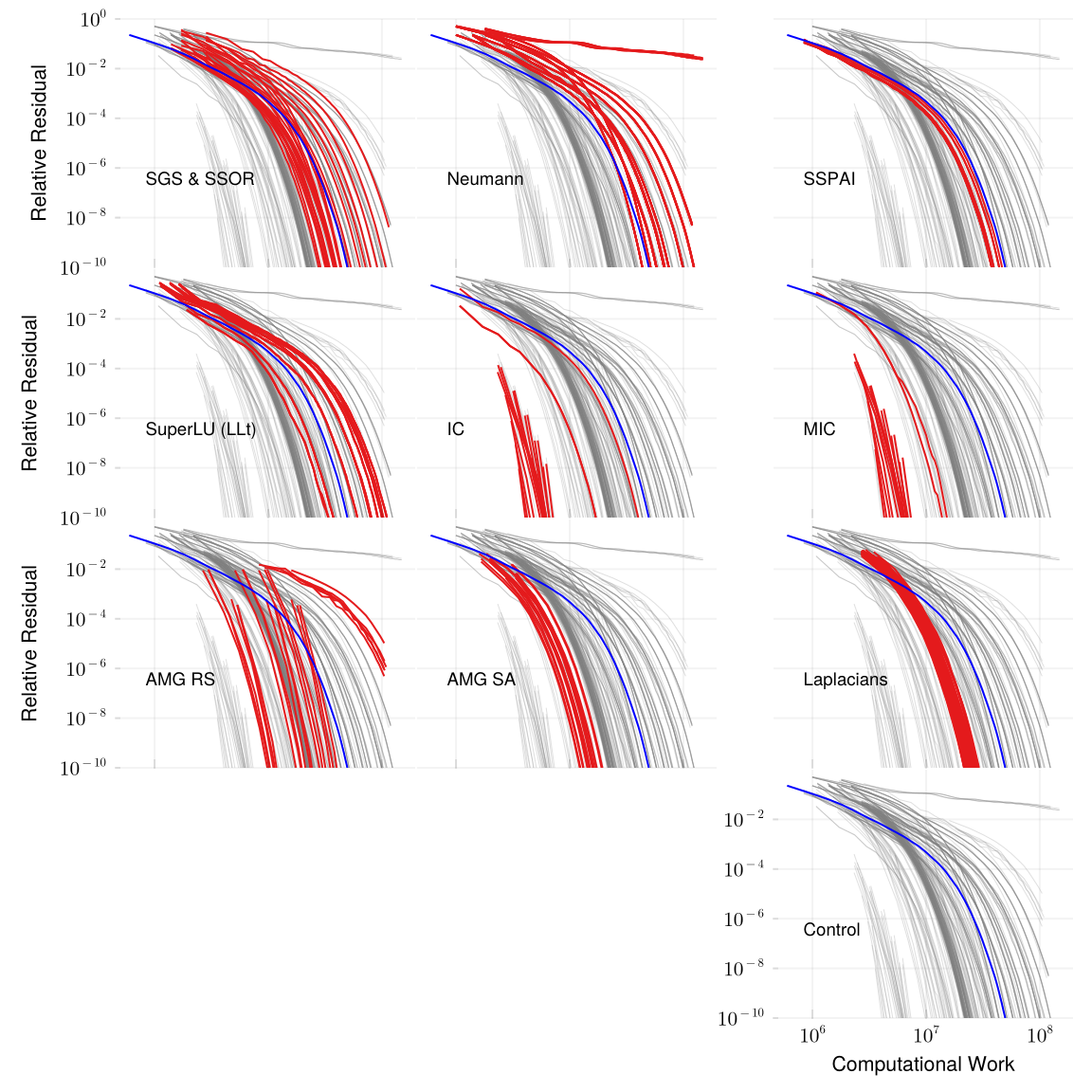}
    \caption{Convergence of the PCG method with various preconditioners applied to the \texttt{jnlbrng1} matrix (40k rows, 199k non-zeros). The plots have a log-log scale.}
    \label{fig:jnlbrng1}
\end{figure}
\clearpage 

\subsection{ldoor}

\begin{figure}[!ht]
    \centering
    \includegraphics[width=\textwidth]{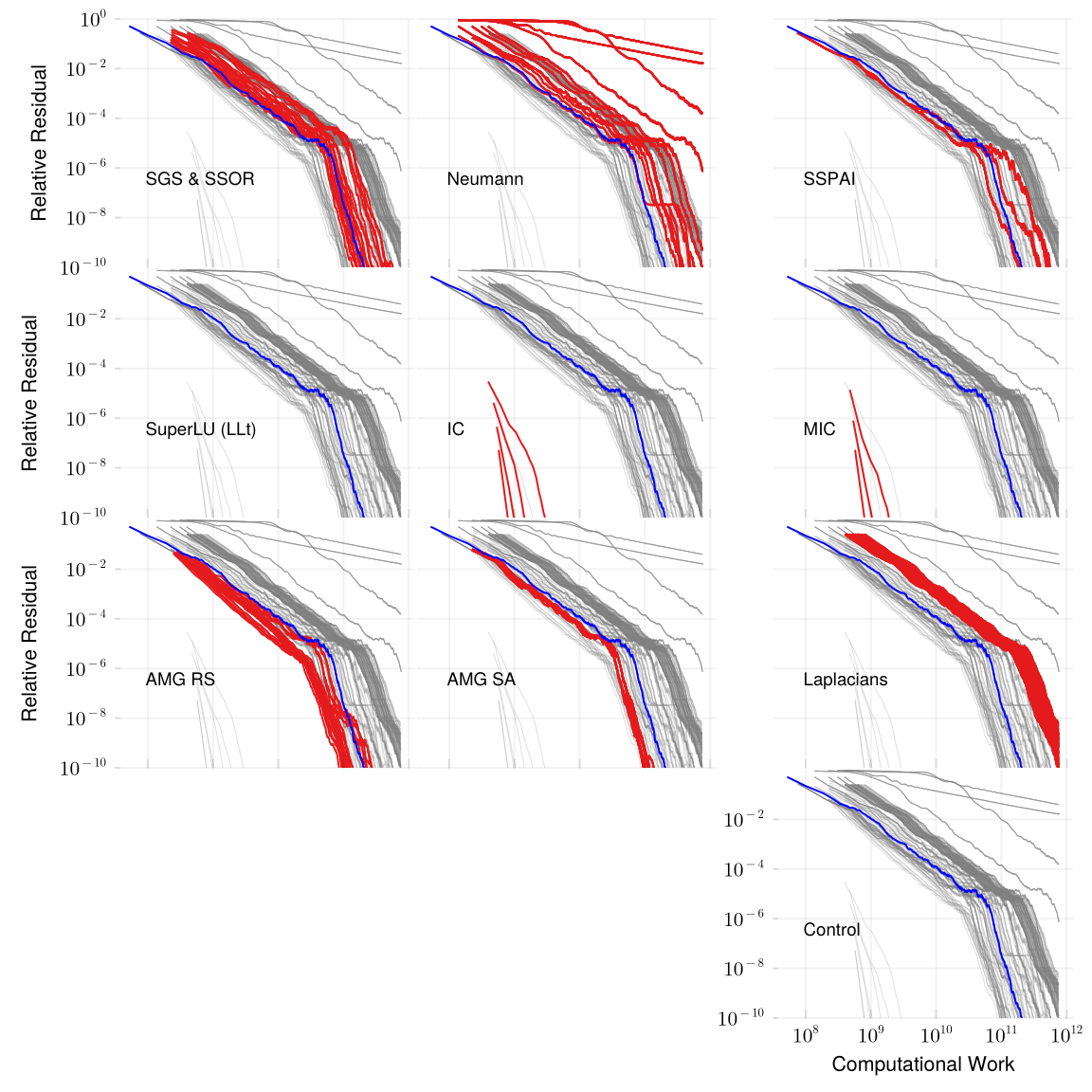}
    \caption{Convergence of the PCG method with various preconditioners applied to the \texttt{ldoor} matrix (952k rows, 42.5m non-zeros). The plots have a log-log scale.}
    \label{fig:ldoor}
\end{figure}
\clearpage 

\subsection{minsurfo}

\begin{figure}[!ht]
    \centering
    \includegraphics[width=\textwidth]{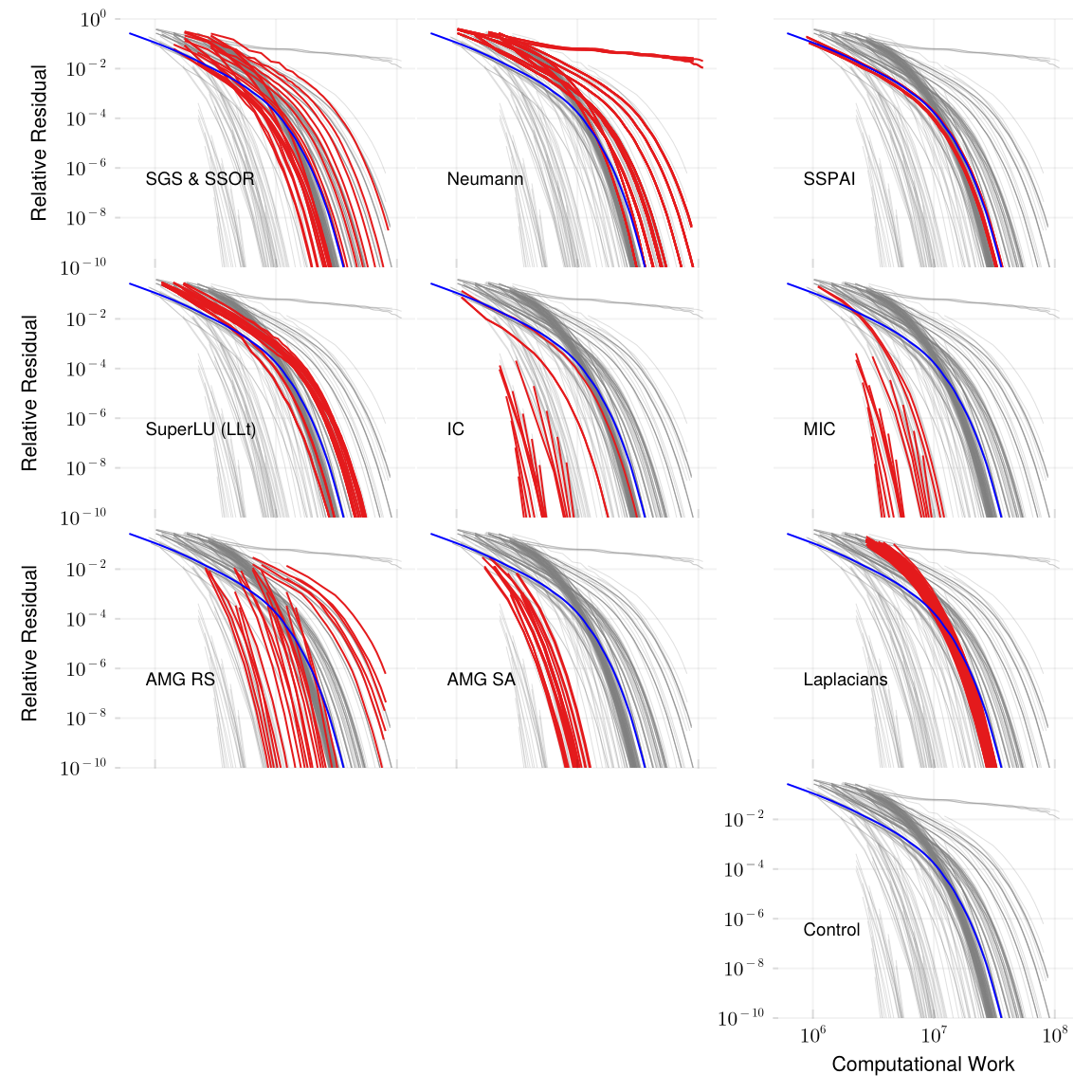}
    \caption{Convergence of the PCG method with various preconditioners applied to the \texttt{minsurfo} matrix (41k rows, 204k non-zeros). The plots have a log-log scale.}
    \label{fig:minsurfo}
\end{figure}
\clearpage 

\subsection{msdoor}

\begin{figure}[!ht]
    \centering
    \includegraphics[width=\textwidth]{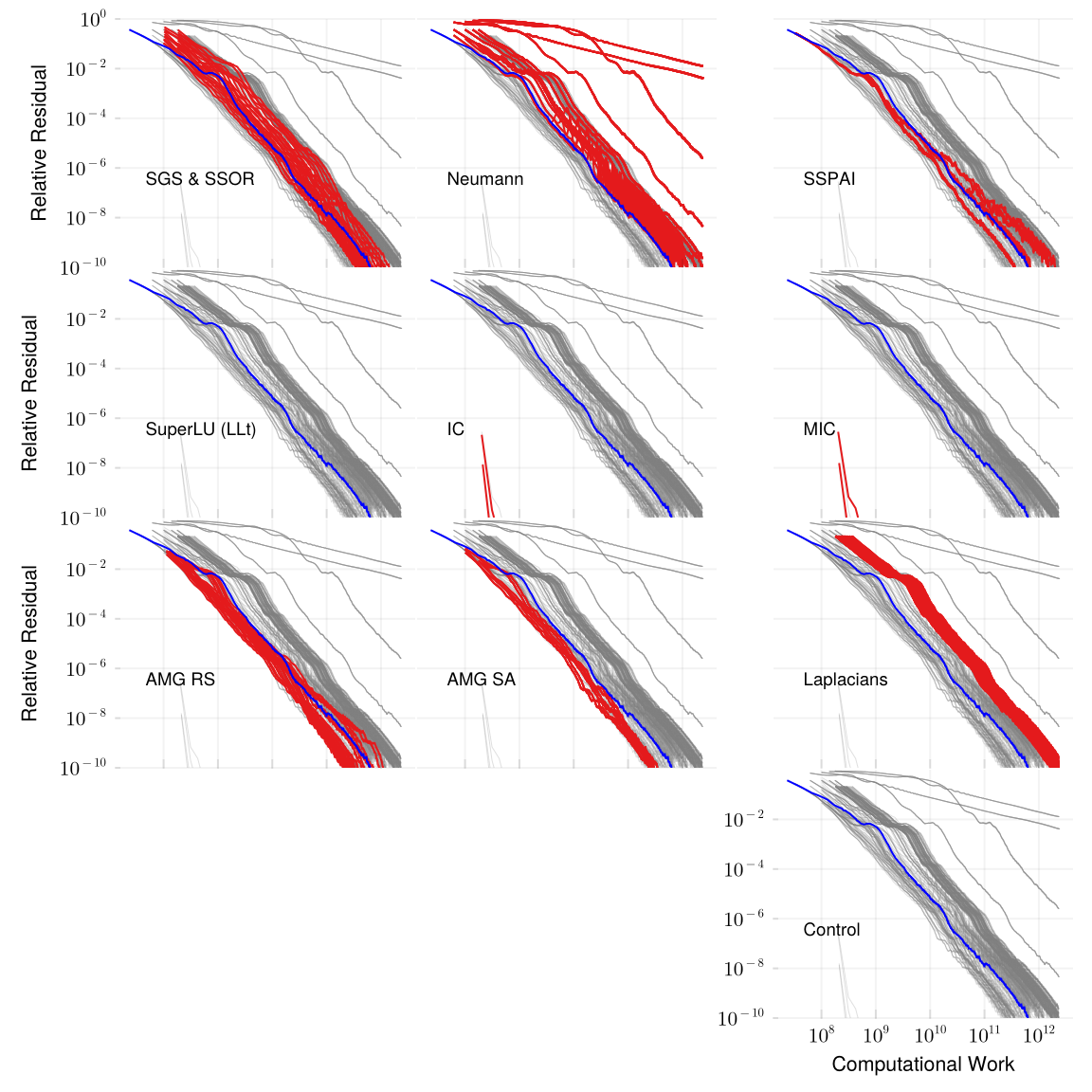}
    \caption{Convergence of the PCG method with various preconditioners applied to the \texttt{msdoor} matrix (416k rows, 19.2m non-zeros). The plots have a log-log scale.}
    \label{fig:msdoor}
\end{figure}
\clearpage 

\subsection{nasasrb}

\begin{figure}[!ht]
    \centering
    \includegraphics[width=\textwidth]{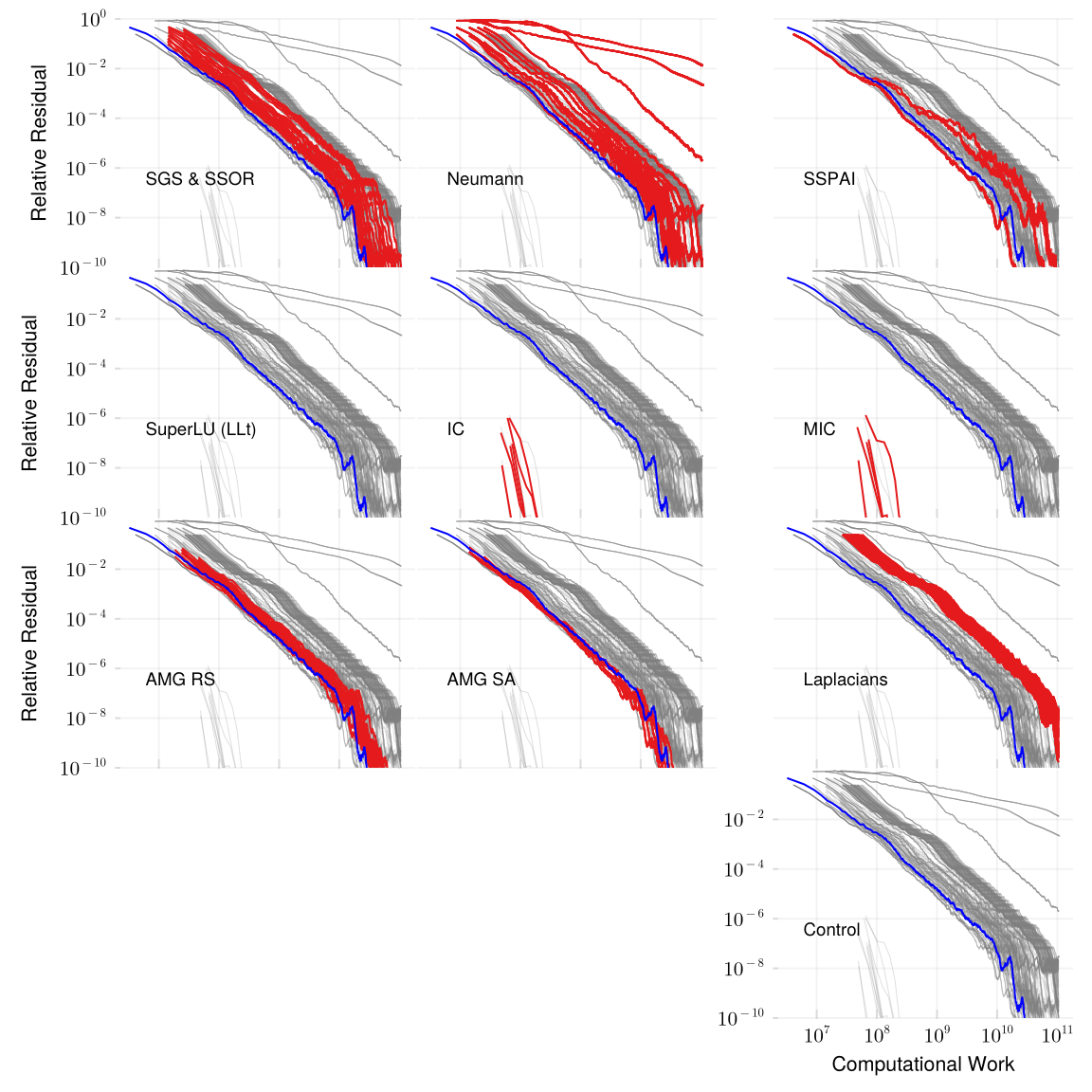}
    \caption{Convergence of the PCG method with various preconditioners applied to the \texttt{nasasrb} matrix (55k rows, 2.7m non-zeros). The plots have a log-log scale.}
    \label{fig:nasasrb}
\end{figure}
\clearpage 

\subsection{obstclae}

\begin{figure}[!ht]
    \centering
    \includegraphics[width=\textwidth]{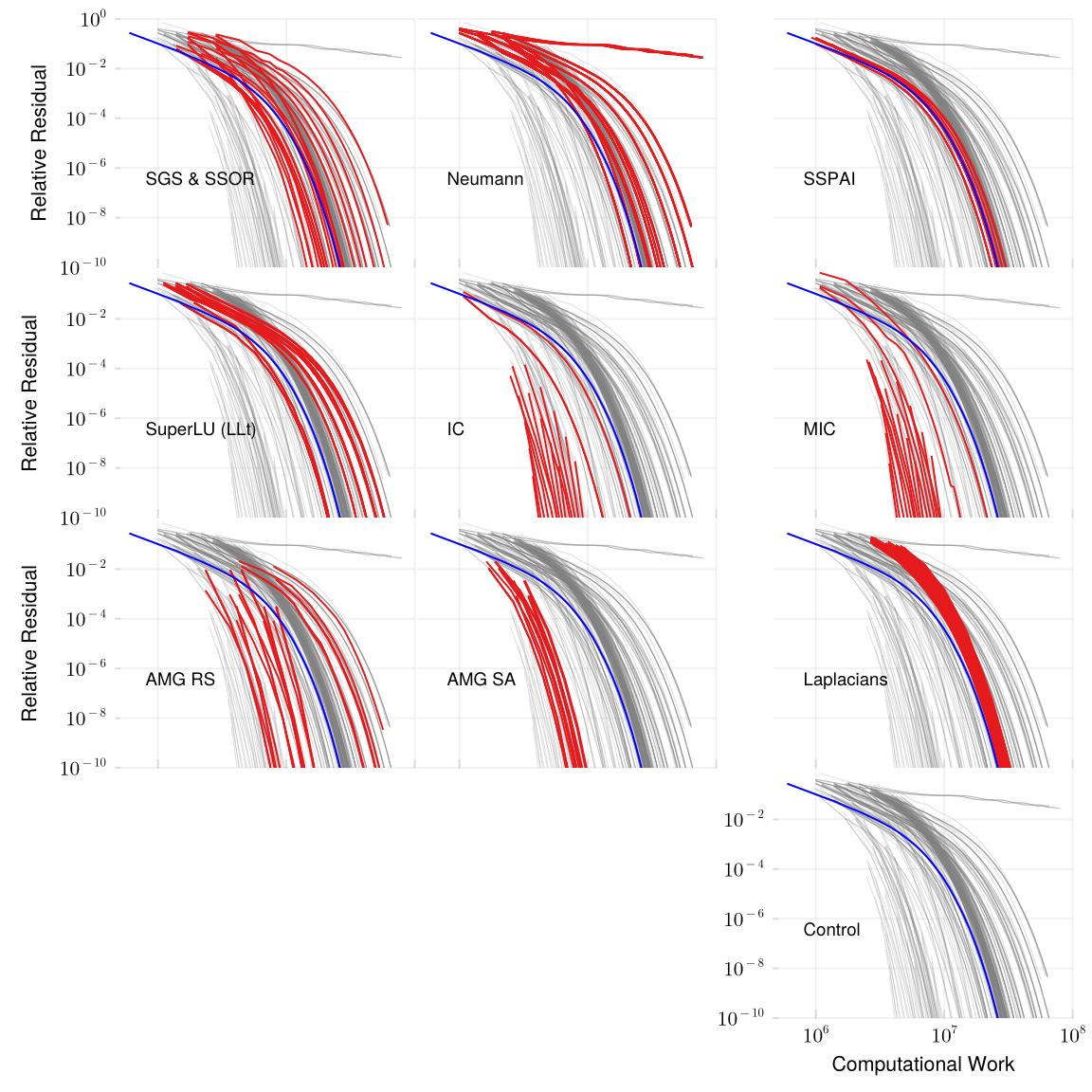}
    \caption{Convergence of the PCG method with various preconditioners applied to the \texttt{obstclae} matrix (40k rows, 198k non-zeros). The plots have a log-log scale.}
    \label{fig:obstclae}
\end{figure}
\clearpage 

\subsection{offshore}

\begin{figure}[!ht]
    \centering
    \includegraphics[width=\textwidth]{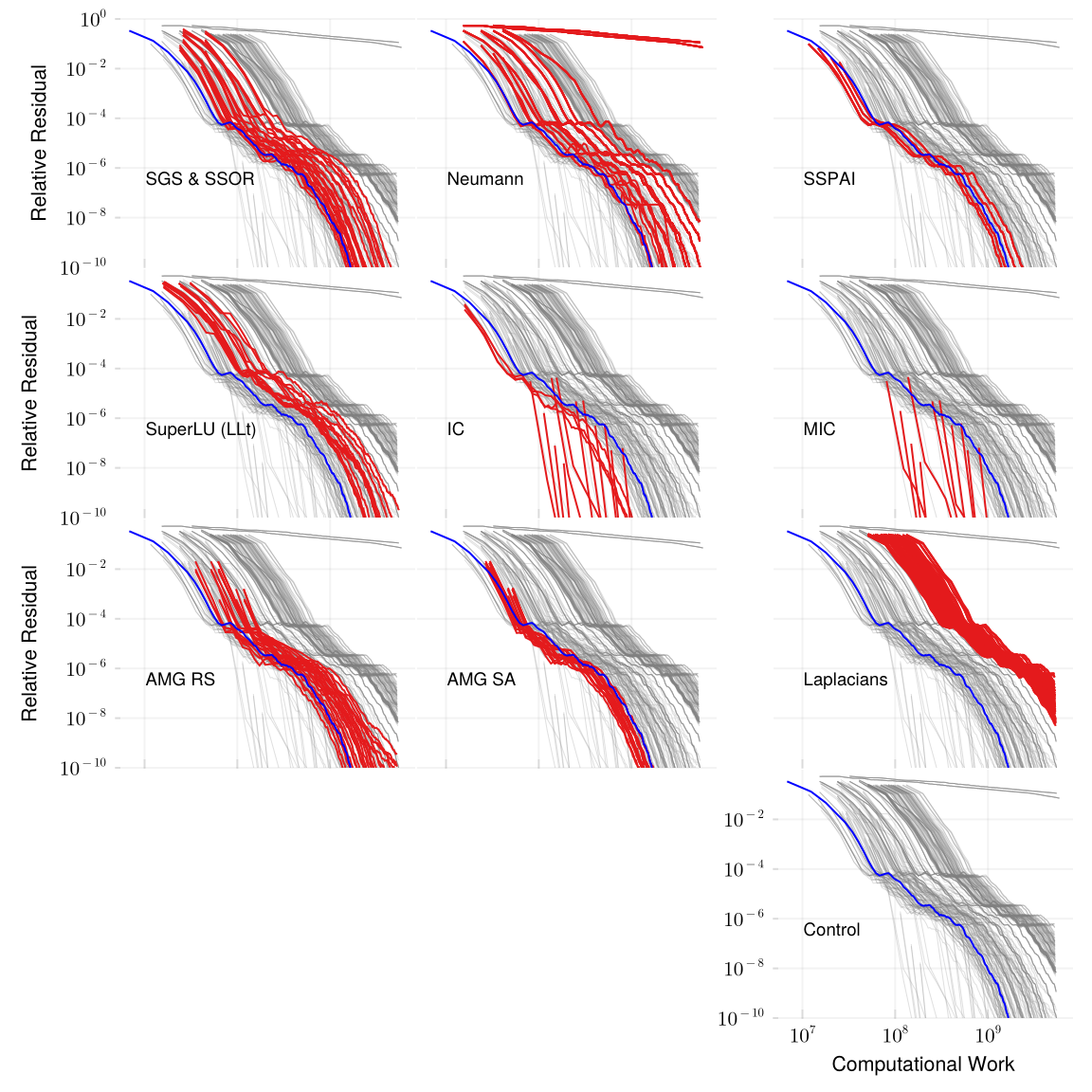}
    \caption{Convergence of the PCG method with various preconditioners applied to the \texttt{offshore} matrix (260k rows, 4.2m non-zeros). The plots have a log-log scale.}
    \label{fig:offshore}
\end{figure}
\clearpage 

\subsection{oilpan}

\begin{figure}[!ht]
    \centering
    \includegraphics[width=\textwidth]{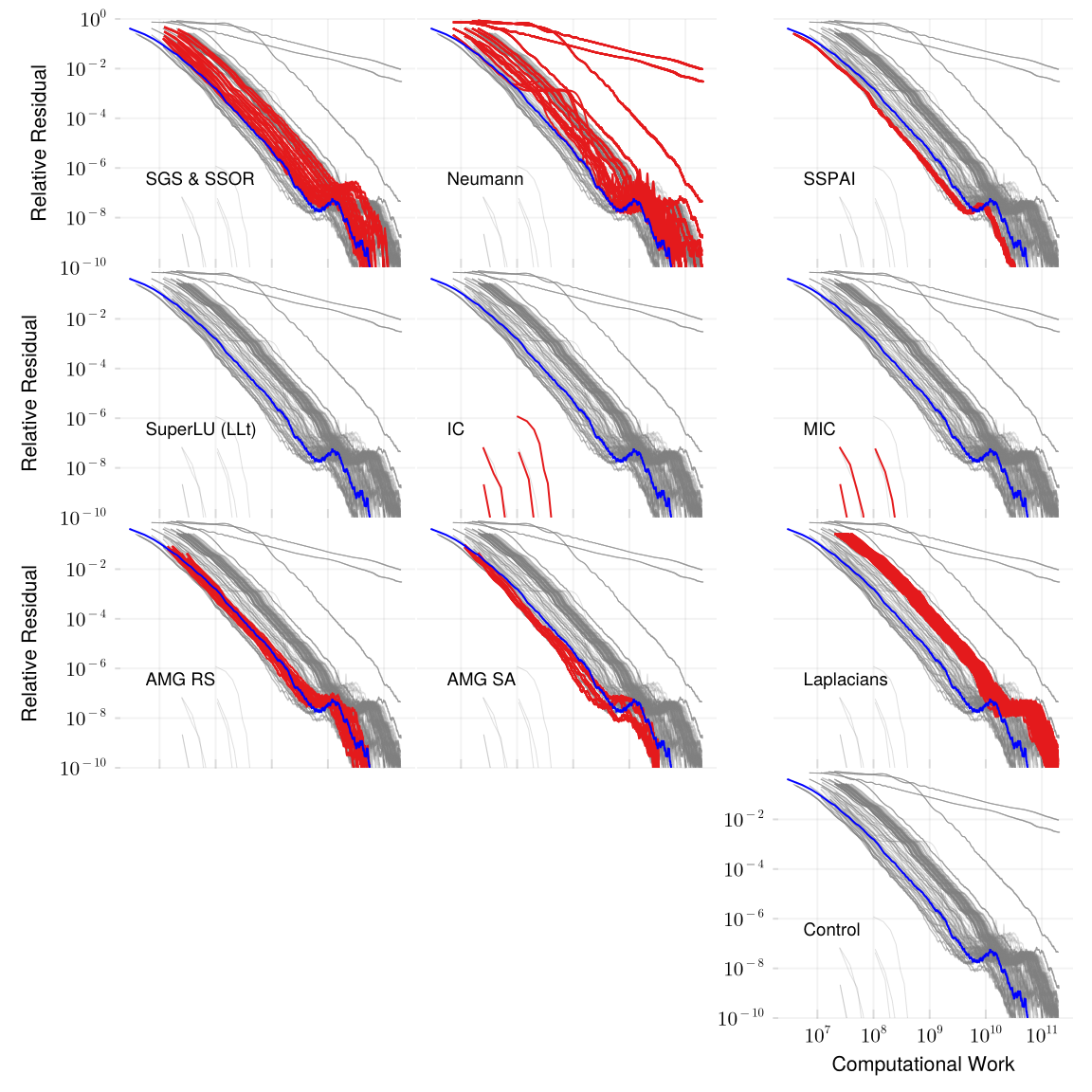}
    \caption{Convergence of the PCG method with various preconditioners applied to the \texttt{oilpan} matrix (74k rows, 2.1m non-zeros). The plots have a log-log scale.}
    \label{fig:oilpan}
\end{figure}
\clearpage 

\subsection{parabolic\_fem}

\begin{figure}[!ht]
    \centering
    \includegraphics[width=\textwidth]{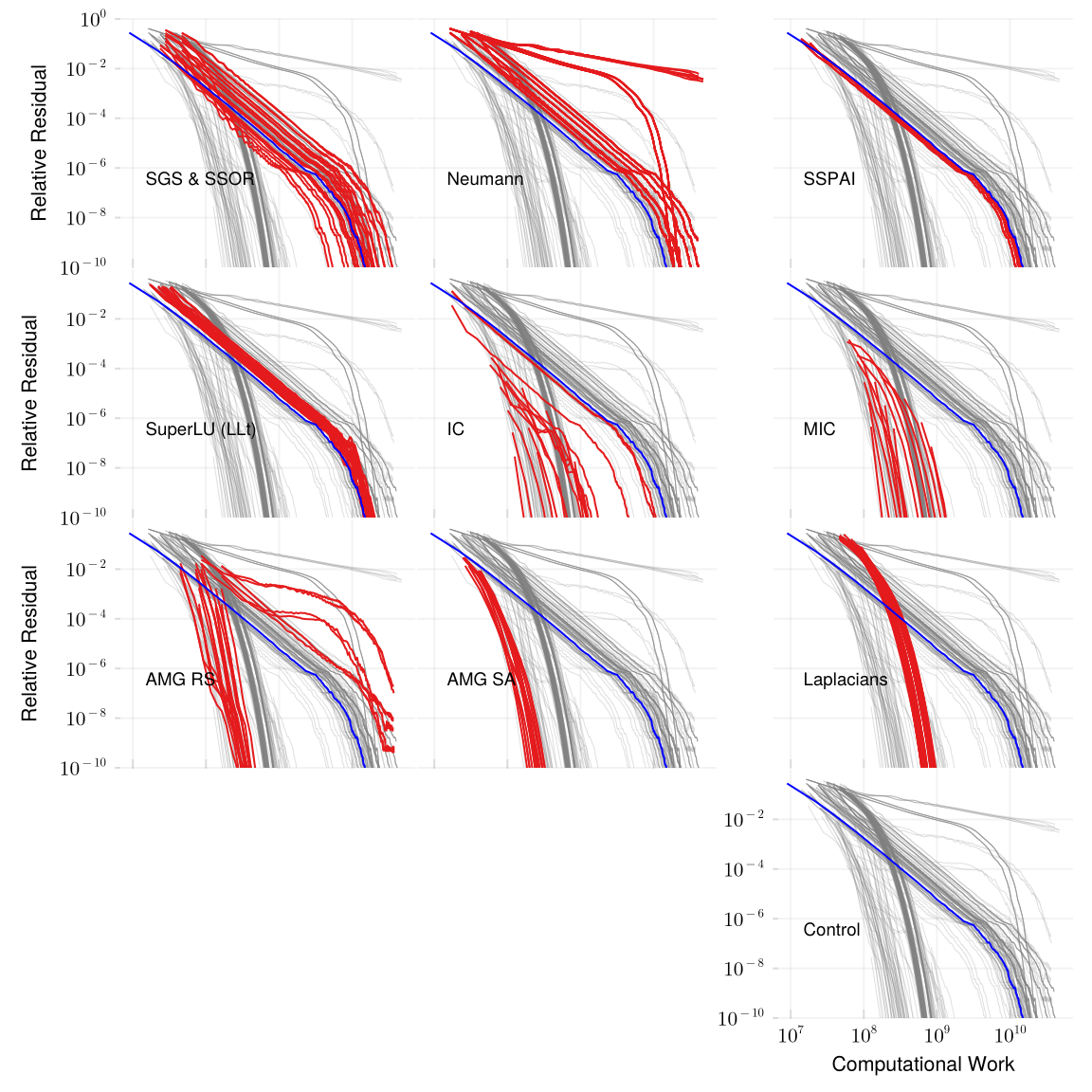}
    \caption{Convergence of the PCG method with various preconditioners applied to the \texttt{parabolic\_fem} matrix (526k rows, 3.7m non-zeros). The plots have a log-log scale.}
    \label{fig:parabolic_fem}
\end{figure}
\clearpage 

\subsection{qa8fm}

\begin{figure}[!ht]
    \centering
    \includegraphics[width=\textwidth]{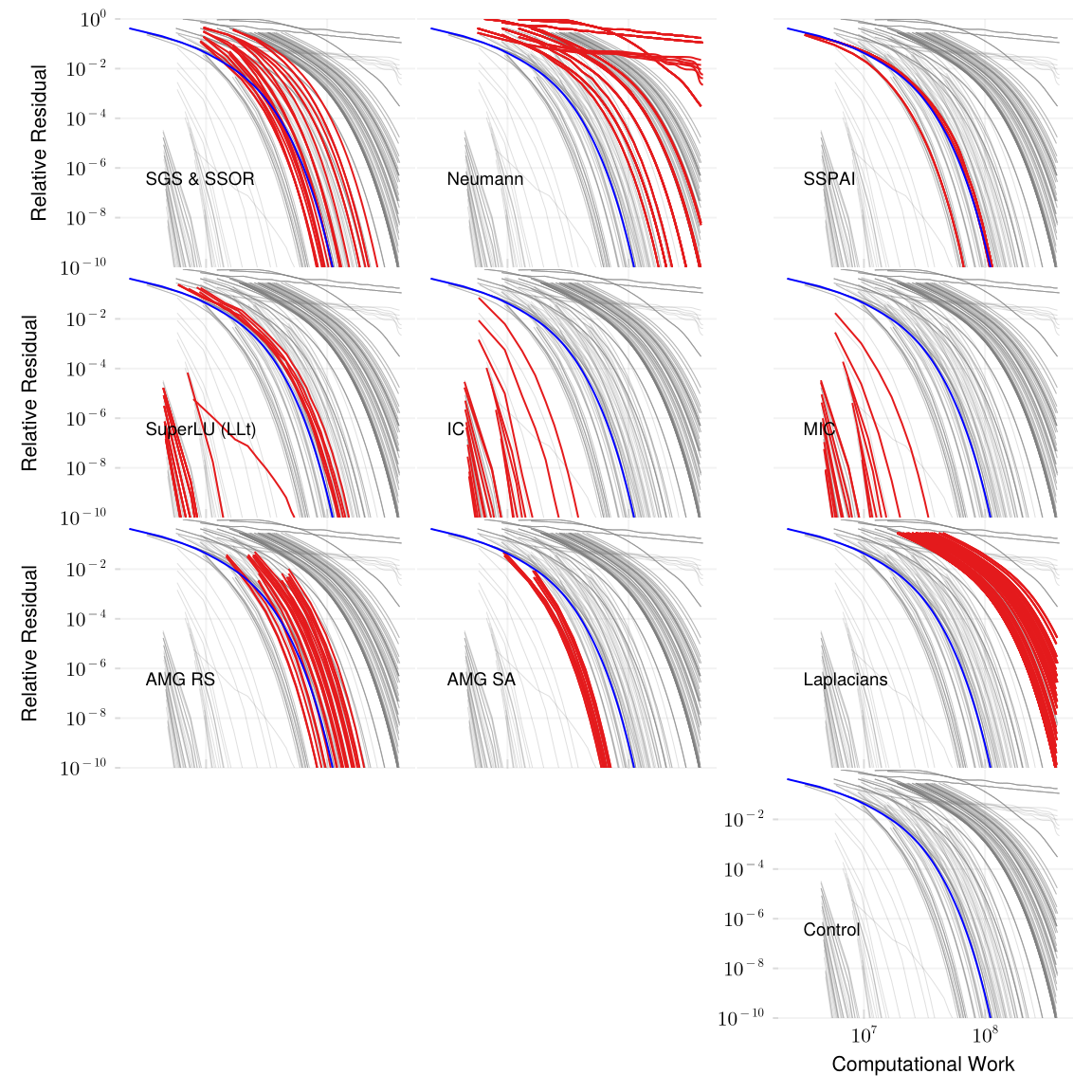}
    \caption{Convergence of the PCG method with various preconditioners applied to the \texttt{qa8fm} matrix (66k rows, 1.7m non-zeros). The plots have a log-log scale.}
    \label{fig:qa8fm}
\end{figure}
\clearpage 

\subsection{shallow\_water1}

\begin{figure}[!ht]
    \centering
    \includegraphics[width=\textwidth]{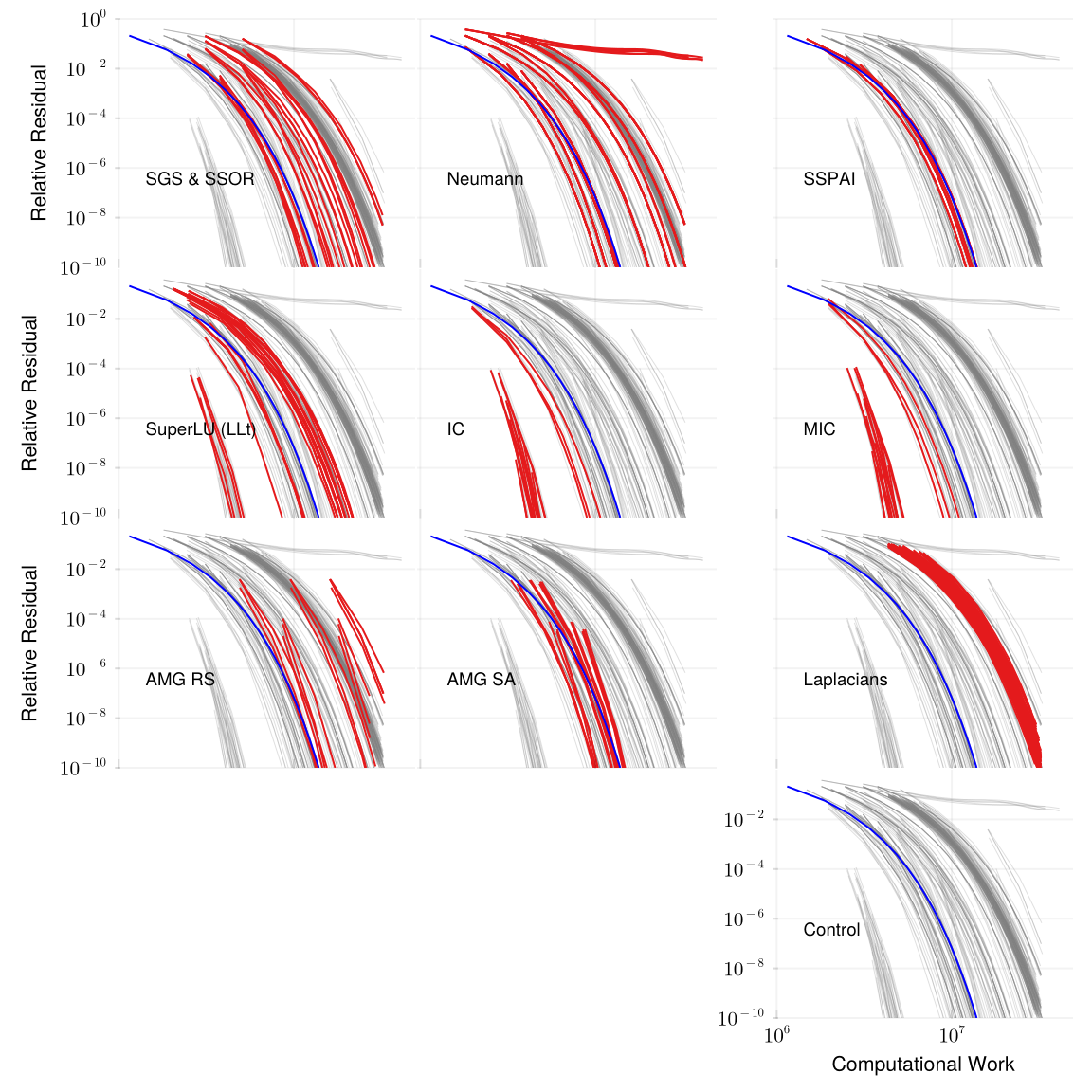}
    \caption{Convergence of the PCG method with various preconditioners applied to the \texttt{shallow\_water1} matrix (82k rows, 328k non-zeros). The plots have a log-log scale.}
    \label{fig:shallow_water1}
\end{figure}
\clearpage 

\subsection{shallow\_water2}

\begin{figure}[!ht]
    \centering
    \includegraphics[width=\textwidth]{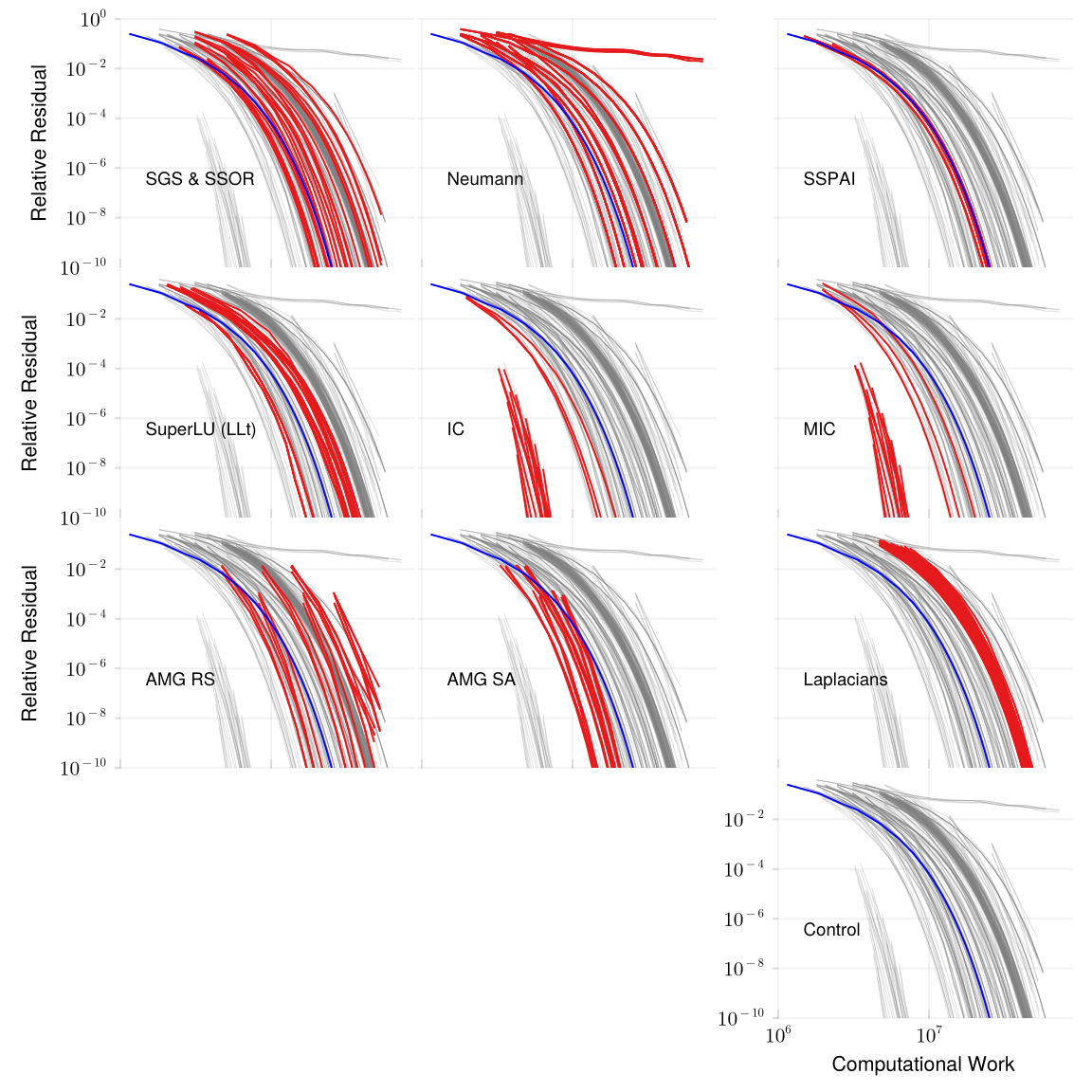}
    \caption{Convergence of the PCG method with various preconditioners applied to the \texttt{shallow\_water2} matrix (82k rows, 328k non-zeros). The plots have a log-log scale.}
    \label{fig:shallow_water2}
\end{figure}
\clearpage 

\subsection{ship\_001}

\begin{figure}[!ht]
    \centering
    \includegraphics[width=\textwidth]{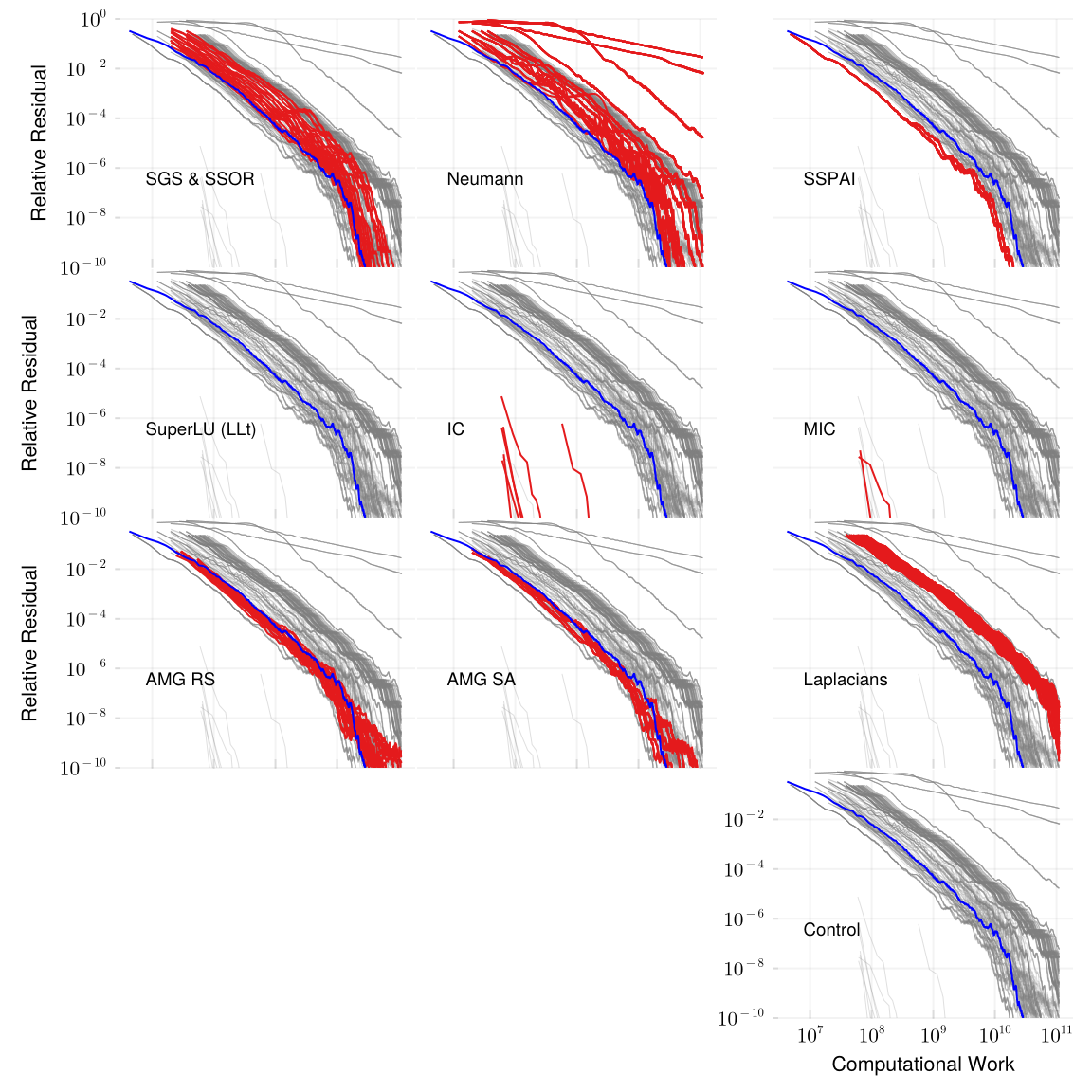}
    \caption{Convergence of the PCG method with various preconditioners applied to the \texttt{ship\_001} matrix (35k rows, 3.9m non-zeros). The plots have a log-log scale.}
    \label{fig:ship_001}
\end{figure}
\clearpage 

\subsection{ship\_003}

\begin{figure}[!ht]
    \centering
    \includegraphics[width=\textwidth]{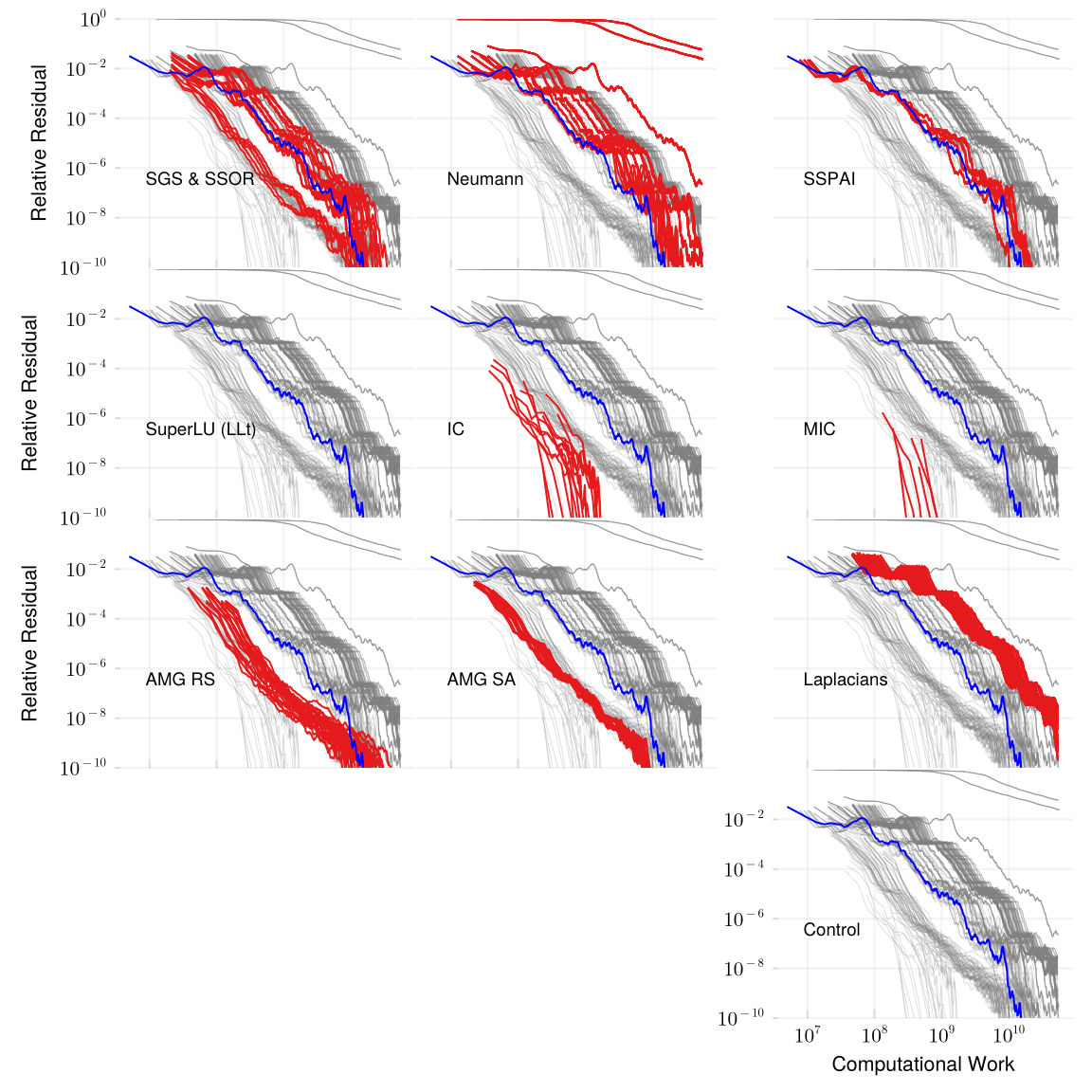}
    \caption{Convergence of the PCG method with various preconditioners applied to the \texttt{ship\_003} matrix (122k rows, 3.8m non-zeros). The plots have a log-log scale.}
    \label{fig:ship_003}
\end{figure}
\clearpage 

\subsection{shipsec1}

\begin{figure}[!ht]
    \centering
    \includegraphics[width=\textwidth]{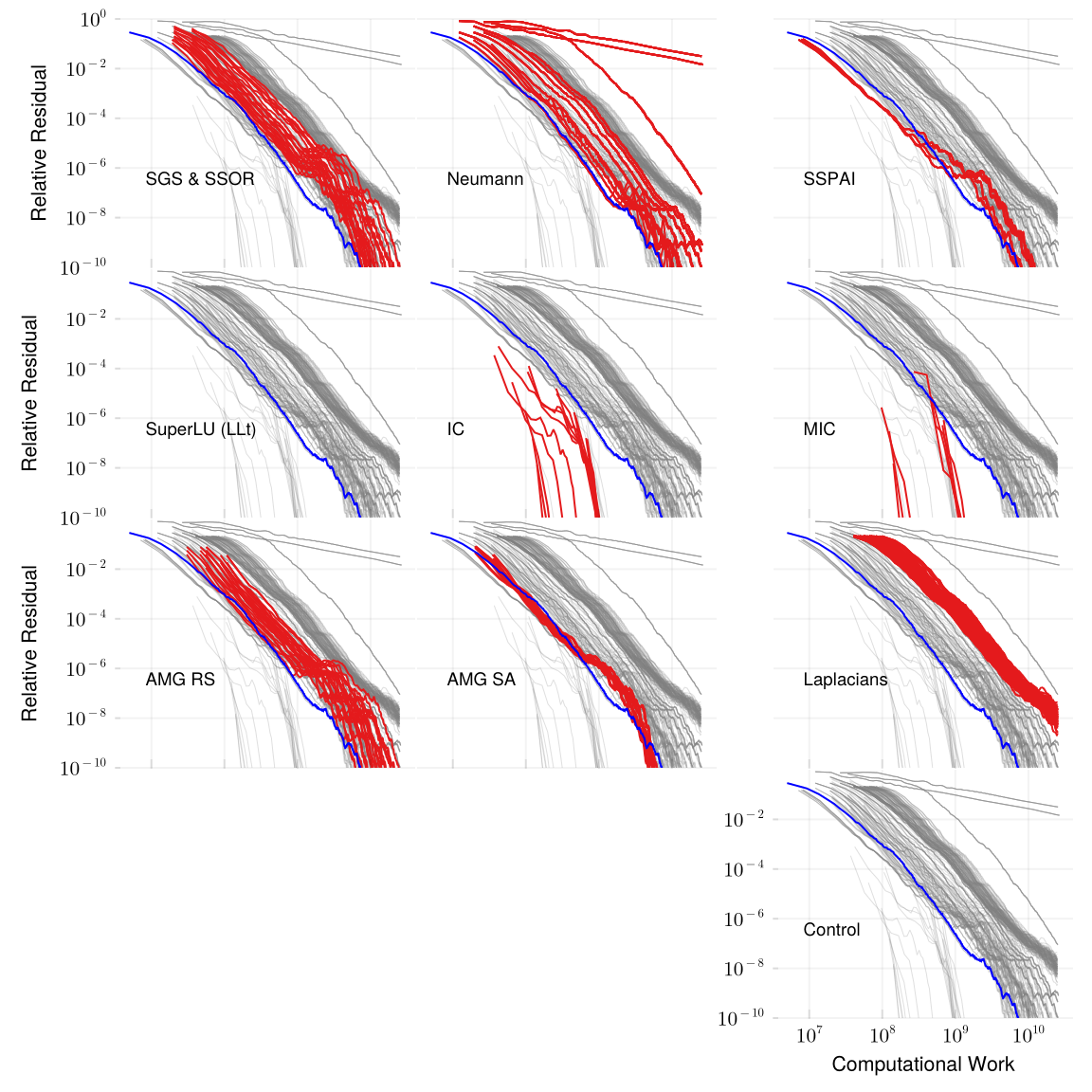}
    \caption{Convergence of the PCG method with various preconditioners applied to the \texttt{shipsec1} matrix (141k rows, 3.6m non-zeros). The plots have a log-log scale.}
    \label{fig:shipsec1}
\end{figure}
\clearpage 

\subsection{shipsec5}

\begin{figure}[!ht]
    \centering
    \includegraphics[width=\textwidth]{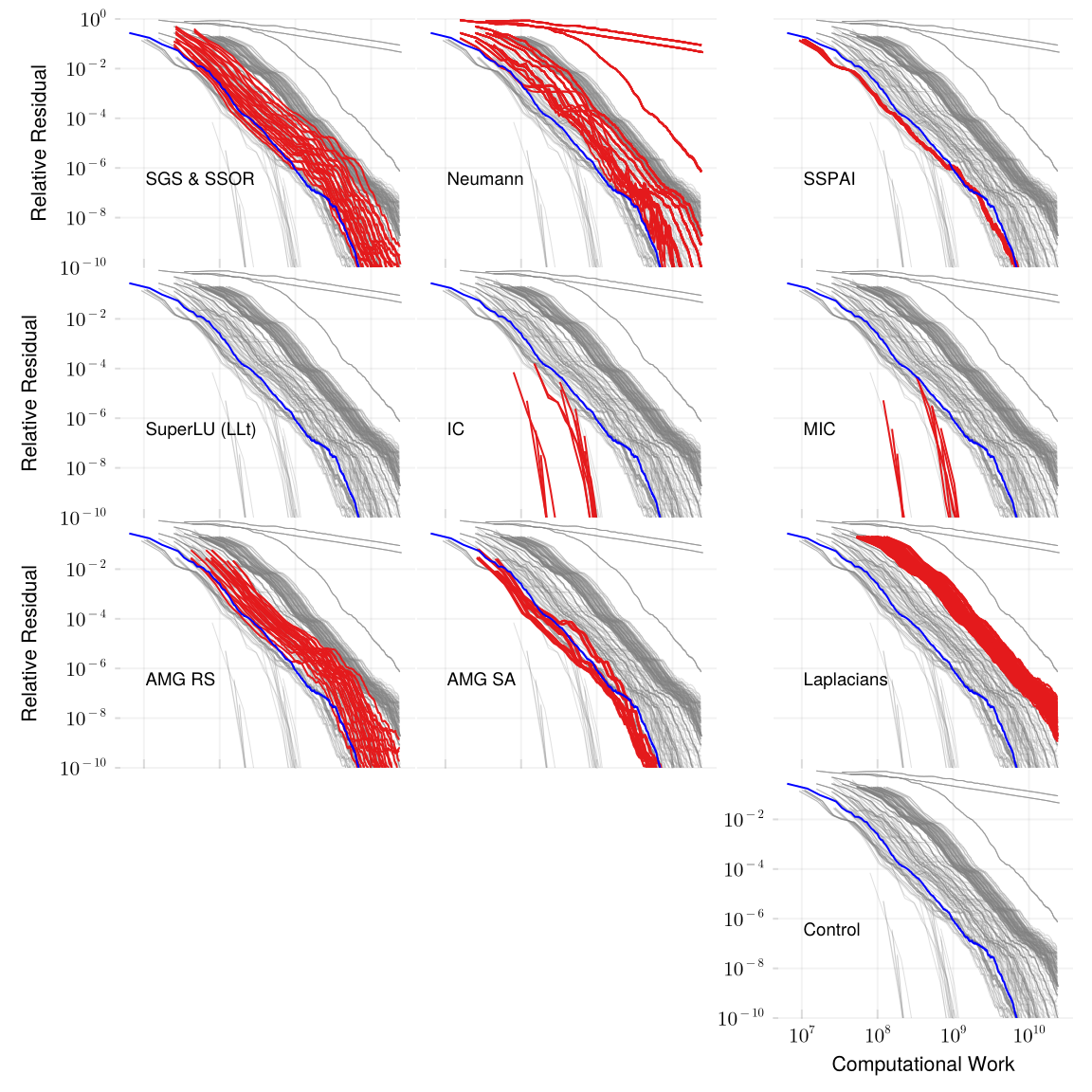}
    \caption{Convergence of the PCG method with various preconditioners applied to the \texttt{shipsec5} matrix (180k rows, 4.6m non-zeros). The plots have a log-log scale.}
    \label{fig:shipsec5}
\end{figure}
\clearpage 

\subsection{shipsec8}

\begin{figure}[!ht]
    \centering
    \includegraphics[width=\textwidth]{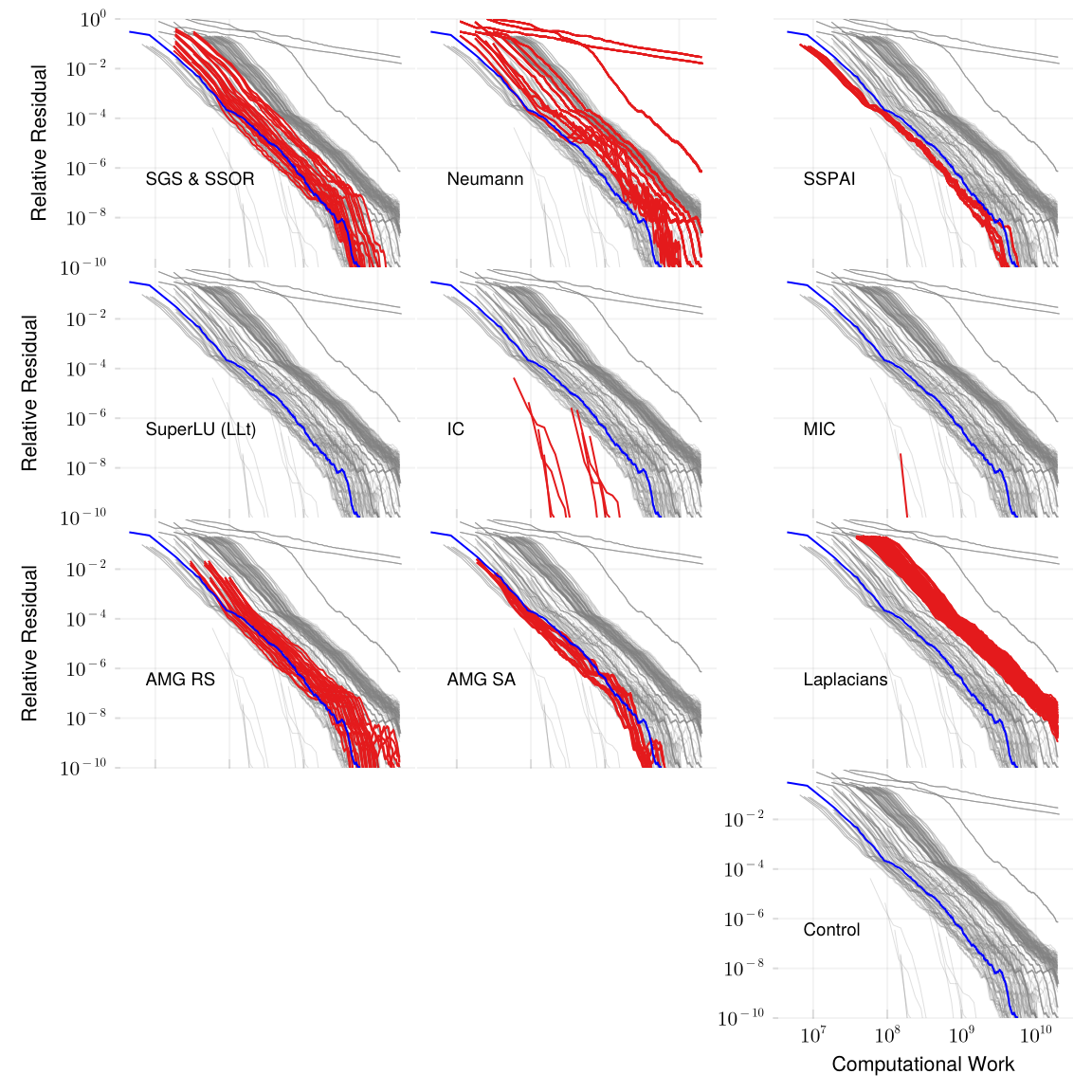}
    \caption{Convergence of the PCG method with various preconditioners applied to the \texttt{shipsec8} matrix (115k rows, 3.3m non-zeros). The plots have a log-log scale.}
    \label{fig:shipsec8}
\end{figure}
\clearpage 

\subsection{smt}

\begin{figure}[!ht]
    \centering
    \includegraphics[width=\textwidth]{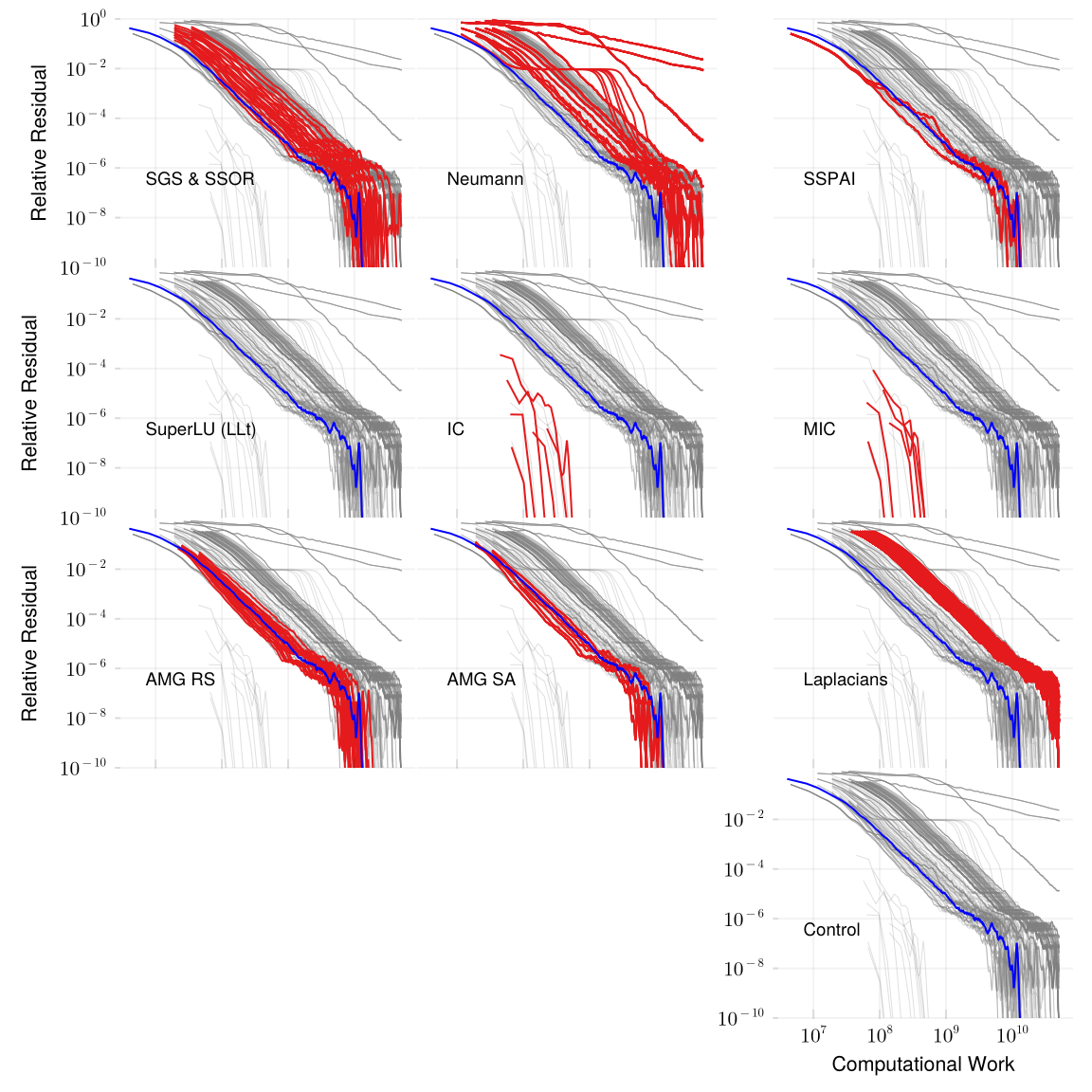}
    \caption{Convergence of the PCG method with various preconditioners applied to the \texttt{smt} matrix (26k rows, 3.7m non-zeros). The plots have a log-log scale.}
    \label{fig:smt}
\end{figure}
\clearpage 

\subsection{t2dah\_e}

\begin{figure}[!ht]
    \centering
    \includegraphics[width=\textwidth]{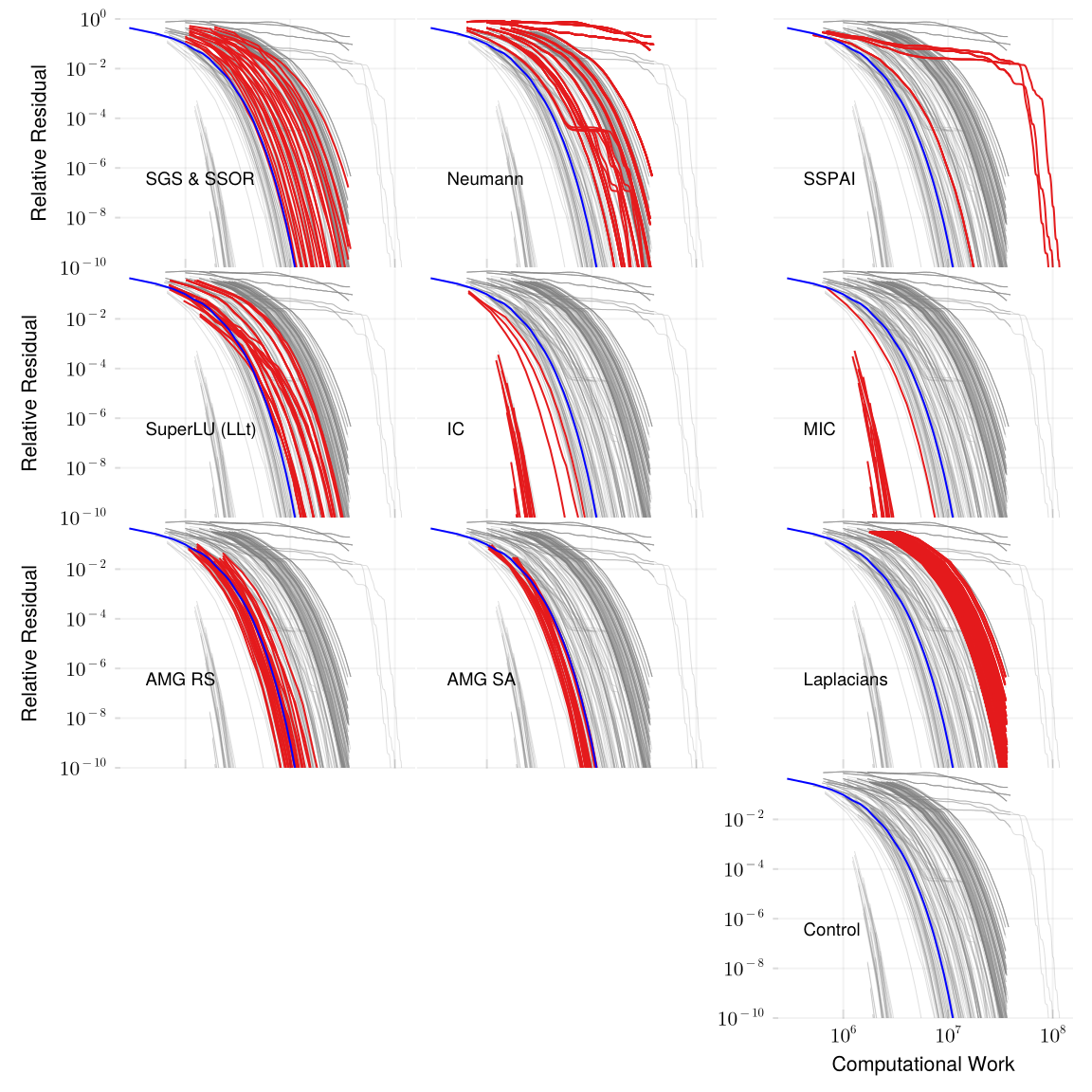}
    \caption{Convergence of the PCG method with various preconditioners applied to the \texttt{t2dah\_e} matrix (11k rows, 176k non-zeros). The plots have a log-log scale.}
    \label{fig:t2dah_e}
\end{figure}
\clearpage 

\subsection{ted\_B}

\begin{figure}[!ht]
    \centering
    \includegraphics[width=\textwidth]{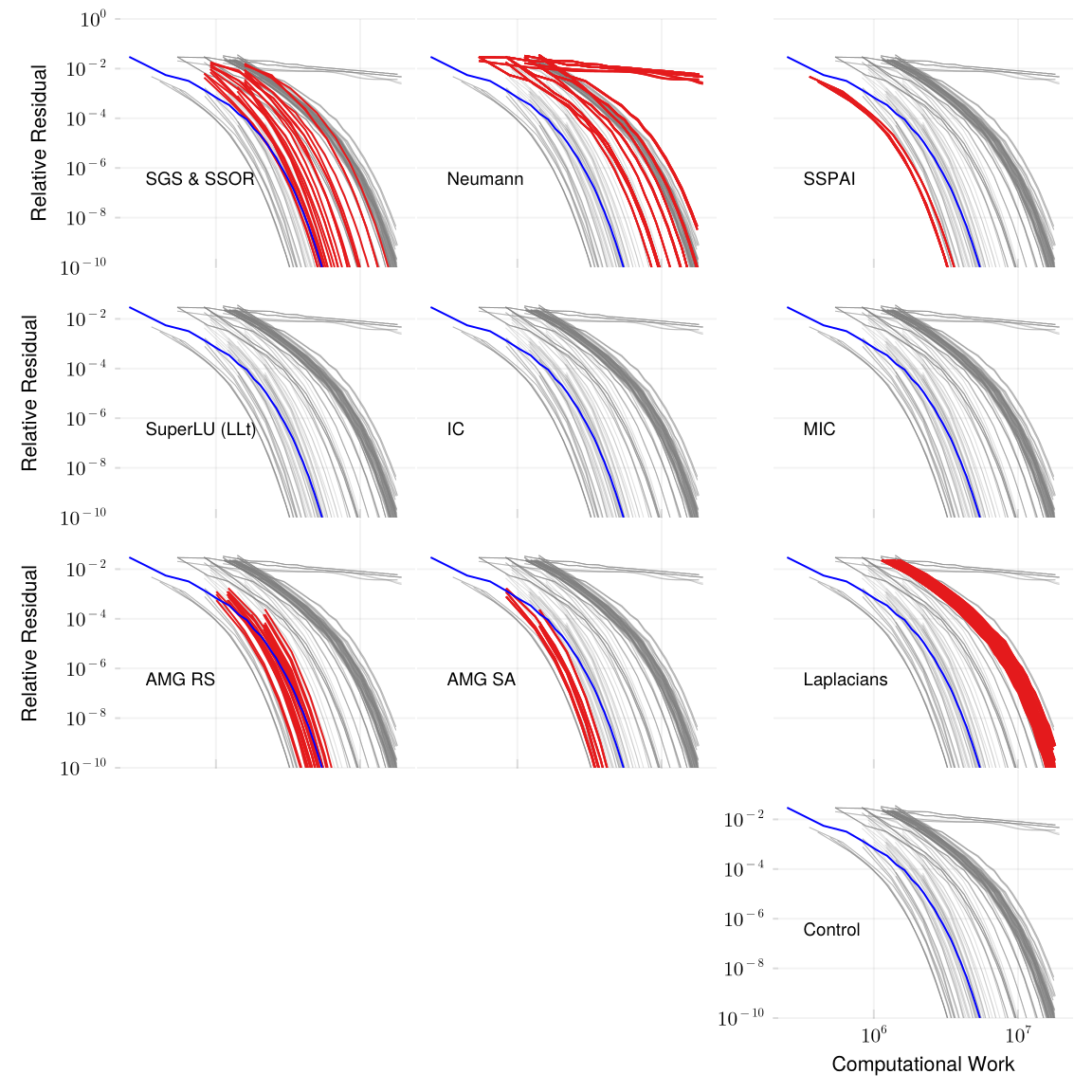}
    \caption{Convergence of the PCG method with various preconditioners applied to the \texttt{ted\_B} matrix (11k rows, 145k non-zeros). The plots have a log-log scale.}
    \label{fig:ted_B}
\end{figure}
\clearpage 

\subsection{thermal1}

\begin{figure}[!ht]
    \centering
    \includegraphics[width=\textwidth]{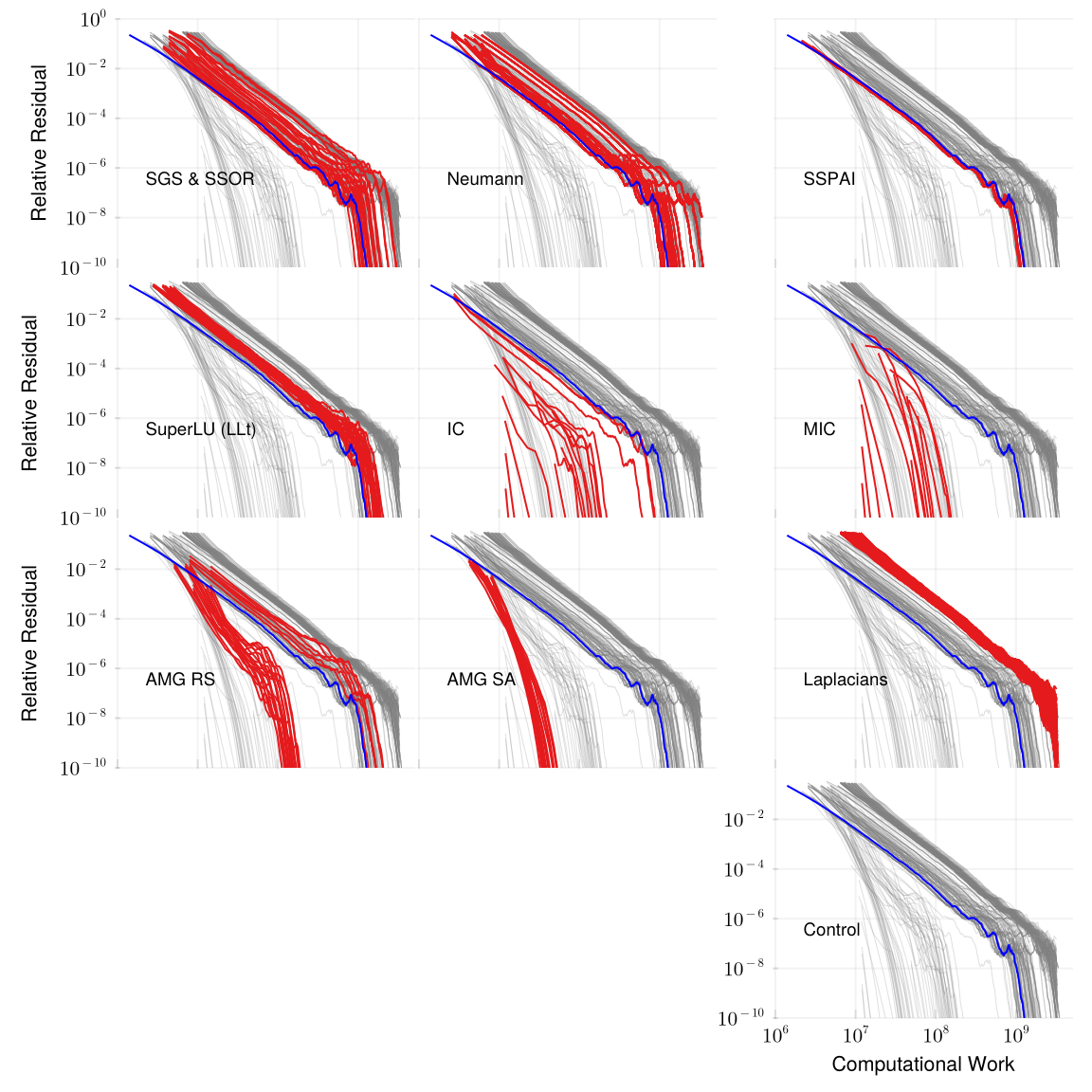}
    \caption{Convergence of the PCG method with various preconditioners applied to the \texttt{thermal1} matrix (83k rows, 574k non-zeros). The plots have a log-log scale.}
    \label{fig:thermal1}
\end{figure}
\clearpage 

\subsection{thermal2}

\begin{figure}[!ht]
    \centering
    \includegraphics[width=\textwidth]{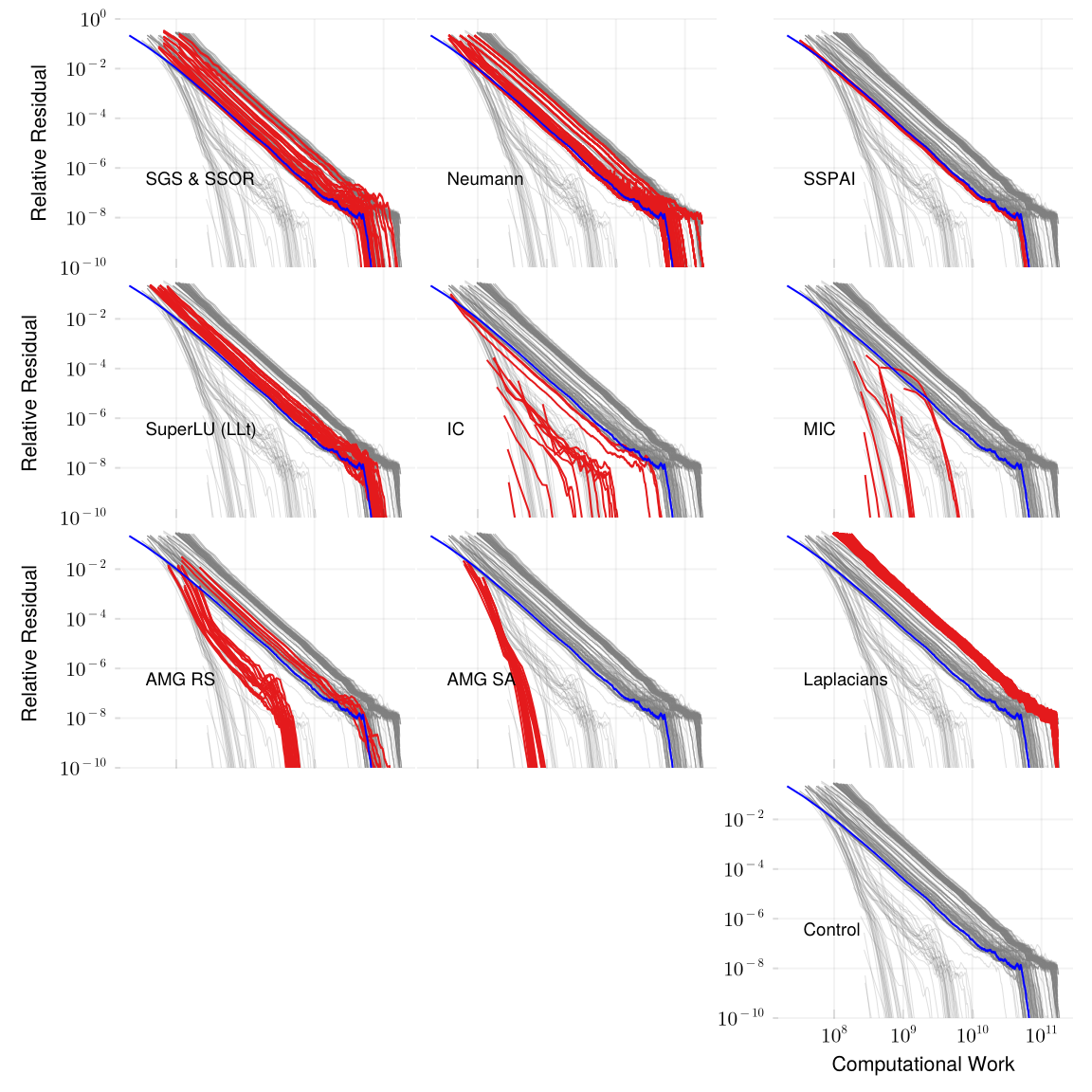}
    \caption{Convergence of the PCG method with various preconditioners applied to the \texttt{thermal2} matrix (1.2m rows, 8.6m non-zeros). The plots have a log-log scale.}
    \label{fig:thermal2}
\end{figure}
\clearpage 

\subsection{thermomech\_dM}

\begin{figure}[!ht]
    \centering
    \includegraphics[width=\textwidth]{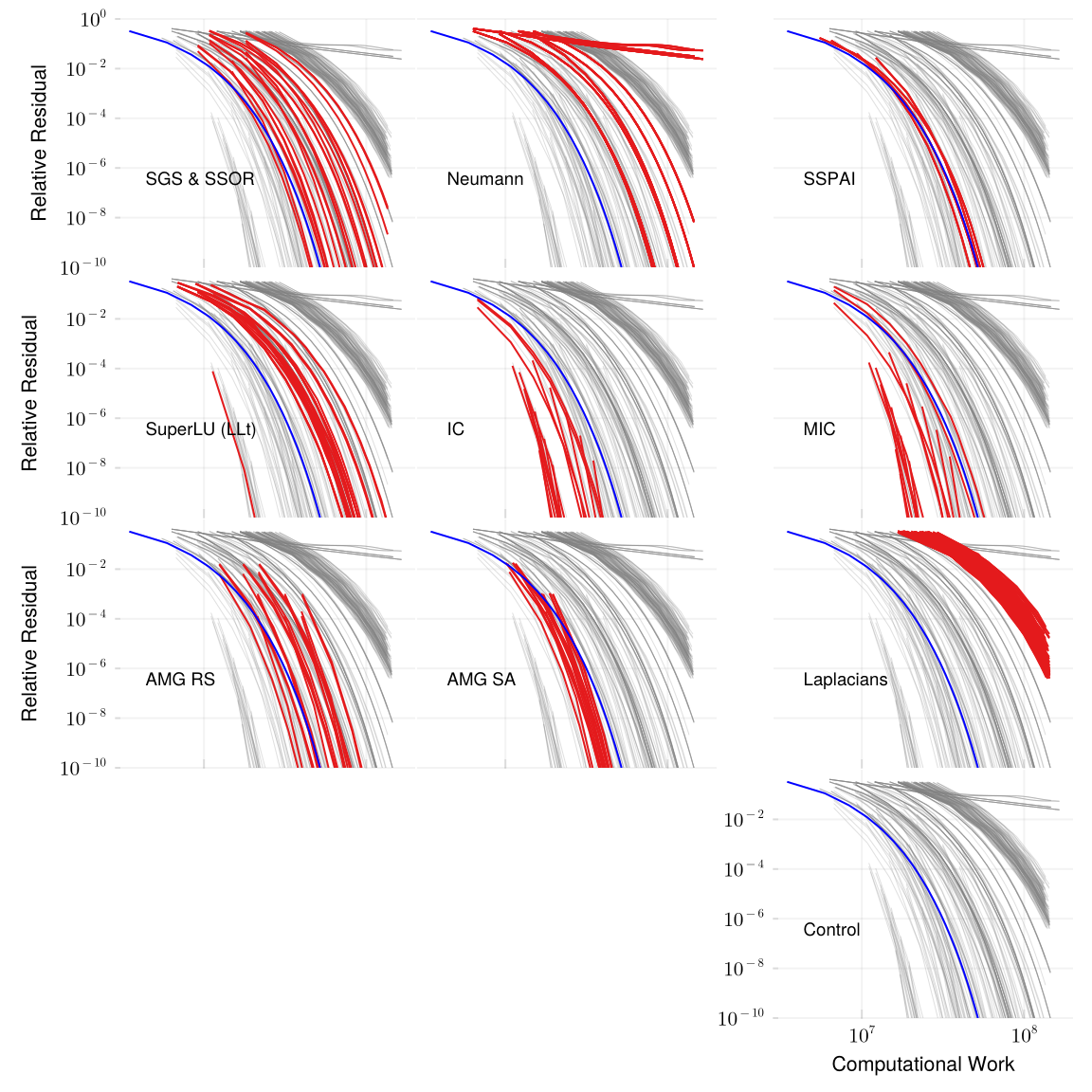}
    \caption{Convergence of the PCG method with various preconditioners applied to the \texttt{thermomech\_dM} matrix (204k rows, 1.4m non-zeros). The plots have a log-log scale.}
    \label{fig:thermomech_dM}
\end{figure}
\clearpage 

\subsection{thermomech\_TC}

\begin{figure}[!ht]
    \centering
    \includegraphics[width=\textwidth]{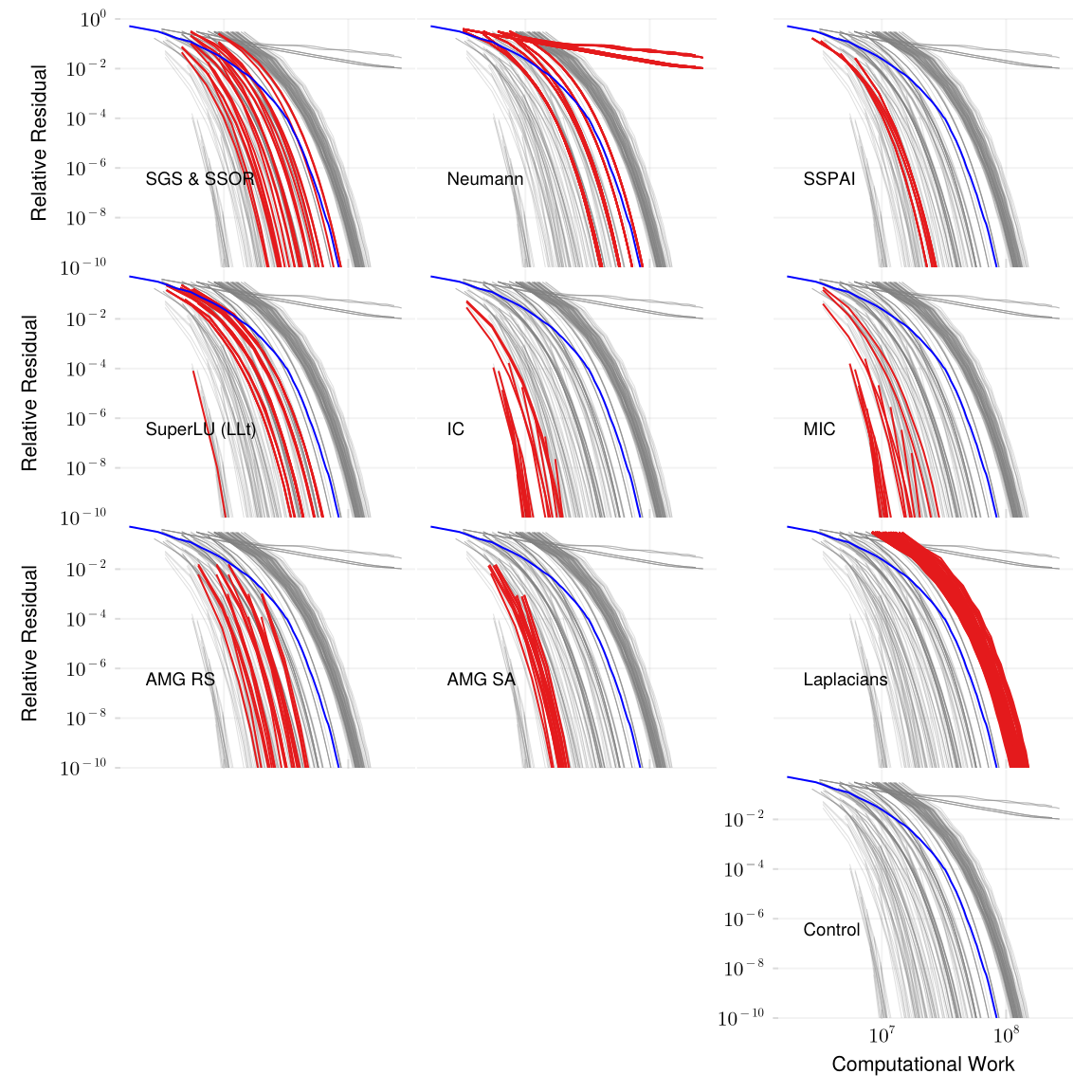}
    \caption{Convergence of the PCG method with various preconditioners applied to the \texttt{thermomech\_TC} matrix (102k rows, 712k non-zeros). The plots have a log-log scale.}
    \label{fig:thermomech_TC}
\end{figure}
\clearpage 

\subsection{thermomech\_TK}

\begin{figure}[!ht]
    \centering
    \includegraphics[width=\textwidth]{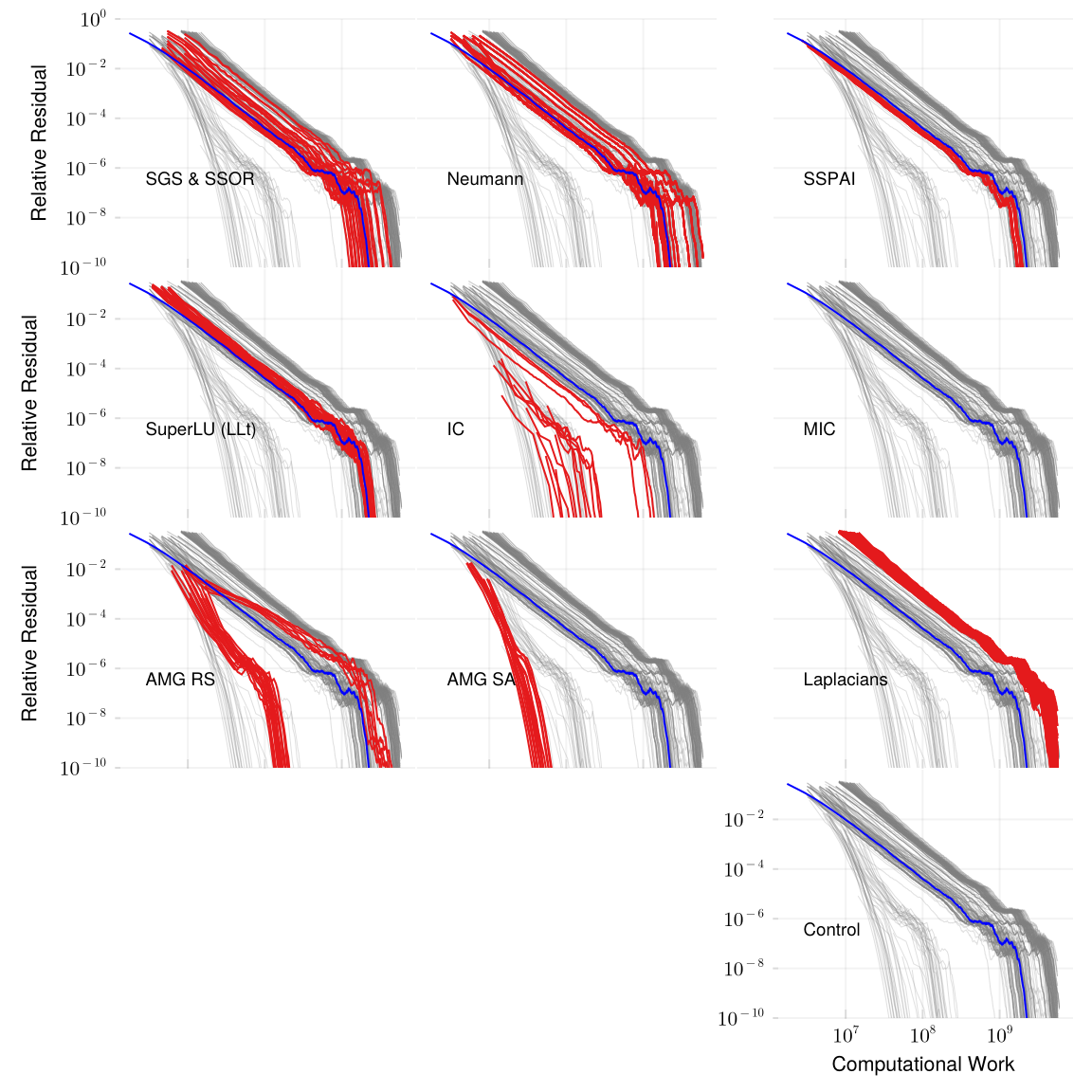}
    \caption{Convergence of the PCG method with various preconditioners applied to the \texttt{thermomech\_TK} matrix (102k rows, 712k non-zeros). The plots have a log-log scale.}
    \label{fig:thermomech_TK}
\end{figure}
\clearpage 

\subsection{tmt\_sym}

\begin{figure}[!ht]
    \centering
    \includegraphics[width=\textwidth]{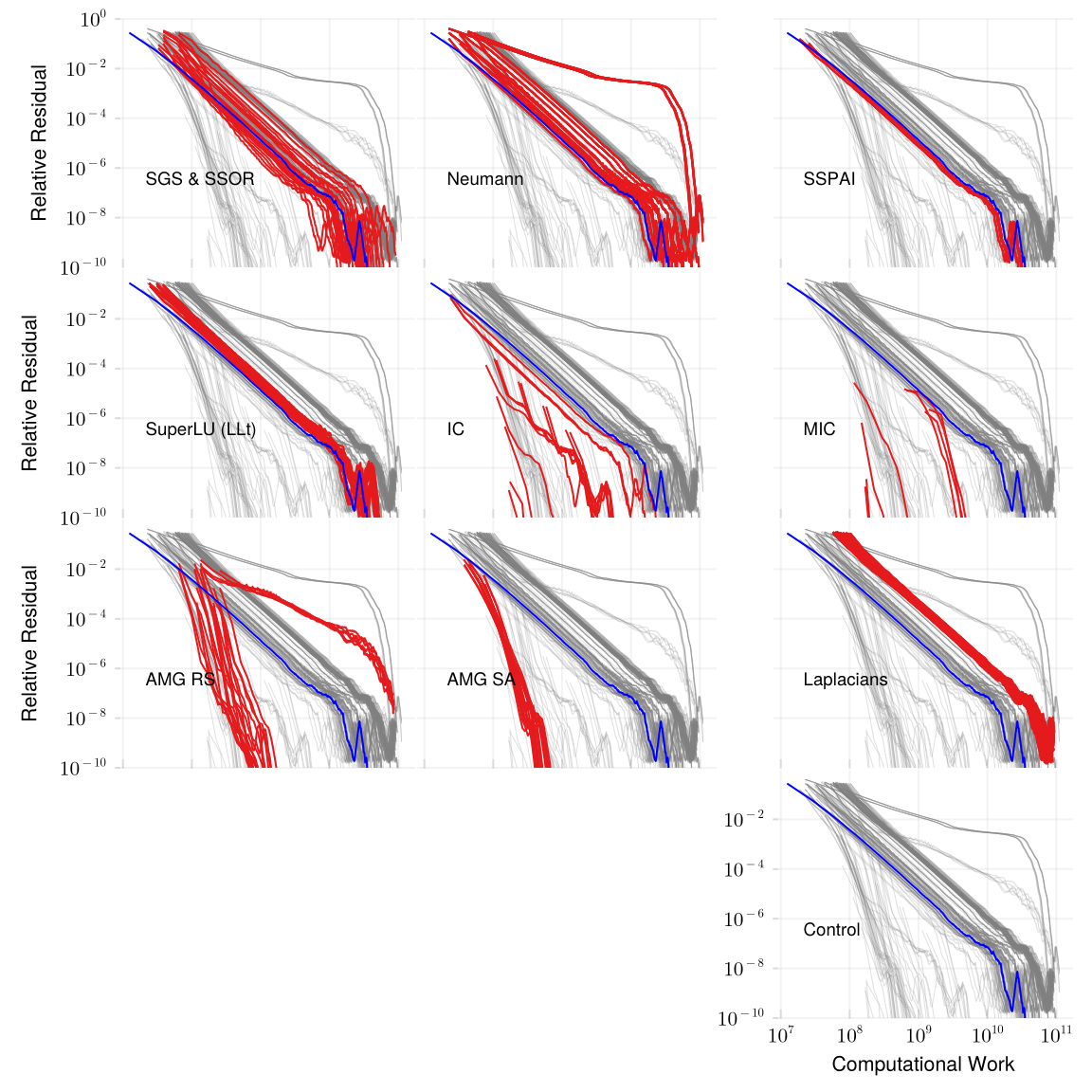}
    \caption{Convergence of the PCG method with various preconditioners applied to the \texttt{tmt\_sym} matrix (727k rows, 5.1m non-zeros). The plots have a log-log scale.}
    \label{fig:tmt_sym}
\end{figure}
\clearpage 

\subsection{wathen100}

\begin{figure}[!ht]
    \centering
    \includegraphics[width=\textwidth]{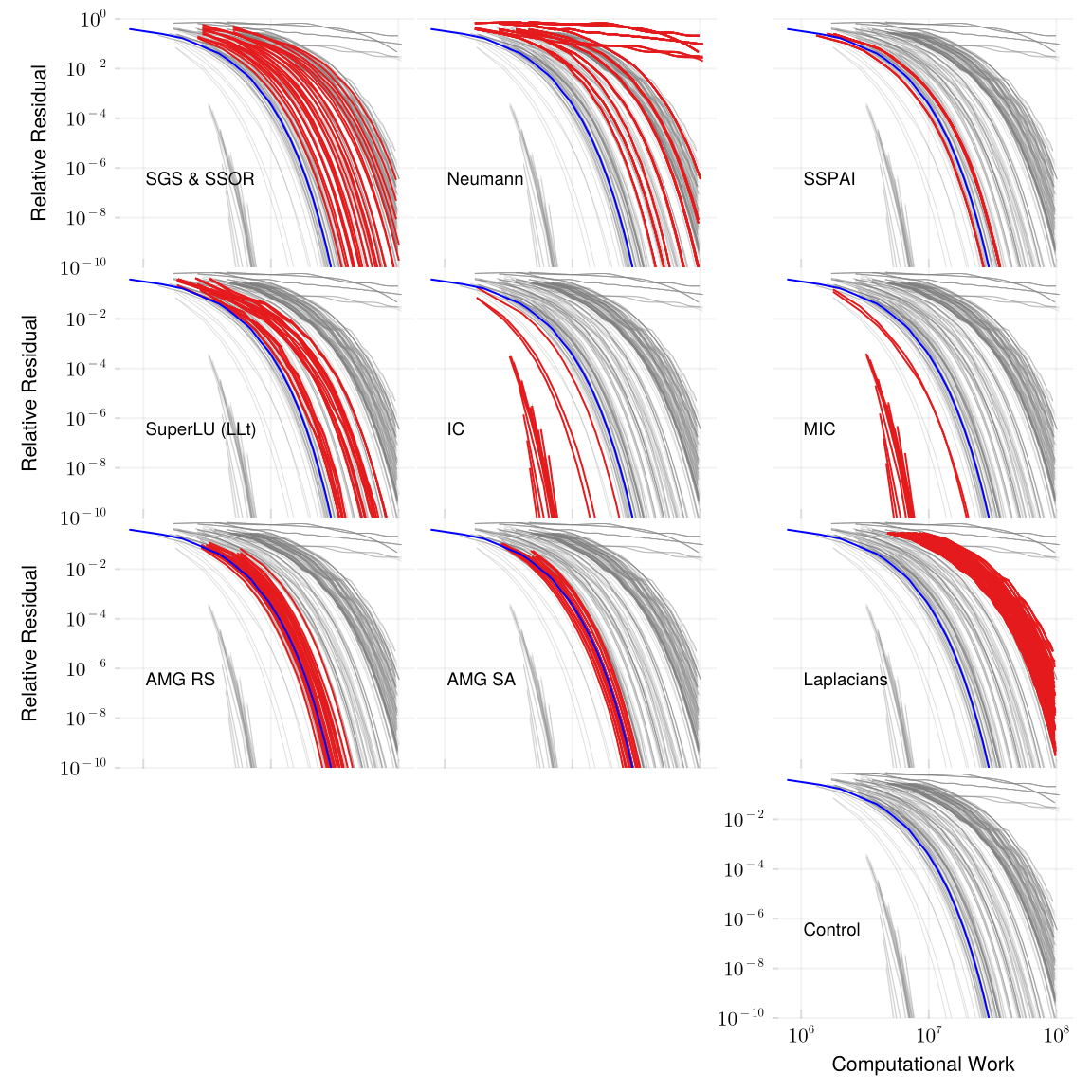}
    \caption{Convergence of the PCG method with various preconditioners applied to the \texttt{wathen100} matrix (30k rows, 472k non-zeros). The plots have a log-log scale.}
    \label{fig:wathen100}
\end{figure}
\clearpage 

\subsection{wathen120}

\begin{figure}[!ht]
    \centering
    \includegraphics[width=\textwidth]{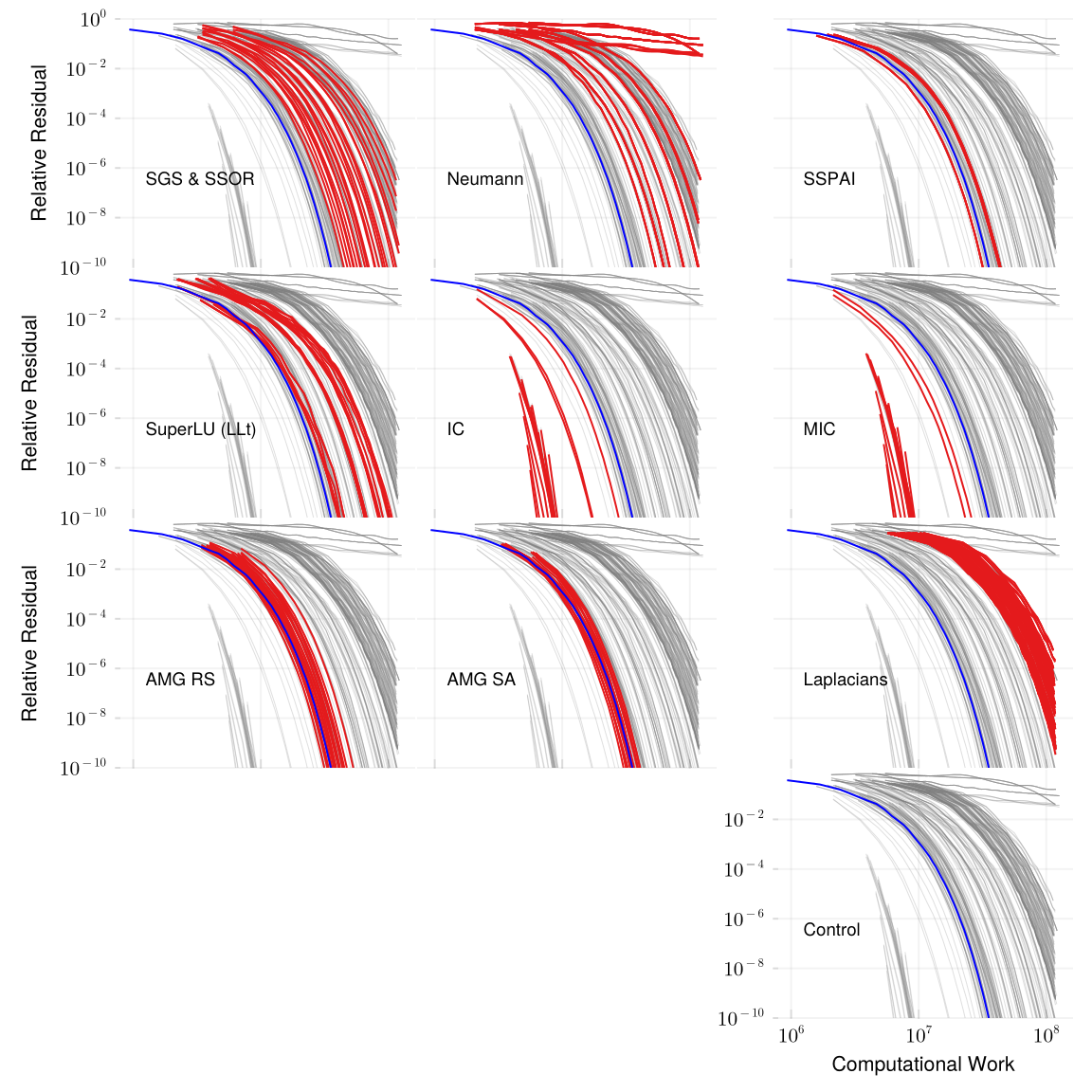}
    \caption{Convergence of the PCG method with various preconditioners applied to the \texttt{wathen120} matrix (36k rows, 566k non-zeros). The plots have a log-log scale.}
    \label{fig:wathen120}
\end{figure}
\clearpage